\documentclass[10pt]{article}
\usepackage{graphics} 
\usepackage{amssymb} 
\usepackage{theorem} 
\usepackage{latexsym} 
 
\makeatletter
  \renewcommand{\section}{\@startsection
  {section}%
  {1}%
  {0mm}%
  {\baselineskip}%
  {0.5\baselineskip}%
  {\bfseries}}%
\makeatother

\makeatletter
  \renewcommand{\subsection}{\@startsection
  {subsection}%
  {2}%
  {0mm}%
  {\baselineskip}%
  {0.5\baselineskip}%
  {\bfseries \itshape}}%
\makeatother

\makeatletter
  \renewcommand{\subsubsection}{\@startsection
  {subsubsection}%
  {3}%
  {0mm}%
  {\baselineskip}%
  {0.5\baselineskip}%
  {\itshape}}%
\makeatother

\makeatletter
  \renewcommand{\@seccntformat}[1]{%
    \normalfont \bfseries \normalsize{\csname the#1\endcsname}.\hspace{0.5em}}
\makeatother

\newcommand{\hidari}{\Bigl\langle\!\!\Bigl\langle} 
\newcommand{\migi}{\Bigr\rangle\!\!\Bigr\rangle}
\newcommand{\hidarii}{\Bigl\langle}
\newcommand{\migii}{\Bigr\rangle} 

\newcommand{\h}{\Bigl(}
\newcommand{\m}{\Bigr)}
\newcommand{\hh}{\Bigl[}
\newcommand{\mm}{\Bigr]} 

\newcommand{\LL}{\mathbf{L}} 
 
\newcommand{\K}{\mathbf{K}}

\newcommand{\noin}{\in\hspace{-0.3cm}/~}
\newcommand{\noeq}{=\hspace{-0.3cm}/~}

\theorembodyfont{\rmfamily}
\newtheorem{definition}{Definition}[section]
\newtheorem{lemma}[definition]{Lemma}
\newtheorem{theorem}{\bfseries Theorem}
\newtheorem{corollary}{\bfseries Corollary}
\newtheorem{fact}[definition]{Fact}
\newtheorem{conjecture}[definition]{Conjecture}
\newtheorem{example}[definition]{Example}

\makeatletter
  
  \@addtoreset{definition}{section}
\makeatother

  

  

\makeatletter
  
  \@addtoreset{equation}{section}
\makeatother

\makeatletter 
  \renewcommand{\fnum@figure}[1]{%
  \footnotesize{Fig. \thefigure.}
  }
\makeatother

\newcommand{\maru}{
\unitlength 0.1in
\begin{picture}(3.40,2.40)(1.5,-3.6)
\special{pn 16}%
\special{ar 320 320 120 120  0.0000000 6.2831853}%
\end{picture}%
}
\newcommand{\marua}{
\unitlength 0.1in
\begin{picture}(3.40,2.40)(1.5,-3.6)
\special{pn 16}%
\special{ar 320 320 120 120  0.0000000 6.2831853}%
\special{pn 4}%
\special{pa 200 320}%
\special{pa 440 320}%
\special{fp}%
\end{picture}%
}
\newcommand{\maruba}{
\unitlength 0.1in
\begin{picture}(3.4,2.42)(1.5,-3.6)
\special{pn 16}%
\special{ar 319 319 121 121  0.0000000 6.2831853}%
\special{pn 4}%
\special{pa 202 289}%
\special{pa 436 289}%
\special{fp}%
\special{pn 4}%
\special{pa 204 349}%
\special{pa 435 349}%
\special{fp}%
\special{pa 435 349}%
\special{pa 436 346}%
\special{fp}%
\end{picture}%
}
\newcommand{\marubb}{
\unitlength 0.1in
\begin{picture}(3.40,2.40)(1.5,-3.6)
\special{pn 16}%
\special{ar 320 320 120 120  1.5625321 6.2831853}%
\special{ar 320 320 120 120  0.0000000 1.5624632}%
\special{pn 4}%
\special{pa 200 320}%
\special{pa 440 320}%
\special{fp}%
\special{pn 4}%
\special{pa 320 200}%
\special{pa 320 440}%
\special{fp}%
\end{picture}%
}

\newcommand{\marubc}{
\unitlength 0.1in
\begin{picture}(3.40,2.40)(1.5,-3.6)
\special{pn 16}%
\special{ar 320 320 120 120  0.0000000 6.2831853}%
\special{pn 4}%
\special{ar 320 320 30 30  0.0000000 6.2831853}%
\special{pn 4}%
\special{pa 200 320}%
\special{pa 290 320}%
\special{fp}%
\special{pn 4}%
\special{pa 350 320}%
\special{pa 440 320}%
\special{fp}%
\end{picture}%
}

\newcommand{\marubd}{
\unitlength 0.1in
\begin{picture}(3.40,2.40)(1.5,-3.6)
\special{pn 16}%
\special{ar 320 320 120 120  0.0000000 6.2831853}%
\special{pn 8}%
\special{pa 200 320}%
\special{pa 440 320}%
\special{dt 0.045}%
\special{pa 440 320}%
\special{pa 439 320}%
\special{dt 0.045}%
\special{pa 440 320}%
\special{pa 440 320}%
\special{dt 0.045}%
\special{pn 8}%
\special{pa 320 200}%
\special{pa 320 440}%
\special{dt 0.035}%
\special{pa 320 440}%
\special{pa 320 439}%
\special{dt 0.035}%
\end{picture}%
}
\newcommand{\marube}{
\unitlength 0.1in
\begin{picture}(3.40,2.40)(1.5,-3.6)
\special{pn 16}%
\special{ar 320 320 120 120  0.0000000 6.2831853}%
\special{pn 4}%
\special{pa 200 320}%
\special{pa 440 320}%
\special{fp}%
\put(3.8000,-2.9000){\makebox(0,0){\tiny{2s}}}%
\end{picture}%
}
\newcommand{\maruca}{
\unitlength 0.1in
\begin{picture}(6.9,2.40)(1.5,-3.60)
\special{pn 16}%
\special{ar 320 320 120 120  0.0000000 6.2831853}%
\special{pn 16}%
\special{ar 680 320 120 120  0.0000000 6.2831853}%
\special{pn 4}%
\special{pa 425 260}%
\special{pa 576 260}%
\special{fp}%
\special{pn 4}%
\special{pa 425 380}%
\special{pa 576 380}%
\special{fp}%
\end{picture}%
}
\newcommand{\marucb}{
\unitlength 0.1in
\begin{picture}(6.9,2.40)(1.5,-3.6)
\special{pn 16}%
\special{ar 320 320 120 120  0.0000000 6.2831853}%
\special{pn 16}%
\special{ar 680 320 120 120  0.0000000 6.2831853}%
\special{pn 4}%
\special{pa 440 320}%
\special{pa 560 320}%
\special{fp}%
\end{picture}%
}
\newcommand{\marucc}{
\unitlength 0.1in
\begin{picture}(6.9,2.40)(1.5,-3.60)
\special{pn 16}%
\special{ar 320 320 120 120  0.0000000 6.2831853}%
\special{pn 16}%
\special{ar 680 320 120 120  0.0000000 6.2831853}%
\special{pn 4}%
\special{pa 440 320}%
\special{pa 560 320}%
\special{fp}%
\special{pn 4}%
\special{pa 320 200}%
\special{pa 320 440}%
\special{fp}%
\end{picture}%
}

\newcommand{\marudd}{
\unitlength 0.1in
\begin{picture}(3.40,2.40)(1.5,-3.6)
\special{pn 16}%
\special{ar 320 320 120 120  0.0000000 6.2831853}%
\special{pn 4}%
\special{pa 320 200}%
\special{pa 320 440}%
\special{fp}%
\special{pn 4}%
\special{pa 209 275}%
\special{pa 431 275}%
\special{fp}%
\special{pn 4}%
\special{pa 209 365}%
\special{pa 429 365}%
\special{fp}%
\end{picture}%
}
\newcommand{\marude}{
\unitlength 0.1in
\begin{picture}(3.40,2.40)(1.5,-3.6)
\special{pn 16}%
\special{ar 320 320 120 120  0.0000000 6.2831853}%
\special{pn 4}%
\special{pa 320 200}%
\special{pa 320 440}%
\special{fp}%
\special{pn 4}%
\special{pa 209 275}%
\special{pa 432 365}%
\special{fp}%
\special{pn 4}%
\special{pa 209 365}%
\special{pa 431 275}%
\special{fp}%
\end{picture}%
}

\newcommand{\marudf}{
\unitlength 0.1in
\begin{picture}(3.40,2.40)(1.50,-3.60)
\special{pn 16}%
\special{ar 320 320 120 120  0.0000000 6.2831853}%
\special{pn 4}%
\special{pa 209 275}%
\special{pa 432 275}%
\special{fp}%
\special{pn 4}%
\special{ar 320 365 30 30  0.0000000 6.2831853}%
\special{pn 4}%
\special{pa 208 366}%
\special{pa 288 366}%
\special{fp}%
\special{pn 4}%
\special{pa 350 366}%
\special{pa 432 366}%
\special{fp}%
\end{picture}%
}

\newcommand{\marudg}{
\unitlength 0.1in
\begin{picture}(3.40,2.40)(1.50,-3.60)
\special{pn 16}%
\special{ar 320 320 120 120  0.0000000 6.2831853}%
\special{pn 4}%
\special{ar 275 320 30 30  0.0000000 6.2831853}%
\special{pn 4}%
\special{ar 365 320 30 30  0.0000000 6.2831853}%
\special{pn 4}%
\special{pa 200 320}%
\special{pa 240 320}%
\special{fp}%
\special{pn 4}%
\special{pa 306 320}%
\special{pa 330 320}%
\special{fp}%
\special{pn 4}%
\special{pa 396 320}%
\special{pa 440 320}%
\special{fp}%
\end{picture}%
}

\newcommand{\marudh}{
\unitlength 0.1in
\begin{picture}(3.40,2.40)(1.5,-3.6)
\special{pn 16}%
\special{ar 320 320 120 120  0.0000000 6.2831853}%
\special{pn 4}%
\special{pa 200 320}%
\special{pa 440 320}%
\special{fp}%
\special{pn 4}%
\special{pa 215 260}%
\special{pa 425 260}%
\special{fp}%
\special{pn 4}%
\special{pa 216 380}%
\special{pa 423 380}%
\special{fp}%
\end{picture}%
}

\newcommand{\marudi}{
\unitlength 0.1in
\begin{picture}(3.40,2.40)(1.5,-3.6)
\special{pn 16}%
\special{ar 320 320 120 120  0.0000000 6.2831853}%
\special{pn 4}%
\special{pa 320 200}%
\special{pa 320 440}%
\special{fp}%
\special{pn 8}%
\special{pa 209 275}%
\special{pa 431 275}%
\special{dt 0.035}%
\special{pa 431 275}%
\special{pa 430 275}%
\special{dt 0.035}%
\special{pn 8}%
\special{pa 209 365}%
\special{pa 431 365}%
\special{dt 0.035}%
\special{pa 431 365}%
\special{pa 430 365}%
\special{dt 0.035}%
\end{picture}%
}

\newcommand{\marudj}{
\unitlength 0.1in
\begin{picture}(3.40,2.40)(1.5,-3.6)
\special{pn 16}%
\special{ar 320 320 120 120  0.0000000 6.2831853}%
\special{pn 8}%
\special{pa 320 200}%
\special{pa 320 440}%
\special{dt 0.04}%
\special{pa 320 440}%
\special{pa 320 439}%
\special{dt 0.04}%
\special{pn 4}%
\special{pa 209 275}%
\special{pa 431 275}%
\special{fp}%
\special{pn 8}%
\special{pa 209 365}%
\special{pa 429 365}%
\special{dt 0.035}%
\special{pa 429 365}%
\special{pa 428 365}%
\special{dt 0.035}%
\end{picture}%
}

\newcommand{\marudk}{
\unitlength 0.1in
\begin{picture}(3.40,2.40)(1.5,-3.6)
\special{pn 16}%
\special{ar 320 320 120 120  0.0000000 6.2831853}%
\special{pn 4}%
\special{pa 320 200}%
\special{pa 320 440}%
\special{fp}%
\special{pn 8}%
\special{pa 215 260}%
\special{pa 423 380}%
\special{dt 0.035}%
\special{pa 423 380}%
\special{pa 422 380}%
\special{dt 0.035}%
\special{pn 8}%
\special{pa 423 260}%
\special{pa 216 380}%
\special{dt 0.035}%
\special{pa 216 380}%
\special{pa 217 380}%
\special{dt 0.035}%
\end{picture}%
}

\newcommand{\marudl}{
\unitlength 0.1in
\begin{picture}(3.4,2.40)(1.5,-3.6)
\special{pn 16}%
\special{ar 320 320 120 120  0.0000000 6.2831853}%
\special{pn 4}%
\special{pa 320 200}%
\special{pa 320 440}%
\special{fp}%
\special{pn 4}%
\special{pa 200 320}%
\special{pa 440 320}%
\special{fp}%
\put(5,-3.2){\makebox(0,0){\tiny{2d}}}%
\end{picture}%
}

\newcommand{\marudm}{
\unitlength 0.1in
\begin{picture}(3.4,2.40)(1.5,-3.6)
\special{pn 16}%
\special{ar 320 320 120 120  0.0000000 6.2831853}%
\special{pn 4}%
\special{pa 320 200}%
\special{pa 320 440}%
\special{fp}%
\special{pn 4}%
\special{pa 200 320}%
\special{pa 440 320}%
\special{fp}%
\put(5,-3.2){\makebox(0,0){\tiny{2s}}}%
\end{picture}%
}

\newcommand{\marudn}{
\unitlength 0.1in
\begin{picture}(3.4,2.40)(1.5,-3.6)
\special{pn 16}%
\special{ar 320 320 120 120  0.0000000 6.2831853}%
\special{pn 4}%
\special{pa 200 320}%
\special{pa 440 320}%
\special{fp}%
\put(3.9,-2.95){\makebox(0,0){\tiny{s}}}%
\end{picture}%
}
\newcommand{\marudo}{
\unitlength 0.1in
\begin{picture}(6.9,2.40)(1.5,-3.6)
\special{pn 16}%
\special{ar 320 320 120 120  0.0000000 6.2831853}%
\special{pn 16}%
\special{ar 680 320 120 120  0.0000000 6.2831853}%
\special{pn 8}%
\special{pa 440 320}%
\special{pa 560 320}%
\special{dt 0.035}%
\special{pa 560 320}%
\special{pa 559 320}%
\special{dt 0.045}%
\end{picture}%
}

\newcommand{\marudp}{
\unitlength 0.1in
\begin{picture}(3.4,2.40)(1.5,-3.6)
\special{pn 16}%
\special{ar 320 320 120 120  0.0000000 6.2831853}%
\special{pn 4}%
\special{pa 320 200}%
\special{pa 320 440}%
\special{fp}%
\special{pn 4}%
\special{pa 200 320}%
\special{pa 440 320}%
\special{fp}%
\put(3.9,-2.95){\makebox(0,0){\tiny{2}}}%
\end{picture}%
}

\newcommand{\maruea}{
\unitlength 0.1in
\begin{picture}(6.9,2.40)(1.5,-3.6)
\special{pn 16}%
\special{ar 320 320 120 120  0.0000000 6.2831853}%
\special{pn 16}%
\special{ar 680 320 120 120  0.0000000 6.2831853}%
\special{pn 4}%
\special{pa 425 260}%
\special{pa 576 260}%
\special{fp}%
\special{pn 4}%
\special{pa 425 380}%
\special{pa 576 380}%
\special{fp}%
\special{pn 4}%
\special{pa 320 200}%
\special{pa 320 440}%
\special{fp}%
\end{picture}%
}

\newcommand{\marueb}{
\unitlength 0.1in
\begin{picture}(6.9,2.40)(1.5,-3.6)
\special{pn 16}%
\special{ar 320 320 120 120  0.0000000 6.2831853}%
\special{pn 16}%
\special{ar 680 320 120 120  0.0000000 6.2831853}%
\special{pn 4}%
\special{pa 425 260}%
\special{pa 576 260}%
\special{fp}%
\special{pn 4}%
\special{pa 425 380}%
\special{pa 576 380}%
\special{fp}%
\special{pn 4}%
\special{pa 200 320}%
\special{pa 440 320}%
\special{fp}%
\end{picture}%
}
\newcommand{\maruec}{
\unitlength 0.1in
\begin{picture}(6.9,2.40)(1.5,-3.6)
\special{pn 16}%
\special{ar 320 320 120 120  0.0000000 6.2831853}%
\special{pn 16}%
\special{ar 680 320 120 120  0.0000000 6.2831853}%
\special{pn 4}%
\special{pa 440 320}%
\special{pa 560 320}%
\special{fp}%
\special{pn 4}%
\special{pa 399 230}%
\special{pa 600 230}%
\special{fp}%
\special{pn 4}%
\special{pa 399 410}%
\special{pa 602 410}%
\special{fp}%
\end{picture}%
}
\newcommand{\marued}{
\unitlength 0.1in
\begin{picture}(6.9,2.40)(1.5,-3.6)
\special{pn 16}%
\special{ar 320 320 120 120  0.0000000 6.2831853}%
\special{pn 16}%
\special{ar 680 320 120 120  0.0000000 6.2831853}%
\special{pn 4}%
\special{pa 401 410}%
\special{pa 600 410}%
\special{fp}%
\special{pn 4}%
\special{pa 399 230}%
\special{pa 563 350}%
\special{fp}%
\special{pn 4}%
\special{pa 599 230}%
\special{pa 437 350}%
\special{fp}%
\end{picture}%
}

\newcommand{\maruef}{
\unitlength 0.1in
\begin{picture}(6.90,2.40)(1.50,-3.60)
\special{pn 16}%
\special{ar 320 320 120 120  0.0000000 6.2831853}%
\special{pn 16}%
\special{ar 680 320 120 120  0.0000000 6.2831853}%
\special{pn 4}%
\special{pa 401 410}%
\special{pa 600 410}%
\special{fp}%
\special{pn 4}%
\special{ar 500 230 30 30  0.0000000 6.2831853}%
\special{pn 4}%
\special{pa 400 230}%
\special{pa 470 230}%
\special{fp}%
\special{pn 4}%
\special{pa 530 230}%
\special{pa 602 230}%
\special{fp}%
\end{picture}%
}

\newcommand{\marueg}{
\unitlength 0.1in
\begin{picture}(6.9,2.40)(1.5,-3.6)
\special{pn 16}%
\special{ar 320 320 120 120  0.0000000 6.2831853}%
\special{pn 16}%
\special{ar 680 320 120 120  0.0000000 6.2831853}%
\special{pn 8}%
\special{pa 425 260}%
\special{pa 576 260}%
\special{dt 0.035}%
\special{pa 576 260}%
\special{pa 575 260}%
\special{dt 0.035}%
\special{pn 8}%
\special{pa 425 380}%
\special{pa 576 380}%
\special{dt 0.035}%
\special{pa 576 380}%
\special{pa 575 380}%
\special{dt 0.035}%
\special{pn 4}%
\special{pa 200 320}%
\special{pa 440 320}%
\special{fp}%
\end{picture}%
}

\newcommand{\marueh}{
\unitlength 0.1in
\begin{picture}(6.9,2.4)(1.5,-3.6)
\special{pn 16}%
\special{ar 320 320 120 120  0.0000000 6.2831853}%
\special{pn 16}%
\special{ar 680 320 120 120  0.0000000 6.2831853}%
\special{pn 8}%
\special{pa 425 260}%
\special{pa 576 260}%
\special{dt 0.035}%
\special{pa 576 260}%
\special{pa 575 260}%
\special{dt 0.035}%
\special{pn 4}%
\special{pa 425 380}%
\special{pa 576 380}%
\special{fp}%
\special{pn 8}%
\special{pa 200 320}%
\special{pa 440 320}%
\special{dt 0.035}%
\special{pa 440 320}%
\special{pa 439 320}%
\special{dt 0.035}%
\end{picture}%
}

\newcommand{\maruei}{
\unitlength 0.1in
\begin{picture}(6.9,2.40)(1.5,-3.6)
\special{pn 16}%
\special{ar 320 320 120 120  0.0000000 6.2831853}%
\special{pn 16}%
\special{ar 680 320 120 120  0.0000000 6.2831853}%
\special{pn 4}%
\special{pa 401 410}%
\special{pa 600 410}%
\special{fp}%
\special{pn 8}%
\special{pa 399 230}%
\special{pa 563 350}%
\special{dt 0.03}%
\special{pa 563 350}%
\special{pa 562 350}%
\special{dt 0.03}%
\special{pn 8}%
\special{pa 599 230}%
\special{pa 437 350}%
\special{dt 0.03}%
\special{pa 437 350}%
\special{pa 438 350}%
\special{dt 0.03}%
\end{picture}%
}

\newcommand{\maruej}{
\unitlength 0.1in
\begin{picture}(6.9,2.4)(1.5,-3.6)
\special{pn 16}%
\special{ar 320 320 120 120  0.0000000 6.2831853}%
\special{pn 16}%
\special{ar 680 320 120 120  0.0000000 6.2831853}%
\special{pn 4}%
\special{pa 425 260}%
\special{pa 576 260}%
\special{fp}%
\special{pn 4}%
\special{pa 425 380}%
\special{pa 576 380}%
\special{fp}%
\put(5.000,-2.3000){\makebox(0,0){\tiny{2s}}}%
\end{picture}%
}

\newcommand{\maruek}{
\unitlength 0.1in
\begin{picture}(6.90,2.4)(1.5,-3.6)
\special{pn 16}%
\special{ar 320 320 120 120  0.0000000 6.2831853}%
\special{pn 16}%
\special{ar 680 320 120 120  0.0000000 6.2831853}%
\special{pn 4}%
\special{pa 440 320}%
\special{pa 560 320}%
\special{fp}%
\put(4.8000,-3.15000){\makebox(0,0)[lb]{\tiny{s}}}%
\end{picture}%
}

\newcommand{\maruel}{
\unitlength 0.1in
\begin{picture}(6.90,2.4)(1.5,-3.6)
\special{pn 16}%
\special{ar 320 320 120 120  0.0000000 6.2831853}%
\special{pn 16}%
\special{ar 680 320 120 120  0.0000000 6.2831853}%
\special{pn 4}%
\special{pa 440 320}%
\special{pa 560 320}%
\special{fp}%
\put(4.8000,-3.15000){\makebox(0,0)[lb]{\tiny{d}}}%
\end{picture}%
}
 
\newcommand{\marufa}{
\unitlength 0.1in
\begin{picture}(3.40,2.40)(1.50,-3.60)
\special{pn 16}%
\special{ar 320 320 120 120  0.0000000 6.2831853}%
\special{pn 4}%
\special{ar 320 320 24 24  0.0000000 6.2831853}%
\special{pn 4}%
\special{ar 250 320 24 24  0.0000000 6.2831853}%
\special{pn 4}%
\special{ar 390 320 24 24  0.0000000 6.2831853}%
\special{pn 4}%
\special{pa 200 320}%
\special{pa 222 320}%
\special{fp}%
\special{pn 4}%
\special{pa 275 320}%
\special{pa 296 320}%
\special{fp}%
\special{pn 4}%
\special{pa 344 320}%
\special{pa 365 320}%
\special{fp}%
\special{pn 4}%
\special{pa 415 320}%
\special{pa 440 320}%
\special{fp}%
\end{picture}%
}

\newcommand{\marufb}{
\unitlength 0.1in
\begin{picture}(3.40,2.40)(1.500,-3.60)
\special{pn 16}%
\special{ar 320 320 120 120  0.0000000 6.2831853}%
\special{pn 4}%
\special{pa 281 281}%
\special{pa 359 281}%
\special{pa 359 359}%
\special{pa 281 359}%
\special{pa 281 281}%
\special{fp}%
\special{pn 4}%
\special{pa 234 234}%
\special{pa 280 280}%
\special{fp}%
\special{pn 4}%
\special{pa 406 234}%
\special{pa 360 280}%
\special{fp}%
\special{pn 4}%
\special{pa 234 406}%
\special{pa 280 360}%
\special{fp}%
\special{pn 4}%
\special{pa 360 360}%
\special{pa 406 406}%
\special{fp}%
\end{picture}%
}

\newcommand{\marufc}{
\unitlength 0.1in
\begin{picture}(4.0,2.4)(3.60,-3.6)
\special{pn 13}%
\special{ar 660 403 60 60  0.0000000 6.2831853}%
\special{pn 13}%
\special{ar 460 403 60 60  1.5541312 1.5707963}%
\special{pn 13}%
\special{ar 460 403 60 60  0.0000000 6.2831853}%
\special{pn 13}%
\special{ar 560 230 60 60  0.0000000 6.2831853}%
\special{pn 4}%
\special{pa 530 282}%
\special{pa 490 352}%
\special{fp}%
\special{pn 4}%
\special{pa 590 282}%
\special{pa 630 354}%
\special{fp}%
\special{pn 4}%
\special{pa 520 404}%
\special{pa 600 404}%
\special{fp}%
\end{picture}%
}

\newcommand{\marufd}{
\unitlength 0.1in
\begin{picture}(6.90,2.40)(1.50,-3.60)
\special{pn 16}%
\special{ar 320 320 120 120  0.0000000 6.2831853}%
\special{pn 16}%
\special{ar 680 320 120 120  0.0000000 6.2831853}%
\special{pn 4}%
\special{pa 437 290}%
\special{pa 564 290}%
\special{fp}%
\special{pn 4}%
\special{pa 435 350}%
\special{pa 564 350}%
\special{fp}%
\special{pn 4}%
\special{pa 399 230}%
\special{pa 600 230}%
\special{fp}%
\special{pn 4}%
\special{pa 399 410}%
\special{pa 602 410}%
\special{fp}%
\end{picture}%
}

\newcommand{\marufe}{
\unitlength 0.1in
\begin{picture}(6.90,2.40)(1.50,-3.60)
\special{pn 16}%
\special{ar 320 320 120 120  0.0000000 6.2831853}%
\special{pn 16}%
\special{ar 680 320 120 120  0.0000000 6.2831853}%
\special{pn 4}%
\special{pa 440 320}%
\special{pa 560 320}%
\special{fp}%
\special{pa 560 320}%
\special{pa 560 320}%
\special{fp}%
\special{pn 4}%
\special{pa 398 410}%
\special{pa 600 410}%
\special{fp}%
\special{pn 4}%
\special{ar 500 230 24 24  0.0000000 6.2831853}%
\special{pn 4}%
\special{pa 400 230}%
\special{pa 476 230}%
\special{fp}%
\special{pn 4}%
\special{pa 524 230}%
\special{pa 600 230}%
\special{fp}%
\end{picture}%
}

\newcommand{\maruff}{
\unitlength 0.1in
\begin{picture}(6.90,2.40)(1.50,-3.60)
\special{pn 16}%
\special{ar 320 320 120 120  0.0000000 6.2831853}%
\special{pn 16}%
\special{ar 680 320 120 120  0.0000000 6.2831853}%
\special{pn 4}%
\special{pa 398 410}%
\special{pa 602 410}%
\special{fp}%
\special{pn 4}%
\special{ar 462 230 24 24  0.0000000 6.2831853}%
\special{pn 4}%
\special{ar 537 230 24 24  0.0000000 6.2831853}%
\special{pn 4}%
\special{pa 400 230}%
\special{pa 438 230}%
\special{fp}%
\special{pn 4}%
\special{pa 486 230}%
\special{pa 514 230}%
\special{fp}%
\special{pn 4}%
\special{pa 562 230}%
\special{pa 600 230}%
\special{fp}%
\end{picture}%
}

\newcommand{\marufg}{
\unitlength 0.1in
\begin{picture}(4.0,2.4)(3.60,-3.6)
\special{pn 13}%
\special{sh 0}%
\special{ar 660 403 60 60  0.0000000 6.2831853}%
\special{pn 13}%
\special{ar 460 403 60 60  1.5541312 1.5707963}%
\special{pn 13}%
\special{sh 0}%
\special{ar 460 403 60 60  0.0000000 6.2831853}%
\special{pn 13}%
\special{sh 0}%
\special{ar 560 230 60 60  0.0000000 6.2831853}%
\special{pn 4}%
\special{pa 540 286}%
\special{pa 500 358}%
\special{fp}%
\special{pn 4}%
\special{pa 580 286}%
\special{pa 620 358}%
\special{fp}%
\special{pn 4}%
\special{pa 620 235}%
\special{pa 680 347}%
\special{fp}%
\special{pn 4}%
\special{pa 500 233}%
\special{pa 440 345}%
\special{fp}%
\end{picture}%
}

\newcommand{\marufh}{
\unitlength 0.1in
\begin{picture}(4.0,2.4)(3.60,-3.6)
\special{pn 4}%
\special{pa 560 230}%
\special{pa 460 403}%
\special{fp}%
\special{pn 4}%
\special{pa 560 230}%
\special{pa 660 403}%
\special{fp}%
\special{pn 13}%
\special{sh 0}%
\special{ar 660 403 60 60  0.0000000 6.2831853}%
\special{pn 13}%
\special{ar 460 403 60 60  1.5541312 1.5707963}%
\special{pn 13}%
\special{sh 0}%
\special{ar 460 403 60 60  0.0000000 6.2831853}%
\special{pn 13}%
\special{sh 0}%
\special{ar 560 230 60 60  0.0000000 6.2831853}%
\special{pn 4}%
\special{pa 516 383}%
\special{pa 604 383}%
\special{fp}%
\special{pn 4}%
\special{pa 516 423}%
\special{pa 603 423}%
\special{fp}%
\end{picture}%
}

\newcommand{\marufi}{
\unitlength 0.1in
\begin{picture}(4.0,2.4)(3.60,-3.6)
\special{pn 13}%
\special{ar 460 403 60 60  0.0000000 6.2831853}%
\special{pn 13}%
\special{ar 660 403 60 60  0.0000000 6.2831853}%
\special{pn 13}%
\special{ar 560 230 60 60  0.0000000 6.2831853}%
\special{pn 4}%
\special{ar 560 403 16 16  0.0000000 6.2831853}%
\special{pn 4}%
\special{pa 530 282}%
\special{pa 490 352}%
\special{fp}%
\special{pn 4}%
\special{pa 592 284}%
\special{pa 632 354}%
\special{fp}%
\special{pa 632 354}%
\special{pa 632 354}%
\special{fp}%
\special{pn 4}%
\special{pa 518 404}%
\special{pa 544 404}%
\special{fp}%
\special{pn 4}%
\special{pa 576 404}%
\special{pa 600 404}%
\special{fp}%
\end{picture}%
}

\newcommand{\marufj}{
\unitlength 0.1in
\begin{picture}(4.20,3.0)(3.50,-6.0)
\special{pn 13}%
\special{ar 460 660 60 60  0.0000000 6.2831853}%
\special{pn 13}%
\special{ar 660 660 60 60  0.0000000 6.2831853}%
\special{pn 13}%
\special{ar 460 460 60 60  0.0000000 6.2831853}%
\special{pn 13}%
\special{ar 660 460 60 60  0.0000000 6.2831853}%
\special{pn 4}%
\special{pa 520 460}%
\special{pa 600 460}%
\special{fp}%
\special{pn 4}%
\special{pa 460 520}%
\special{pa 460 600}%
\special{fp}%
\special{pn 4}%
\special{pa 660 520}%
\special{pa 660 600}%
\special{fp}%
\special{pn 4}%
\special{pa 520 660}%
\special{pa 600 660}%
\special{fp}%
\end{picture}%
}

\newcommand{\marufk}{
\unitlength 0.1in
\begin{picture}(4.20,3.0)(3.50,-6.0)
\special{pn 13}%
\special{ar 460 660 60 60  0.0000000 6.2831853}%
\special{pn 13}%
\special{ar 660 660 60 60  0.0000000 6.2831853}%
\special{pn 13}%
\special{ar 460 460 60 60  0.0000000 6.2831853}%
\special{pn 13}%
\special{ar 660 460 60 60  0.0000000 6.2831853}%
\special{pn 4}%
\special{pa 506 420}%
\special{pa 615 420}%
\special{fp}%
\special{pn 4}%
\special{pa 505 500}%
\special{pa 616 500}%
\special{fp}%
\special{pn 4}%
\special{pa 505 620}%
\special{pa 616 620}%
\special{fp}%
\special{pn 4}%
\special{pa 505 700}%
\special{pa 616 700}%
\special{fp}%
\end{picture}%
} 
\newcommand{\maruga}{
\unitlength 0.1in
\begin{picture}(3.40,2.40)(1.50,-3.60)
\special{pn 16}%
\special{ar 320 320 120 120  0.0000000 6.2831853}%
\special{pn 4}%
\special{pa 320 200}%
\special{pa 320 440}%
\special{fp}%
\special{pn 4}%
\special{pa 200 320}%
\special{pa 440 320}%
\special{fp}%
\special{pn 4}%
\special{pa 216 260}%
\special{pa 423 260}%
\special{fp}%
\special{pn 4}%
\special{pa 216 380}%
\special{pa 423 380}%
\special{fp}%
\end{picture}%
}

\newcommand{\marugb}{
\unitlength 0.1in
\begin{picture}(3.40,2.40)(1.50,-3.60)
\special{pn 16}%
\special{ar 320 320 120 120  0.0000000 6.2831853}%
\special{pn 4}%
\special{pa 380 216}%
\special{pa 230 396}%
\special{fp}%
\special{pn 4}%
\special{pa 260 218}%
\special{pa 410 396}%
\special{fp}%
\special{pn 4}%
\special{ar 440 410 75 75  2.8577985 4.7123890}%
\special{pn 4}%
\special{ar 200 410 75 75  4.7257215 6.2831853}%
\special{ar 200 410 75 75  0.0000000 0.2801080}%
\end{picture}%
}

\newcommand{\marugc}{
\unitlength 0.1in
\begin{picture}(3.40,2.40)(1.50,-3.60)
\special{pn 16}%
\special{ar 320 320 120 120  0.0000000 6.2831853}%
\special{pn 4}%
\special{pa 320 200}%
\special{pa 320 440}%
\special{fp}%
\special{pn 4}%
\special{pa 200 320}%
\special{pa 399 230}%
\special{fp}%
\special{pn 4}%
\special{pa 440 320}%
\special{pa 240 230}%
\special{fp}%
\special{pn 4}%
\special{pa 216 380}%
\special{pa 423 380}%
\special{fp}%
\end{picture}%
}

\newcommand{\marugd}{
\unitlength 0.1in
\begin{picture}(3.40,2.40)(1.50,-3.60)
\special{pn 16}%
\special{ar 320 320 120 120  0.0000000 6.2831853}%
\special{pn 4}%
\special{pa 275 207}%
\special{pa 275 431}%
\special{fp}%
\special{pn 4}%
\special{pa 365 209}%
\special{pa 365 431}%
\special{fp}%
\special{pn 4}%
\special{pa 209 275}%
\special{pa 431 275}%
\special{fp}%
\special{pn 4}%
\special{pa 209 365}%
\special{pa 431 365}%
\special{fp}%
\end{picture}%
}

\newcommand{\maruge}{
\unitlength 0.1in
\begin{picture}(3.40,2.40)(1.50,-3.60)
\special{pn 16}%
\special{ar 320 320 120 120  0.0000000 6.2831853}%
\special{pn 4}%
\special{pa 209 275}%
\special{pa 431 275}%
\special{fp}%
\special{pn 4}%
\special{pa 207 365}%
\special{pa 431 365}%
\special{fp}%
\special{pn 4}%
\special{pa 365 209}%
\special{pa 275 431}%
\special{fp}%
\special{pn 4}%
\special{pa 275 207}%
\special{pa 365 429}%
\special{fp}%
\end{picture}%
}

\newcommand{\marugf}{
\unitlength 0.1in
\begin{picture}(3.40,2.40)(1.50,-3.60)
\special{pn 16}%
\special{ar 320 320 120 120  0.0000000 6.2831853}%
\special{pn 4}%
\special{pa 200 320}%
\special{pa 440 320}%
\special{fp}%
\special{pn 4}%
\special{pa 320 200}%
\special{pa 320 440}%
\special{fp}%
\special{pn 4}%
\special{pa 234 234}%
\special{pa 405 405}%
\special{fp}%
\special{pn 4}%
\special{pa 405 234}%
\special{pa 234 405}%
\special{fp}%
\end{picture}%
}

\newcommand{\marugg}{
\unitlength 0.1in
\begin{picture}(3.40,2.4)(1.5,-3.60)
\special{pn 16}%
\special{ar 320 320 120 120  0.0000000 6.2831853}%
\special{pn 4}%
\special{pa 320 200}%
\special{pa 320 440}%
\special{fp}%
\special{pn 4}%
\special{pa 209 275}%
\special{pa 431 275}%
\special{fp}%
\put(4.700,-2.7600){\makebox(0,0){\tiny 2}}%
\special{pn 4}%
\special{pa 209 365}%
\special{pa 431 365}%
\special{fp}%
\end{picture}%
}

\newcommand{\marugh}{
\unitlength 0.1in
\begin{picture}(3.40,2.40)(1.50,-2)
\special{pn 16}%
\special{ar 320 163 120 120  0.0000000 6.2831853}%
\special{pn 4}%
\special{pa 320 43}%
\special{pa 320 283}%
\special{fp}%
\special{pn 4}%
\special{pa 209 118}%
\special{pa 431 118}%
\special{fp}%
\special{pn 4}%
\special{pa 209 208}%
\special{pa 431 208}%
\special{fp}%
\put(3.2000,-0.050){\makebox(0,0){\tiny 2}}%
\end{picture}%
}

\newcommand{\marugi}{
\unitlength 0.1in
\begin{picture}(3.40,2.40)(1.50,-3.6)
\special{pn 16}%
\special{ar 320 320 120 120  0.0000000 6.2831853}%
\special{pn 4}%
\special{pa 320 200}%
\special{pa 320 440}%
\special{fp}%
\special{pn 4}%
\special{pa 216 260}%
\special{pa 423 380}%
\special{fp}%
\special{pn 4}%
\special{pa 423 260}%
\special{pa 216 380}%
\special{fp}%
\put(3.2000,-1.620){\makebox(0,0){\tiny 2}}%
\end{picture}%
}   
\newcommand{\maruha}{
\unitlength 0.1in
\begin{picture}(6.90,2.40)(1.50,-3.60)
\special{pn 16}%
\special{ar 320 320 120 120  0.0000000 6.2831853}%
\special{pn 16}%
\special{ar 680 320 120 120  0.0000000 6.2831853}%
\special{pn 4}%
\special{pa 200 320}%
\special{pa 440 320}%
\special{fp}%
\special{pn 4}%
\special{pa 320 200}%
\special{pa 320 440}%
\special{fp}%
\special{pa 320 440}%
\special{pa 320 440}%
\special{fp}%
\special{pn 4}%
\special{pa 413 245}%
\special{pa 585 245}%
\special{fp}%
\special{pn 4}%
\special{pa 413 395}%
\special{pa 585 395}%
\special{fp}%
\end{picture}%
}

\newcommand{\maruhb}{
\unitlength 0.1in
\begin{picture}(6.90,2.40)(1.50,-3.60)
\special{pn 16}%
\special{ar 320 320 120 120  0.0000000 6.2831853}%
\special{pn 16}%
\special{ar 680 320 120 120  0.0000000 6.2831853}%
\special{pn 4}%
\special{pa 401 230}%
\special{pa 600 230}%
\special{fp}%
\special{pn 4}%
\special{pa 399 410}%
\special{pa 602 410}%
\special{fp}%
\special{pn 4}%
\special{pa 204 290}%
\special{pa 437 290}%
\special{fp}%
\special{pn 4}%
\special{pa 204 350}%
\special{pa 435 350}%
\special{fp}%
\end{picture}%
}

\newcommand{\maruhc}{
\unitlength 0.1in
\begin{picture}(6.90,2.40)(1.50,-3.60)
\special{pn 16}%
\special{ar 320 320 120 120  0.0000000 6.2831853}%
\special{pn 16}%
\special{ar 680 320 120 120  0.0000000 6.2831853}%
\special{pn 4}%
\special{pa 401 230}%
\special{pa 600 230}%
\special{fp}%
\special{pn 4}%
\special{pa 399 410}%
\special{pa 602 410}%
\special{fp}%
\special{pn 4}%
\special{pa 431 275}%
\special{pa 209 365}%
\special{fp}%
\special{pn 4}%
\special{pa 209 275}%
\special{pa 429 365}%
\special{fp}%
\end{picture}%
}

\newcommand{\maruhd}{
\unitlength 0.1in
\begin{picture}(6.90,2.40)(1.50,-3.60)
\special{pn 16}%
\special{ar 320 320 120 120  0.0000000 6.2831853}%
\special{pn 16}%
\special{ar 680 320 120 120  0.0000000 6.2831853}%
\special{pn 4}%
\special{pa 401 230}%
\special{pa 600 230}%
\special{fp}%
\special{pn 4}%
\special{pa 437 290}%
\special{pa 599 410}%
\special{fp}%
\special{pn 4}%
\special{pa 564 290}%
\special{pa 399 410}%
\special{fp}%
\special{pn 4}%
\special{pa 201 335}%
\special{pa 438 335}%
\special{fp}%
\end{picture}%
}

\newcommand{\maruhe}{
\unitlength 0.1in
\begin{picture}(6.90,2.40)(1.50,-3.60)
\special{pn 16}%
\special{ar 320 320 120 120  0.0000000 6.2831853}%
\special{pn 16}%
\special{ar 680 320 120 120  0.0000000 6.2831853}%
\special{pn 4}%
\special{pa 398 230}%
\special{pa 564 290}%
\special{fp}%
\special{pn 4}%
\special{pa 600 230}%
\special{pa 437 290}%
\special{fp}%
\special{pn 4}%
\special{pa 437 350}%
\special{pa 602 410}%
\special{fp}%
\special{pn 4}%
\special{pa 563 350}%
\special{pa 401 410}%
\special{fp}%
\end{picture}%
}

\newcommand{\maruhf}{
\unitlength 0.1in
\begin{picture}(6.90,2.40)(1.50,-3.60)
\special{pn 16}%
\special{ar 320 320 120 120  0.0000000 6.2831853}%
\special{pn 16}%
\special{ar 680 320 120 120  0.0000000 6.2831853}%
\special{pn 4}%
\special{pa 401 230}%
\special{pa 600 230}%
\special{fp}%
\special{pn 4}%
\special{pa 399 410}%
\special{pa 602 410}%
\special{fp}%
\special{pn 4}%
\special{pa 200 320}%
\special{pa 440 320}%
\special{fp}%
\special{pn 4}%
\special{pa 560 320}%
\special{pa 800 320}%
\special{fp}%
\end{picture}%
}

\newcommand{\maruhg}{
\unitlength 0.1in
\begin{picture}(6.90,2.40)(1.50,-3.60)
\special{pn 16}%
\special{ar 320 320 120 120  0.0000000 6.2831853}%
\special{pn 16}%
\special{ar 680 320 120 120  0.0000000 6.2831853}%
\special{pn 4}%
\special{pa 401 230}%
\special{pa 600 230}%
\special{fp}%
\special{pn 4}%
\special{pa 399 410}%
\special{pa 602 410}%
\special{fp}%
\special{pn 4}%
\special{pa 200 320}%
\special{pa 440 320}%
\special{fp}%
\put(3.2000,-2.90){\makebox(0,0){\tiny 2}}%
\end{picture}%
}

\newcommand{\maruhh}{
\unitlength 0.1in
\begin{picture}(6.90,2.4)(1.50,-3.60)
\special{pn 16}%
\special{ar 320 320 120 120  0.0000000 6.2831853}%
\special{pn 16}%
\special{ar 680 320 120 120  0.0000000 6.2831853}%
\special{pn 4}%
\special{pa 401 230}%
\special{pa 600 230}%
\special{fp}%
\special{pn 4}%
\special{pa 437 290}%
\special{pa 599 410}%
\special{fp}%
\special{pn 4}%
\special{pa 564 290}%
\special{pa 399 410}%
\special{fp}%
\put(5.0000,-2.0000){\makebox(0,0){\tiny 2}}%
\end{picture}%
}

\newcommand{\maruhi}{
\unitlength 0.1in
\begin{picture}(6.90,2.40)(1.50,-3.60)
\special{pn 16}%
\special{ar 320 320 120 120  0.0000000 6.2831853}%
\special{pn 16}%
\special{ar 680 320 120 120  0.0000000 6.2831853}%
\special{pn 4}%
\special{pa 440 320}%
\special{pa 560 320}%
\special{fp}%
\special{pn 4}%
\special{pa 320 200}%
\special{pa 320 440}%
\special{fp}%
\end{picture}%
}

\newcommand{\maruhj}{
\unitlength 0.1in
\begin{picture}(4.0,2.4)(3.60,-3.6)
\special{pn 13}%
\special{ar 560 230 60 60  0.0000000 6.2831853}%
\special{pn 13}%
\special{ar 460 403 60 60  0.0000000 6.2831853}%
\special{pn 13}%
\special{ar 660 403 60 60  0.0000000 6.2831853}%
\special{pn 4}%
\special{pa 500 234}%
\special{pa 433 350}%
\special{fp}%
\special{pn 4}%
\special{pa 586 285}%
\special{pa 520 398}%
\special{fp}%
\special{pn 4}%
\special{pa 535 286}%
\special{pa 600 398}%
\special{fp}%
\special{pn 4}%
\special{pa 621 232}%
\special{pa 688 350}%
\special{fp}%
\end{picture}%
}

\newcommand{\maruhk}{
\unitlength 0.1in
\begin{picture}(4.0,2.4)(3.60,-3.6)
\special{pn 13}%
\special{ar 560 230 60 60  0.0000000 6.2831853}%
\special{pn 13}%
\special{ar 460 403 60 60  0.0000000 6.2831853}%
\special{pn 13}%
\special{ar 660 403 60 60  0.0000000 6.2831853}%
\special{pn 4}%
\special{pa 520 403}%
\special{pa 600 403}%
\special{fp}%
\special{pn 4}%
\special{pa 560 170}%
\special{pa 560 290}%
\special{fp}%
\special{pn 4}%
\special{pa 508 259}%
\special{pa 459 342}%
\special{fp}%
\special{pn 4}%
\special{pa 612 261}%
\special{pa 660 344}%
\special{fp}%
\end{picture}%
}

\newcommand{\maruhl}{
\unitlength 0.1in
\begin{picture}(4.0,2.4)(3.60,-3.6)
\special{pn 13}%
\special{ar 560 230 60 60  0.0000000 6.2831853}%
\special{pn 13}%
\special{ar 460 403 60 60  0.0000000 6.2831853}%
\special{pn 13}%
\special{ar 660 403 60 60  0.0000000 6.2831853}%
\special{pn 4}%
\special{pa 512 373}%
\special{pa 609 433}%
\special{fp}%
\special{pn 4}%
\special{pa 607 373}%
\special{pa 513 433}%
\special{fp}%
\special{pn 4}%
\special{pa 515 269}%
\special{pa 472 343}%
\special{fp}%
\special{pn 4}%
\special{pa 606 270}%
\special{pa 650 343}%
\special{fp}%
\end{picture}%
}

\newcommand{\maruhm}{
\unitlength 0.1in
\begin{picture}(5.12,1.28)(2.00,-3)
\special{pn 13}%
\special{ar 264 264 64 64  0.0000000 6.2831853}%
\special{pn 13}%
\special{ar 456 264 64 64  0.0000000 6.2831853}%
\special{pn 13}%
\special{ar 648 264 64 64  0.0000000 6.2831853}%
\special{pn 4}%
\special{pa 328 264}%
\special{pa 392 264}%
\special{fp}%
\special{pn 4}%
\special{pa 520 264}%
\special{pa 584 264}%
\special{fp}%
\end{picture}%
}
  
\newcommand{\maruia}{
\unitlength 0.1in
\begin{picture}(3.40,2.40)(1.50,-3.60)
\special{pn 16}%
\special{ar 320 320 120 120  0.0000000 6.2831853}%
\special{pn 4}%
\special{pa 320 200}%
\special{pa 320 440}%
\special{fp}%
\special{pn 8}%
\special{pa 200 320}%
\special{pa 440 320}%
\special{dt 0.035}%
\special{pa 440 320}%
\special{pa 439 320}%
\special{dt 0.035}%
\special{pn 8}%
\special{pa 260 215}%
\special{pa 260 423}%
\special{dt 0.035}%
\special{pa 260 423}%
\special{pa 260 422}%
\special{dt 0.035}%
\end{picture}%
}

\newcommand{\maruib}{
\unitlength 0.1in
\begin{picture}(3.40,2.40)(1.50,-3.60)
\special{pn 16}%
\special{ar 320 320 120 120  0.0000000 6.2831853}%
\special{pn 4}%
\special{pa 320 200}%
\special{pa 320 440}%
\special{fp}%
\special{pn 8}%
\special{pa 275 209}%
\special{pa 264 258}%
\special{fp}%
\special{pn 8}%
\special{pa 230 242}%
\special{pa 239 279}%
\special{fp}%
\special{pn 8}%
\special{pa 243 296}%
\special{pa 254 339}%
\special{fp}%
\special{pn 8}%
\special{pa 258 360}%
\special{pa 269 402}%
\special{fp}%
\special{pn 8}%
\special{pa 258 278}%
\special{pa 240 350}%
\special{fp}%
\special{pn 8}%
\special{pa 236 371}%
\special{pa 230 398}%
\special{fp}%
\special{pn 8}%
\special{pa 273 417}%
\special{pa 276 431}%
\special{fp}%
\end{picture}%
}

\newcommand{\maruic}{
\unitlength 0.1in
\begin{picture}(3.40,2.40)(1.50,-3.60)
\special{pn 16}%
\special{ar 320 320 120 120  0.0000000 6.2831853}%
\special{pn 4}%
\special{pa 320 200}%
\special{pa 320 440}%
\special{fp}%
\special{pn 8}%
\special{pa 200 320}%
\special{pa 440 320}%
\special{dt 0.035}%
\special{pa 440 320}%
\special{pa 439 320}%
\special{dt 0.035}%
\end{picture}%
}

\newcommand{\maruid}{
\unitlength 0.1in
\begin{picture}(6.90,2.40)(1.50,-3.60)
\special{pn 16}%
\special{ar 320 320 120 120  0.0000000 6.2831853}%
\special{pn 16}%
\special{ar 680 320 120 120  0.0000000 6.2831853}%
\special{pn 4}%
\special{pa 440 320}%
\special{pa 560 320}%
\special{fp}%
\special{pn 8}%
\special{pa 215 260}%
\special{pa 575 260}%
\special{dt 0.035}%
\special{pa 575 260}%
\special{pa 574 260}%
\special{dt 0.035}%
\special{pn 8}%
\special{pa 320 200}%
\special{pa 320 440}%
\special{dt 0.035}%
\special{pa 320 440}%
\special{pa 320 439}%
\special{dt 0.035}%
\end{picture}%
}

\newcommand{\maruie}{
\unitlength 0.1in
\begin{picture}(6.90,2.40)(1.50,-3.60)
\special{pn 16}%
\special{ar 320 320 120 120  0.0000000 6.2831853}%
\special{pn 16}%
\special{ar 680 320 120 120  0.0000000 6.2831853}%
\special{pn 4}%
\special{pa 440 320}%
\special{pa 560 320}%
\special{fp}%
\special{pn 8}%
\special{pa 599 230}%
\special{pa 566 260}%
\special{fp}%
\special{pn 8}%
\special{pa 548 275}%
\special{pa 530 293}%
\special{fp}%
\special{pn 8}%
\special{pa 515 305}%
\special{pa 485 335}%
\special{fp}%
\special{pn 8}%
\special{pa 468 348}%
\special{pa 447 368}%
\special{fp}%
\special{pn 8}%
\special{pa 428 386}%
\special{pa 402 410}%
\special{fp}%
\special{pn 8}%
\special{pa 401 230}%
\special{pa 425 251}%
\special{fp}%
\special{pn 8}%
\special{pa 437 263}%
\special{pa 462 285}%
\special{fp}%
\special{pn 8}%
\special{pa 482 302}%
\special{pa 519 338}%
\special{fp}%
\special{pn 8}%
\special{pa 539 353}%
\special{pa 555 368}%
\special{fp}%
\special{pn 8}%
\special{pa 566 377}%
\special{pa 600 410}%
\special{fp}%
\end{picture}%
}

\newcommand{\maruja}{
\unitlength 0.1in
\begin{picture}(3.40,2.4)(1.5,-3.60)
\special{pn 16}%
\special{ar 320 320 120 120  0.0000000 6.2831853}%
\special{pn 4}%
\special{pa 320 200}%
\special{pa 320 440}%
\special{fp}%
\special{pn 4}%
\special{pa 209 275}%
\special{pa 431 275}%
\special{fp}%
\put(4.900,-2.7600){\makebox(0,0){\tiny 2s}}%
\special{pn 4}%
\special{pa 209 365}%
\special{pa 431 365}%
\special{fp}%
\end{picture}%
}

\newcommand{\marujb}{
\unitlength 0.1in
\begin{picture}(3.40,2.40)(1.50,-2)
\special{pn 16}%
\special{ar 320 163 120 120  0.0000000 6.2831853}%
\special{pn 4}%
\special{pa 320 43}%
\special{pa 320 283}%
\special{fp}%
\special{pn 4}%
\special{pa 209 118}%
\special{pa 431 118}%
\special{fp}%
\special{pn 4}%
\special{pa 209 208}%
\special{pa 431 208}%
\special{fp}%
\put(3.2000,-0.050){\makebox(0,0){\tiny 2s}}%
\end{picture}%
}

\newcommand{\marujc}{
\unitlength 0.1in
\begin{picture}(3.40,2.40)(1.50,-3.6)
\special{pn 16}%
\special{ar 320 320 120 120  0.0000000 6.2831853}%
\special{pn 4}%
\special{pa 320 200}%
\special{pa 320 440}%
\special{fp}%
\special{pn 4}%
\special{pa 216 260}%
\special{pa 423 380}%
\special{fp}%
\special{pn 4}%
\special{pa 423 260}%
\special{pa 216 380}%
\special{fp}%
\put(3.2000,-1.620){\makebox(0,0){\tiny 2s}}%
\end{picture}%
}

\newcommand{\marujd}{
\unitlength 0.1in
\begin{picture}(3.4,2.40)(1.5,-3.6)
\special{pn 16}%
\special{ar 320 320 120 120  0.0000000 6.2831853}%
\special{pn 4}%
\special{pa 320 200}%
\special{pa 320 440}%
\special{fp}%
\special{pn 4}%
\special{pa 200 320}%
\special{pa 440 320}%
\special{fp}%
\put(3.9,-2.95){\makebox(0,0){\tiny{s}}}%
\end{picture}%
}

\newcommand{\maruje}{
\unitlength 0.1in
\begin{picture}(3.4,2.40)(1.5,-3.6)
\special{pn 16}%
\special{ar 320 320 120 120  0.0000000 6.2831853}%
\special{pn 4}%
\special{pa 320 200}%
\special{pa 320 440}%
\special{fp}%
\special{pn 4}%
\special{pa 200 320}%
\special{pa 440 320}%
\special{fp}%
\put(3.9,-2.9){\makebox(0,0){\tiny{d}}}%
\end{picture}%
}

\newcommand{\marujf}{
\unitlength 0.1in
\begin{picture}(6.90,2.40)(1.50,-3.60)
\special{pn 16}%
\special{ar 320 320 120 120  0.0000000 6.2831853}%
\special{pn 16}%
\special{ar 680 320 120 120  0.0000000 6.2831853}%
\special{pn 4}%
\special{pa 401 230}%
\special{pa 600 230}%
\special{fp}%
\special{pn 4}%
\special{pa 399 410}%
\special{pa 602 410}%
\special{fp}%
\special{pn 4}%
\special{pa 200 320}%
\special{pa 440 320}%
\special{fp}%
\put(3.2000,-2.90){\makebox(0,0){\tiny 2s}}%
\end{picture}%
}

\newcommand{\marujg}{
\unitlength 0.1in
\begin{picture}(6.90,2.4)(1.50,-3.60)
\special{pn 16}%
\special{ar 320 320 120 120  0.0000000 6.2831853}%
\special{pn 16}%
\special{ar 680 320 120 120  0.0000000 6.2831853}%
\special{pn 4}%
\special{pa 401 230}%
\special{pa 600 230}%
\special{fp}%
\special{pn 4}%
\special{pa 437 290}%
\special{pa 599 410}%
\special{fp}%
\special{pn 4}%
\special{pa 564 290}%
\special{pa 399 410}%
\special{fp}%
\put(5.0000,-2.0000){\makebox(0,0){\tiny 2d}}%
\end{picture}%
}

\newcommand{\marujh}{
\unitlength 0.1in
\begin{picture}(6.9,2.4)(1.5,-3.6)
\special{pn 16}%
\special{ar 320 320 120 120  0.0000000 6.2831853}%
\special{pn 16}%
\special{ar 680 320 120 120  0.0000000 6.2831853}%
\special{pn 4}%
\special{pa 425 260}%
\special{pa 576 260}%
\special{fp}%
\special{pn 4}%
\special{pa 425 380}%
\special{pa 576 380}%
\special{fp}%
\put(5.000,-2.3000){\makebox(0,0){\tiny{s}}}%
\end{picture}%
}

\newcommand{\maruji}{
\unitlength 0.1in
\begin{picture}(6.90,2.40)(1.50,-3.60)
\special{pn 16}%
\special{ar 320 320 120 120  0.0000000 6.2831853}%
\special{pn 16}%
\special{ar 680 320 120 120  0.0000000 6.2831853}%
\special{pn 4}%
\special{pa 401 230}%
\special{pa 600 230}%
\special{fp}%
\special{pn 4}%
\special{pa 399 410}%
\special{pa 602 410}%
\special{fp}%
\special{pn 4}%
\special{pa 200 320}%
\special{pa 440 320}%
\special{fp}%
\put(5.000,-2.0){\makebox(0,0){\tiny 2s}}%
\end{picture}%
}

\newcommand{\marujj}{
\unitlength 0.1in
\begin{picture}(4.0,2.4)(3.60,-3.6)
\special{pn 13}%
\special{ar 660 403 60 60  0.0000000 6.2831853}%
\special{pn 13}%
\special{ar 460 403 60 60  1.5541312 1.5707963}%
\special{pn 13}%
\special{ar 460 403 60 60  0.0000000 6.2831853}%
\special{pn 13}%
\special{ar 560 230 60 60  0.0000000 6.2831853}%
\special{pn 4}%
\special{pa 530 282}%
\special{pa 490 352}%
\special{fp}%
\special{pn 4}%
\special{pa 590 282}%
\special{pa 630 354}%
\special{fp}%
\special{pn 4}%
\special{pa 520 404}%
\special{pa 600 404}%
\special{fp}%
\put(5.6100,-4.5){\makebox(0,0){\tiny 2s}}%
\end{picture}%
}

\newcommand{\marujk}{
\unitlength 0.1in
\begin{picture}(3.40,2.40)(11.50,-3.60)
\special{pn 16}%
\special{ar 1320 320 120 120  0.0000000 6.2831853}%
\special{pn 4}%
\special{ar 1470 260 120 120  1.9295670 3.6052403}%
\special{pn 4}%
\special{ar 1470 380 120 120  2.7000003 4.3440401}%
\special{pn 4}%
\special{ar 1170 260 120 120  5.8304068 6.2831853}%
\special{ar 1170 260 120 120  0.0000000 1.2277724}%
\special{pn 4}%
\special{ar 1170 380 120 120  5.0564683 6.2831853}%
\special{ar 1170 380 120 120  0.0000000 0.4411795}%
\end{picture}%
}

\newcommand{\marujl}{
\unitlength 0.1in
\begin{picture}(3.40,2.40)(1.50,-3.60)
\special{pn 16}%
\special{ar 320 320 120 120  0.0000000 6.2831853}%
\special{pn 4}%
\special{ar 320 275 30 30  0.0000000 6.2831853}%
\special{pn 4}%
\special{ar 320 365 30 30  0.0000000 6.2831853}%
\special{pn 4}%
\special{pa 209 275}%
\special{pa 290 275}%
\special{fp}%
\special{pn 4}%
\special{pa 350 275}%
\special{pa 431 275}%
\special{fp}%
\special{pn 4}%
\special{pa 209 365}%
\special{pa 290 365}%
\special{fp}%
\special{pn 4}%
\special{pa 350 365}%
\special{pa 431 365}%
\special{fp}%
\end{picture}%
}

\newcommand{\marukd}{
\unitlength 0.1in
\begin{picture}(3.40,2.40)(11.50,-3.60)
\special{pn 16}%
\special{ar 1320 320 120 120  0.0000000 6.2831853}%
\special{pn 4}%
\special{pa 1209 275}%
\special{pa 1431 275}%
\special{fp}%
\special{pn 4}%
\special{ar 1275 365 30 30  0.0000000 6.2831853}%
\special{pn 4}%
\special{ar 1365 365 30 30  0.0000000 6.2831853}%
\special{pn 4}%
\special{pa 1209 365}%
\special{pa 1243 365}%
\special{fp}%
\special{pn 4}%
\special{pa 1306 365}%
\special{pa 1335 365}%
\special{fp}%
\special{pn 4}%
\special{pa 1395 365}%
\special{pa 1431 365}%
\special{fp}%
\end{picture}%
}

\newcommand{\maruke}{
\unitlength 0.1in
\begin{picture}(3.40,2.40)(1.50,-3.60)
\special{pn 16}%
\special{ar 320 320 120 120  0.0000000 6.2831853}%
\special{pn 4}%
\special{pa 200 320}%
\special{pa 440 320}%
\special{fp}%
\special{pn 4}%
\special{pa 216 260}%
\special{pa 423 260}%
\special{fp}%
\special{pn 4}%
\special{ar 320 380 30 30  0.0000000 6.2831853}%
\special{pn 4}%
\special{pa 216 380}%
\special{pa 290 380}%
\special{fp}%
\special{pn 4}%
\special{pa 350 380}%
\special{pa 422 380}%
\special{fp}%
\end{picture}%
}

\newcommand{\marukf}{
\unitlength 0.1in
\begin{picture}(3.40,2.40)(1.50,-3.60)
\special{pn 16}%
\special{ar 320 320 120 120  0.0000000 6.2831853}%
\special{pn 4}%
\special{pa 225 245}%
\special{pa 413 245}%
\special{fp}%
\special{pn 4}%
\special{pa 227 395}%
\special{pa 414 395}%
\special{fp}%
\special{pn 4}%
\special{pa 203 294}%
\special{pa 437 294}%
\special{fp}%
\special{pn 4}%
\special{pa 203 345}%
\special{pa 437 345}%
\special{fp}%
\end{picture}%
}

\newcommand{\marukg}{
\unitlength 0.1in
\begin{picture}(6.90,2.40)(1.50,-3.60)
\special{pn 16}%
\special{ar 320 320 120 120  0.0000000 6.2831853}%
\special{pn 16}%
\special{ar 680 320 120 120  0.0000000 6.2831853}%
\special{pn 4}%
\special{pa 260 216}%
\special{pa 380 423}%
\special{fp}%
\special{pn 4}%
\special{pa 380 216}%
\special{pa 260 423}%
\special{fp}%
\special{pn 4}%
\special{pa 425 260}%
\special{pa 575 260}%
\special{fp}%
\special{pn 4}%
\special{pa 425 380}%
\special{pa 575 380}%
\special{fp}%
\end{picture}%
}

\newcommand{\marukh}{
\unitlength 0.1in
\begin{picture}(6.90,2.40)(1.50,-3.60)
\special{pn 16}%
\special{ar 320 320 120 120  0.0000000 6.2831853}%
\special{pn 16}%
\special{ar 680 320 120 120  0.0000000 6.2831853}%
\special{pn 4}%
\special{pa 290 203}%
\special{pa 290 437}%
\special{fp}%
\special{pn 4}%
\special{pa 350 203}%
\special{pa 350 437}%
\special{fp}%
\special{pn 4}%
\special{pa 425 260}%
\special{pa 575 260}%
\special{fp}%
\special{pn 4}%
\special{pa 425 380}%
\special{pa 575 380}%
\special{fp}%
\end{picture}%
}

\newcommand{\maruki}{
\unitlength 0.1in
\begin{picture}(6.90,2.40)(1.50,-3.60)
\special{pn 16}%
\special{ar 320 320 120 120  0.0000000 6.2831853}%
\special{pn 16}%
\special{ar 680 320 120 120  0.0000000 6.2831853}%
\special{pn 4}%
\special{ar 320 320 30 30  0.0000000 6.2831853}%
\special{pn 4}%
\special{pa 320 200}%
\special{pa 320 290}%
\special{fp}%
\special{pn 4}%
\special{pa 320 350}%
\special{pa 320 440}%
\special{fp}%
\special{pn 4}%
\special{pa 425 260}%
\special{pa 575 260}%
\special{fp}%
\special{pn 4}%
\special{pa 425 380}%
\special{pa 575 380}%
\special{fp}%
\end{picture}%
}

\newcommand{\marukk}{
\unitlength 0.1in
\begin{picture}(6.90,2.40)(1.50,-3.60)
\special{pn 16}%
\special{ar 320 320 120 120  0.0000000 6.2831853}%
\special{pn 16}%
\special{ar 680 320 120 120  0.0000000 6.2831853}%
\special{pn 4}%
\special{pa 320 200}%
\special{pa 320 440}%
\special{fp}%
\special{pn 4}%
\special{pa 416 395}%
\special{pa 587 395}%
\special{fp}%
\special{pn 4}%
\special{ar 500 245 30 30  0.0000000 6.2831853}%
\special{pn 4}%
\special{pa 414 245}%
\special{pa 470 245}%
\special{fp}%
\special{pn 4}%
\special{pa 530 245}%
\special{pa 585 245}%
\special{fp}%
\end{picture}%
}

\newcommand{\marukm}{
\unitlength 0.1in
\begin{picture}(6.90,2.40)(1.50,-3.60)
\special{pn 16}%
\special{ar 320 320 120 120  0.0000000 6.2831853}%
\special{pn 16}%
\special{ar 680 320 120 120  0.0000000 6.2831853}%
\special{pn 4}%
\special{pa 320 200}%
\special{pa 320 440}%
\special{fp}%
\special{pn 4}%
\special{pa 440 320}%
\special{pa 560 320}%
\special{fp}%
\special{pn 4}%
\special{pa 413 245}%
\special{pa 588 245}%
\special{fp}%
\special{pn 4}%
\special{pa 413 395}%
\special{pa 588 395}%
\special{fp}%
\end{picture}%
}

\newcommand{\marukn}{
\unitlength 0.1in
\begin{picture}(6.90,2.40)(1.500,-3.60)
\special{pn 16}%
\special{ar 320 320 120 120  0.0000000 6.2831853}%
\special{pn 16}%
\special{ar 680 320 120 120  0.0000000 6.2831853}%
\special{pn 4}%
\special{pa 399 230}%
\special{pa 602 230}%
\special{fp}%
\special{pn 4}%
\special{pa 437 290}%
\special{pa 563 290}%
\special{fp}%
\special{pn 4}%
\special{pa 401 410}%
\special{pa 602 410}%
\special{fp}%
\special{pn 4}%
\special{pa 204 350}%
\special{pa 437 350}%
\special{fp}%
\end{picture}%
}

\newcommand{\marukp}{
\unitlength 0.1in
\begin{picture}(6.90,2.40)(1.50,-3.60)
\special{pn 16}%
\special{ar 320 320 120 120  0.0000000 6.2831853}%
\special{pn 16}%
\special{ar 680 320 120 120  0.0000000 6.2831853}%
\special{pn 4}%
\special{pa 320 200}%
\special{pa 320 440}%
\special{fp}%
\special{pn 4}%
\special{pa 680 200}%
\special{pa 680 440}%
\special{fp}%
\special{pn 4}%
\special{pa 414 245}%
\special{pa 587 245}%
\special{fp}%
\special{pn 4}%
\special{pa 416 395}%
\special{pa 587 395}%
\special{fp}%
\end{picture}%
}

\newcommand{\marukq}{
\unitlength 0.1in
\begin{picture}(6.90,2.40)(1.500,-3.60)
\special{pn 16}%
\special{ar 320 320 120 120  0.0000000 6.2831853}%
\special{pn 16}%
\special{ar 680 320 120 120  0.0000000 6.2831853}%
\special{pn 4}%
\special{pa 401 230}%
\special{pa 599 230}%
\special{fp}%
\special{pn 4}%
\special{pa 431 275}%
\special{pa 567 275}%
\special{fp}%
\special{pn 4}%
\special{pa 401 410}%
\special{pa 560 320}%
\special{fp}%
\special{pn 4}%
\special{pa 440 320}%
\special{pa 602 410}%
\special{fp}%
\end{picture}%
}

\newcommand{\marukr}{
\unitlength 0.1in
\begin{picture}(4.0,2.4)(3.60,-3.6)
\special{pn 13}%
\special{ar 560 230 60 60  0.0000000 6.2831853}%
\special{pn 13}%
\special{ar 460 403 60 60  0.0000000 6.2831853}%
\special{pn 13}%
\special{ar 660 403 60 60  0.0000000 6.2831853}%
\special{pn 4}%
\special{pa 530 281}%
\special{pa 490 351}%
\special{fp}%
\special{pn 4}%
\special{pa 520 403}%
\special{pa 600 403}%
\special{fp}%
\special{pn 4}%
\special{pa 590 283}%
\special{pa 630 350}%
\special{fp}%
\special{pn 4}%
\special{pa 500 230}%
\special{pa 620 230}%
\special{fp}%
\end{picture}%
}

\newcommand{\pice}{
\unitlength 0.1in
\begin{picture}(7.00,4.00)(6.00,-8.4)
\special{pn 4}%
\special{pa 1000 800}%
\special{pa 1300 800}%
\special{fp}% 
\special{pn 4}%
\special{ar 800 800 200 200  4.7123890 6.2831853}%
\special{ar 800 800 200 200  0.0000000 1.5707963}%
\special{pn 4}%
\special{ar 800 600 200 400  1.5707963 3.1415927}%
\special{pn 4}%
\special{ar 800 1000 200 400  3.1415927 4.7123890}%
\end{picture}%
}
\newcommand{\picf}{
\unitlength 0.1in
\begin{picture}(6.00,4.00)(8.00,-10.400)
\special{pn 4}%
\special{pa 1000 1000}%
\special{pa 1400 1000}%
\special{fp}%
\special{pn 4}%
\special{pa 1000 1000}%
\special{pa 800 800}%
\special{fp}%
\special{pn 4}%
\special{pa 1000 1000}%
\special{pa 800 1200}%
\special{fp}%
\end{picture}%
}
\newcommand{\picg}{
\unitlength 0.1in
\begin{picture}(4.00,4.00)(4.00,-6.40)
\special{pn 4}%
\special{pa 400 400}%
\special{pa 800 400}%
\special{fp}%
\special{pa 400 800}%
\special{pa 800 800}%
\special{fp}%
\special{pa 600 400}%
\special{pa 600 800}%
\special{fp}%
\end{picture}%
}

\newcommand{\pich}{
\unitlength 0.1in
\begin{picture}(4.00,4.00)(4.00,-6.40) 
\special{pn 4}%
\special{pa 400 400}%
\special{pa 400 800}%
\special{fp}%
\special{pa 800 400}%
\special{pa 800 800}%
\special{fp}%
\special{pa 400 600}%
\special{pa 800 600}%
\special{fp}%
\end{picture}%
}

\newcommand{\pici}{
\unitlength 0.1in
\begin{picture}(4.00,4.00)(4.00,-6.40)
\special{pn 4}%
\special{pa 400 400}%
\special{pa 800 800}%
\special{fp}%
\special{pa 800 400}%
\special{pa 400 800}%
\special{fp}%
\special{pn 4}%
\special{pa 500 700}%
\special{pa 700 700}%
\special{fp}%
\end{picture}%
}

\newcommand{\picj}{
\unitlength 0.1in
\begin{picture}(5.00,4.05)(4.00,-6.05)
\special{pn 16}%
\special{pa 400 800}%
\special{pa 880 800}%
\special{fp}%
\special{pn 4}%
\special{sh 1}%
\special{pa 900 800}%
\special{pa 820 760}%
\special{pa 847 800}%
\special{pa 820 840}%
\special{pa 900 800}%
\special{fp}%
\special{pn 4}%
\special{pa 500 400}%
\special{pa 700 800}%
\special{fp}%
\special{pn 4}%
\special{pa 700 400}%
\special{pa 500 800}%
\special{fp}%
\end{picture}%
}

\newcommand{\pick}{%WinTpicVersion2.15
\unitlength 0.1in
\begin{picture}(5.00,4.05)(4.00,-6.05)
\special{pn 16}%
\special{pa 400 800}%
\special{pa 880 800}%
\special{fp}%
\special{pn 4}%
\special{sh 1}%
\special{pa 900 800}%
\special{pa 820 760}%
\special{pa 847 800}%
\special{pa 820 840}%
\special{pa 900 800}%
\special{fp}%
\special{pn 4}%
\special{pa 500 400}%
\special{pa 500 800}%
\special{fp}%
\special{pn 4}%
\special{pa 700 400}%
\special{pa 700 800}%
\special{fp}%
\end{picture}%
}

\newcommand{\picl}{
\unitlength 0.1in
\begin{picture}(5.00,4.05)(4.00,-6.05)
\special{pn 16}%
\special{pa 400 800}%
\special{pa 880 800}%
\special{fp}%
\special{pn 4}%
\special{sh 1}%
\special{pa 900 800}%
\special{pa 820 760}%
\special{pa 847 800}%
\special{pa 820 840}%
\special{pa 900 800}%
\special{fp}%
\special{pn 4}%
\special{pa 500 400}%
\special{pa 600 600}%
\special{fp}%
\special{pn 4}%
\special{pa 700 400}%
\special{pa 600 600}%
\special{fp}%
\special{pn 4}%
\special{pa 600 600}%
\special{pa 600 800}%
\special{fp}%
\end{picture}%
}

\newcommand{\picm}{
\unitlength 0.1in
\begin{picture}(8.35,6.00)(2.65,-9.15)
\special{pn 4}%
\special{pa 800 800}%
\special{pa 500 500}%
\special{fp}%
\special{pn 4}%
\special{pa 800 800}%
\special{pa 1100 500}%
\special{fp}%
\special{pn 4}%
\special{pa 800 800}%
\special{pa 800 1200}%
\special{fp}%
\put(12.0000,-4.0000){\makebox(0,0){$a$}}%
\put(4.0000,-4.0000){\makebox(0,0){$b$}}%
\put(8.0000,-13.0000){\makebox(0,0){$c$}}%
\end{picture}%
}

\newcommand{\picn}{
\unitlength 0.1in
\begin{picture}(8.00,8.90)(4.00,-9.15)
\special{pn 16}%
\special{pa 400 1200}%
\special{pa 1180 1200}%
\special{fp}%
\special{pn 4}%
\special{sh 1}%
\special{pa 1200 1200}%
\special{pa 1120 1160}%
\special{pa 1147 1200}%
\special{pa 1120 1240}%
\special{pa 1200 1200}%
\special{fp}%
\special{pn 4}%
\special{pa 800 600}%
\special{pa 800 1200}%
\special{fp}%
\put(8.0000,-5.0000){\makebox(0,0){$a$}}%
\put(8.0000,-13.9000){\makebox(0,0){$T^a$}}%
\end{picture}%
}

\newcommand{\pico}{
\unitlength 0.1in
\begin{picture}(14.19,5.05)(1.61,-4.80)
\special{pn 16}%
\special{ar 440 440 240 240  0.0000000 6.2831853}%
\special{pn 16}%
\special{ar 1340 440 240 240  0.0000000 6.2831853}%
\special{pn 4}%
\special{pa 626 290}%
\special{pa 1154 290}%
\special{fp}%
\special{pn 4}%
\special{pa 629 590}%
\special{pa 1154 590}%
\special{fp}%
\special{pn 4}%
\special{pa 440 440}%
\special{pa 680 440}%
\special{fp}%
\special{pn 4}%
\special{pa 440 440}%
\special{pa 269 269}%
\special{fp}%
\special{pn 4}%
\special{pa 440 440}%
\special{pa 272 608}%
\special{fp}%
\special{pn 4}%
\special{pa 1100 440}%
\special{pa 1580 440}%
\special{fp}%
\put(5.6000,-4.0700){\makebox(0,0){{\scriptsize $a$}}}%
\put(2.9600,-3.5000){\makebox(0,0){{\scriptsize $b$}}}%
\put(3.4700,-5.9000){\makebox(0,0){}}%
\put(3.5300,-5.9000){\makebox(0,0){{\scriptsize $c$}}}%
\put(8.9000,-2.4000){\makebox(0,0){{\scriptsize $d$}}}%
\put(8.9000,-6.3000){\makebox(0,0){{\scriptsize $e$}}}%
\put(13.4300,-3.90){\makebox(0,0){{\scriptsize $f$}}}%
\end{picture}%
}

\newcommand{\picca}{
\unitlength 0.1in
\begin{picture}(3.30,4.90)(1.95,-4)
\special{pn 13}%
\special{ar 360 360 160 160  0.0000000 6.2831853}%
\special{pn 4}%
\special{pa 360 200}%
\special{pa 360 520}%
\special{fp}%
\special{pn 8}%
\special{pa 518 370}%
\special{pa 520 360}%
\special{fp}%
\special{sh 1}%
\special{pa 520 360}%
\special{pa 487 421}%
\special{pa 510 412}%
\special{pa 527 429}%
\special{pa 520 360}%
\special{fp}%
\special{pn 8}%
\special{pa 202 350}%
\special{pa 200 360}%
\special{fp}%
\special{sh 1}%
\special{pa 200 360}%
\special{pa 233 299}%
\special{pa 210 308}%
\special{pa 193 291}%
\special{pa 200 360}%
\special{fp}%
\put(3.7000,-1.6000){\makebox(0,0){\scriptsize $a$}}%
\put(3.6000,-6.8000){\makebox(0,0){$\{K:a\}$}}%
\end{picture}%
}

\newcommand{\piccb}{
\unitlength 0.1in
\begin{picture}(4.51,4.55)(4.43,-4.0)
\special{pn 8}%
\special{ar 648 360 160 160  1.5707963 4.7123890}%
\special{pn 8}%
\special{pa 648 520}%
\special{pa 648 200}%
\special{fp}%
\special{sh 1}%
\special{pa 648 200}%
\special{pa 628 267}%
\special{pa 648 253}%
\special{pa 668 267}%
\special{pa 648 200}%
\special{fp}%
\special{pn 8}%
\special{ar 729 360 160 160  4.7123890 6.2831853}%
\special{ar 729 360 160 160  0.0000000 1.5707963}%
\special{pn 8}%
\special{pa 729 200}%
\special{pa 729 520}%
\special{fp}%
\special{sh 1}%
\special{pa 729 520}%
\special{pa 749 453}%
\special{pa 729 467}%
\special{pa 709 453}%
\special{pa 729 520}%
\special{fp}%
\put(8.8800,-6.7000){\makebox(0,0){$K^{[1]}_-$}}%
\put(4.8800,-6.7000){\makebox(0,0){$K^{[1]}_+$}}%
\special{pn 8}%
\special{pa 887 370}%
\special{pa 889 360}%
\special{fp}%
\special{sh 1}%
\special{pa 889 360}%
\special{pa 856 421}%
\special{pa 879 412}%
\special{pa 896 429}%
\special{pa 889 360}%
\special{fp}%
\special{pn 8}%
\special{pa 489 350}%
\special{pa 487 360}%
\special{fp}%
\special{sh 1}%
\special{pa 487 360}%
\special{pa 520 299}%
\special{pa 497 308}%
\special{pa 480 291}%
\special{pa 487 360}%
\special{fp}%
\put(6.9000,-1.5000){\makebox(0,0){\scriptsize $a$}}%
\end{picture}%
}

\newcommand{\piccc}{
\unitlength 0.1in
\begin{picture}(8.10,3.20)(1.95,-4.0)
\special{pn 13}%
\special{ar 360 360 160 160  0.0000000 6.2831853}%
\special{pn 13}%
\special{ar 840 360 160 160  0.0000000 6.2831853}%
\special{pn 4}%
\special{pa 520 360}%
\special{pa 680 360}%
\special{fp}%
\special{pn 4}%
\special{pa 998 370}%
\special{pa 1000 360}%
\special{fp}%
\special{sh 1}%
\special{pa 1000 360}%
\special{pa 967 421}%
\special{pa 990 412}%
\special{pa 1007 429}%
\special{pa 1000 360}%
\special{fp}%
\special{pn 4}%
\special{pa 202 350}%
\special{pa 200 360}%
\special{fp}%
\special{sh 1}%
\special{pa 200 360}%
\special{pa 233 299}%
\special{pa 210 308}%
\special{pa 193 291}%
\special{pa 200 360}%
\special{fp}%
\put(6.0000,-3.1000){\makebox(0,0){\scriptsize $a$}}%
\put(6.0000,-6.5000){\makebox(0,0){$\{K_1,K_2:a\}$}}%
\end{picture}%
}

\newcommand{\piccd}{
\unitlength 0.1in
\begin{picture}(8.10,3.45)(1.95,-4.0)
\special{pn 8}%
\special{ar 360 360 160 160  0.2449787 6.0382066}%
\special{pn 8}%
\special{ar 840 360 160 160  3.3865713 6.2831853}%
\special{ar 840 360 160 160  0.0000000 2.8966140}%
\special{pn 8}%
\special{pa 686 320}%
\special{pa 514 320}%
\special{fp}%
\special{sh 1}%
\special{pa 514 320}%
\special{pa 581 340}%
\special{pa 567 320}%
\special{pa 581 300}%
\special{pa 514 320}%
\special{fp}%
\special{pn 8}%
\special{pa 514 400}%
\special{pa 686 400}%
\special{fp}%
\special{sh 1}%
\special{pa 686 400}%
\special{pa 619 380}%
\special{pa 633 400}%
\special{pa 619 420}%
\special{pa 686 400}%
\special{fp}%
\put(6.0000,-6.5000){\makebox(0,0){$K^{[2]}$}}%
\special{pn 8}%
\special{pa 998 370}%
\special{pa 1000 360}%
\special{fp}%
\special{sh 1}%
\special{pa 1000 360}%
\special{pa 967 421}%
\special{pa 990 412}%
\special{pa 1007 429}%
\special{pa 1000 360}%
\special{fp}%
\special{pn 8}%
\special{pa 202 350}%
\special{pa 200 360}%
\special{fp}%
\special{sh 1}%
\special{pa 200 360}%
\special{pa 233 299}%
\special{pa 210 308}%
\special{pa 193 291}%
\special{pa 200 360}%
\special{fp}%
\put(6.1000,-2.6000){\makebox(0,0){\scriptsize $a$}}% 
\end{picture}%
}

\newcommand{\picce}{
\unitlength 0.1in
\begin{picture}(4.40,4.45)(1.05,-4.0)
\special{pn 13}%
\special{ar 380 360 160 160  0.0000000 6.2831853}%
\special{pn 4}%
\special{pa 380 200}%
\special{pa 380 520}%
\special{fp}%
\special{pn 4}%
\special{pa 538 370}%
\special{pa 540 360}%
\special{fp}%
\special{sh 1}%
\special{pa 540 360}%
\special{pa 507 421}%
\special{pa 530 412}%
\special{pa 547 429}%
\special{pa 540 360}%
\special{fp}%
\special{pn 4}%
\special{pa 222 350}%
\special{pa 220 360}%
\special{fp}%
\special{sh 1}%
\special{pa 220 360}%
\special{pa 253 299}%
\special{pa 230 308}%
\special{pa 207 291}%
\special{pa 220 360}%
\special{fp}%
\special{pn 4}%
\special{pa 220 360}%
\special{pa 540 360}%
\special{fp}%
\put(3.8000,-1.6000){\makebox(0,0){\scriptsize $a_1$}}%
\put(1.5000,-3.6000){\makebox(0,0){\scriptsize $a_2$}}%
\put(3.8000,-6.5000){\makebox(0,0){$\{K:a_1,a_2\}$}}%
\end{picture}%
}

\newcommand{\piccf}{
\unitlength 0.1in
\begin{picture}(4.51,4.65)(4.43,-4.0)
\special{pn 8}%
\special{ar 648 360 160 160  1.5707963 4.7123890}%
\special{pn 8}%
\special{pa 648 520}%
\special{pa 648 200}%
\special{fp}%
\special{sh 1}%
\special{pa 648 200}%
\special{pa 628 267}%
\special{pa 648 253}%
\special{pa 668 267}%
\special{pa 648 200}%
\special{fp}%
\special{pn 8}%
\special{ar 729 360 160 160  4.7123890 6.2831853}%
\special{ar 729 360 160 160  0.0000000 1.5707963}%
\special{pn 8}%
\special{pa 729 200}%
\special{pa 729 520}%
\special{fp}%
\special{sh 1}%
\special{pa 729 520}%
\special{pa 749 453}%
\special{pa 729 467}%
\special{pa 709 453}%
\special{pa 729 520}%
\special{fp}%
\put(8.8800,-6.5000){\makebox(0,0){$K^{[3]}_-$}}%
\put(4.8800,-6.5000){\makebox(0,0){$K^{[3]}_+$}}%
\special{pn 8}%
\special{pa 887 370}%
\special{pa 889 360}%
\special{fp}%
\special{sh 1}%
\special{pa 889 360}%
\special{pa 856 421}%
\special{pa 879 412}%
\special{pa 896 429}%
\special{pa 889 360}%
\special{fp}%
\special{pn 8}%
\special{pa 489 350}%
\special{pa 487 360}%
\special{fp}%
\special{sh 1}%
\special{pa 487 360}%
\special{pa 520 299}%
\special{pa 497 308}%
\special{pa 480 291}%
\special{pa 487 360}%
\special{fp}%
\put(6.9000,-1.4000){\makebox(0,0){\scriptsize $a_1$}}%
\end{picture}%
}

\newcommand{\piccg}{
\unitlength 0.1in
\begin{picture}(6.06,4.06)(0.19,-4.3)
\special{pn 8}%
\special{ar 244 360 160 160  3.1415927 6.2831853}%
\special{pn 8}%
\special{pa 84 360}%
\special{pa 404 360}%
\special{fp}%
\special{sh 1}%
\special{pa 404 360}%
\special{pa 337 340}%
\special{pa 351 360}%
\special{pa 337 380}%
\special{pa 404 360}%
\special{fp}%
\special{pn 8}%
\special{pa 254 202}%
\special{pa 244 200}%
\special{fp}%
\special{sh 1}%
\special{pa 244 200}%
\special{pa 305 233}%
\special{pa 296 210}%
\special{pa 313 193}%
\special{pa 244 200}%
\special{fp}%
\special{pn 8}%
\special{ar 244 400 160 160  6.2831853 6.2831853}%
\special{ar 244 400 160 160  0.0000000 3.1415927}%
\special{pn 8}%
\special{pa 404 400}%
\special{pa 84 400}%
\special{fp}%
\special{sh 1}%
\special{pa 84 400}%
\special{pa 151 420}%
\special{pa 137 400}%
\special{pa 151 380}%
\special{pa 84 400}%
\special{fp}%
\special{pn 8}%
\special{pa 234 558}%
\special{pa 244 560}%
\special{fp}%
\special{sh 1}%
\special{pa 244 560}%
\special{pa 183 527}%
\special{pa 192 550}%
\special{pa 175 567}%
\special{pa 244 560}%
\special{fp}%
\put(4.2400,-3.8000){\makebox(0,0){}}%
\put(0.20,-3.8000){\makebox(0,0){\scriptsize $a_2$}}%
\put(5.6600,-5.4000){\makebox(0,0){$K^{[4]}_-$}}%
\put(5.7000,-2.4400){\makebox(0,0){$K^{[4]}_+$}}%
\end{picture}%
}

\newcommand{\picch}{
\unitlength 0.1in
\begin{picture}(4.15,4.75)(1.25,-4.0)
\special{pn 8}%
\special{ar 380 380 160 160  4.7123890 6.2831853}%
\special{pn 8}%
\special{ar 380 380 160 160  1.5707963 3.1415927}%
\special{pn 8}%
\special{ar 380 380 160 160  3.1415927 3.2165927}%
\special{ar 380 380 160 160  3.4415927 3.5165927}%
\special{ar 380 380 160 160  3.7415927 3.8165927}%
\special{ar 380 380 160 160  4.0415927 4.1165927}%
\special{ar 380 380 160 160  4.3415927 4.4165927}%
\special{ar 380 380 160 160  4.6415927 4.7123890}%
\special{pn 8}%
\special{pa 380 220}%
\special{pa 380 540}%
\special{fp}%
\special{sh 1}%
\special{pa 380 540}%
\special{pa 400 473}%
\special{pa 380 487}%
\special{pa 360 473}%
\special{pa 380 540}%
\special{fp}%
\special{pn 8}%
\special{pa 220 380}%
\special{pa 540 380}%
\special{fp}%
\special{sh 1}%
\special{pa 540 380}%
\special{pa 473 360}%
\special{pa 487 380}%
\special{pa 473 400}%
\special{pa 540 380}%
\special{fp}%
\special{pn 8}%
\special{ar 380 380 160 160  6.2831853 6.3581853}%
\special{ar 380 380 160 160  6.5831853 6.6581853}%
\special{ar 380 380 160 160  6.8831853 6.9581853}%
\special{ar 380 380 160 160  7.1831853 7.2581853}%
\special{ar 380 380 160 160  7.4831853 7.5581853}%
\special{ar 380 380 160 160  7.7831853 7.8539816}%
\special{pn 8}%
\special{pa 390 222}%
\special{pa 380 220}%
\special{fp}%
\special{sh 1}%
\special{pa 380 220}%
\special{pa 441 253}%
\special{pa 432 230}%
\special{pa 449 213}%
\special{pa 380 220}%
\special{fp}%
\special{pn 8}%
\special{pa 222 390}%
\special{pa 220 380}%
\special{fp}%
\special{sh 1}%
\special{pa 220 380}%
\special{pa 213 449}%
\special{pa 230 432}%
\special{pa 253 441}%
\special{pa 220 380}%
\special{fp}%
\put(3.8000,-6.5000){\makebox(0,0){$K^{[5]}_+$}}%
\put(3.8000,-1.5000){\makebox(0,0){\scriptsize $a_1$}}%
\put(1.600,-3.8000){\makebox(0,0){\scriptsize $a_2$}}%
\end{picture}%
}

\newcommand{\picci}{
\unitlength 0.1in
\begin{picture}(3.85,4.30)(0.19,-2.3)
\special{pn 8}%
\special{ar 244 212 160 160  3.1415927 4.7123890}%
\special{pn 8}%
\special{ar 244 212 160 160  6.2831853 6.2831853}%
\special{ar 244 212 160 160  0.0000000 1.5707963}%
\special{pn 8}%
\special{ar 244 212 160 160  4.7123890 4.7873890}%
\special{ar 244 212 160 160  5.0123890 5.0873890}%
\special{ar 244 212 160 160  5.3123890 5.3873890}%
\special{ar 244 212 160 160  5.6123890 5.6873890}%
\special{ar 244 212 160 160  5.9123890 5.9873890}%
\special{ar 244 212 160 160  6.2123890 6.2831853}%
\special{pn 8}%
\special{ar 244 212 160 160  1.5707963 1.6457963}%
\special{ar 244 212 160 160  1.8707963 1.9457963}%
\special{ar 244 212 160 160  2.1707963 2.2457963}%
\special{ar 244 212 160 160  2.4707963 2.5457963}%
\special{ar 244 212 160 160  2.7707963 2.8457963}%
\special{ar 244 212 160 160  3.0707963 3.1415927}%
\special{pn 8}%
\special{pa 84 212}%
\special{pa 404 212}%
\special{fp}%
\special{sh 1}%
\special{pa 404 212}%
\special{pa 337 192}%
\special{pa 351 212}%
\special{pa 337 232}%
\special{pa 404 212}%
\special{fp}%
\special{pn 8}%
\special{pa 244 372}%
\special{pa 244 52}%
\special{fp}%
\special{sh 1}%
\special{pa 244 52}%
\special{pa 224 119}%
\special{pa 244 105}%
\special{pa 264 119}%
\special{pa 244 52}%
\special{fp}%
\special{pn 8}%
\special{pa 86 202}%
\special{pa 84 212}%
\special{fp}%
\special{sh 1}%
\special{pa 84 212}%
\special{pa 117 151}%
\special{pa 94 160}%
\special{pa 77 143}%
\special{pa 84 212}%
\special{fp}%
\put(0.2500,-2.1200){\makebox(0,0){\scriptsize $a_2$}}%
\put(2.4400,-0.1500){\makebox(0,0){\scriptsize $a_1$}}%
\put(2.4400,-5.00){\makebox(0,0){$K^{[5]}_-$}}%
\special{pn 8}%
\special{pa 254 370}%
\special{pa 244 372}%
\special{fp}%
\special{sh 1}%
\special{pa 244 372}%
\special{pa 313 379}%
\special{pa 296 362}%
\special{pa 305 339}%
\special{pa 244 372}%
\special{fp}%
\end{picture}%
}

\newcommand{\piccj}{
\unitlength 0.1in
\begin{picture}(3.30,4.21)(1.95,-3.3)
\special{pn 13}%
\special{ar 360 328 160 160  0.0000000 6.2831853}%
\special{pn 4}%
\special{pa 420 180}%
\special{pa 420 474}%
\special{fp}%
\special{pn 4}%
\special{pa 300 180}%
\special{pa 300 478}%
\special{fp}%
\put(4.6000,-1.30){\makebox(0,0){\scriptsize $a_2$}}%
\put(3.0000,-1.30){\makebox(0,0){\scriptsize $a_1$}}%
\put(3.6000,-6.00){\makebox(0,0){$\{K:a_1,a_2\}$}}%
\special{pn 8}%
\special{pa 518 338}%
\special{pa 520 328}%
\special{dt 0.045}%
\special{sh 1}%
\special{pa 520 328}%
\special{pa 487 389}%
\special{pa 510 380}%
\special{pa 527 397}%
\special{pa 520 328}%
\special{fp}%
\special{pn 8}%
\special{pa 202 318}%
\special{pa 200 328}%
\special{dt 0.045}%
\special{sh 1}%
\special{pa 200 328}%
\special{pa 233 267}%
\special{pa 210 276}%
\special{pa 193 259}%
\special{pa 200 328}%
\special{fp}%
\end{picture}%
}

\newcommand{\picck}{
\unitlength 0.1in
\begin{picture}(3.90,4.19)(1.95,-2.5)
\special{pn 8}%
\special{ar 360 236 160 160  1.9513027 4.3272315}%
\special{pn 8}%
\special{ar 420 236 160 160  4.3318826 6.2831853}%
\special{ar 420 236 160 160  0.0000000 1.9513027}%
\special{pn 8}%
\special{pa 360 86}%
\special{pa 360 386}%
\special{fp}%
\special{sh 1}%
\special{pa 360 386}%
\special{pa 380 319}%
\special{pa 360 333}%
\special{pa 340 319}%
\special{pa 360 386}%
\special{fp}%
\special{pn 8}%
\special{pa 300 386}%
\special{pa 300 86}%
\special{fp}%
\special{sh 1}%
\special{pa 300 86}%
\special{pa 280 153}%
\special{pa 300 139}%
\special{pa 320 153}%
\special{pa 300 86}%
\special{fp}%
\special{pn 8}%
\special{pa 480 86}%
\special{pa 480 384}%
\special{dt 0.045}%
\special{pa 480 384}%
\special{pa 480 383}%
\special{dt 0.045}%
\special{pn 8}%
\special{pa 578 246}%
\special{pa 580 236}%
\special{dt 0.045}%
\special{sh 1}%
\special{pa 580 236}%
\special{pa 547 297}%
\special{pa 570 288}%
\special{pa 587 305}%
\special{pa 580 236}%
\special{fp}%
\special{pn 8}%
\special{pa 202 226}%
\special{pa 200 236}%
\special{dt 0.045}%
\special{sh 1}%
\special{pa 200 236}%
\special{pa 233 175}%
\special{pa 210 184}%
\special{pa 193 167}%
\special{pa 200 236}%
\special{fp}%
\put(3.3200,-0.200){\makebox(0,0){\scriptsize $a_1$}}%
\put(2.4000,-5.5000){\makebox(0,0){$K^{[6]}_+$}}%
\put(5.50,-5.5000){\makebox(0,0){$K^{[6]}_-$}}%
\end{picture}%
}

\newcommand{\piccl}{
\unitlength 0.1in
\begin{picture}(3.90,4.19)(1.95,-2.5)
\special{pn 8}%
\special{ar 420 236 160 160  5.0975465 6.2831853}%
\special{ar 420 236 160 160  0.0000000 1.1808785}%
\special{pn 8}%
\special{ar 360 236 160 160  1.1808785 5.0975465}%
\special{pn 8}%
\special{pa 480 88}%
\special{pa 480 386}%
\special{fp}%
\special{sh 1}%
\special{pa 480 386}%
\special{pa 500 319}%
\special{pa 480 333}%
\special{pa 460 319}%
\special{pa 480 386}%
\special{fp}%
\special{pn 8}%
\special{pa 420 386}%
\special{pa 420 88}%
\special{fp}%
\special{sh 1}%
\special{pa 420 88}%
\special{pa 400 155}%
\special{pa 420 141}%
\special{pa 440 155}%
\special{pa 420 88}%
\special{fp}%
\special{pn 8}%
\special{pa 578 246}%
\special{pa 580 236}%
\special{fp}%
\special{sh 1}%
\special{pa 580 236}%
\special{pa 547 297}%
\special{pa 570 288}%
\special{pa 587 305}%
\special{pa 580 236}%
\special{fp}%
\special{pn 8}%
\special{pa 202 226}%
\special{pa 200 236}%
\special{fp}%
\special{sh 1}%
\special{pa 200 236}%
\special{pa 233 175}%
\special{pa 210 184}%
\special{pa 193 167}%
\special{pa 200 236}%
\special{fp}%
\special{pn 8}%
\special{pa 300 86}%
\special{pa 300 386}%
\special{dt 0.045}%
\special{pa 300 386}%
\special{pa 300 385}%
\special{dt 0.045}%
\put(4.5000,-0.200){\makebox(0,0){\scriptsize $a_2$}}%
\put(2.9000,-5.5000){\makebox(0,0){$K^{[7]}_+$}}%
\put(6.000,-5.5000){\makebox(0,0){$K^{[7]}_-$}}%
\end{picture}%
}

\newcommand{\piccm}{
\unitlength 0.1in
\begin{picture}(5.44,4.15)(3.33,-7.4)
\special{pn 8}%
\special{ar 712 748 160 160  5.0928954 6.2831853}%
\special{ar 712 748 160 160  0.0000000 1.1902899}%
\special{pn 8}%
\special{ar 612 744 160 160  4.3272315 5.0975465}%
\special{pn 8}%
\special{ar 612 744 160 160  1.1856388 1.9559538}%
\special{pn 8}%
\special{ar 512 744 160 160  1.9559538 4.3364280}%
\special{pn 8}%
\special{pa 772 602}%
\special{pa 772 894}%
\special{fp}%
\special{sh 1}%
\special{pa 772 894}%
\special{pa 792 827}%
\special{pa 772 841}%
\special{pa 752 827}%
\special{pa 772 894}%
\special{fp}%
\special{pn 8}%
\special{pa 672 892}%
\special{pa 672 594}%
\special{fp}%
\special{sh 1}%
\special{pa 672 594}%
\special{pa 652 661}%
\special{pa 672 647}%
\special{pa 692 661}%
\special{pa 672 594}%
\special{fp}%
\special{pn 8}%
\special{pa 552 596}%
\special{pa 552 894}%
\special{fp}%
\special{sh 1}%
\special{pa 552 894}%
\special{pa 572 827}%
\special{pa 552 841}%
\special{pa 532 827}%
\special{pa 552 894}%
\special{fp}%
\special{pn 8}%
\special{pa 452 892}%
\special{pa 452 592}%
\special{fp}%
\special{sh 1}%
\special{pa 452 592}%
\special{pa 432 659}%
\special{pa 452 645}%
\special{pa 472 659}%
\special{pa 452 592}%
\special{fp}%
\special{pn 8}%
\special{pa 870 754}%
\special{pa 872 744}%
\special{fp}%
\special{sh 1}%
\special{pa 872 744}%
\special{pa 839 805}%
\special{pa 862 796}%
\special{pa 879 813}%
\special{pa 872 744}%
\special{fp}%
\special{pn 8}%
\special{pa 354 734}%
\special{pa 352 744}%
\special{fp}%
\special{sh 1}%
\special{pa 352 744}%
\special{pa 385 683}%
\special{pa 362 692}%
\special{pa 345 675}%
\special{pa 352 744}%
\special{fp}%
\put(5.0200,-5.20){\makebox(0,0){\scriptsize $a_1$}}%
\put(7.3200,-5.20){\makebox(0,0){\scriptsize $a_2$}}%
\put(3.30,-10.50){\makebox(0,0){$K^{[8]}_+$}}%
\put(6.40,-10.50){\makebox(0,0){$K^{[8]}_0$}}% 
\put(9.50,-10.500){\makebox(0,0){$K^{[8]}_-$}}%
\special{pn 8}%
\special{pa 562 594}%
\special{pa 552 596}%
\special{fp}%
\special{sh 1}%
\special{pa 552 596}%
\special{pa 621 603}%
\special{pa 604 586}%
\special{pa 613 563}%
\special{pa 552 596}%
\special{fp}%
\special{pn 8}%
\special{pa 662 894}%
\special{pa 672 892}%
\special{fp}%
\special{sh 1}%
\special{pa 672 892}%
\special{pa 603 885}%
\special{pa 620 902}%
\special{pa 611 925}%
\special{pa 672 892}%
\special{fp}%
\end{picture}%
}

\newcommand{\piccn}{
\unitlength 0.1in
\begin{picture}(8.10,3.73)(1.95,-6)
\special{pn 13}%
\special{ar 360 572 160 160  0.0000000 6.2831853}%
\special{pn 13}%
\special{ar 840 572 160 160  0.0000000 6.2831853}%
\special{pn 4}%
\special{pa 486 472}%
\special{pa 714 472}%
\special{fp}%
\special{pn 4}%
\special{pa 484 672}%
\special{pa 714 672}%
\special{fp}%
\put(5.9600,-4.300){\makebox(0,0){\scriptsize $a_1$}}%
\put(6.100,-6.300){\makebox(0,0){\scriptsize $a_2$}}%
\put(6.0000,-8.500){\makebox(0,0){$\{K_1,K_2:a_1,a_2\}$}}%
\special{pn 8}%
\special{pa 998 580}%
\special{pa 1000 570}%
\special{dt 0.045}%
\special{sh 1}%
\special{pa 1000 570}%
\special{pa 967 631}%
\special{pa 990 622}%
\special{pa 1007 639}%
\special{pa 1000 570}%
\special{fp}%
\special{pn 8}%
\special{pa 682 560}%
\special{pa 680 570}%
\special{dt 0.045}%
\special{sh 1}%
\special{pa 680 570}%
\special{pa 713 509}%
\special{pa 690 518}%
\special{pa 673 501}%
\special{pa 680 570}%
\special{fp}%
\special{pn 8}%
\special{pa 518 580}%
\special{pa 520 570}%
\special{dt 0.045}%
\special{sh 1}%
\special{pa 520 570}%
\special{pa 487 631}%
\special{pa 510 622}%
\special{pa 527 639}%
\special{pa 520 570}%
\special{fp}%
\special{pn 8}%
\special{pa 202 560}%
\special{pa 200 570}%
\special{dt 0.045}%
\special{sh 1}%
\special{pa 200 570}%
\special{pa 233 509}%
\special{pa 210 518}%
\special{pa 193 501}%
\special{pa 200 570}%
\special{fp}%
\end{picture}%
}

\newcommand{\picco}{
\unitlength 0.1in
\begin{picture}(8.05,3.67)(2.35,-7.5)
\special{pn 8}%
\special{ar 400 736 160 160  0.6708575 5.6045228}%
\special{pn 8}%
\special{ar 880 736 160 160  2.4629302 3.8202551}%
\special{pn 8}%
\special{ar 880 736 160 160  3.8202551 3.8952551}%
\special{ar 880 736 160 160  4.1202551 4.1952551}%
\special{ar 880 736 160 160  4.4202551 4.4952551}%
\special{ar 880 736 160 160  4.7202551 4.7952551}%
\special{ar 880 736 160 160  5.0202551 5.0952551}%
\special{ar 880 736 160 160  5.3202551 5.3952551}%
\special{ar 880 736 160 160  5.6202551 5.6952551}%
\special{ar 880 736 160 160  5.9202551 5.9952551}%
\special{ar 880 736 160 160  6.2202551 6.2952551}%
\special{ar 880 736 160 160  6.5202551 6.5952551}%
\special{ar 880 736 160 160  6.8202551 6.8952551}%
\special{ar 880 736 160 160  7.1202551 7.1952551}%
\special{ar 880 736 160 160  7.4202551 7.4952551}%
\special{ar 880 736 160 160  7.7202551 7.7952551}%
\special{ar 880 736 160 160  8.0202551 8.0952551}%
\special{ar 880 736 160 160  8.3202551 8.3952551}%
\special{ar 880 736 160 160  8.6202551 8.6952551}%
\special{pn 8}%
\special{pa 756 636}%
\special{pa 526 636}%
\special{fp}%
\special{sh 1}%
\special{pa 526 636}%
\special{pa 593 656}%
\special{pa 579 636}%
\special{pa 593 616}%
\special{pa 526 636}%
\special{fp}%
\special{pn 8}%
\special{ar 400 736 160 160  5.6123278 5.6873278}%
\special{ar 400 736 160 160  5.9123278 5.9873278}%
\special{ar 400 736 160 160  6.2123278 6.2873278}%
\special{ar 400 736 160 160  6.5123278 6.5873278}%
\special{ar 400 736 160 160  6.8123278 6.8873278}%
\special{pn 8}%
\special{pa 526 836}%
\special{pa 756 836}%
\special{fp}%
\special{sh 1}%
\special{pa 756 836}%
\special{pa 689 816}%
\special{pa 703 836}%
\special{pa 689 856}%
\special{pa 756 836}%
\special{fp}%
\special{pn 8}%
\special{pa 242 726}%
\special{pa 240 736}%
\special{dt 0.045}%
\special{sh 1}%
\special{pa 240 736}%
\special{pa 273 675}%
\special{pa 250 684}%
\special{pa 233 667}%
\special{pa 240 736}%
\special{fp}%
\put(6.700,-6.00){\makebox(0,0){\scriptsize $a_1$}}%
\put(6.4000,-7.900){\makebox(0,0){\scriptsize $a_2$}}%
\put(4.7200,-10.00){\makebox(0,0){$K^{[9]}_+$}}%
\special{pn 8}%
\special{pa 728 681}%
\special{pa 730 671}%
\special{dt 0.045}%
\special{sh 1}%
\special{pa 730 671}%
\special{pa 697 732}%
\special{pa 720 723}%
\special{pa 737 740}%
\special{pa 730 671}%
\special{fp}%
\end{picture}%
}

\newcommand{\piccp}{
\unitlength 0.1in
\begin{picture}(8.05,3.65)(2.00,-7.3)
\special{pn 8}%
\special{ar 360 720 160 160  0.6708575 0.7458575}%
\special{ar 360 720 160 160  0.9708575 1.0458575}%
\special{ar 360 720 160 160  1.2708575 1.3458575}%
\special{ar 360 720 160 160  1.5708575 1.6458575}%
\special{ar 360 720 160 160  1.8708575 1.9458575}%
\special{ar 360 720 160 160  2.1708575 2.2458575}%
\special{ar 360 720 160 160  2.4708575 2.5458575}%
\special{ar 360 720 160 160  2.7708575 2.8458575}%
\special{ar 360 720 160 160  3.0708575 3.1458575}%
\special{ar 360 720 160 160  3.3708575 3.4458575}%
\special{ar 360 720 160 160  3.6708575 3.7458575}%
\special{ar 360 720 160 160  3.9708575 4.0458575}%
\special{ar 360 720 160 160  4.2708575 4.3458575}%
\special{ar 360 720 160 160  4.5708575 4.6458575}%
\special{ar 360 720 160 160  4.8708575 4.9458575}%
\special{ar 360 720 160 160  5.1708575 5.2458575}%
\special{ar 360 720 160 160  5.4708575 5.5458575}%
\special{pn 8}%
\special{ar 360 720 160 160  5.6123278 6.2831853}%
\special{ar 360 720 160 160  0.0000000 0.6708575}%
\special{pn 8}%
\special{ar 840 720 160 160  3.8202551 6.2831853}%
\special{ar 840 720 160 160  0.0000000 2.4629302}%
\special{pn 8}%
\special{pa 716 620}%
\special{pa 486 620}%
\special{fp}%
\special{sh 1}%
\special{pa 486 620}%
\special{pa 553 640}%
\special{pa 539 620}%
\special{pa 553 600}%
\special{pa 486 620}%
\special{fp}%
\special{pn 8}%
\special{pa 486 820}%
\special{pa 716 820}%
\special{fp}%
\special{sh 1}%
\special{pa 716 820}%
\special{pa 649 800}%
\special{pa 663 820}%
\special{pa 649 840}%
\special{pa 716 820}%
\special{fp}%
\special{pn 8}%
\special{pa 998 730}%
\special{pa 1000 720}%
\special{fp}%
\special{sh 1}%
\special{pa 1000 720}%
\special{pa 967 781}%
\special{pa 990 772}%
\special{pa 1007 789}%
\special{pa 1000 720}%
\special{fp}%
\put(6.0800,-5.8000){\makebox(0,0){\scriptsize $a_1$}}%
\put(6.1000,-7.700){\makebox(0,0){\scriptsize $a_2$}}%
\put(7.2600,-10.00){\makebox(0,0){$K^{[9]}_-$}}%
\special{pn 8}%
\special{pa 514 764}%
\special{pa 512 774}%
\special{fp}%
\special{sh 1}%
\special{pa 512 774}%
\special{pa 545 713}%
\special{pa 522 722}%
\special{pa 505 705}%
\special{pa 512 774}%
\special{fp}%
\special{pn 8}%
\special{ar 840 720 160 160  2.4783897 2.5533897}%
\special{ar 840 720 160 160  2.7783897 2.8533897}%
\special{ar 840 720 160 160  3.0783897 3.1533897}%
\special{ar 840 720 160 160  3.3783897 3.4533897}%
\special{ar 840 720 160 160  3.6783897 3.7533897}%
\end{picture}%
}

\newcommand{\piccq}{
\unitlength 0.1in
\begin{picture}(8.10,4.25)(2.15,-2.5)
\special{pn 13}%
\special{ar 380 232 160 160  0.0000000 6.2831853}%
\special{pn 13}%
\special{ar 860 232 160 160  0.0000000 6.2831853}%
\special{pn 4}%
\special{pa 540 232}%
\special{pa 700 232}%
\special{fp}%
\special{pn 4}%
\special{pa 380 72}%
\special{pa 380 392}%
\special{fp}%
\put(3.8000,-0.300){\makebox(0,0){\scriptsize $a_1$}}%
\put(6.2200,-1.800){\makebox(0,0){\scriptsize $a_2$}}%
\put(6.2000,-5.00){\makebox(0,0){$\{K_1,K_2:a_1,a_2\}$}}%
\special{pn 8}%
\special{pa 1018 242}%
\special{pa 1020 232}%
\special{fp}%
\special{sh 1}%
\special{pa 1020 232}%
\special{pa 987 293}%
\special{pa 1010 284}%
\special{pa 1027 301}%
\special{pa 1020 232}%
\special{fp}%
\special{pn 8}%
\special{pa 222 222}%
\special{pa 220 232}%
\special{fp}%
\special{sh 1}%
\special{pa 220 232}%
\special{pa 253 171}%
\special{pa 230 180}%
\special{pa 213 163}%
\special{pa 220 232}%
\special{fp}%
\end{picture}%
}

\newcommand{\piccr}{
\unitlength 0.1in
\begin{picture}(8.10,3.20)(1.95,-5.20)
\special{pn 8}%
\special{ar 360 360 160 160  0.2510027 6.0321826}%
\special{pn 8}%
\special{ar 840 360 160 160  3.3957174 6.2831853}%
\special{ar 840 360 160 160  0.0000000 2.8874679}%
\special{pn 8}%
\special{pa 686 320}%
\special{pa 516 320}%
\special{fp}%
\special{sh 1}%
\special{pa 516 320}%
\special{pa 583 340}%
\special{pa 569 320}%
\special{pa 583 300}%
\special{pa 516 320}%
\special{fp}%
\special{pn 8}%
\special{pa 516 400}%
\special{pa 686 400}%
\special{fp}%
\special{sh 1}%
\special{pa 686 400}%
\special{pa 619 380}%
\special{pa 633 400}%
\special{pa 619 420}%
\special{pa 686 400}%
\special{fp}%
\special{pn 8}%
\special{pa 998 370}%
\special{pa 1000 360}%
\special{fp}%
\special{sh 1}%
\special{pa 1000 360}%
\special{pa 967 421}%
\special{pa 990 412}%
\special{pa 1007 429}%
\special{pa 1000 360}%
\special{fp}%
\special{pn 8}%
\special{pa 202 350}%
\special{pa 200 360}%
\special{fp}%
\special{sh 1}%
\special{pa 200 360}%
\special{pa 233 299}%
\special{pa 210 308}%
\special{pa 193 291}%
\special{pa 200 360}%
\special{fp}%
\special{pn 8}%
\special{pa 360 200}%
\special{pa 360 520}%
\special{dt 0.045}%
\special{pa 360 520}%
\special{pa 360 519}%
\special{dt 0.045}%
\put(6.0000,-2.500){\makebox(0,0){\scriptsize $a_2$}}%
\put(6.0200,-6.00){\makebox(0,0){$K^{[10]}$}}%
\end{picture}%
}

\newcommand{\piccs}{
\unitlength 0.1in
\begin{picture}(9.14,4.25)(1.71,-3.72)
\special{pn 8}%
\special{ar 360 212 160 160  1.5707963 4.7123890}%
\special{pn 8}%
\special{ar 440 212 160 160  0.2510027 1.5707963}%
\special{pn 8}%
\special{ar 440 212 160 160  4.7123890 6.0290606}%
\special{pn 8}%
\special{ar 920 212 160 160  3.3957174 6.2831853}%
\special{ar 920 212 160 160  0.0000000 2.8966140}%
\special{pn 8}%
\special{pa 764 172}%
\special{pa 594 172}%
\special{fp}%
\special{sh 1}%
\special{pa 594 172}%
\special{pa 661 192}%
\special{pa 647 172}%
\special{pa 661 152}%
\special{pa 594 172}%
\special{fp}%
\special{pn 8}%
\special{pa 440 52}%
\special{pa 440 372}%
\special{fp}%
\special{sh 1}%
\special{pa 440 372}%
\special{pa 460 305}%
\special{pa 440 319}%
\special{pa 420 305}%
\special{pa 440 372}%
\special{fp}%
\special{pn 8}%
\special{pa 594 252}%
\special{pa 766 252}%
\special{fp}%
\special{sh 1}%
\special{pa 766 252}%
\special{pa 699 232}%
\special{pa 713 252}%
\special{pa 699 272}%
\special{pa 766 252}%
\special{fp}%
\special{pn 8}%
\special{pa 360 372}%
\special{pa 360 52}%
\special{fp}%
\special{sh 1}%
\special{pa 360 52}%
\special{pa 340 119}%
\special{pa 360 105}%
\special{pa 380 119}%
\special{pa 360 52}%
\special{fp}%
\special{pn 8}%
\special{pa 1078 222}%
\special{pa 1080 212}%
\special{fp}%
\special{sh 1}%
\special{pa 1080 212}%
\special{pa 1047 273}%
\special{pa 1070 264}%
\special{pa 1087 281}%
\special{pa 1080 212}%
\special{fp}%
\special{pn 8}%
\special{pa 202 202}%
\special{pa 200 212}%
\special{fp}%
\special{sh 1}%
\special{pa 200 212}%
\special{pa 233 151}%
\special{pa 210 160}%
\special{pa 193 143}%
\special{pa 200 212}%
\special{fp}%
\put(4.100,-0.100){\makebox(0,0){\scriptsize $a_1$}}%
\put(7.000,-1.200){\makebox(0,0){\scriptsize $a_2$}}%
\put(2.1600,-5.00){\makebox(0,0){$K^{[11]}_+$}}%
\put(7.5000,-5.00){\makebox(0,0){$K^{[11]}_-$}}%
\end{picture}%
}

\newcommand{\picct}{
\unitlength 0.1in
\begin{picture}(8.10,3.20)(1.95,-5.20)
\special{pn 8}%
\special{ar 360 360 160 160  0.0000000 6.2831853}%
\special{pn 8}%
\special{ar 840 360 160 160  0.0000000 6.2831853}%
\special{pn 8}%
\special{pa 360 200}%
\special{pa 360 520}%
\special{dt 0.045}%
\special{pa 360 520}%
\special{pa 360 519}%
\special{dt 0.045}%
\special{pn 8}%
\special{pa 520 360}%
\special{pa 680 360}%
\special{dt 0.045}%
\special{pa 680 360}%
\special{pa 679 360}%
\special{dt 0.045}%
\put(3.6000,-6.500){\makebox(0,0){$K^{[12]}_+$}}%
\put(8.800,-6.500){\makebox(0,0){$K^{[12]}_-$}}%
\special{pn 8}%
\special{pa 998 370}%
\special{pa 1000 360}%
\special{fp}%
\special{sh 1}%
\special{pa 1000 360}%
\special{pa 967 421}%
\special{pa 990 412}%
\special{pa 1007 429}%
\special{pa 1000 360}%
\special{fp}%
\special{pn 8}%
\special{pa 202 350}%
\special{pa 200 360}%
\special{fp}%
\special{sh 1}%
\special{pa 200 360}%
\special{pa 233 299}%
\special{pa 210 308}%
\special{pa 193 291}%
\special{pa 200 360}%
\special{fp}%
\end{picture}%
}

\newcommand{\piccu}{
\unitlength 0.1in
\begin{picture}(8.90,4.25)(1.95,-3.72)
\special{pn 8}%
\special{ar 360 212 160 160  1.5707963 4.7123890}%
\special{pn 8}%
\special{ar 440 212 160 160  4.7123890 6.2831853}%
\special{ar 440 212 160 160  0.0000000 1.5707963}%
\special{pn 8}%
\special{ar 920 212 160 160  0.0000000 6.2831853}%
\special{pn 8}%
\special{pa 360 372}%
\special{pa 360 52}%
\special{fp}%
\special{sh 1}%
\special{pa 360 52}%
\special{pa 340 119}%
\special{pa 360 105}%
\special{pa 380 119}%
\special{pa 360 52}%
\special{fp}%
\special{pn 8}%
\special{pa 440 52}%
\special{pa 440 372}%
\special{fp}%
\special{sh 1}%
\special{pa 440 372}%
\special{pa 460 305}%
\special{pa 440 319}%
\special{pa 420 305}%
\special{pa 440 372}%
\special{fp}%
\special{pn 8}%
\special{pa 600 212}%
\special{pa 760 212}%
\special{dt 0.045}%
\special{pa 760 212}%
\special{pa 759 212}%
\special{dt 0.045}%
\special{pn 8}%
\special{pa 1078 222}%
\special{pa 1080 212}%
\special{fp}%
\special{sh 1}%
\special{pa 1080 212}%
\special{pa 1047 273}%
\special{pa 1070 264}%
\special{pa 1087 281}%
\special{pa 1080 212}%
\special{fp}%
\special{pn 8}%
\special{pa 202 202}%
\special{pa 200 212}%
\special{fp}%
\special{sh 1}%
\special{pa 200 212}%
\special{pa 233 151}%
\special{pa 210 160}%
\special{pa 193 143}%
\special{pa 200 212}%
\special{fp}%
\special{pn 8}%
\special{pa 450 54}%
\special{pa 440 52}%
\special{fp}%
\special{sh 1}%
\special{pa 440 52}%
\special{pa 501 85}%
\special{pa 492 62}%
\special{pa 509 45}%
\special{pa 440 52}%
\special{fp}%
\put(4.0000,-0.100){\makebox(0,0){\scriptsize $a_1$}}%
\put(6.8200,-1.4000){\makebox(0,0){\scriptsize $a_2$}}%
\put(2.4200,-5.00){\makebox(0,0){$K^{[13]}_+$}}%
\put(6.000,-5.00){\makebox(0,0){$K^{[13]}_0$}}%
\put(10.500,-5.000){\makebox(0,0){$K^{[13]}_-$}}%
\end{picture}%
}

\newcommand{\piccw}{
\unitlength 0.1in
\begin{picture}(12.90,3.30)(1.95,-4)
\special{pn 13}%
\special{ar 360 360 160 160  0.0000000 6.2831853}%
\special{pn 13}%
\special{ar 840 360 160 160  0.0000000 6.2831853}%
\special{pn 13}%
\special{ar 1320 360 160 160  0.0000000 6.2831853}%
\special{pn 4}%
\special{pa 520 360}%
\special{pa 680 360}%
\special{fp}%
\special{pn 4}%
\special{pa 1000 360}%
\special{pa 1160 360}%
\special{fp}%
\put(6.0200,-3.000){\makebox(0,0){\scriptsize $a_1$}}%
\put(10.8000,-3.00){\makebox(0,0){\scriptsize $a_2$}}%
\put(8.4200,-6.5000){\makebox(0,0){$\{K_1,K_2,K_3:a_1,a_2\}$}}%
\special{pn 8}%
\special{pa 1478 370}%
\special{pa 1480 360}%
\special{fp}%
\special{sh 1}%
\special{pa 1480 360}%
\special{pa 1447 421}%
\special{pa 1470 412}%
\special{pa 1487 429}%
\special{pa 1480 360}%
\special{fp}%
\special{pn 8}%
\special{pa 850 202}%
\special{pa 840 200}%
\special{fp}%
\special{sh 1}%
\special{pa 840 200}%
\special{pa 901 233}%
\special{pa 892 210}%
\special{pa 909 193}%
\special{pa 840 200}%
\special{fp}%
\special{pn 8}%
\special{pa 202 350}%
\special{pa 200 360}%
\special{fp}%
\special{sh 1}%
\special{pa 200 360}%
\special{pa 233 299}%
\special{pa 210 308}%
\special{pa 193 291}%
\special{pa 200 360}%
\special{fp}%
\special{pn 8}%
\special{pa 830 518}%
\special{pa 840 520}%
\special{fp}%
\special{sh 1}%
\special{pa 840 520}%
\special{pa 779 487}%
\special{pa 788 510}%
\special{pa 771 527}%
\special{pa 840 520}%
\special{fp}%
\end{picture}%
}

\newcommand{\piccx}{
\unitlength 0.1in
\begin{picture}(12.90,4.15)(2.50,-17.25)
\special{pn 8}%
\special{ar 415 1560 160 160  0.2494697 6.0290606}%
\special{pn 8}%
\special{ar 895 1560 160 160  3.3910623 6.0337157}%
\special{pn 8}%
\special{ar 895 1560 160 160  0.2541247 2.8874679}%
\special{pn 8}%
\special{ar 1375 1560 160 160  3.3957174 6.2831853}%
\special{ar 1375 1560 160 160  0.0000000 2.8921230}%
\special{pn 8}%
\special{pa 1221 1520}%
\special{pa 1049 1520}%
\special{fp}%
\special{sh 1}%
\special{pa 1049 1520}%
\special{pa 1116 1540}%
\special{pa 1102 1520}%
\special{pa 1116 1500}%
\special{pa 1049 1520}%
\special{fp}%
\special{pn 8}%
\special{pa 738 1520}%
\special{pa 569 1520}%
\special{fp}%
\special{sh 1}%
\special{pa 569 1520}%
\special{pa 636 1540}%
\special{pa 622 1520}%
\special{pa 636 1500}%
\special{pa 569 1520}%
\special{fp}%
\special{pn 8}%
\special{pa 572 1600}%
\special{pa 738 1600}%
\special{fp}%
\special{sh 1}%
\special{pa 738 1600}%
\special{pa 671 1580}%
\special{pa 685 1600}%
\special{pa 671 1620}%
\special{pa 738 1600}%
\special{fp}%
\special{pn 8}%
\special{pa 1049 1600}%
\special{pa 1218 1600}%
\special{fp}%
\special{sh 1}%
\special{pa 1218 1600}%
\special{pa 1151 1580}%
\special{pa 1165 1600}%
\special{pa 1151 1620}%
\special{pa 1218 1600}%
\special{fp}%
\special{pn 8}%
\special{pa 1533 1570}%
\special{pa 1535 1560}%
\special{fp}%
\special{sh 1}%
\special{pa 1535 1560}%
\special{pa 1502 1621}%
\special{pa 1525 1612}%
\special{pa 1542 1629}%
\special{pa 1535 1560}%
\special{fp}%
\special{pn 8}%
\special{pa 905 1402}%
\special{pa 895 1400}%
\special{fp}%
\special{sh 1}%
\special{pa 895 1400}%
\special{pa 956 1433}%
\special{pa 947 1410}%
\special{pa 964 1393}%
\special{pa 895 1400}%
\special{fp}%
\special{pn 8}%
\special{pa 257 1550}%
\special{pa 255 1560}%
\special{fp}%
\special{sh 1}%
\special{pa 255 1560}%
\special{pa 288 1499}%
\special{pa 265 1508}%
\special{pa 248 1491}%
\special{pa 255 1560}%
\special{fp}%
\special{pn 8}%
\special{pa 885 1718}%
\special{pa 895 1720}%
\special{fp}%
\special{sh 1}%
\special{pa 895 1720}%
\special{pa 834 1687}%
\special{pa 843 1710}%
\special{pa 826 1727}%
\special{pa 895 1720}%
\special{fp}%
\put(6.5500,-14.6000){\makebox(0,0){\scriptsize $a_1$}}%
\put(11.3700,-14.5000){\makebox(0,0){\scriptsize $a_2$}}%
\put(9.200,-18.5000){\makebox(0,0){$K^{[14]}$}}%
\end{picture}%
}

\newcommand{\piccy}{
\unitlength 0.1in
\begin{picture}(12.90,3.30)(1.95,-5.25)
\special{pn 8}%
\special{ar 360 360 160 160  0.0000000 6.2831853}%
\special{pn 8}%
\special{ar 840 360 160 160  0.2541247 6.0290606}%
\special{pn 8}%
\special{ar 1320 360 160 160  3.3957174 6.2831853}%
\special{ar 1320 360 160 160  0.0000000 2.8874679}%
\special{pn 8}%
\special{pa 1166 320}%
\special{pa 994 320}%
\special{fp}%
\special{sh 1}%
\special{pa 994 320}%
\special{pa 1061 340}%
\special{pa 1047 320}%
\special{pa 1061 300}%
\special{pa 994 320}%
\special{fp}%
\special{pn 8}%
\special{pa 994 400}%
\special{pa 1166 400}%
\special{fp}%
\special{sh 1}%
\special{pa 1166 400}%
\special{pa 1099 380}%
\special{pa 1113 400}%
\special{pa 1099 420}%
\special{pa 1166 400}%
\special{fp}%
\special{pn 8}%
\special{pa 520 360}%
\special{pa 680 360}%
\special{dt 0.045}%
\special{pa 680 360}%
\special{pa 679 360}%
\special{dt 0.045}%
\special{pn 8}%
\special{pa 1478 370}%
\special{pa 1480 360}%
\special{fp}%
\special{sh 1}%
\special{pa 1480 360}%
\special{pa 1447 421}%
\special{pa 1470 412}%
\special{pa 1487 429}%
\special{pa 1480 360}%
\special{fp}%
\special{pn 8}%
\special{pa 850 202}%
\special{pa 840 200}%
\special{fp}%
\special{sh 1}%
\special{pa 840 200}%
\special{pa 901 233}%
\special{pa 892 210}%
\special{pa 909 193}%
\special{pa 840 200}%
\special{fp}%
\special{pn 8}%
\special{pa 830 518}%
\special{pa 840 520}%
\special{fp}%
\special{sh 1}%
\special{pa 840 520}%
\special{pa 779 487}%
\special{pa 788 510}%
\special{pa 771 527}%
\special{pa 840 520}%
\special{fp}%
\special{pn 8}%
\special{pa 202 350}%
\special{pa 200 360}%
\special{fp}%
\special{sh 1}%
\special{pa 200 360}%
\special{pa 233 299}%
\special{pa 210 308}%
\special{pa 193 291}%
\special{pa 200 360}%
\special{fp}%
\put(6.0000,-3.2200){\makebox(0,0){}}%
\put(10.80,-2.6000){\makebox(0,0){\scriptsize $a_2$}}%
\put(4.00,-6.3000){\makebox(0,0){$K^{[15]}_+$}}%
\put(11.00,-6.3000){\makebox(0,0){$K^{[15]}_-$}}%
\end{picture}%
}

\newcommand{\piccz}{
\unitlength 0.1in
\begin{picture}(12.90,3.30)(1.95,-5.25)
\special{pn 8}%
\special{ar 360 360 160 160  0.2541247 6.0290606}%
\special{pn 8}%
\special{ar 840 360 160 160  3.3957174 6.2831853}%
\special{ar 840 360 160 160  0.0000000 2.8874679}%
\special{pn 8}%
\special{ar 1320 360 160 160  0.0000000 6.2831853}%
\special{pn 8}%
\special{pa 686 320}%
\special{pa 514 320}%
\special{fp}%
\special{sh 1}%
\special{pa 514 320}%
\special{pa 581 340}%
\special{pa 567 320}%
\special{pa 581 300}%
\special{pa 514 320}%
\special{fp}%
\special{pn 8}%
\special{pa 514 400}%
\special{pa 686 400}%
\special{fp}%
\special{sh 1}%
\special{pa 686 400}%
\special{pa 619 380}%
\special{pa 633 400}%
\special{pa 619 420}%
\special{pa 686 400}%
\special{fp}%
\special{pn 8}%
\special{pa 1000 360}%
\special{pa 1160 360}%
\special{dt 0.045}%
\special{pa 1160 360}%
\special{pa 1159 360}%
\special{dt 0.045}%
\special{pn 8}%
\special{pa 1478 370}%
\special{pa 1480 360}%
\special{dt 0.045}%
\special{sh 1}%
\special{pa 1480 360}%
\special{pa 1447 421}%
\special{pa 1470 412}%
\special{pa 1487 429}%
\special{pa 1480 360}%
\special{fp}%
\special{pn 8}%
\special{pa 850 202}%
\special{pa 840 200}%
\special{dt 0.045}%
\special{sh 1}%
\special{pa 840 200}%
\special{pa 901 233}%
\special{pa 892 210}%
\special{pa 909 193}%
\special{pa 840 200}%
\special{fp}%
\special{pn 8}%
\special{pa 202 350}%
\special{pa 200 360}%
\special{dt 0.045}%
\special{sh 1}%
\special{pa 200 360}%
\special{pa 233 299}%
\special{pa 210 308}%
\special{pa 193 291}%
\special{pa 200 360}%
\special{fp}%
\special{pn 8}%
\special{pa 830 518}%
\special{pa 840 520}%
\special{dt 0.045}%
\special{sh 1}%
\special{pa 840 520}%
\special{pa 779 487}%
\special{pa 788 510}%
\special{pa 771 527}%
\special{pa 840 520}%
\special{fp}%
\put(6.200,-2.8000){\makebox(0,0){\scriptsize $a_1$}}%
\put(10.8000,-3.2400){\makebox(0,0){}}%
\put(6.5000,-6.5000){\makebox(0,0){$K^{[16]}_+$}}%
\put(13.800,-6.500){\makebox(0,0){$K^{[16]}_-$}}%
\end{picture}%
}

\newcommand{\picdb}{
\unitlength 0.1in
\begin{picture}(20.99,5)(0.95,-7)
\special{pn 16}%
\special{ar 595 657 400 400  0.0000000 6.2831853}%
\special{pn 16}%
\special{ar 1794 657 400 400  0.0000000 6.2831853}%
\special{pn 4}%
\special{pa 196 657}%
\special{pa 994 657}%
\special{fp}%
\special{pn 4}%
\special{pa 964 507}%
\special{pa 1424 507}%
\special{fp}%
\special{pn 4}%
\special{pa 860 357}%
\special{pa 1529 357}%
\special{fp}%
\special{pn 4}%
\special{pa 860 962}%
\special{pa 1424 807}%
\special{fp}%
\special{pn 4}%
\special{pa 970 807}%
\special{pa 1534 957}%
\special{fp}%
\put(0.800,-6.5700){\makebox(0,0){\scriptsize $(-)$}}%
\put(16.70,-3.700){\makebox(0,0){\scriptsize $(+)$}}%
\put(15.500,-5.0700){\makebox(0,0){\scriptsize $(+)$}}%
\put(15.400,-8.1200){\makebox(0,0){\scriptsize $(-)$}}%
\put(16.500,-9.400){\makebox(0,0){\scriptsize $(-)$}}%
\put(12.1000,-12.500){}%
\end{picture}%
}

\newcommand{\picdc}{
\unitlength 0.1in
\begin{picture}(6.00,0.00)(4.00,-4.00)
\special{pn 8}%
\special{pa 400 400}%
\special{pa 1000 400}%
\special{fp}%
\end{picture}%
}

\newcommand{\picdd}{
\unitlength 0.1in
\begin{picture}(6.00,1.65)(4.00,-4.00)
\special{pn 8}%
\special{pa 400 400}%
\special{pa 1000 400}%
\special{fp}%
\put(7.0000,-3.4000){\makebox(0,0){\scriptsize$2$}}%
\end{picture}%
}

\newcommand{\picde}{
\unitlength 0.1in
\begin{picture}(4.50,3.38)(1.7,-2)
\special{pn 16}%
\special{ar 320 171 120 120  0.0000000 6.2831853}%
\special{pn 4}%
\special{pa 320 51}%
\special{pa 320 291}%
\special{fp}%
\special{pn 4}%
\special{pa 209 126}%
\special{pa 431 126}%
\special{fp}%
\special{pn 4}%
\special{pa 209 216}%
\special{pa 432 216}%
\special{fp}%
\put(3.3000,0.00){\makebox(0,0){\scriptsize$\epsilon_1$}}%
\put(5.10,-1.300){\makebox(0,0){\scriptsize$\epsilon_2$}}%
\put(5.10,-2.300){\makebox(0,0){\scriptsize$\epsilon_3$}}%
\end{picture}%
}

\newcommand{\picdj}{
\unitlength 0.1in
\begin{picture}(7.0,2.40)(1.50,-3.60)
\special{pn 16}%
\special{ar 320 320 120 120  0.0000000 6.2831853}%
\special{pn 16}%
\special{ar 680 320 120 120  0.0000000 6.2831853}%
\special{pn 4}%
\special{pa 401 230}%
\special{pa 600 230}%
\special{fp}%
\special{pn 4}%
\special{pa 399 410}%
\special{pa 602 410}%
\special{fp}%
\special{pn 4}%
\special{pa 200 320}%
\special{pa 440 320}%
\special{fp}%
\put(5.000,-2.00){\makebox(0,0){\tiny 2}}%
\end{picture}%
}

\newcommand{\picdp}{
\unitlength 0.1in
\begin{picture}(8.5,4.40)(-0.4,-2)
\special{pn 16}%
\special{ar 188 163 120 120  0.0000000 6.2831853}%
\special{pn 16}%
\special{ar 548 163 120 120  0.0000000 6.2831853}%
\special{pn 4}%
\special{pa 268 73}%
\special{pa 469 73}%
\special{fp}%
\special{pn 4}%
\special{pa 268 253}%
\special{pa 469 253}%
\special{fp}%
\special{pn 4}%
\special{pa 77 118}%
\special{pa 299 118}%
\special{fp}%
\special{pn 4}%
\special{pa 77 208}%
\special{pa 299 208}%
\special{fp}%
\special{pn 4}%
\special{pa 428 163}%
\special{pa 668 163}%
\special{fp}%
\special{pn 4}%
\special{pa 578 47}%
\special{pa 578 278}%
\special{fp}%
\special{pn 4}%
\special{pa 158 46}%
\special{pa 158 278}%
\special{fp}%
\put(0.20,-1.1800){\makebox(0,0){\scriptsize$\epsilon_1$}}%
\put(0.200,-2.0800){\makebox(0,0){\scriptsize$\epsilon_2$}}%
\put(1.5800,0.00){\makebox(0,0){\scriptsize$\epsilon_3$}}%
\put(3.80,-0.3){\makebox(0,0){\scriptsize$\epsilon_4$}}%
\put(3.800,-3.00){\makebox(0,0){\scriptsize$\epsilon_5$}}%
\put(7.3000,-1.70){\makebox(0,0){\scriptsize$\epsilon_7$}}%
\put(6.00,-0.100){\makebox(0,0){\scriptsize$\epsilon_6$}}%
\end{picture}%
}

\newcommand{\picea}{
\unitlength 0.1in
\begin{picture}(6.00,1.65)(4.00,-4.00)
\special{pn 8}%
\special{pa 400 400}%
\special{pa 1000 400}%
\special{fp}%
\put(7.0000,-3.4000){\makebox(0,0){\small s}}%
\end{picture}%
}

\newcommand{\piceb}{
\unitlength 0.1in
\begin{picture}(6.00,1.65)(4.00,-4.00)
\special{pn 8}%
\special{pa 400 400}%
\special{pa 1000 400}%
\special{fp}%
\put(7.0000,-3.4000){\makebox(0,0){\small d}}%
\end{picture}%
}

\newcommand{\picec}{
\unitlength 0.1in
\begin{picture}(6.00,1.65)(4.00,-4.00)
\special{pn 8}%
\special{pa 400 400}%
\special{pa 1000 400}%
\special{fp}%
\put(7.0000,-3.4000){\makebox(0,0){\small$2$s}}%
\end{picture}%
}

\newcommand{\piced}{
\unitlength 0.1in
\begin{picture}(6.00,1.65)(4.00,-4.00)
\special{pn 8}%
\special{pa 400 400}%
\special{pa 1000 400}%
\special{fp}%
\put(7.0000,-3.4000){\makebox(0,0){\small$2$d}}%
\end{picture}%
}

\newcommand{\picee}{
\unitlength 0.1in
\begin{picture}(6.00,2.00)(4.00,-4.500)
\special{pn 8}%
\special{pa 400 400}%
\special{pa 1000 400}%
\special{dt 0.045}%
\special{pa 1000 400}%
\special{pa 999 400}%
\special{dt 0.045}%
\end{picture}%
}

\newcommand{\picef}{
\unitlength 0.1in
\begin{picture}(3.20,3.20)(2.00,-5.20)
\special{pn 8}%
\special{ar 360 360 160 160  0.0000000 0.0750000}%
\special{ar 360 360 160 160  0.3000000 0.3750000}%
\special{ar 360 360 160 160  0.6000000 0.6750000}%
\special{ar 360 360 160 160  0.9000000 0.9750000}%
\special{ar 360 360 160 160  1.2000000 1.2750000}%
\special{ar 360 360 160 160  1.5000000 1.5750000}%
\special{ar 360 360 160 160  1.8000000 1.8750000}%
\special{ar 360 360 160 160  2.1000000 2.1750000}%
\special{ar 360 360 160 160  2.4000000 2.4750000}%
\special{ar 360 360 160 160  2.7000000 2.7750000}%
\special{ar 360 360 160 160  3.0000000 3.0750000}%
\special{ar 360 360 160 160  3.3000000 3.3750000}%
\special{ar 360 360 160 160  3.6000000 3.6750000}%
\special{ar 360 360 160 160  3.9000000 3.9750000}%
\special{ar 360 360 160 160  4.2000000 4.2750000}%
\special{ar 360 360 160 160  4.5000000 4.5750000}%
\special{ar 360 360 160 160  4.8000000 4.8750000}%
\special{ar 360 360 160 160  5.1000000 5.1750000}%
\special{ar 360 360 160 160  5.4000000 5.4750000}%
\special{ar 360 360 160 160  5.7000000 5.7750000}%
\special{ar 360 360 160 160  6.0000000 6.0750000}%
\special{pn 8}%
\special{pa 248 472}%
\special{pa 472 248}%
\special{fp}%
\special{sh 1}%
\special{pa 472 248}%
\special{pa 411 281}%
\special{pa 434 286}%
\special{pa 439 309}%
\special{pa 472 248}%
\special{fp}%
\special{pn 8}%
\special{pa 330 330}%
\special{pa 248 248}%
\special{fp}%
\special{sh 1}%
\special{pa 248 248}%
\special{pa 281 309}%
\special{pa 286 286}%
\special{pa 309 281}%
\special{pa 248 248}%
\special{fp}%
\special{pn 8}%
\special{pa 472 472}%
\special{pa 390 390}%
\special{fp}%
\put(3.5000,-6.0000){\makebox(0,0){($+$,s)}}%
\end{picture}%
}

\newcommand{\piceg}{
\unitlength 0.1in
\begin{picture}(3.20,3.20)(2.00,-5.20)
\special{pn 8}%
\special{ar 360 360 160 160  0.0000000 0.0750000}%
\special{ar 360 360 160 160  0.3000000 0.3750000}%
\special{ar 360 360 160 160  0.6000000 0.6750000}%
\special{ar 360 360 160 160  0.9000000 0.9750000}%
\special{ar 360 360 160 160  1.2000000 1.2750000}%
\special{ar 360 360 160 160  1.5000000 1.5750000}%
\special{ar 360 360 160 160  1.8000000 1.8750000}%
\special{ar 360 360 160 160  2.1000000 2.1750000}%
\special{ar 360 360 160 160  2.4000000 2.4750000}%
\special{ar 360 360 160 160  2.7000000 2.7750000}%
\special{ar 360 360 160 160  3.0000000 3.0750000}%
\special{ar 360 360 160 160  3.3000000 3.3750000}%
\special{ar 360 360 160 160  3.6000000 3.6750000}%
\special{ar 360 360 160 160  3.9000000 3.9750000}%
\special{ar 360 360 160 160  4.2000000 4.2750000}%
\special{ar 360 360 160 160  4.5000000 4.5750000}%
\special{ar 360 360 160 160  4.8000000 4.8750000}%
\special{ar 360 360 160 160  5.1000000 5.1750000}%
\special{ar 360 360 160 160  5.4000000 5.4750000}%
\special{ar 360 360 160 160  5.7000000 5.7750000}%
\special{ar 360 360 160 160  6.0000000 6.0750000}%
\special{pn 8}%
\special{pa 472 248}%
\special{pa 248 472}%
\special{fp}%
\special{sh 1}%
\special{pa 248 472}%
\special{pa 309 439}%
\special{pa 286 434}%
\special{pa 281 411}%
\special{pa 248 472}%
\special{fp}%
\special{pn 8}%
\special{pa 248 248}%
\special{pa 330 330}%
\special{fp}%
\special{pn 8}%
\special{pa 390 390}%
\special{pa 472 472}%
\special{fp}%
\special{sh 1}%
\special{pa 472 472}%
\special{pa 439 411}%
\special{pa 434 434}%
\special{pa 411 439}%
\special{pa 472 472}%
\special{fp}%
\put(3.5000,-6.0000){\makebox(0,0){($+$,s)}}%
\end{picture}%
}

\newcommand{\piceh}{
\unitlength 0.1in
\begin{picture}(3.20,3.20)(2.00,-5.20)
\special{pn 8}%
\special{ar 360 360 160 160  0.0000000 0.0750000}%
\special{ar 360 360 160 160  0.3000000 0.3750000}%
\special{ar 360 360 160 160  0.6000000 0.6750000}%
\special{ar 360 360 160 160  0.9000000 0.9750000}%
\special{ar 360 360 160 160  1.2000000 1.2750000}%
\special{ar 360 360 160 160  1.5000000 1.5750000}%
\special{ar 360 360 160 160  1.8000000 1.8750000}%
\special{ar 360 360 160 160  2.1000000 2.1750000}%
\special{ar 360 360 160 160  2.4000000 2.4750000}%
\special{ar 360 360 160 160  2.7000000 2.7750000}%
\special{ar 360 360 160 160  3.0000000 3.0750000}%
\special{ar 360 360 160 160  3.3000000 3.3750000}%
\special{ar 360 360 160 160  3.6000000 3.6750000}%
\special{ar 360 360 160 160  3.9000000 3.9750000}%
\special{ar 360 360 160 160  4.2000000 4.2750000}%
\special{ar 360 360 160 160  4.5000000 4.5750000}%
\special{ar 360 360 160 160  4.8000000 4.8750000}%
\special{ar 360 360 160 160  5.1000000 5.1750000}%
\special{ar 360 360 160 160  5.4000000 5.4750000}%
\special{ar 360 360 160 160  5.7000000 5.7750000}%
\special{ar 360 360 160 160  6.0000000 6.0750000}%
\special{pn 8}%
\special{pa 248 248}%
\special{pa 472 472}%
\special{fp}%
\special{sh 1}%
\special{pa 472 472}%
\special{pa 439 411}%
\special{pa 434 434}%
\special{pa 411 439}%
\special{pa 472 472}%
\special{fp}%
\special{pn 8}%
\special{pa 248 472}%
\special{pa 330 390}%
\special{fp}%
\special{pn 8}%
\special{pa 390 330}%
\special{pa 472 248}%
\special{fp}%
\special{sh 1}%
\special{pa 472 248}%
\special{pa 411 281}%
\special{pa 434 286}%
\special{pa 439 309}%
\special{pa 472 248}%
\special{fp}%
\put(3.6000,-6.0000){\makebox(0,0){($+$,d)}}%
\end{picture}%
}

\newcommand{\picei}{
\unitlength 0.1in
\begin{picture}(3.20,3.20)(2.00,-5.20)
\special{pn 8}%
\special{ar 360 360 160 160  0.0000000 0.0750000}%
\special{ar 360 360 160 160  0.3000000 0.3750000}%
\special{ar 360 360 160 160  0.6000000 0.6750000}%
\special{ar 360 360 160 160  0.9000000 0.9750000}%
\special{ar 360 360 160 160  1.2000000 1.2750000}%
\special{ar 360 360 160 160  1.5000000 1.5750000}%
\special{ar 360 360 160 160  1.8000000 1.8750000}%
\special{ar 360 360 160 160  2.1000000 2.1750000}%
\special{ar 360 360 160 160  2.4000000 2.4750000}%
\special{ar 360 360 160 160  2.7000000 2.7750000}%
\special{ar 360 360 160 160  3.0000000 3.0750000}%
\special{ar 360 360 160 160  3.3000000 3.3750000}%
\special{ar 360 360 160 160  3.6000000 3.6750000}%
\special{ar 360 360 160 160  3.9000000 3.9750000}%
\special{ar 360 360 160 160  4.2000000 4.2750000}%
\special{ar 360 360 160 160  4.5000000 4.5750000}%
\special{ar 360 360 160 160  4.8000000 4.8750000}%
\special{ar 360 360 160 160  5.1000000 5.1750000}%
\special{ar 360 360 160 160  5.4000000 5.4750000}%
\special{ar 360 360 160 160  5.7000000 5.7750000}%
\special{ar 360 360 160 160  6.0000000 6.0750000}%
\special{pn 8}%
\special{pa 472 472}%
\special{pa 248 248}%
\special{fp}%
\special{sh 1}%
\special{pa 248 248}%
\special{pa 281 309}%
\special{pa 286 286}%
\special{pa 309 281}%
\special{pa 248 248}%
\special{fp}%
\special{pn 8}%
\special{pa 472 248}%
\special{pa 390 330}%
\special{fp}%
\special{pn 8}%
\special{pa 330 390}%
\special{pa 248 472}%
\special{fp}%
\special{sh 1}%
\special{pa 248 472}%
\special{pa 309 439}%
\special{pa 286 434}%
\special{pa 281 411}%
\special{pa 248 472}%
\special{fp}%
\put(3.6000,-6.0000){\makebox(0,0){($+$,d)}}%
\end{picture}%
}

\newcommand{\picej}{
\unitlength 0.1in
\begin{picture}(3.20,3.20)(2.00,-5.20)
\special{pn 8}%
\special{ar 360 360 160 160  0.0000000 0.0750000}%
\special{ar 360 360 160 160  0.3000000 0.3750000}%
\special{ar 360 360 160 160  0.6000000 0.6750000}%
\special{ar 360 360 160 160  0.9000000 0.9750000}%
\special{ar 360 360 160 160  1.2000000 1.2750000}%
\special{ar 360 360 160 160  1.5000000 1.5750000}%
\special{ar 360 360 160 160  1.8000000 1.8750000}%
\special{ar 360 360 160 160  2.1000000 2.1750000}%
\special{ar 360 360 160 160  2.4000000 2.4750000}%
\special{ar 360 360 160 160  2.7000000 2.7750000}%
\special{ar 360 360 160 160  3.0000000 3.0750000}%
\special{ar 360 360 160 160  3.3000000 3.3750000}%
\special{ar 360 360 160 160  3.6000000 3.6750000}%
\special{ar 360 360 160 160  3.9000000 3.9750000}%
\special{ar 360 360 160 160  4.2000000 4.2750000}%
\special{ar 360 360 160 160  4.5000000 4.5750000}%
\special{ar 360 360 160 160  4.8000000 4.8750000}%
\special{ar 360 360 160 160  5.1000000 5.1750000}%
\special{ar 360 360 160 160  5.4000000 5.4750000}%
\special{ar 360 360 160 160  5.7000000 5.7750000}%
\special{ar 360 360 160 160  6.0000000 6.0750000}%
\special{pn 8}%
\special{pa 472 472}%
\special{pa 248 248}%
\special{fp}%
\special{sh 1}%
\special{pa 248 248}%
\special{pa 281 309}%
\special{pa 286 286}%
\special{pa 309 281}%
\special{pa 248 248}%
\special{fp}%
\special{pn 8}%
\special{pa 248 472}%
\special{pa 330 390}%
\special{fp}%
\special{pn 8}%
\special{pa 390 330}%
\special{pa 472 248}%
\special{fp}%
\special{sh 1}%
\special{pa 472 248}%
\special{pa 411 281}%
\special{pa 434 286}%
\special{pa 439 309}%
\special{pa 472 248}%
\special{fp}%
\put(3.5000,-6.0000){\makebox(0,0){($-$,s)}}%
\end{picture}%
}

\newcommand{\picek}{
\unitlength 0.1in
\begin{picture}(3.20,3.20)(2.00,-5.20)
\special{pn 8}%
\special{ar 360 360 160 160  0.0000000 0.0750000}%
\special{ar 360 360 160 160  0.3000000 0.3750000}%
\special{ar 360 360 160 160  0.6000000 0.6750000}%
\special{ar 360 360 160 160  0.9000000 0.9750000}%
\special{ar 360 360 160 160  1.2000000 1.2750000}%
\special{ar 360 360 160 160  1.5000000 1.5750000}%
\special{ar 360 360 160 160  1.8000000 1.8750000}%
\special{ar 360 360 160 160  2.1000000 2.1750000}%
\special{ar 360 360 160 160  2.4000000 2.4750000}%
\special{ar 360 360 160 160  2.7000000 2.7750000}%
\special{ar 360 360 160 160  3.0000000 3.0750000}%
\special{ar 360 360 160 160  3.3000000 3.3750000}%
\special{ar 360 360 160 160  3.6000000 3.6750000}%
\special{ar 360 360 160 160  3.9000000 3.9750000}%
\special{ar 360 360 160 160  4.2000000 4.2750000}%
\special{ar 360 360 160 160  4.5000000 4.5750000}%
\special{ar 360 360 160 160  4.8000000 4.8750000}%
\special{ar 360 360 160 160  5.1000000 5.1750000}%
\special{ar 360 360 160 160  5.4000000 5.4750000}%
\special{ar 360 360 160 160  5.7000000 5.7750000}%
\special{ar 360 360 160 160  6.0000000 6.0750000}%
\special{pn 8}%
\special{pa 248 248}%
\special{pa 472 472}%
\special{fp}%
\special{sh 1}%
\special{pa 472 472}%
\special{pa 439 411}%
\special{pa 434 434}%
\special{pa 411 439}%
\special{pa 472 472}%
\special{fp}%
\special{pn 8}%
\special{pa 472 248}%
\special{pa 390 330}%
\special{fp}%
\special{pn 8}%
\special{pa 330 390}%
\special{pa 248 472}%
\special{fp}%
\special{sh 1}%
\special{pa 248 472}%
\special{pa 309 439}%
\special{pa 286 434}%
\special{pa 281 411}%
\special{pa 248 472}%
\special{fp}%
\put(3.5000,-6.0000){\makebox(0,0){($-$,s)}}%
\end{picture}%
}

\newcommand{\picel}{
\unitlength 0.1in
\begin{picture}(3.20,3.20)(2.00,-5.20)
\special{pn 8}%
\special{ar 360 360 160 160  0.0000000 0.0750000}%
\special{ar 360 360 160 160  0.3000000 0.3750000}%
\special{ar 360 360 160 160  0.6000000 0.6750000}%
\special{ar 360 360 160 160  0.9000000 0.9750000}%
\special{ar 360 360 160 160  1.2000000 1.2750000}%
\special{ar 360 360 160 160  1.5000000 1.5750000}%
\special{ar 360 360 160 160  1.8000000 1.8750000}%
\special{ar 360 360 160 160  2.1000000 2.1750000}%
\special{ar 360 360 160 160  2.4000000 2.4750000}%
\special{ar 360 360 160 160  2.7000000 2.7750000}%
\special{ar 360 360 160 160  3.0000000 3.0750000}%
\special{ar 360 360 160 160  3.3000000 3.3750000}%
\special{ar 360 360 160 160  3.6000000 3.6750000}%
\special{ar 360 360 160 160  3.9000000 3.9750000}%
\special{ar 360 360 160 160  4.2000000 4.2750000}%
\special{ar 360 360 160 160  4.5000000 4.5750000}%
\special{ar 360 360 160 160  4.8000000 4.8750000}%
\special{ar 360 360 160 160  5.1000000 5.1750000}%
\special{ar 360 360 160 160  5.4000000 5.4750000}%
\special{ar 360 360 160 160  5.7000000 5.7750000}%
\special{ar 360 360 160 160  6.0000000 6.0750000}%
\special{pn 8}%
\special{pa 248 472}%
\special{pa 472 248}%
\special{fp}%
\special{sh 1}%
\special{pa 472 248}%
\special{pa 411 281}%
\special{pa 434 286}%
\special{pa 439 309}%
\special{pa 472 248}%
\special{fp}%
\special{pn 8}%
\special{pa 248 248}%
\special{pa 330 330}%
\special{fp}%
\special{pn 8}%
\special{pa 390 390}%
\special{pa 472 472}%
\special{fp}%
\special{sh 1}%
\special{pa 472 472}%
\special{pa 439 411}%
\special{pa 434 434}%
\special{pa 411 439}%
\special{pa 472 472}%
\special{fp}%
\put(3.5000,-6.0000){\makebox(0,0){($-$,d)}}%
\end{picture}%
}

\newcommand{\picem}{
\unitlength 0.1in
\begin{picture}(3.20,3.20)(2.00,-5.20)
\special{pn 8}%
\special{ar 360 360 160 160  0.0000000 0.0750000}%
\special{ar 360 360 160 160  0.3000000 0.3750000}%
\special{ar 360 360 160 160  0.6000000 0.6750000}%
\special{ar 360 360 160 160  0.9000000 0.9750000}%
\special{ar 360 360 160 160  1.2000000 1.2750000}%
\special{ar 360 360 160 160  1.5000000 1.5750000}%
\special{ar 360 360 160 160  1.8000000 1.8750000}%
\special{ar 360 360 160 160  2.1000000 2.1750000}%
\special{ar 360 360 160 160  2.4000000 2.4750000}%
\special{ar 360 360 160 160  2.7000000 2.7750000}%
\special{ar 360 360 160 160  3.0000000 3.0750000}%
\special{ar 360 360 160 160  3.3000000 3.3750000}%
\special{ar 360 360 160 160  3.6000000 3.6750000}%
\special{ar 360 360 160 160  3.9000000 3.9750000}%
\special{ar 360 360 160 160  4.2000000 4.2750000}%
\special{ar 360 360 160 160  4.5000000 4.5750000}%
\special{ar 360 360 160 160  4.8000000 4.8750000}%
\special{ar 360 360 160 160  5.1000000 5.1750000}%
\special{ar 360 360 160 160  5.4000000 5.4750000}%
\special{ar 360 360 160 160  5.7000000 5.7750000}%
\special{ar 360 360 160 160  6.0000000 6.0750000}%
\special{pn 8}%
\special{pa 472 248}%
\special{pa 248 472}%
\special{fp}%
\special{sh 1}%
\special{pa 248 472}%
\special{pa 309 439}%
\special{pa 286 434}%
\special{pa 281 411}%
\special{pa 248 472}%
\special{fp}%
\special{pn 8}%
\special{pa 472 472}%
\special{pa 390 390}%
\special{fp}%
\special{pn 8}%
\special{pa 330 330}%
\special{pa 248 248}%
\special{fp}%
\special{sh 1}%
\special{pa 248 248}%
\special{pa 281 309}%
\special{pa 286 286}%
\special{pa 309 281}%
\special{pa 248 248}%
\special{fp}%
\put(3.5000,-6.0000){\makebox(0,0){($-$,d)}}%
\end{picture}%
}

\newcommand{\picen}{
\unitlength 0.1in
\begin{picture}(1.50,6)(2.55,-8.00)
\special{pn 8}%
\special{pa 400 1000}%
\special{pa 400 200}%
\special{fp}%
\special{sh 1}%
\special{pa 400 200}%
\special{pa 380 267}%
\special{pa 400 253}%
\special{pa 420 267}%
\special{pa 400 200}%
\special{fp}%
\put(3.0000,-2.0000){\makebox(0,0){t}}%
\end{picture}%
}

\newcommand{\piceo}{
\unitlength 0.1in
\begin{picture}(6.40,8.51)(2.00,-12)
\special{pn 8}%
\special{ar 520 409 320 320  0.0000000 6.2831853}%
\special{pn 8}%
\special{pa 520 89}%
\special{pa 520 729}%
\special{fp}%
\put(5.2000,-0.00){\makebox(0,0){($-$,d)}}%
\special{pn 8}%
\special{pa 244 249}%
\special{pa 796 569}%
\special{fp}%
\special{pn 8}%
\special{pa 796 249}%
\special{pa 244 569}%
\special{fp}%
\put(10.00,-2.200){\makebox(0,0){($-$,d)}}%
\put(10.00,-5.8500){\makebox(0,0){($-$,s)}}%
\put(4.600,-6.7700){\makebox(0,0){c}}%
\put(2.6400,-4.9000){\makebox(0,0){b}}%
\put(3.200,-2.300){\makebox(0,0){a}}%
\put(5.2000,-9.0000){\makebox(0,0){$G^e(K)$}}%
\end{picture}%
}

\newcommand{\picfc}{
\unitlength 0.1in
\begin{picture}(2.50,3.00)(8.00,-14.00)
\special{pn 8}%
\special{pa 1050 650}%
\special{pa 900 500}%
\special{fp}%
\special{pn 8}%
\special{pa 1050 650}%
\special{pa 900 800}%
\special{fp}%
\special{pn 8}%
\special{pa 800 600}%
\special{pa 1000 600}%
\special{fp}%
\special{pn 8}%
\special{pa 800 700}%
\special{pa 1000 700}%
\special{fp}%
\end{picture}%
}

\newcommand{\picfd}{
\unitlength 0.1in
\begin{picture}(7.25,7.24)(0.28,-8.78)
\special{pn 8}%
\special{ar 433 403 320 320  0.0000000 6.2831853}%
\special{pn 8}%
\special{pa 117 443}%
\special{pa 749 443}%
\special{fp}%
\special{pn 8}%
\special{pa 353 95}%
\special{pa 353 711}%
\special{fp}%
\special{pn 8}%
\special{pa 513 99}%
\special{pa 513 711}%
\special{fp}%
\special{pn 8}%
\special{pa 121 323}%
\special{pa 681 603}%
\special{fp}%
\special{pn 8}%
\special{pa 745 323}%
\special{pa 181 603}%
\special{fp}%
\special{pn 8}%
\special{pa 185 203}%
\special{pa 685 203}%
\special{fp}%
\put(3.500,-7.900){\makebox(0,0){\footnotesize a1}}%
\put(5.400,-7.900){\makebox(0,0){\footnotesize a2}}%
\put(7.500,-6.3500){\makebox(0,0){\footnotesize a3}}%
\put(8.300,-4.4300){\makebox(0,0){\footnotesize a4}}%
\put(8.300,-2.9100){\makebox(0,0){\footnotesize a5}}%
\put(7.700,-1.800){\makebox(0,0){\footnotesize a6}}%
\put(5.1300,-0.3500){\makebox(0,0){$+$}}%
\put(3.5300,-0.3500){\makebox(0,0){$+$}}%
\put(1.100,-2.0300){\makebox(0,0){$-$}}%
\put(0.500,-3.0700){\makebox(0,0){$-$}}%
\put(0.400,-4.4300){\makebox(0,0){$-$}}%
\put(1.00,-6.3100){\makebox(0,0){$-$}}%
\put(4.3700,-9.6300){\makebox(0,0){$G(K)$}}%
\end{picture}%
}

\newcommand{\picfe}{
\unitlength 0.1in
\begin{picture}(7.25,7.24)(0.28,-8.78)
\special{pn 8}%
\special{ar 433 403 320 320  0.0000000 6.2831853}%
\special{pn 8}%
\special{pa 117 443}%
\special{pa 749 443}%
\special{fp}%
\special{pn 8}%
\special{pa 353 95}%
\special{pa 353 711}%
\special{fp}%
\special{pn 8}%
\special{pa 513 99}%
\special{pa 513 711}%
\special{fp}%
\special{pn 8}%
\special{pa 121 323}%
\special{pa 681 603}%
\special{fp}%
\special{pn 8}%
\special{pa 745 323}%
\special{pa 181 603}%
\special{fp}%
\special{pn 8}%
\special{pa 185 203}%
\special{pa 685 203}%
\special{fp}%
\put(3.500,-7.900){\makebox(0,0){\footnotesize a1}}%
\put(5.400,-7.900){\makebox(0,0){\footnotesize a2}}%
\put(7.500,-6.3500){\makebox(0,0){\footnotesize a3}}%
\put(8.300,-4.4300){\makebox(0,0){\footnotesize a4}}%
\put(8.300,-2.9100){\makebox(0,0){\footnotesize a5}}%
\put(7.700,-1.800){\makebox(0,0){\footnotesize a6}}%
\put(5.1300,-0.3500){\makebox(0,0){$-$}}%
\put(3.5300,-0.3500){\makebox(0,0){$+$}}%
\put(1.100,-2.0300){\makebox(0,0){$-$}}%
\put(0.500,-3.0700){\makebox(0,0){$-$}}%
\put(0.400,-4.4300){\makebox(0,0){$+$}}%
\put(1.00,-6.3100){\makebox(0,0){$-$}}%
\put(4.3700,-9.6300){\makebox(0,0){$G(\alpha(K))$}}%
\end{picture}%
}

\newcommand{\picff}{
\unitlength 0.1in
\begin{picture}(3.20,3.20)(2.00,-5.20)
\special{pn 8}%
\special{ar 360 360 160 160  0.0000000 0.0750000}%
\special{ar 360 360 160 160  0.3000000 0.3750000}%
\special{ar 360 360 160 160  0.6000000 0.6750000}%
\special{ar 360 360 160 160  0.9000000 0.9750000}%
\special{ar 360 360 160 160  1.2000000 1.2750000}%
\special{ar 360 360 160 160  1.5000000 1.5750000}%
\special{ar 360 360 160 160  1.8000000 1.8750000}%
\special{ar 360 360 160 160  2.1000000 2.1750000}%
\special{ar 360 360 160 160  2.4000000 2.4750000}%
\special{ar 360 360 160 160  2.7000000 2.7750000}%
\special{ar 360 360 160 160  3.0000000 3.0750000}%
\special{ar 360 360 160 160  3.3000000 3.3750000}%
\special{ar 360 360 160 160  3.6000000 3.6750000}%
\special{ar 360 360 160 160  3.9000000 3.9750000}%
\special{ar 360 360 160 160  4.2000000 4.2750000}%
\special{ar 360 360 160 160  4.5000000 4.5750000}%
\special{ar 360 360 160 160  4.8000000 4.8750000}%
\special{ar 360 360 160 160  5.1000000 5.1750000}%
\special{ar 360 360 160 160  5.4000000 5.4750000}%
\special{ar 360 360 160 160  5.7000000 5.7750000}%
\special{ar 360 360 160 160  6.0000000 6.0750000}%
\special{pn 8}%
\special{pa 248 472}%
\special{pa 472 248}%
\special{fp}%
\special{sh 1}%
\special{pa 472 248}%
\special{pa 411 281}%
\special{pa 434 286}%
\special{pa 439 309}%
\special{pa 472 248}%
\special{fp}%
\special{pn 8}%
\special{pa 330 330}%
\special{pa 248 248}%
\special{fp}%
\special{sh 1}%
\special{pa 248 248}%
\special{pa 281 309}%
\special{pa 286 286}%
\special{pa 309 281}%
\special{pa 248 248}%
\special{fp}%
\special{pn 8}%
\special{pa 472 472}%
\special{pa 390 390}%
\special{fp}%
\put(3.5000,-6.0000){\makebox(0,0){$+1$}}%
\end{picture}%
}

\newcommand{\picfg}{
\unitlength 0.1in
\begin{picture}(3.20,3.20)(2.00,-5.20)
\special{pn 8}%
\special{ar 360 360 160 160  0.0000000 0.0750000}%
\special{ar 360 360 160 160  0.3000000 0.3750000}%
\special{ar 360 360 160 160  0.6000000 0.6750000}%
\special{ar 360 360 160 160  0.9000000 0.9750000}%
\special{ar 360 360 160 160  1.2000000 1.2750000}%
\special{ar 360 360 160 160  1.5000000 1.5750000}%
\special{ar 360 360 160 160  1.8000000 1.8750000}%
\special{ar 360 360 160 160  2.1000000 2.1750000}%
\special{ar 360 360 160 160  2.4000000 2.4750000}%
\special{ar 360 360 160 160  2.7000000 2.7750000}%
\special{ar 360 360 160 160  3.0000000 3.0750000}%
\special{ar 360 360 160 160  3.3000000 3.3750000}%
\special{ar 360 360 160 160  3.6000000 3.6750000}%
\special{ar 360 360 160 160  3.9000000 3.9750000}%
\special{ar 360 360 160 160  4.2000000 4.2750000}%
\special{ar 360 360 160 160  4.5000000 4.5750000}%
\special{ar 360 360 160 160  4.8000000 4.8750000}%
\special{ar 360 360 160 160  5.1000000 5.1750000}%
\special{ar 360 360 160 160  5.4000000 5.4750000}%
\special{ar 360 360 160 160  5.7000000 5.7750000}%
\special{ar 360 360 160 160  6.0000000 6.0750000}%
\special{pn 8}%
\special{pa 472 472}%
\special{pa 248 248}%
\special{fp}%
\special{sh 1}%
\special{pa 248 248}%
\special{pa 281 309}%
\special{pa 286 286}%
\special{pa 309 281}%
\special{pa 248 248}%
\special{fp}%
\special{pn 8}%
\special{pa 248 472}%
\special{pa 330 390}%
\special{fp}%
\special{pn 8}%
\special{pa 390 330}%
\special{pa 472 248}%
\special{fp}%
\special{sh 1}%
\special{pa 472 248}%
\special{pa 411 281}%
\special{pa 434 286}%
\special{pa 439 309}%
\special{pa 472 248}%
\special{fp}%
\put(3.5000,-6.0000){\makebox(0,0){$-1$}}%
\end{picture}%
}

\newcommand{\picfh}{
\unitlength 0.1in
\begin{picture}(8.96,10.50)(8.60,-14.50)
\special{pn 13}%
\special{ar 958 680 700 490  6.2831853 6.2831853}%
\special{ar 958 680 700 490  0.0000000 1.5707963}%
\special{pn 4}%
\special{ar 1308 680 350 140  0.7954988 6.2831853}%
\special{ar 1308 680 350 140  0.0000000 0.7853982}%
\special{pn 4}%
\special{pa 958 680}%
\special{pa 958 1170}%
\special{fp}%
\special{pn 4}%
\special{pa 1658 680}%
\special{pa 1658 1170}%
\special{fp}%
\special{pn 13}%
\special{ar 958 1170 700 490  4.7123890 4.7325570}%
\special{ar 958 1170 700 490  4.7930612 4.8132293}%
\special{ar 958 1170 700 490  4.8737335 4.8939016}%
\special{ar 958 1170 700 490  4.9544058 4.9745739}%
\special{ar 958 1170 700 490  5.0350781 5.0552461}%
\special{ar 958 1170 700 490  5.1157503 5.1359184}%
\special{ar 958 1170 700 490  5.1964226 5.2165907}%
\special{ar 958 1170 700 490  5.2770949 5.2972629}%
\special{ar 958 1170 700 490  5.3577671 5.3779352}%
\special{ar 958 1170 700 490  5.4384394 5.4586075}%
\special{ar 958 1170 700 490  5.5191117 5.5392797}%
\special{ar 958 1170 700 490  5.5997839 5.6199520}%
\special{ar 958 1170 700 490  5.6804562 5.7006243}%
\special{ar 958 1170 700 490  5.7611285 5.7812965}%
\special{ar 958 1170 700 490  5.8418007 5.8619688}%
\special{ar 958 1170 700 490  5.9224730 5.9426411}%
\special{ar 958 1170 700 490  6.0031453 6.0233134}%
\special{ar 958 1170 700 490  6.0838176 6.1039856}%
\special{ar 958 1170 700 490  6.1644898 6.1846579}%
\special{ar 958 1170 700 490  6.2451621 6.2653302}%
\special{pn 8}%
\special{ar 1308 1170 350 140  3.1415927 3.1905722}%
\special{ar 1308 1170 350 140  3.3375110 3.3864906}%
\special{ar 1308 1170 350 140  3.5334294 3.5824090}%
\special{ar 1308 1170 350 140  3.7293478 3.7783273}%
\special{ar 1308 1170 350 140  3.9252661 3.9742457}%
\special{ar 1308 1170 350 140  4.1211845 4.1701641}%
\special{ar 1308 1170 350 140  4.3171029 4.3660824}%
\special{ar 1308 1170 350 140  4.5130212 4.5620008}%
\special{ar 1308 1170 350 140  4.7089396 4.7579192}%
\special{ar 1308 1170 350 140  4.9048580 4.9538376}%
\special{ar 1308 1170 350 140  5.1007763 5.1497559}%
\special{ar 1308 1170 350 140  5.2966947 5.3456743}%
\special{ar 1308 1170 350 140  5.4926131 5.5415927}%
\special{ar 1308 1170 350 140  5.6885314 5.7375110}%
\special{ar 1308 1170 350 140  5.8844498 5.9334294}%
\special{ar 1308 1170 350 140  6.0803682 6.1293478}%
\special{ar 1308 1170 350 140  6.2762865 6.2831853}%
\special{pn 4}%
\special{ar 1308 1170 350 140  6.2831853 6.2831853}%
\special{ar 1308 1170 350 140  0.0000000 3.1415927}%
\special{pn 13}%
\special{ar 1308 400 448 448  6.2831853 6.2831853}%
\special{ar 1308 400 448 448  0.0000000 0.6747409}%
\special{pn 13}%
\special{ar 1308 400 448 448  2.4668517 3.1415927}%
\special{pn 13}%
\special{ar 1308 1450 448 448  5.6084444 6.2831853}%
\special{pn 13}%
\special{ar 1308 1450 448 448  3.1415927 3.8163336}%
\end{picture}%
}

\newcommand{\picfi}{
\unitlength 0.1in
\begin{picture}(8.96,10.50)(8.60,-12.50)
\special{pn 8}%
\special{pa 958 480}%
\special{pa 958 970}%
\special{dt 0.04}%
\special{pa 958 970}%
\special{pa 958 969}%
\special{dt 0.04}%
\special{pn 8}%
\special{pa 1658 480}%
\special{pa 1658 970}%
\special{dt 0.04}%
\special{pa 1658 970}%
\special{pa 1658 969}%
\special{dt 0.04}%
\special{pn 8}%
\special{pa 958 480}%
\special{pa 1658 480}%
\special{dt 0.04}%
\special{pa 1658 480}%
\special{pa 1657 480}%
\special{dt 0.04}%
\special{pn 8}%
\special{pa 958 970}%
\special{pa 1658 970}%
\special{dt 0.04}%
\special{pa 1658 970}%
\special{pa 1657 970}%
\special{dt 0.04}%
\special{pn 13}%
\special{pa 958 970}%
\special{pa 1658 480}%
\special{fp}%
\special{pn 13}%
\special{pa 1658 970}%
\special{pa 1378 774}%
\special{fp}%
\special{pn 13}%
\special{pa 1238 676}%
\special{pa 958 480}%
\special{fp}%
\special{pn 13}%
\special{ar 1308 200 448 448  6.2831853 6.2831853}%
\special{ar 1308 200 448 448  0.0000000 0.6747409}%
\special{pn 13}%
\special{ar 1308 200 448 448  2.4668517 3.1415927}%
\special{pn 13}%
\special{ar 1308 1250 448 448  5.6084444 6.2831853}%
\special{pn 13}%
\special{ar 1308 1250 448 448  3.1415927 3.8163336}%
\end{picture}%
}

\newcommand{\picfj}{
\unitlength 0.1in
\begin{picture}(19.25,19.25)(5.96,-21.25)
\special{pn 4}%
\special{pa 601 2120}%
\special{pa 601 200}%
\special{fp}%
\special{sh 1}%
\special{pa 601 200}%
\special{pa 581 267}%
\special{pa 601 253}%
\special{pa 621 267}%
\special{pa 601 200}%
\special{fp}%
\special{pn 8}%
\special{pa 601 2120}%
\special{pa 2521 2120}%
\special{fp}%
\special{sh 1}%
\special{pa 2521 2120}%
\special{pa 2454 2100}%
\special{pa 2468 2120}%
\special{pa 2454 2140}%
\special{pa 2521 2120}%
\special{fp}%
\special{pn 8}%
\special{pa 601 2120}%
\special{pa 841 1160}%
\special{fp}%
\special{sh 1}%
\special{pa 841 1160}%
\special{pa 805 1220}%
\special{pa 828 1212}%
\special{pa 844 1230}%
\special{pa 841 1160}%
\special{fp}%
\end{picture}%
}

\newcommand{\picfk}{
\unitlength 0.1in
\begin{picture}(8.96,10.50)(7.19,-18.50)
\special{pn 4}%
\special{ar 1167 1080 350 140  0.0000000 6.2831853}%
\special{pn 4}%
\special{ar 1167 1570 350 140  6.2831853 6.2831853}%
\special{ar 1167 1570 350 140  0.0000000 3.1415927}%
\special{pn 8}%
\special{ar 1167 1570 350 140  3.1415927 3.1905722}%
\special{ar 1167 1570 350 140  3.3375110 3.3864906}%
\special{ar 1167 1570 350 140  3.5334294 3.5824090}%
\special{ar 1167 1570 350 140  3.7293478 3.7783273}%
\special{ar 1167 1570 350 140  3.9252661 3.9742457}%
\special{ar 1167 1570 350 140  4.1211845 4.1701641}%
\special{ar 1167 1570 350 140  4.3171029 4.3660824}%
\special{ar 1167 1570 350 140  4.5130212 4.5620008}%
\special{ar 1167 1570 350 140  4.7089396 4.7579192}%
\special{ar 1167 1570 350 140  4.9048580 4.9538376}%
\special{ar 1167 1570 350 140  5.1007763 5.1497559}%
\special{ar 1167 1570 350 140  5.2966947 5.3456743}%
\special{ar 1167 1570 350 140  5.4926131 5.5415927}%
\special{ar 1167 1570 350 140  5.6885314 5.7375110}%
\special{ar 1167 1570 350 140  5.8844498 5.9334294}%
\special{ar 1167 1570 350 140  6.0803682 6.1293478}%
\special{ar 1167 1570 350 140  6.2762865 6.2831853}%
\special{pn 4}%
\special{pa 817 1080}%
\special{pa 817 1570}%
\special{fp}%
\special{pn 4}%
\special{pa 1517 1080}%
\special{pa 1517 1570}%
\special{fp}%
\special{pn 13}%
\special{ar 1517 1080 700 490  1.5707963 3.1415927}%
\special{pn 13}%
\special{ar 1517 1570 700 490  3.1415927 3.1617607}%
\special{ar 1517 1570 700 490  3.2222649 3.2424330}%
\special{ar 1517 1570 700 490  3.3029372 3.3231053}%
\special{ar 1517 1570 700 490  3.3836095 3.4037775}%
\special{ar 1517 1570 700 490  3.4642817 3.4844498}%
\special{ar 1517 1570 700 490  3.5449540 3.5651221}%
\special{ar 1517 1570 700 490  3.6256263 3.6457943}%
\special{ar 1517 1570 700 490  3.7062985 3.7264666}%
\special{ar 1517 1570 700 490  3.7869708 3.8071389}%
\special{ar 1517 1570 700 490  3.8676431 3.8878111}%
\special{ar 1517 1570 700 490  3.9483153 3.9684834}%
\special{ar 1517 1570 700 490  4.0289876 4.0491557}%
\special{ar 1517 1570 700 490  4.1096599 4.1298279}%
\special{ar 1517 1570 700 490  4.1903321 4.2105002}%
\special{ar 1517 1570 700 490  4.2710044 4.2911725}%
\special{ar 1517 1570 700 490  4.3516767 4.3718448}%
\special{ar 1517 1570 700 490  4.4323490 4.4525170}%
\special{ar 1517 1570 700 490  4.5130212 4.5331893}%
\special{ar 1517 1570 700 490  4.5936935 4.6138616}%
\special{ar 1517 1570 700 490  4.6743658 4.6945338}%
\special{pn 13}%
\special{ar 1167 800 448 448  6.2831853 6.2831853}%
\special{ar 1167 800 448 448  0.0000000 0.6747409}%
\special{pn 13}%
\special{ar 1167 800 448 448  2.4668517 3.1415927}%
\special{pn 13}%
\special{ar 1167 1850 448 448  3.1415927 3.8163336}%
\special{pn 13}%
\special{ar 1167 1850 448 448  5.6084444 6.2831853}%
\end{picture}%
}

\newcommand{\picfl}{
\unitlength 0.1in
\begin{picture}(8.96,10.50)(8.60,-12.50)
\special{pn 8}%
\special{pa 958 480}%
\special{pa 958 970}%
\special{dt 0.04}%
\special{pa 958 970}%
\special{pa 958 969}%
\special{dt 0.04}%
\special{pn 8}%
\special{pa 1658 480}%
\special{pa 1658 970}%
\special{dt 0.04}%
\special{pa 1658 970}%
\special{pa 1658 969}%
\special{dt 0.04}%
\special{pn 8}%
\special{pa 958 480}%
\special{pa 1658 480}%
\special{dt 0.04}%
\special{pa 1658 480}%
\special{pa 1657 480}%
\special{dt 0.04}%
\special{pn 8}%
\special{pa 958 970}%
\special{pa 1658 970}%
\special{dt 0.04}%
\special{pa 1658 970}%
\special{pa 1657 970}%
\special{dt 0.04}%
\special{pn 13}%
\special{ar 1308 200 448 448  6.2831853 6.2831853}%
\special{ar 1308 200 448 448  0.0000000 0.6747409}%
\special{pn 13}%
\special{ar 1308 200 448 448  2.4668517 3.1415927}%
\special{pn 13}%
\special{ar 1308 1250 448 448  5.6084444 6.2831853}%
\special{pn 13}%
\special{ar 1308 1250 448 448  3.1415927 3.8163336}%
\special{pn 13}%
\special{pa 958 480}%
\special{pa 1658 970}%
\special{fp}%
\special{pn 13}%
\special{pa 958 970}%
\special{pa 1238 774}%
\special{fp}%
\special{pa 1238 774}%
\special{pa 1238 774}%
\special{fp}%
\special{pn 13}%
\special{pa 1378 676}%
\special{pa 1658 480}%
\special{fp}%
\end{picture}%
}

\newcommand{\picfm}{
\unitlength 0.1in
\begin{picture}(3.20,3.35)(12.00,-5.35)
\special{pn 8}%
\special{ar 1360 360 160 160  0.0000000 0.0750000}%
\special{ar 1360 360 160 160  0.3000000 0.3750000}%
\special{ar 1360 360 160 160  0.6000000 0.6750000}%
\special{ar 1360 360 160 160  0.9000000 0.9750000}%
\special{ar 1360 360 160 160  1.2000000 1.2750000}%
\special{ar 1360 360 160 160  1.5000000 1.5750000}%
\special{ar 1360 360 160 160  1.8000000 1.8750000}%
\special{ar 1360 360 160 160  2.1000000 2.1750000}%
\special{ar 1360 360 160 160  2.4000000 2.4750000}%
\special{ar 1360 360 160 160  2.7000000 2.7750000}%
\special{ar 1360 360 160 160  3.0000000 3.0750000}%
\special{ar 1360 360 160 160  3.3000000 3.3750000}%
\special{ar 1360 360 160 160  3.6000000 3.6750000}%
\special{ar 1360 360 160 160  3.9000000 3.9750000}%
\special{ar 1360 360 160 160  4.2000000 4.2750000}%
\special{ar 1360 360 160 160  4.5000000 4.5750000}%
\special{ar 1360 360 160 160  4.8000000 4.8750000}%
\special{ar 1360 360 160 160  5.1000000 5.1750000}%
\special{ar 1360 360 160 160  5.4000000 5.4750000}%
\special{ar 1360 360 160 160  5.7000000 5.7750000}%
\special{ar 1360 360 160 160  6.0000000 6.0750000}%
\special{pn 8}%
\special{ar 1040 360 260 260  5.7678596 6.2831853}%
\special{ar 1040 360 260 260  0.0000000 0.5153257}%
\special{pn 8}%
\special{ar 1680 360 260 260  2.6262670 3.6635932}%
\special{pn 8}%
\special{pa 1448 240}%
\special{pa 1452 230}%
\special{fp}%
\special{sh 1}%
\special{pa 1452 230}%
\special{pa 1409 284}%
\special{pa 1432 280}%
\special{pa 1446 299}%
\special{pa 1452 230}%
\special{fp}%
\special{pn 8}%
\special{pa 1272 240}%
\special{pa 1268 230}%
\special{fp}%
\special{sh 1}%
\special{pa 1268 230}%
\special{pa 1274 299}%
\special{pa 1288 280}%
\special{pa 1311 284}%
\special{pa 1268 230}%
\special{fp}%
\put(13.7000,-6.2000){\makebox(0,0){$L_0$}}%
\end{picture}%
}

\newcommand{\picfn}{
\unitlength 0.1in
\begin{picture}(3.20,3.20)(2.00,-5.20)
\special{pn 8}%
\special{ar 360 360 160 160  0.0000000 0.0750000}%
\special{ar 360 360 160 160  0.3000000 0.3750000}%
\special{ar 360 360 160 160  0.6000000 0.6750000}%
\special{ar 360 360 160 160  0.9000000 0.9750000}%
\special{ar 360 360 160 160  1.2000000 1.2750000}%
\special{ar 360 360 160 160  1.5000000 1.5750000}%
\special{ar 360 360 160 160  1.8000000 1.8750000}%
\special{ar 360 360 160 160  2.1000000 2.1750000}%
\special{ar 360 360 160 160  2.4000000 2.4750000}%
\special{ar 360 360 160 160  2.7000000 2.7750000}%
\special{ar 360 360 160 160  3.0000000 3.0750000}%
\special{ar 360 360 160 160  3.3000000 3.3750000}%
\special{ar 360 360 160 160  3.6000000 3.6750000}%
\special{ar 360 360 160 160  3.9000000 3.9750000}%
\special{ar 360 360 160 160  4.2000000 4.2750000}%
\special{ar 360 360 160 160  4.5000000 4.5750000}%
\special{ar 360 360 160 160  4.8000000 4.8750000}%
\special{ar 360 360 160 160  5.1000000 5.1750000}%
\special{ar 360 360 160 160  5.4000000 5.4750000}%
\special{ar 360 360 160 160  5.7000000 5.7750000}%
\special{ar 360 360 160 160  6.0000000 6.0750000}%
\special{pn 8}%
\special{pa 248 472}%
\special{pa 472 248}%
\special{fp}%
\special{sh 1}%
\special{pa 472 248}%
\special{pa 411 281}%
\special{pa 434 286}%
\special{pa 439 309}%
\special{pa 472 248}%
\special{fp}%
\special{pn 8}%
\special{pa 330 330}%
\special{pa 248 248}%
\special{fp}%
\special{sh 1}%
\special{pa 248 248}%
\special{pa 281 309}%
\special{pa 286 286}%
\special{pa 309 281}%
\special{pa 248 248}%
\special{fp}%
\special{pn 8}%
\special{pa 472 472}%
\special{pa 390 390}%
\special{fp}%
\put(3.5000,-6.0000){\makebox(0,0){$L_+$}}%
\end{picture}%
}

\newcommand{\picfo}{
\unitlength 0.1in
\begin{picture}(3.20,3.20)(2.00,-5.20)
\special{pn 8}%
\special{ar 360 360 160 160  0.0000000 0.0750000}%
\special{ar 360 360 160 160  0.3000000 0.3750000}%
\special{ar 360 360 160 160  0.6000000 0.6750000}%
\special{ar 360 360 160 160  0.9000000 0.9750000}%
\special{ar 360 360 160 160  1.2000000 1.2750000}%
\special{ar 360 360 160 160  1.5000000 1.5750000}%
\special{ar 360 360 160 160  1.8000000 1.8750000}%
\special{ar 360 360 160 160  2.1000000 2.1750000}%
\special{ar 360 360 160 160  2.4000000 2.4750000}%
\special{ar 360 360 160 160  2.7000000 2.7750000}%
\special{ar 360 360 160 160  3.0000000 3.0750000}%
\special{ar 360 360 160 160  3.3000000 3.3750000}%
\special{ar 360 360 160 160  3.6000000 3.6750000}%
\special{ar 360 360 160 160  3.9000000 3.9750000}%
\special{ar 360 360 160 160  4.2000000 4.2750000}%
\special{ar 360 360 160 160  4.5000000 4.5750000}%
\special{ar 360 360 160 160  4.8000000 4.8750000}%
\special{ar 360 360 160 160  5.1000000 5.1750000}%
\special{ar 360 360 160 160  5.4000000 5.4750000}%
\special{ar 360 360 160 160  5.7000000 5.7750000}%
\special{ar 360 360 160 160  6.0000000 6.0750000}%
\special{pn 8}%
\special{pa 472 472}%
\special{pa 248 248}%
\special{fp}%
\special{sh 1}%
\special{pa 248 248}%
\special{pa 281 309}%
\special{pa 286 286}%
\special{pa 309 281}%
\special{pa 248 248}%
\special{fp}%
\special{pn 8}%
\special{pa 248 472}%
\special{pa 330 390}%
\special{fp}%
\special{pn 8}%
\special{pa 390 330}%
\special{pa 472 248}%
\special{fp}%
\special{sh 1}%
\special{pa 472 248}%
\special{pa 411 281}%
\special{pa 434 286}%
\special{pa 439 309}%
\special{pa 472 248}%
\special{fp}%
\put(3.5000,-6.0000){\makebox(0,0){$L_-$}}%
\end{picture}%
}

\newcommand{\picfq}{
\unitlength 0.1in
\begin{picture}(39.85,13.35)(4.15,-16.00)
\special{pn 8}%
\special{pa 600 1600}%
\special{pa 600 400}%
\special{fp}%
\special{sh 1}%
\special{pa 600 400}%
\special{pa 580 467}%
\special{pa 600 453}%
\special{pa 620 467}%
\special{pa 600 400}%
\special{fp}%
\special{pn 8}%
\special{pa 600 1200}%
\special{pa 4400 1200}%
\special{dt 0.045}%
\special{pa 4400 1200}%
\special{pa 4399 1200}%
\special{dt 0.045}%
\special{pn 8}%
\special{pa 600 800}%
\special{pa 4400 800}%
\special{dt 0.045}%
\special{pa 4400 800}%
\special{pa 4399 800}%
\special{dt 0.045}%
\special{pn 13}%
\special{pa 1200 600}%
\special{pa 1200 1400}%
\special{fp}%
\special{pn 13}%
\special{pa 1600 600}%
\special{pa 1600 1400}%
\special{fp}%
\special{pn 13}%
\special{pa 3600 600}%
\special{pa 3600 1400}%
\special{fp}%
\special{pn 13}%
\special{pa 4000 600}%
\special{pa 4000 1400}%
\special{fp}%
\special{pn 13}%
\special{pa 2400 600}%
\special{pa 2400 800}%
\special{fp}%
\special{pn 13}%
\special{pa 2800 600}%
\special{pa 2800 800}%
\special{fp}%
\special{pn 13}%
\special{pa 2400 1200}%
\special{pa 2400 1400}%
\special{fp}%
\special{pn 13}%
\special{pa 2800 1200}%
\special{pa 2800 1400}%
\special{fp}%
\put(4.6000,-3.5000){\makebox(0,0){$t$}}%
\put(4.7000,-10.1000){\makebox(0,0){$I_a$}}%
\put(20.0000,-10.0000){\makebox(0,0){$\cdots$}}%
\put(32.0000,-10.0000){\makebox(0,0){$\cdots$}}%
\put(26.0000,-14.8000){\makebox(0,0){$a$}}%
\special{pn 4}%
\special{ar 2600 800 200 100  0.0000000 6.2831853}%
\special{pn 8}%
\special{ar 2600 1200 200 100  3.1415927 3.2215927}%
\special{ar 2600 1200 200 100  3.4615927 3.5415927}%
\special{ar 2600 1200 200 100  3.7815927 3.8615927}%
\special{ar 2600 1200 200 100  4.1015927 4.1815927}%
\special{ar 2600 1200 200 100  4.4215927 4.5015927}%
\special{ar 2600 1200 200 100  4.7415927 4.8215927}%
\special{ar 2600 1200 200 100  5.0615927 5.1415927}%
\special{ar 2600 1200 200 100  5.3815927 5.4615927}%
\special{ar 2600 1200 200 100  5.7015927 5.7815927}%
\special{ar 2600 1200 200 100  6.0215927 6.1015927}%
\special{pn 4}%
\special{ar 2600 1200 200 100  6.2831853 6.2831853}%
\special{ar 2600 1200 200 100  0.0000000 3.1415927}%
\special{pn 4}%
\special{pa 2400 800}%
\special{pa 2400 1200}%
\special{fp}%
\special{pn 4}%
\special{pa 2800 800}%
\special{pa 2800 1200}%
\special{fp}%
\special{pn 13}%
\special{ar 2400 800 400 400  6.2831853 6.2831853}%
\special{ar 2400 800 400 400  0.0000000 1.5707963}%
\special{pn 13}%
\special{ar 2400 1200 400 400  4.7123890 4.7423890}%
\special{ar 2400 1200 400 400  4.8323890 4.8623890}%
\special{ar 2400 1200 400 400  4.9523890 4.9823890}%
\special{ar 2400 1200 400 400  5.0723890 5.1023890}%
\special{ar 2400 1200 400 400  5.1923890 5.2223890}%
\special{ar 2400 1200 400 400  5.3123890 5.3423890}%
\special{ar 2400 1200 400 400  5.4323890 5.4623890}%
\special{ar 2400 1200 400 400  5.5523890 5.5823890}%
\special{ar 2400 1200 400 400  5.6723890 5.7023890}%
\special{ar 2400 1200 400 400  5.7923890 5.8223890}%
\special{ar 2400 1200 400 400  5.9123890 5.9423890}%
\special{ar 2400 1200 400 400  6.0323890 6.0623890}%
\special{ar 2400 1200 400 400  6.1523890 6.1823890}%
\special{ar 2400 1200 400 400  6.2723890 6.2831853}%
\end{picture}%
}

\newcommand{\picfr}{
\unitlength 0.1in
\begin{picture}(8.20,4.35)(9.50,-15.85)
\special{pn 8}%
\special{ar 1360 1320 160 160  0.0000000 0.0750000}%
\special{ar 1360 1320 160 160  0.3000000 0.3750000}%
\special{ar 1360 1320 160 160  0.6000000 0.6750000}%
\special{ar 1360 1320 160 160  0.9000000 0.9750000}%
\special{ar 1360 1320 160 160  1.2000000 1.2750000}%
\special{ar 1360 1320 160 160  1.5000000 1.5750000}%
\special{ar 1360 1320 160 160  1.8000000 1.8750000}%
\special{ar 1360 1320 160 160  2.1000000 2.1750000}%
\special{ar 1360 1320 160 160  2.4000000 2.4750000}%
\special{ar 1360 1320 160 160  2.7000000 2.7750000}%
\special{ar 1360 1320 160 160  3.0000000 3.0750000}%
\special{ar 1360 1320 160 160  3.3000000 3.3750000}%
\special{ar 1360 1320 160 160  3.6000000 3.6750000}%
\special{ar 1360 1320 160 160  3.9000000 3.9750000}%
\special{ar 1360 1320 160 160  4.2000000 4.2750000}%
\special{ar 1360 1320 160 160  4.5000000 4.5750000}%
\special{ar 1360 1320 160 160  4.8000000 4.8750000}%
\special{ar 1360 1320 160 160  5.1000000 5.1750000}%
\special{ar 1360 1320 160 160  5.4000000 5.4750000}%
\special{ar 1360 1320 160 160  5.7000000 5.7750000}%
\special{ar 1360 1320 160 160  6.0000000 6.0750000}%
\special{pn 8}%
\special{pa 1248 1432}%
\special{pa 1472 1208}%
\special{fp}%
\special{sh 1}%
\special{pa 1472 1208}%
\special{pa 1411 1241}%
\special{pa 1434 1246}%
\special{pa 1439 1269}%
\special{pa 1472 1208}%
\special{fp}%
\special{pn 8}%
\special{pa 1330 1290}%
\special{pa 1248 1208}%
\special{fp}%
\special{sh 1}%
\special{pa 1248 1208}%
\special{pa 1281 1269}%
\special{pa 1286 1246}%
\special{pa 1309 1241}%
\special{pa 1248 1208}%
\special{fp}%
\special{pn 8}%
\special{pa 1472 1432}%
\special{pa 1390 1350}%
\special{fp}%
\special{pn 8}%
\special{ar 1600 1320 170 170  3.8604227 6.2831853}%
\special{ar 1600 1320 170 170  0.0000000 2.4227627}%
\special{pn 8}%
\special{ar 1120 1320 170 170  0.7188300 5.5643553}%
\put(13.7000,-16.7000){\makebox(0,0){$K_n^+$}}%
\end{picture}%
}

\newcommand{\picfs}{
\unitlength 0.1in
\begin{picture}(8.20,4.35)(9.50,-15.85)
\special{pn 8}%
\special{ar 1360 1320 160 160  0.0000000 0.0750000}%
\special{ar 1360 1320 160 160  0.3000000 0.3750000}%
\special{ar 1360 1320 160 160  0.6000000 0.6750000}%
\special{ar 1360 1320 160 160  0.9000000 0.9750000}%
\special{ar 1360 1320 160 160  1.2000000 1.2750000}%
\special{ar 1360 1320 160 160  1.5000000 1.5750000}%
\special{ar 1360 1320 160 160  1.8000000 1.8750000}%
\special{ar 1360 1320 160 160  2.1000000 2.1750000}%
\special{ar 1360 1320 160 160  2.4000000 2.4750000}%
\special{ar 1360 1320 160 160  2.7000000 2.7750000}%
\special{ar 1360 1320 160 160  3.0000000 3.0750000}%
\special{ar 1360 1320 160 160  3.3000000 3.3750000}%
\special{ar 1360 1320 160 160  3.6000000 3.6750000}%
\special{ar 1360 1320 160 160  3.9000000 3.9750000}%
\special{ar 1360 1320 160 160  4.2000000 4.2750000}%
\special{ar 1360 1320 160 160  4.5000000 4.5750000}%
\special{ar 1360 1320 160 160  4.8000000 4.8750000}%
\special{ar 1360 1320 160 160  5.1000000 5.1750000}%
\special{ar 1360 1320 160 160  5.4000000 5.4750000}%
\special{ar 1360 1320 160 160  5.7000000 5.7750000}%
\special{ar 1360 1320 160 160  6.0000000 6.0750000}%
\special{pn 8}%
\special{ar 1600 1320 170 170  3.8604227 6.2831853}%
\special{ar 1600 1320 170 170  0.0000000 2.4227627}%
\special{pn 8}%
\special{ar 1120 1320 170 170  0.7188300 5.5643553}%
\special{pn 8}%
\special{pa 1472 1432}%
\special{pa 1248 1208}%
\special{fp}%
\special{sh 1}%
\special{pa 1248 1208}%
\special{pa 1281 1269}%
\special{pa 1286 1246}%
\special{pa 1309 1241}%
\special{pa 1248 1208}%
\special{fp}%
\special{pn 8}%
\special{pa 1390 1290}%
\special{pa 1472 1208}%
\special{fp}%
\special{sh 1}%
\special{pa 1472 1208}%
\special{pa 1411 1241}%
\special{pa 1434 1246}%
\special{pa 1439 1269}%
\special{pa 1472 1208}%
\special{fp}%
\special{pn 8}%
\special{pa 1248 1432}%
\special{pa 1330 1350}%
\special{fp}%
\put(13.6000,-16.7000){\makebox(0,0){$K_n^-$}}%
\end{picture}%
}

\newcommand{\picft}{
\unitlength 0.1in
\begin{picture}(8.20,3.95)(9.50,-15.45)
\special{pn 8}%
\special{ar 1360 1320 160 160  0.0000000 0.0750000}%
\special{ar 1360 1320 160 160  0.3000000 0.3750000}%
\special{ar 1360 1320 160 160  0.6000000 0.6750000}%
\special{ar 1360 1320 160 160  0.9000000 0.9750000}%
\special{ar 1360 1320 160 160  1.2000000 1.2750000}%
\special{ar 1360 1320 160 160  1.5000000 1.5750000}%
\special{ar 1360 1320 160 160  1.8000000 1.8750000}%
\special{ar 1360 1320 160 160  2.1000000 2.1750000}%
\special{ar 1360 1320 160 160  2.4000000 2.4750000}%
\special{ar 1360 1320 160 160  2.7000000 2.7750000}%
\special{ar 1360 1320 160 160  3.0000000 3.0750000}%
\special{ar 1360 1320 160 160  3.3000000 3.3750000}%
\special{ar 1360 1320 160 160  3.6000000 3.6750000}%
\special{ar 1360 1320 160 160  3.9000000 3.9750000}%
\special{ar 1360 1320 160 160  4.2000000 4.2750000}%
\special{ar 1360 1320 160 160  4.5000000 4.5750000}%
\special{ar 1360 1320 160 160  4.8000000 4.8750000}%
\special{ar 1360 1320 160 160  5.1000000 5.1750000}%
\special{ar 1360 1320 160 160  5.4000000 5.4750000}%
\special{ar 1360 1320 160 160  5.7000000 5.7750000}%
\special{ar 1360 1320 160 160  6.0000000 6.0750000}%
\special{pn 8}%
\special{ar 1600 1320 170 170  0.0000000 6.2831853}%
\special{pn 8}%
\special{ar 1120 1320 170 170  0.0000000 6.2831853}%
\special{pn 8}%
\special{pa 1464 1220}%
\special{pa 1472 1208}%
\special{fp}%
\special{sh 1}%
\special{pa 1472 1208}%
\special{pa 1418 1252}%
\special{pa 1442 1252}%
\special{pa 1452 1275}%
\special{pa 1472 1208}%
\special{fp}%
\special{pn 8}%
\special{pa 1256 1220}%
\special{pa 1248 1208}%
\special{fp}%
\special{sh 1}%
\special{pa 1248 1208}%
\special{pa 1268 1275}%
\special{pa 1278 1252}%
\special{pa 1302 1252}%
\special{pa 1248 1208}%
\special{fp}%
\put(11.0000,-16.8000){\makebox(0,0){$K_n^0$}}%
\put(16.5000,-16.8000){\makebox(0,0){$K_{n+1}^0$}}%
\end{picture}%
}

\newcommand{\picfu}{
\unitlength 0.1in
\begin{picture}(6.72,5.65)(13.04,-17.65)
\special{pn 8}%
\special{ar 1640 1360 160 160  0.0000000 0.0750000}%
\special{ar 1640 1360 160 160  0.3000000 0.3750000}%
\special{ar 1640 1360 160 160  0.6000000 0.6750000}%
\special{ar 1640 1360 160 160  0.9000000 0.9750000}%
\special{ar 1640 1360 160 160  1.2000000 1.2750000}%
\special{ar 1640 1360 160 160  1.5000000 1.5750000}%
\special{ar 1640 1360 160 160  1.8000000 1.8750000}%
\special{ar 1640 1360 160 160  2.1000000 2.1750000}%
\special{ar 1640 1360 160 160  2.4000000 2.4750000}%
\special{ar 1640 1360 160 160  2.7000000 2.7750000}%
\special{ar 1640 1360 160 160  3.0000000 3.0750000}%
\special{ar 1640 1360 160 160  3.3000000 3.3750000}%
\special{ar 1640 1360 160 160  3.6000000 3.6750000}%
\special{ar 1640 1360 160 160  3.9000000 3.9750000}%
\special{ar 1640 1360 160 160  4.2000000 4.2750000}%
\special{ar 1640 1360 160 160  4.5000000 4.5750000}%
\special{ar 1640 1360 160 160  4.8000000 4.8750000}%
\special{ar 1640 1360 160 160  5.1000000 5.1750000}%
\special{ar 1640 1360 160 160  5.4000000 5.4750000}%
\special{ar 1640 1360 160 160  5.7000000 5.7750000}%
\special{ar 1640 1360 160 160  6.0000000 6.0750000}%
\special{pn 8}%
\special{pa 1528 1472}%
\special{pa 1752 1248}%
\special{fp}%
\special{sh 1}%
\special{pa 1752 1248}%
\special{pa 1691 1281}%
\special{pa 1714 1286}%
\special{pa 1719 1309}%
\special{pa 1752 1248}%
\special{fp}%
\special{pn 8}%
\special{pa 1610 1330}%
\special{pa 1528 1248}%
\special{fp}%
\special{sh 1}%
\special{pa 1528 1248}%
\special{pa 1561 1309}%
\special{pa 1566 1286}%
\special{pa 1589 1281}%
\special{pa 1528 1248}%
\special{fp}%
\special{pn 8}%
\special{pa 1752 1472}%
\special{pa 1670 1390}%
\special{fp}%
\special{pn 8}%
\special{ar 1752 1472 224 224  4.7123890 6.2831853}%
\special{ar 1752 1472 224 224  0.0000000 3.1415927}%
\special{pn 8}%
\special{ar 1528 1472 224 224  1.2433961 4.7123890}%
\special{pn 8}%
\special{ar 1528 1472 224 224  6.2831853 6.2831853}%
\special{ar 1528 1472 224 224  0.0000000 0.8233546}%
\put(18.9000,-18.5000){\makebox(0,0){$K_n^+$}}%
\put(14.1000,-18.5000){\makebox(0,0){$K_{n+1}^+$}}%
\end{picture}%
}

\newcommand{\picfw}{
\unitlength 0.1in
\begin{picture}(6.76,5.75)(13.00,-17.75)
\special{pn 8}%
\special{ar 1640 1360 160 160  0.0000000 0.0750000}%
\special{ar 1640 1360 160 160  0.3000000 0.3750000}%
\special{ar 1640 1360 160 160  0.6000000 0.6750000}%
\special{ar 1640 1360 160 160  0.9000000 0.9750000}%
\special{ar 1640 1360 160 160  1.2000000 1.2750000}%
\special{ar 1640 1360 160 160  1.5000000 1.5750000}%
\special{ar 1640 1360 160 160  1.8000000 1.8750000}%
\special{ar 1640 1360 160 160  2.1000000 2.1750000}%
\special{ar 1640 1360 160 160  2.4000000 2.4750000}%
\special{ar 1640 1360 160 160  2.7000000 2.7750000}%
\special{ar 1640 1360 160 160  3.0000000 3.0750000}%
\special{ar 1640 1360 160 160  3.3000000 3.3750000}%
\special{ar 1640 1360 160 160  3.6000000 3.6750000}%
\special{ar 1640 1360 160 160  3.9000000 3.9750000}%
\special{ar 1640 1360 160 160  4.2000000 4.2750000}%
\special{ar 1640 1360 160 160  4.5000000 4.5750000}%
\special{ar 1640 1360 160 160  4.8000000 4.8750000}%
\special{ar 1640 1360 160 160  5.1000000 5.1750000}%
\special{ar 1640 1360 160 160  5.4000000 5.4750000}%
\special{ar 1640 1360 160 160  5.7000000 5.7750000}%
\special{ar 1640 1360 160 160  6.0000000 6.0750000}%
\special{pn 8}%
\special{ar 1752 1472 224 224  4.7123890 6.2831853}%
\special{ar 1752 1472 224 224  0.0000000 3.1415927}%
\special{pn 8}%
\special{ar 1528 1472 224 224  1.2433961 4.7123890}%
\special{pn 8}%
\special{ar 1528 1472 224 224  6.2831853 6.2831853}%
\special{ar 1528 1472 224 224  0.0000000 0.8233546}%
\special{pn 8}%
\special{pa 1752 1472}%
\special{pa 1528 1248}%
\special{fp}%
\special{sh 1}%
\special{pa 1528 1248}%
\special{pa 1561 1309}%
\special{pa 1566 1286}%
\special{pa 1589 1281}%
\special{pa 1528 1248}%
\special{fp}%
\special{pn 8}%
\special{pa 1528 1472}%
\special{pa 1610 1390}%
\special{fp}%
\special{pn 8}%
\special{pa 1670 1330}%
\special{pa 1752 1248}%
\special{fp}%
\special{sh 1}%
\special{pa 1752 1248}%
\special{pa 1691 1281}%
\special{pa 1714 1286}%
\special{pa 1719 1309}%
\special{pa 1752 1248}%
\special{fp}%
\put(19.0000,-18.6000){\makebox(0,0){$K_n^-$}}%
\put(13.9000,-18.5000){\makebox(0,0){$K_{n+1}^-$}}%
\end{picture}%
}

\newcommand{\picfx}{
\unitlength 0.1in
\begin{picture}(6.72,5.55)(13.04,-17.55)
\special{pn 8}%
\special{ar 1640 1360 160 160  0.0000000 0.0750000}%
\special{ar 1640 1360 160 160  0.3000000 0.3750000}%
\special{ar 1640 1360 160 160  0.6000000 0.6750000}%
\special{ar 1640 1360 160 160  0.9000000 0.9750000}%
\special{ar 1640 1360 160 160  1.2000000 1.2750000}%
\special{ar 1640 1360 160 160  1.5000000 1.5750000}%
\special{ar 1640 1360 160 160  1.8000000 1.8750000}%
\special{ar 1640 1360 160 160  2.1000000 2.1750000}%
\special{ar 1640 1360 160 160  2.4000000 2.4750000}%
\special{ar 1640 1360 160 160  2.7000000 2.7750000}%
\special{ar 1640 1360 160 160  3.0000000 3.0750000}%
\special{ar 1640 1360 160 160  3.3000000 3.3750000}%
\special{ar 1640 1360 160 160  3.6000000 3.6750000}%
\special{ar 1640 1360 160 160  3.9000000 3.9750000}%
\special{ar 1640 1360 160 160  4.2000000 4.2750000}%
\special{ar 1640 1360 160 160  4.5000000 4.5750000}%
\special{ar 1640 1360 160 160  4.8000000 4.8750000}%
\special{ar 1640 1360 160 160  5.1000000 5.1750000}%
\special{ar 1640 1360 160 160  5.4000000 5.4750000}%
\special{ar 1640 1360 160 160  5.7000000 5.7750000}%
\special{ar 1640 1360 160 160  6.0000000 6.0750000}%
\special{pn 8}%
\special{ar 1752 1472 224 224  4.7123890 6.2831853}%
\special{ar 1752 1472 224 224  0.0000000 3.1415927}%
\special{pn 8}%
\special{ar 1528 1472 224 224  1.2433961 4.7123890}%
\special{pn 8}%
\special{ar 1528 1472 224 224  6.2831853 6.2831853}%
\special{ar 1528 1472 224 224  0.0000000 0.8233546}%
\special{pn 8}%
\special{ar 1880 1360 170 170  2.4227627 3.8604227}%
\special{pn 8}%
\special{ar 1400 1360 170 170  5.5643553 6.2831853}%
\special{ar 1400 1360 170 170  0.0000000 0.7188300}%
\special{pn 8}%
\special{pa 1744 1260}%
\special{pa 1752 1248}%
\special{fp}%
\special{sh 1}%
\special{pa 1752 1248}%
\special{pa 1698 1292}%
\special{pa 1722 1292}%
\special{pa 1732 1315}%
\special{pa 1752 1248}%
\special{fp}%
\special{pn 8}%
\special{pa 1536 1260}%
\special{pa 1528 1248}%
\special{fp}%
\special{sh 1}%
\special{pa 1528 1248}%
\special{pa 1548 1315}%
\special{pa 1558 1292}%
\special{pa 1582 1292}%
\special{pa 1528 1248}%
\special{fp}%
\put(16.5000,-18.4000){\makebox(0,0){$K_n^0$}}%
\end{picture}%
}
\newcommand{\picga}{
\unitlength 0.1in
\begin{picture}(10.00,4.00)(2.00,-6.00)
\special{pn 16}%
\special{ar 400 400 200 200  0.0000000 6.2831853}%
\special{pn 16}%
\special{ar 1000 400 200 200  0.0000000 6.2831853}%
\special{pn 4}%
\special{pa 400 200}%
\special{pa 400 600}%
\special{fp}%
\special{pn 4}%
\special{pa 200 400}%
\special{pa 400 400}%
\special{fp}%
\special{pn 4}%
\special{ar 1000 400 50 50  0.0000000 6.2831853}%
\special{pn 4}%
\special{pa 1000 200}%
\special{pa 1000 350}%
\special{fp}%
\special{pn 4}%
\special{pa 1000 450}%
\special{pa 1000 600}%
\special{fp}%
\special{pn 4}%
\special{pa 532 250}%
\special{pa 868 250}%
\special{fp}%
\special{pn 4}%
\special{pa 532 550}%
\special{pa 868 550}%
\special{fp}%
\special{pn 4}%
\special{pa 600 400}%
\special{pa 700 250}%
\special{fp}%
\end{picture}%
}

\newcommand{\picgb}{
\unitlength 0.1in
\begin{picture}(8.96,10.50)(8.60,-14.50)
\special{pn 13}%
\special{ar 958 680 700 490  6.2831853 6.2831853}%
\special{ar 958 680 700 490  0.0000000 1.5707963}%
\special{pn 4}%
\special{ar 1308 680 350 140  0.7954988 6.2831853}%
\special{ar 1308 680 350 140  0.0000000 0.7853982}%
\special{pn 4}%
\special{pa 958 680}%
\special{pa 958 1170}%
\special{fp}%
\special{pn 4}%
\special{pa 1658 680}%
\special{pa 1658 1170}%
\special{fp}%
\special{pn 13}%
\special{ar 958 1170 700 490  4.7123890 4.7325570}%
\special{ar 958 1170 700 490  4.7930612 4.8132293}%
\special{ar 958 1170 700 490  4.8737335 4.8939016}%
\special{ar 958 1170 700 490  4.9544058 4.9745739}%
\special{ar 958 1170 700 490  5.0350781 5.0552461}%
\special{ar 958 1170 700 490  5.1157503 5.1359184}%
\special{ar 958 1170 700 490  5.1964226 5.2165907}%
\special{ar 958 1170 700 490  5.2770949 5.2972629}%
\special{ar 958 1170 700 490  5.3577671 5.3779352}%
\special{ar 958 1170 700 490  5.4384394 5.4586075}%
\special{ar 958 1170 700 490  5.5191117 5.5392797}%
\special{ar 958 1170 700 490  5.5997839 5.6199520}%
\special{ar 958 1170 700 490  5.6804562 5.7006243}%
\special{ar 958 1170 700 490  5.7611285 5.7812965}%
\special{ar 958 1170 700 490  5.8418007 5.8619688}%
\special{ar 958 1170 700 490  5.9224730 5.9426411}%
\special{ar 958 1170 700 490  6.0031453 6.0233134}%
\special{ar 958 1170 700 490  6.0838176 6.1039856}%
\special{ar 958 1170 700 490  6.1644898 6.1846579}%
\special{ar 958 1170 700 490  6.2451621 6.2653302}%
\special{pn 8}%
\special{ar 1308 1170 350 140  3.1415927 3.1905722}%
\special{ar 1308 1170 350 140  3.3375110 3.3864906}%
\special{ar 1308 1170 350 140  3.5334294 3.5824090}%
\special{ar 1308 1170 350 140  3.7293478 3.7783273}%
\special{ar 1308 1170 350 140  3.9252661 3.9742457}%
\special{ar 1308 1170 350 140  4.1211845 4.1701641}%
\special{ar 1308 1170 350 140  4.3171029 4.3660824}%
\special{ar 1308 1170 350 140  4.5130212 4.5620008}%
\special{ar 1308 1170 350 140  4.7089396 4.7579192}%
\special{ar 1308 1170 350 140  4.9048580 4.9538376}%
\special{ar 1308 1170 350 140  5.1007763 5.1497559}%
\special{ar 1308 1170 350 140  5.2966947 5.3456743}%
\special{ar 1308 1170 350 140  5.4926131 5.5415927}%
\special{ar 1308 1170 350 140  5.6885314 5.7375110}%
\special{ar 1308 1170 350 140  5.8844498 5.9334294}%
\special{ar 1308 1170 350 140  6.0803682 6.1293478}%
\special{ar 1308 1170 350 140  6.2762865 6.2831853}%
\special{pn 4}%
\special{ar 1308 1170 350 140  6.2831853 6.2831853}%
\special{ar 1308 1170 350 140  0.0000000 3.1415927}%
\special{pn 13}%
\special{ar 1308 400 448 448  6.2831853 6.2831853}%
\special{ar 1308 400 448 448  0.0000000 0.6747409}%
\special{pn 13}%
\special{ar 1308 400 448 448  2.4668517 3.1415927}%
\special{pn 13}%
\special{ar 1308 1450 448 448  5.6084444 6.2831853}%
\special{pn 13}%
\special{ar 1308 1450 448 448  3.1415927 3.8163336}%
\special{pn 30}%
\special{sh 1}%
\special{ar 1308 1100 10 10 0  6.28318530717959E+0000}%
\special{sh 1}%
\special{ar 1308 1100 10 10 0  6.28318530717959E+0000}%
\special{pn 30}%
\special{sh 1}%
\special{ar 1308 750 10 10 0  6.28318530717959E+0000}%
\special{sh 1}%
\special{ar 1308 750 10 10 0  6.28318530717959E+0000}%
\put(13.300,-12.100){\makebox(0,0){$d_k$}}%
\put(13.300,-6.4000){\makebox(0,0){$d_k^\prime$}}%
\end{picture}%
}

\newcommand{\picgc}{
\unitlength 0.1in
\begin{picture}(4.00,4.00)(12.00,-4.20)
\special{pn 8}%
\special{ar 1400 400 150 150  0.0000000 0.0800000}%
\special{ar 1400 400 150 150  0.3200000 0.4000000}%
\special{ar 1400 400 150 150  0.6400000 0.7200000}%
\special{ar 1400 400 150 150  0.9600000 1.0400000}%
\special{ar 1400 400 150 150  1.2800000 1.3600000}%
\special{ar 1400 400 150 150  1.6000000 1.6800000}%
\special{ar 1400 400 150 150  1.9200000 2.0000000}%
\special{ar 1400 400 150 150  2.2400000 2.3200000}%
\special{ar 1400 400 150 150  2.5600000 2.6400000}%
\special{ar 1400 400 150 150  2.8800000 2.9600000}%
\special{ar 1400 400 150 150  3.2000000 3.2800000}%
\special{ar 1400 400 150 150  3.5200000 3.6000000}%
\special{ar 1400 400 150 150  3.8400000 3.9200000}%
\special{ar 1400 400 150 150  4.1600000 4.2400000}%
\special{ar 1400 400 150 150  4.4800000 4.5600000}%
\special{ar 1400 400 150 150  4.8000000 4.8800000}%
\special{ar 1400 400 150 150  5.1200000 5.2000000}%
\special{ar 1400 400 150 150  5.4400000 5.5200000}%
\special{ar 1400 400 150 150  5.7600000 5.8400000}%
\special{ar 1400 400 150 150  6.0800000 6.1600000}%
\special{pn 8}%
\special{pa 1200 600}%
\special{pa 1600 200}%
\special{fp}%
\special{sh 1}%
\special{pa 1600 200}%
\special{pa 1539 233}%
\special{pa 1562 238}%
\special{pa 1567 261}%
\special{pa 1600 200}%
\special{fp}%
\special{pn 8}%
\special{pa 1600 600}%
\special{pa 1200 200}%
\special{fp}%
\special{sh 1}%
\special{pa 1200 200}%
\special{pa 1233 261}%
\special{pa 1238 238}%
\special{pa 1261 233}%
\special{pa 1200 200}%
\special{fp}%
\put(13.4,-7){$a_i$}
\end{picture}%
}

\newcommand{\picgd}{
\unitlength 0.1in
\begin{picture}(4.00,4.00)(12.00,-4.20)
\special{pn 8}%
\special{ar 1400 400 150 150  0.0000000 0.0800000}%
\special{ar 1400 400 150 150  0.3200000 0.4000000}%
\special{ar 1400 400 150 150  0.6400000 0.7200000}%
\special{ar 1400 400 150 150  0.9600000 1.0400000}%
\special{ar 1400 400 150 150  1.2800000 1.3600000}%
\special{ar 1400 400 150 150  1.6000000 1.6800000}%
\special{ar 1400 400 150 150  1.9200000 2.0000000}%
\special{ar 1400 400 150 150  2.2400000 2.3200000}%
\special{ar 1400 400 150 150  2.5600000 2.6400000}%
\special{ar 1400 400 150 150  2.8800000 2.9600000}%
\special{ar 1400 400 150 150  3.2000000 3.2800000}%
\special{ar 1400 400 150 150  3.5200000 3.6000000}%
\special{ar 1400 400 150 150  3.8400000 3.9200000}%
\special{ar 1400 400 150 150  4.1600000 4.2400000}%
\special{ar 1400 400 150 150  4.4800000 4.5600000}%
\special{ar 1400 400 150 150  4.8000000 4.8800000}%
\special{ar 1400 400 150 150  5.1200000 5.2000000}%
\special{ar 1400 400 150 150  5.4400000 5.5200000}%
\special{ar 1400 400 150 150  5.7600000 5.8400000}%
\special{ar 1400 400 150 150  6.0800000 6.1600000}%
\special{pn 8}%
\special{pa 1200 600}%
\special{pa 1295 505}%
\special{fp}%
\special{pn 8}%
\special{pa 1505 295}%
\special{pa 1600 200}%
\special{fp}%
\special{pn 8}%
\special{pa 1200 200}%
\special{pa 1295 295}%
\special{fp}%
\special{pn 8}%
\special{pa 1505 505}%
\special{pa 1600 600}%
\special{fp}%
\special{pn 8}%
\special{ar 1600 400 142 142  2.3062361 3.9769492}%
\special{pn 8}%
\special{ar 1200 400 142 142  5.4478287 6.2831853}%
\special{ar 1200 400 142 142  0.0000000 0.8353566}%
\end{picture}%
}

\newcommand{\picge}{
\unitlength 0.1in
\begin{picture}(4.00,4.00)(12.00,-4.20)
\special{pn 8}%
\special{ar 1400 400 150 150  0.0000000 0.0800000}%
\special{ar 1400 400 150 150  0.3200000 0.4000000}%
\special{ar 1400 400 150 150  0.6400000 0.7200000}%
\special{ar 1400 400 150 150  0.9600000 1.0400000}%
\special{ar 1400 400 150 150  1.2800000 1.3600000}%
\special{ar 1400 400 150 150  1.6000000 1.6800000}%
\special{ar 1400 400 150 150  1.9200000 2.0000000}%
\special{ar 1400 400 150 150  2.2400000 2.3200000}%
\special{ar 1400 400 150 150  2.5600000 2.6400000}%
\special{ar 1400 400 150 150  2.8800000 2.9600000}%
\special{ar 1400 400 150 150  3.2000000 3.2800000}%
\special{ar 1400 400 150 150  3.5200000 3.6000000}%
\special{ar 1400 400 150 150  3.8400000 3.9200000}%
\special{ar 1400 400 150 150  4.1600000 4.2400000}%
\special{ar 1400 400 150 150  4.4800000 4.5600000}%
\special{ar 1400 400 150 150  4.8000000 4.8800000}%
\special{ar 1400 400 150 150  5.1200000 5.2000000}%
\special{ar 1400 400 150 150  5.4400000 5.5200000}%
\special{ar 1400 400 150 150  5.7600000 5.8400000}%
\special{ar 1400 400 150 150  6.0800000 6.1600000}%
\special{pn 8}%
\special{pa 1200 600}%
\special{pa 1295 505}%
\special{fp}%
\special{pn 8}%
\special{pa 1505 295}%
\special{pa 1600 200}%
\special{fp}%
\special{pn 8}%
\special{pa 1200 200}%
\special{pa 1295 295}%
\special{fp}%
\special{pn 8}%
\special{pa 1600 600}%
\special{pa 1505 505}%
\special{fp}%
\special{pn 8}%
\special{ar 1400 200 142 142  0.7354398 2.4061529}%
\special{pn 8}%
\special{ar 1400 600 142 142  3.8770324 5.5477455}%
\end{picture}%
}

\newcommand{\picgf}{
\unitlength 0.1in
\begin{picture}(4.00,4.00)(12.00,-4.20)
\special{pn 8}%
\special{ar 1400 400 150 150  0.0000000 0.0800000}%
\special{ar 1400 400 150 150  0.3200000 0.4000000}%
\special{ar 1400 400 150 150  0.6400000 0.7200000}%
\special{ar 1400 400 150 150  0.9600000 1.0400000}%
\special{ar 1400 400 150 150  1.2800000 1.3600000}%
\special{ar 1400 400 150 150  1.6000000 1.6800000}%
\special{ar 1400 400 150 150  1.9200000 2.0000000}%
\special{ar 1400 400 150 150  2.2400000 2.3200000}%
\special{ar 1400 400 150 150  2.5600000 2.6400000}%
\special{ar 1400 400 150 150  2.8800000 2.9600000}%
\special{ar 1400 400 150 150  3.2000000 3.2800000}%
\special{ar 1400 400 150 150  3.5200000 3.6000000}%
\special{ar 1400 400 150 150  3.8400000 3.9200000}%
\special{ar 1400 400 150 150  4.1600000 4.2400000}%
\special{ar 1400 400 150 150  4.4800000 4.5600000}%
\special{ar 1400 400 150 150  4.8000000 4.8800000}%
\special{ar 1400 400 150 150  5.1200000 5.2000000}%
\special{ar 1400 400 150 150  5.4400000 5.5200000}%
\special{ar 1400 400 150 150  5.7600000 5.8400000}%
\special{ar 1400 400 150 150  6.0800000 6.1600000}%
\special{pn 8}%
\special{pa 1200 600}%
\special{pa 1600 200}%
\special{fp}%
\special{pn 8}%
\special{pa 1600 600}%
\special{pa 1200 200}%
\special{fp}%
\end{picture}%
}

\newcommand{\picgg}{
\unitlength 0.1in
\begin{picture}(20.05,12.00)(1.95,-14.00)
\special{pn 13}%
\special{pa 1000 400}%
\special{pa 1000 600}%
\special{fp}%
\special{pn 13}%
\special{pa 1000 1000}%
\special{pa 1000 1200}%
\special{fp}%
\special{pn 13}%
\special{pa 1400 400}%
\special{pa 1400 600}%
\special{fp}%
\special{pn 13}%
\special{pa 1400 1000}%
\special{pa 1400 1200}%
\special{fp}%
\special{pn 13}%
\special{pa 1000 600}%
\special{pa 1400 1000}%
\special{fp}%
\special{pn 13}%
\special{pa 1400 600}%
\special{pa 1000 1000}%
\special{fp}%
\special{pn 13}%
\special{pa 2200 400}%
\special{pa 2200 1200}%
\special{fp}%
\special{pn 8}%
\special{pa 1100 900}%
\special{pa 2200 900}%
\special{dt 0.045}%
\special{pa 2200 900}%
\special{pa 2199 900}%
\special{dt 0.045}%
\special{pn 8}%
\special{pa 1100 700}%
\special{pa 1300 700}%
\special{dt 0.045}%
\special{pa 1300 700}%
\special{pa 1299 700}%
\special{dt 0.045}%
\special{pn 8}%
\special{pa 200 1400}%
\special{pa 200 200}%
\special{fp}%
\special{sh 1}%
\special{pa 200 200}%
\special{pa 180 267}%
\special{pa 200 253}%
\special{pa 220 267}%
\special{pa 200 200}%
\special{fp}%
\end{picture}%
}

\newcommand{\picgh}{
\unitlength 0.1in
\begin{picture}(4.40,4.45)(1.05,-4.0)
\special{pn 13}%
\special{ar 380 360 160 160  0.0000000 6.2831853}%
\special{pn 4}%
\special{pa 538 370}%
\special{pa 540 360}%
\special{fp}%
\special{sh 1}%
\special{pa 540 360}%
\special{pa 507 421}%
\special{pa 530 412}%
\special{pa 547 429}%
\special{pa 540 360}%
\special{fp}%
\special{pn 4}%
\special{pa 222 350}%
\special{pa 220 360}%
\special{fp}%
\special{sh 1}%
\special{pa 220 360}%
\special{pa 253 299}%
\special{pa 230 308}%
\special{pa 207 291}%
\special{pa 220 360}%
\special{fp}%
\special{pn 4}%
\special{pa 220 360}%
\special{pa 540 360}%
\special{fp}%
\put(1.5000,-3.6000){\makebox(0,0){\scriptsize $a$}}%
\put(6.500,-6.5000){\makebox(0,0){$\{K_1,K_2:a\}$}}%
\end{picture}%
}

\newcommand{\picgi}{
\unitlength 0.1in
\begin{picture}(4.40,4.45)(1.05,-4.0)
\special{pn 8}%
\special{ar 380 360 160 160  0.0000000 6.2831853}%
\special{pn 4}%
\special{pa 538 370}%
\special{pa 540 360}%
\special{fp}%
\special{sh 1}%
\special{pa 540 360}%
\special{pa 507 421}%
\special{pa 530 412}%
\special{pa 547 429}%
\special{pa 540 360}%
\special{fp}%
\special{pn 4}%
\special{pa 222 350}%
\special{pa 220 360}%
\special{fp}%
\special{sh 1}%
\special{pa 220 360}%
\special{pa 253 299}%
\special{pa 230 308}%
\special{pa 207 291}%
\special{pa 220 360}%
\special{fp}%
\put(1.5000,-3.6000){\makebox(0,0){}}%
\put(4.500,-6.5000){\makebox(0,0){$K^{[17]}_0$}}%
\end{picture}%
}

\newcommand{\picgj}{
\unitlength 0.1in
\begin{picture}(4.40,4.45)(1.05,-4.0)
\special{pn 13}%
\special{ar 380 360 160 160  0.0000000 6.2831853}%
\special{pn 4}%
\special{pa 538 370}%
\special{pa 540 360}%
\special{fp}%
\special{sh 1}%
\special{pa 540 360}%
\special{pa 507 421}%
\special{pa 530 412}%
\special{pa 547 429}%
\special{pa 540 360}%
\special{fp}%
\special{pn 4}%
\special{pa 222 350}%
\special{pa 220 360}%
\special{fp}%
\special{sh 1}%
\special{pa 220 360}%
\special{pa 253 299}%
\special{pa 230 308}%
\special{pa 207 291}%
\special{pa 220 360}%
\special{fp}%
\put(1.5000,-3.6000){\makebox(0,0){}}%
\put(5.000,-6.5000){\makebox(0,0){}}%
\end{picture}%
}

\newcommand{\picgk}{
\unitlength 0.1in
\begin{picture}(6.06,4.06)(0.19,-4.3)
\special{pn 8}%
\special{ar 244 360 160 160  3.1415927 6.2831853}%
\special{pn 8}%
\special{pa 84 360}%
\special{pa 404 360}%
\special{fp}%
\special{sh 1}%
\special{pa 404 360}%
\special{pa 337 340}%
\special{pa 351 360}%
\special{pa 337 380}%
\special{pa 404 360}%
\special{fp}%
\special{pn 8}%
\special{pa 254 202}%
\special{pa 244 200}%
\special{fp}%
\special{sh 1}%
\special{pa 244 200}%
\special{pa 305 233}%
\special{pa 296 210}%
\special{pa 313 193}%
\special{pa 244 200}%
\special{fp}%
\special{pn 8}%
\special{ar 244 400 160 160  6.2831853 6.2831853}%
\special{ar 244 400 160 160  0.0000000 3.1415927}%
\special{pn 8}%
\special{pa 404 400}%
\special{pa 84 400}%
\special{fp}%
\special{sh 1}%
\special{pa 84 400}%
\special{pa 151 420}%
\special{pa 137 400}%
\special{pa 151 380}%
\special{pa 84 400}%
\special{fp}%
\special{pn 8}%
\special{pa 234 558}%
\special{pa 244 560}%
\special{fp}%
\special{sh 1}%
\special{pa 244 560}%
\special{pa 183 527}%
\special{pa 192 550}%
\special{pa 175 567}%
\special{pa 244 560}%
\special{fp}%
\put(4.2400,-3.8000){\makebox(0,0){}}%
\put(0.20,-3.8000){\makebox(0,0){\scriptsize $a$}}%
\put(5.6600,-5.4000){\makebox(0,0){$K^{[17]}_-$}}%
\put(5.7000,-2.4400){\makebox(0,0){$K^{[17]}_+$}}%
\end{picture}%
}

\newcommand{\pichc}{
\unitlength 0.1in
\begin{picture}(20.99,5)(0.95,-7)
\special{pn 16}%
\special{ar 595 657 400 400  0.0000000 6.2831853}%
\special{pn 16}%
\special{ar 1794 657 400 400  0.0000000 6.2831853}%
\special{pn 4}%
\special{pa 964 507}%
\special{pa 1424 507}%
\special{fp}%
\special{pn 4}%
\special{pa 860 357}%
\special{pa 1529 357}%
\special{fp}%
\put(16.70,-3.700){\makebox(0,0){\scriptsize $(+)$}}%
\put(15.500,-5.0700){\makebox(0,0){\scriptsize $(+)$}}%
\put(12.1000,-12.500){}%
\end{picture}%
}

\newcommand{\pichd}{
\unitlength 0.1in
\begin{picture}(9,6)(0,-5)
\special{pn 8}%
\special{ar 433 403 320 320  0.0000000 6.2831853}%
\special{pn 8}%
\special{pa 117 443}%
\special{pa 749 443}%
\special{fp}%
\special{pn 8}%
\special{pa 121 323}%
\special{pa 681 603}%
\special{fp}%
\special{pn 8}%
\special{pa 745 323}%
\special{pa 181 603}%
\special{fp}%
\put(3.500,-7.900){\makebox(0,0){}}%
\put(5.400,-7.900){\makebox(0,0){}}%
\put(7.500,-6.3500){\makebox(0,0){\footnotesize a3}}%
\put(8.300,-4.4300){\makebox(0,0){\footnotesize a4}}%
\put(8.300,-2.9100){\makebox(0,0){\footnotesize a5}}%
\put(7.700,-1.800){\makebox(0,0){}}%
\put(5.1300,-0.3500){\makebox(0,0){}}%
\put(3.5300,-0.3500){\makebox(0,0){}}%
\put(1.100,-2.0300){\makebox(0,0){}}%
\put(0.500,-3.0700){\makebox(0,0){$-$}}%
\put(0.400,-4.4300){\makebox(0,0){$+$}}%
\put(1.00,-6.3100){\makebox(0,0){$-$}}%
\put(4.3700,-9.6300){\makebox(0,0){}}%
\end{picture}%
}

\newcommand{\piche}{
\unitlength 0.1in
\begin{picture}(9,6)(0,-5)
\special{pn 8}%
\special{ar 433 403 320 320  0.0000000 6.2831853}%
\special{pn 8}%
\special{pa 117 443}%
\special{pa 749 443}%
\special{fp}%
\special{pn 8}%
\special{pa 353 95}%
\special{pa 353 711}%
\special{fp}%
\special{pn 8}%
\special{pa 513 99}%
\special{pa 513 711}%
\special{fp}%
\special{pn 8}%
\special{pa 121 323}%
\special{pa 681 603}%
\special{fp}%
\special{pn 8}%
\special{pa 745 323}%
\special{pa 181 603}%
\special{fp}%
\special{pn 8}%
\special{pa 185 203}%
\special{pa 685 203}%
\special{fp}%
\put(3.500,-7.900){\makebox(0,0){\footnotesize a1}}%
\put(5.400,-7.900){\makebox(0,0){\footnotesize a2}}%
\put(7.500,-6.3500){\makebox(0,0){\footnotesize a3}}%
\put(8.300,-4.4300){\makebox(0,0){\footnotesize a4}}%
\put(8.300,-2.9100){\makebox(0,0){\footnotesize a5}}%
\put(7.700,-1.800){\makebox(0,0){\footnotesize a6}}%
\put(5.1300,-0.3500){\makebox(0,0){$-$}}%
\put(3.5300,-0.3500){\makebox(0,0){$+$}}%
\put(1.100,-2.0300){\makebox(0,0){$-$}}%
\put(0.500,-3.0700){\makebox(0,0){$-$}}%
\put(0.400,-4.4300){\makebox(0,0){$+$}}%
\put(1.00,-6.3100){\makebox(0,0){$-$}}%
\put(4.3700,-9.6300){\makebox(0,0){}}%
\end{picture}%
}

\newcommand{\pichf}{
\unitlength 0.1in
\begin{picture}(9,6)(0,-5)
\special{pn 8}%
\special{ar 433 403 320 320  0.0000000 6.2831853}%
\end{picture}%
}

\begin{document}
\vspace*{0.8in}
\begin{center}
{\bfseries
THE COMBINATORIAL GAUSS DIAGRAM FORMULA \newline
FOR KONTSEVICH INTEGRAL}
\end{center}

\begin{center}
 TOMOSHIRO OCHIAI\\[0.2cm]
 {\footnotesize \itshape Department of Physics,\\
  University of Tokyo,\\
   Tokyo 113, Japan}\\
   e-mail: ochiai@hep-th.phys.s.u-tokyo.ac.jp\\
  {\footnotesize June, 2000}
\end{center}

\vspace{0.7in}

\begin{center}
  ABSTRACT
\end{center}

\begin{center}
  \begin{minipage}{0.9\textwidth}
  {\footnotesize \setlength{\baselineskip}{10pt}
   \hspace{15pt} In this paper, we shall give an explicit Gauss diagram formula for the Kontsevich integral of links up to degree four. This practical formula enables us to actually compute the Kontsevich integral in a combinatorial way.
    }
  \end{minipage}
\end{center}

\begin{center}
  \begin{minipage}{0.9\textwidth}
  {\footnotesize \setlength{\baselineskip}{10pt}
   {\it Keywords}:~~ Kontsevich Integral, Gauss Diagram, Combinatorial, 
   Vassiliev Invariant
    }
\end{minipage}
\end{center}

\section{Introduction}

There are several types of formulas for Vassiliev invariants. However most of them are not suited for actual computations. So we provide more practical formulas for them.

Kontsevich \cite{Ko} defined the famous link invariant (Kontsevich integral) using iterated integrals. The Kontsevich integral is a universal Vassiliev invariant of links which dominates all the other Vassiliev invariants. We give an explicit Gauss diagram formula for the Kontsevich integral up to degree four which is useful for actual computations.

In this paper we shall show the following results. 
We prove that the Kontsevich integral of links up to degree four can be expressed by some link invariants $v_{1},v_{2},v_{3.1},v_{3.2},v_{4.1},v_{4.2},v_{4.3},v_{4.4}$ (See Theorem \ref{th:Kontsevich}). We give an explicit Gauss diagram formula for these link invariants $v_{1},v_{2},v_{3.1},v_{3.2},v_{4.1},v_{4.2},v_{4.3},v_{4.4}$ in terms of Gauss diagrams (See Theorem \ref{th:mainresult}). This formula is obtained by evaluating Kontsevich integral using inductive argument.
As a corollary, we obtain an explicit Gauss diagram formula for the power series expansion of the Homfly polynomial up to degree four, since the Kontsevich integral and the weight system of $su(N)$ gives the Homfly polynomial
(See Corollary \ref{co:Homfly}).

Witten \cite{Wi} showed that the Chern-Simons quantum field theory gives a link invariant. We believe that the theory in this paper is the mathematical counter part of \cite{HS},\cite{HS2},\cite{HS3},\cite{La}, in which the quantum field theoretical method is used. In fact, these Gauss diagram formula for $v_{2},v_{3.1},v_{3.2},v_{4.1}$ and $v_{4.2}$ coincide with those obtained by different 
methods in  \cite{HS,La}, \cite{HS2,La}, \cite{HS3}, \cite{La} and \cite{La} respectively. The Gauss diagram formulas for $v_{4.3}$ and $v_{4.4}$ are completely new.     

The present paper is organized as follows. In section 2 we review the Kontsevich integral and discuss its property. In section 3 we give the Gauss diagram formula. In section 4  we discuss the relation to Homfly polynomial and give an example of the Gauss diagram formula. In section 5 we derive the Gauss diagram formula from Kontsevich integral. In section 6 we make a consistency check for the Gauss diagram formula.  
\newline

\section{Kontsevich Integral}

\subsection{Weight system}

We shall review the weight system  as in \cite{LMO} and fix our notation. (see also \cite{BN},\cite{Ko},\cite{Oh})
\begin{definition}\hspace{-5pt}{\bf.}\hspace{6pt}(CC Diagram)
A {\it uni-trivalent graph} is a graph every vertex of which is either univalent or trivalent. A uni-trivalent graph is said to be vertex-oriented if at each trivalent vertex a cyclic order of edges is specified. Let $X=\cup_{i=1}^nS^1_i$ be $n$-oriented circles and G a vertex-oriented uni-trivalent graph. {\it A Chinese Character Diagram} (CC Diagram) is the pair $\{X,G\}$ where all the univalent vertices of $G$ are on X. 
In all figures in the sequel, the component of $X$ will be drawed by thick circles and the edges of graph $G$ by thin lines. By convention, we set the orientation of each component of $X$ and the orientation of each trivalent vertex counterclockwise, unless otherwise stated. 

Two CC diagrams $D=\{X,G\}$, $D^\prime=\{X^\prime,G^\prime\}$ are regarded as equal if there is a homeomorphism $F:D\to D^\prime$ such that $F|_X$ is a homeomorphism from $X$ to $X^\prime$ which preserves orientation and $F|_G$ preserves the vertex orientation at each trivalent vertices. The degree of a CC diagram is defined to be half the number of vertices of the CC diagram.

For example, one of CC diagrams of degree 7 is  
\begin{eqnarray*}
  \picga.
\end{eqnarray*}
\end{definition}

\begin{definition}\hspace{-5pt}{\bf.}\hspace{6pt}(Chord Diagram)
A CC diagram is called a {\it Chord Diagram} if all the vertices are univalent.
An edge of a chord diagram is called a {\it chord}. Then it is clear that the degree of a chord diagram is equal to the number of the chords. 
We give an example of chord diagrams of degree 2 with 2-circles : 
\begin{eqnarray*}
  \Bigl\{\maruba\maru,\marubb\maru,\marucc,\maruca,\marua\marua\Bigr\}.
\end{eqnarray*} 
\end{definition}

\begin{definition}\hspace{-5pt}{\bf.}\hspace{6pt}
Let $\mathfrak{D}^t$ be the set of all CC diagrams.
We define the vector space $\mathcal{A}$ by
\begin{eqnarray*}
  \mathcal{A}=\mbox{span}(\mathfrak{D}^t)/\mbox{AS,IHX,STU}
\end{eqnarray*}
where the $AS,IHX$ and $STU$ relations are shown below: 
\begin{eqnarray*}
  &&AS:\quad\pice=-\picf\\
  &&IHX:\quad\picg=\quad\pich\quad-\pici\\
  &&STU:\quad\picj=\quad\pick\quad-\picl.
\end{eqnarray*}	
\newline	
$\square$
\end{definition}

Although we restrict our consideration to the natural (fundamental) representation of $su(N)$, all the argument in the sequel is valid for any simple Lie algebra and its  irreducible representation.
The matrix basis $\{T^a\}_{a\in I}$ of the natural (fundamental) representation of $su(N)$ are normalized as follows:
\begin{eqnarray*}
  [T^a,T^b]=if^{abc}T^c, \quad \mbox{Tr}(T^a T^b)=\frac{1}{2}\delta^{ab}
\end{eqnarray*}
with the structure constant $if^{abc}$. 
\begin{definition}\hspace{-5pt}{\bf.}\hspace{6pt}\label{Def:weightsystem}
(Weight System)
We define a map $W_{su(N)}:\mathcal{A}\to \mathbb{C}$ which is called the weight system for the natural (fundamental) representaion of $su(N)$. Let $D=\{X,G\}$ be a CC diagram and $E(G)$ the set of all the edges of the graph $G$. By a labelling of $D$, we mean a map $\rho:E(G)\to I$. For each labelling, we assign the structure constant $if^{abc}$ to each trivalent vertex where the three edges around the vertex are labeled by $a,b,c$ along its orientation. We assign the basis $T^a$ to each univalent vertex where the edge emanating from the vertex is labeled by $a$.  
\begin{center}
\setlength{\unitlength}{1cm}
\begin{picture}(15,4)(0,-2)
\put(2,0){\picm} 
\put(3.1,-2){$if_{abc}$}
\put(8,0){\picn}
\end{picture}
\end{center}

Define $W_{su(N)}(D)$ as follows. For each labelling, make the product of all the assigned structure constants $if^{abc}$ and all the traces of the product of the basis $T^a$ along each circle of $X$. Define $W_{su(N)}(D)$ to be the sum of these products where the sum is taken over all the labelling:
\begin{eqnarray*}
  W_{su(N)}(D)=\frac{x^m}{N^n}\sum_{a,b,c,\cdots=1}^{N^2-1}\{\mbox{product of } (if^{abc})  \}\prod^{n}
  \{\mbox{Trace}(\mbox{product of }T^a)\},
\end{eqnarray*}
where $m$ denotes the degree of $D$ and $n$ the number of the circles.
For example,
\begin{eqnarray*} 
  &&W_{su(N)}\biggl(\pico\biggr)\\
  &&\\
  &&\hspace{2cm}=\frac{x^5}{N^2}\sum_{a,b,c,d,e,f=1}^{N^2-1}
  if^{abc}\mbox{Tr}(T^eT^cT^bT^dT^a)\mbox{Tr}(T^fT^eT^fT^d).
\end{eqnarray*}
\end{definition}

\subsection{Kontsevich integral}

\begin{definition}\hspace{-5pt}{\bf.}\hspace{6pt}(Kontsevich Integral)
Let $\hat{\mathcal{A}}$ be the quotient of $\mathcal{A}$ by the framing independence relation. 
We shall define the Kontsevich integral on $\hat{\mathcal{A}}$. For more detail, see \cite{BN},\cite{Ko}. Let $X=\cup_{i=1}^nS^1_i$ be $n$-oriented circles and $\vec{x}:X\rightarrow\mathbb{R}^3$ an imbedding. An $n$-component oriented link $\LL$ is its image 
$\mathbf{L}=\{\K_1,\cdots,\K_n\}~(\K_i=\vec{x}(S^1_i))$ with the natural orientation. Let us introduce coordinates 
$z,t$ in $\mathbb{R}^3$ by $z=x_1+ix_2\in\mathbb{C}$, $t=x_3\in\mathbb{R}$. Let $t_{\mathrm{min}}$ (resp. $t_{\mathrm{max}}$) be the minimum (resp. maximum) value of $t$ on $\LL$. 
We consider $m$-planes $t=t_k, (k=1,\cdots,m)$ where $t_{\mathrm{min}}<t_1<\cdots<t_m<t_{\mathrm{max}}$.
Define a height function $\pi$ on $X$ by $\pi(s)=t(\vec{x}(s)),~(s\in X)$.
The inverse function $\pi^{-1}(t)$ is a multi-valued function on $\mathbb{R}$.
So set $(\pi^{-1})(t_k)=\{s_k^1,\cdots,s_k^{n(t_k)}\}$, where $n(t_k)$ denotes the number of points on the section $t=t_k$ of the link $\LL$.
For $1\le i\le j\le n(t_k)$, set $z_{ij}(t_k)=z\bigl\{
\vec{x}(s^i_k)-\vec{x}(s^j_k)\bigr\}$. Define the collection of all the pairings by $P=\{(i_1,j_1),(i_2,j_2),\cdots,(i_m,j_m):1\le i_k< j_k\le n(t_k)~~(k=1,\cdots,m)\}$. For a pairing $p\in P$, write $D_p$ for the chord diagram of degree $m$ obtained by joining $s_k^{i_k}$ and $s_k^{j_k}$ by chords on $X$ $(k=1,\cdots,m)$. It is to be regarded as an element of the quotient $\hat{\mathcal{A}}$.
The Kontsevich integral is defined as follows:
\begin{eqnarray}\label{eqn:KontsevichIntegral}
  Z(\LL)=\sum_{m=0}^\infty \frac{1}{(i\pi)^m}
  \int_{t_{\mathrm{max}}>t_1>\cdots>t_m>t_{\mathrm{min}}}  \sum_{p\in P} 
  D_p \prod_{k=1}^m\{\epsilon~d\log(z_{i_kj_k}(t_k))\},
\end{eqnarray}
where the signature $\epsilon$ in front of $d\log(z_{i_kj_k}(t_k))$ is $+1$ if the two orientations of $\LL$ at $\vec{x}(s_k^{i_k})$ and $\vec{x}(s_k^{j_k})$ are the same with respect to $t$-axis and $-1$ if they are different. 
Notice we have used the slightly different normalization from \cite{BN},\cite{Ko}.
\end{definition}

\begin{figure}[h]
\setlength{\unitlength}{1cm}
\begin{picture}(15,1.5)(1.4,0)
\put(6.7,-0.2){\includegraphics{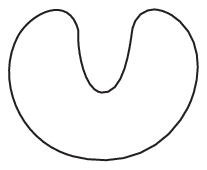}}
\end{picture}
\caption{$U_0$}
\label{fig:Unknot}
\end{figure}

\begin{definition}\hspace{-5pt}{\bf.}\hspace{6pt}\label{def:modifiedKontsevich}
(Modified Kontsevich Integral)
Define $Z_W(\LL)$ by
\begin{eqnarray}\label{eqn:ModifiedKontsevichIntegral1}
  Z_W(\LL)=\hat{W}_{su(N)}(Z(\LL)),
\end{eqnarray}
where $\hat{W}_{su(N)}$ is the renormalized version of $W_{su(N)}$ to be compatible with the framing independence (see \cite{BN} page 426).
It is known that the Kontsevich integral is invariant under only horizontal deformation of $\LL$. So we define the modified Kontsevich integral by
\begin{eqnarray}\label{eqn:ModifiedKontsevichIntegral2}
  \hat{Z}_W(\LL)=Z_W(\LL)Z_W(U_0)^{-m(\LL)},
\end{eqnarray}
where $m(\LL)$ denotes the number of maximal points of link $\LL$ and $U_0$ is a knot given in Figure \ref{fig:Unknot}. It is known that $\hat{Z}_W(\LL)$ is invariant under arbitrary deformations of the link $\LL$. We remark $\hat{Z}_W(\LL)$ is a formal power series with respect to $x$.
\end{definition} 

\subsection{The Kontsevich integral up to degree four}

In this section, we shall prove that the modified Kontsevich integral of links up to degree four can be expressed by some link invariants $v_{1},v_{2},v_{3.1},v_{3.2},v_{4.1},v_{4.2},v_{4.3}$ and $v_{4.4}$ ({Theorem \ref{th:Kontsevich}}).
\begin{definition}\hspace{-5pt}{\bf.}\hspace{6pt}
Let $f(x)=\sum_{n=0}^\infty a_n x^n$ be a power series of $x$.
Define $\bigl[f(x)\bigr]^{(k)}$ by
\begin{eqnarray*}
  \bigl[f(x)\bigr]^{(k)}=\sum_{n=0}^k a_n x^n.
\end{eqnarray*}
\end{definition}
\begin{definition}\hspace{-5pt}{\bf.}\hspace{6pt}\label{df:Kontsevich}
Let $\LL$ be a link and $D$ a chord diagram of degree $m$ without any isolated 
chord (namely, $D$ cannot be decomposed as the product of \marua~ and a chord diagram). Define $\hidari \LL,D\migi$ by
\begin{eqnarray*}
  \hidari \LL,D\migi=\frac{1}{(i\pi)^m}
  \int_{t_{\mathrm{max}}>t_1>\cdots>t_m>t_{\mathrm{min}}}
  \sum_{p\in P}\prod_{k=1}^m\bigl\{\epsilon~d\log(z_{i_kj_k}(t_k))\bigr\}
  \Theta(D_p,D),
\end{eqnarray*}
where $D_p$ denotes the chord diagram corresponding to the pairing $p\in P$. The sum is taken over all the pairings $p\in P$. $\Theta(D_p,D)$ is defined by
\begin{eqnarray*}
  \Theta(D_p,D)=
  \left\{\begin{array}{cl}
  1 & \mbox{if}~~D_p=D \\
  0 & \mbox{if}~~D_p\ne D 
  \end{array}\right..
\end{eqnarray*}
More generally, for a formal linear combinaiton of chord diagrams 
$\displaystyle \sum_{i}b_i D_i$ ($b_i\in {\mathbb C}$, $D_i$ is a chord diagram without any isolated chord), set
\begin{eqnarray*}
  \hidari \LL,\sum_{i}b_i D_i\migi=\sum_{i}b_i \hidari \LL,D_i\migi.
\end{eqnarray*}
\end{definition}

\begin{theorem}\hspace{-5pt}{\bf.}\hspace{6pt}\label{th:Kontsevich}
Let $\LL=\{\K_1,\cdots,\K_n\}$ be a link where $\K_i$ denotes each component of the link $\LL$ $(i=1,\cdots,n)$.
The modified Kontsevich integral up to degree four $\bigl[\hat{Z}_W(\LL)\bigr]^{(4)}$ can be expressed as
\begin{eqnarray}\label{eqn:FormKontsevichIntegral}
  \bigl[\hat{Z}_W(\LL)\bigr]^{(4)}
  =W_{su(N)}^{(4)}\Biggl(\biggl\{\exp\Bigl(
  \sum_{D\in\mathfrak{D}_K}D~w(D:\LL)
  \Bigr)\biggr\}\biggl\{\sum_{D\in\mathfrak{D}_L}D~w(D:\LL) 
  \biggr\}\Biggr),
\end{eqnarray}
where $W_{su(N)}^{(4)}(D)=\bigl[W_{su(N)}(D)\bigr]^{(4)}$ , the sum is taken over the following CC diagrams:
\begin{eqnarray}\label{eq:CC diagrams}
  &&\hspace{-1.1cm}
  \mathfrak{D}_K=\Bigl\{\marubc,\marudg,\marufa,\marufb\Bigr\},
  \nonumber\\
  &&\hspace{-1.1cm}
  \mathfrak{D}_L=\Bigl\{\maru,\maruca,\maruec,\maruef,\marufc,\marufd
  \nonumber\\
  &&\hspace{0.5cm}\marufe,\maruff,\marufg,\marufh,\marufi,\marufk,\marufj\Bigr\}
\end{eqnarray}
and
\begin{eqnarray}\label{eq:w(D,L)}
  &&\hspace{-0.5cm}
  \bullet~~w\Bigl(\marubc:\LL\Bigr)=(-\frac{1}{2})\sum_{i=1}^nv_2(\K_i),
   \hspace{0.5cm}
  \bullet~~w\Bigl(\marudg:\LL\Bigr)=(-\frac{1}{2})^2\sum_{i=1}^nv_{3.1}(\K_i),
  \nonumber\\
  &&\hspace{-0.5cm}
  \bullet~~w\Bigl(\marufa:\LL\Bigr)=(-\frac{1}{2})^3\sum_{i=1}^n
  v_{4.1}(\K_i), \hspace{0.5cm}
  \bullet~~w\Bigl(\marufb:\LL\Bigr)=\sum_{i=1}^nv_{4.2}(\K_i),\nonumber\\
  &&\hspace{-0.5cm}\bullet~~w\Bigl(\maru:\LL\Bigr)=1, \hspace{0.8cm}
  \bullet~~w\Bigl(\maruca:\LL\Bigr)
  =\sum_{1\le i< j\le n}\frac{1}{2}\bigl(v_{1}(\{\K_i,\K_j\})\bigr)^2 ,
  \nonumber\\
  &&\hspace{-0.5cm}\bullet~~w\Bigl(\maruec:\LL\Bigr)
  =\sum_{1\le i< j\le n}\frac{1}{3!}
  \bigl(v_{1}(\{\K_i,\K_j\})\bigr)^3,\nonumber\\
  &&\hspace{-0.5cm}\bullet~~w\Bigl(\maruef:\LL\Bigr)
  =(-\frac{1}{2})\sum_{1\le i< j\le n}
  v_{3.2}(\{\K_i,\K_j\}),\nonumber\\
  &&\hspace{-0.5cm}\bullet~~w\Bigl(\marufc:\LL\Bigr)
  =\sum_{1\le i<j<k\le n}v_{1}(\{\K_i,\K_j\})
  v_{1}(\{\K_j,\K_k\})
  v_{1}(\{\K_k,\K_i\}),
  \nonumber\\
  &&\hspace{-0.5cm}\bullet~~w\Bigl(\marufd:\LL\Bigr)
  =\sum_{1\le i< j\le n}\frac{1}{4!}
  \bigl(v_{1}(\{\K_i,\K_j\})\bigr)^4,\nonumber\\
  &&\hspace{-0.5cm}\bullet~~w\Bigl(\marufe:\LL\Bigr)
  =(-\frac{1}{2})\sum_{1\le i< j\le n}\frac{1}{2}
  v_{1}(\{\K_i,\K_j\})
  v_{3.2}(\{\K_i,\K_j\}),\nonumber\\
  &&\hspace{-0.5cm}\bullet~~w\Bigl(\maruff:\LL\Bigr)
  =(-\frac{1}{2})^2\sum_{1\le i< j\le n}
  v_{4.3}(\{\K_i,\K_j\}),\nonumber\\
  &&\hspace{-0.5cm}\bullet~~w\Bigl(\marufg:\LL\Bigr)
  =\sum_{
  {{\scriptstyle 1\le i< j<k\le n}\atop
  {\scriptstyle 1\le j<i<k\le n}}\atop
  {\scriptstyle 1\le j<k<i\le n}
  }
  \frac{1}{2}\bigl(v_{1}(\{\K_i,\K_j\})\bigr)^2
  \frac{1}{2}\bigl(v_{1}(\{\K_i,\K_k\})\bigr)^2,\nonumber\\
  &&\hspace{-0.5cm}\bullet~~w\Bigl(\marufh:\LL\Bigr)
  =\sum_{
  {{\scriptstyle 1\le i< j<k\le n}\atop
  {\scriptstyle 1\le j<i<k\le n}}\atop
  {\scriptstyle 1\le j<k<i\le n}
  }v_{1}(\{\K_i,\K_j\})
  v_{1}(\{\K_i,\K_k\})
  \frac{1}{2}\bigl(v_{1}(\{\K_j,\K_k\})\bigr)^2,\nonumber\\
  &&\hspace{-0.5cm}\bullet~~w\Bigl(\marufi:\LL\Bigr)
  =(-\frac{1}{2})\sum_{1\le i<j<k\le n}v_{4.4}(\{\K_i,\K_j,\K_k\}),
  \nonumber\\
  &&\hspace{-0.5cm}\bullet~~w\Bigl(\marufk:\LL\Bigr)
  =\sum_{
  {{\scriptstyle 1\le i< j<k<l\le n}\atop
  {\scriptstyle 1\le i<k<j<l\le n}}\atop
  {\scriptstyle 1\le i<k<l<j\le n}
  }
  \frac{1}{2}(v_{1}(\{\K_i,\K_j\}))^2
  \frac{1}{2}(v_{1}(\{\K_k,\K_l\}))^2,\nonumber\\
  &&\hspace{-0.5cm}\bullet~~w\Bigl(\marufj:\LL\Bigr)\nonumber\\
  &&\nonumber\\
  &&\hspace{0.5cm}=\sum_{
  {{\scriptstyle 1\le i< k<j<l\le n}\atop
  {\scriptstyle 1\le i<j<k<l\le n}}\atop
  {\scriptstyle 1\le i<j<l<k\le n}
  }
  v_{1}(\{\K_i,\K_j\})
  v_{1}(\{\K_j,\K_k\})
  v_{1}(\{\K_k,\K_l\})
  v_{1}(\{\K_l,\K_i\}),\nonumber\\
\end{eqnarray}
and $v_{1},v_{2},v_{3.1},v_{3.2},v_{4.1},v_{4.2},v_{4.3},v_{4.4}$ are given as follows:
\begin{eqnarray}
  &&\bullet~~v_{1}(\{\K_i,\K_j\})
  =\hidari \{\K_i,\K_j\},\marucb\migi,\label{eq:v12}\\
  &&\bullet~~v_2(\K_i)
  =\hidari \K_i,\marubb\migi-\frac{1}{6}m(\K_i),\label{eq:v2}\\
  &&\bullet~~v_{3.1}(\K_i)=\hidari \K_i,\marudd+2\marude\migi,\label{eq:v31}\\
  &&\bullet~~v_{3,2}(\{\K_i,\K_j\})=\hidari \{\K_i,\K_j\}
  ,\marueb+\marued\migi,\label{eq:v32}\\
  &&\bullet~~v_{4.1}(\K_i)=\hidari \K_i,\maruga+\marugb
  +2\marugc+4\marugd+5\maruge+7\marugf\migi\nonumber\\
  &&\hspace{9.5cm}
  +\frac{1}{360}m(\K_i),\nonumber\\
  &&\bullet~~v_{4.2}(\K_i)=\hidari \K_i,\marugd+\maruge+\marugf\migi
  -\frac{1}{360}m(\K_i),\nonumber\\
  &&\bullet~~
  v_{4.3}(\{\K_i,\K_j\})=\hidari \{\K_i,\K_j\},\maruha+\maruhb
  +2\maruhc \nonumber\\
  &&\hspace{5cm}+\maruhd+\maruhe+\maruhf\migi,\nonumber\\
  &&\bullet~~v_{4.4}(\{\K_i,\K_j,\K_k\})
  =\hidari \{\K_i,\K_j,\K_k\},\maruhj+\maruhk+\maruhl\migi,\label{eq:v44}
\end{eqnarray}
where $m(\K_i)$ is the number of maximal points of $\K_i$. 

Moreover 
$v_{1},v_{2},v_{3.1},v_{3.2},v_{4.1},v_{4.2},v_{4.3},v_{4.4}$ are link invariants.
\end{theorem}
{\bfseries Proof of Theorem \ref{th:Kontsevich}.} The computation is long but straightforward. Kontsevich integral (\ref{eqn:KontsevichIntegral}) can be rewritten in the following form: 
\begin{eqnarray*}
  Z(\LL)=\sum_{m=0}^\infty\sum_{D\in {\mathfrak D}_m}D
  \hidari \LL,D\migi,
\end{eqnarray*}
where ${\mathfrak D}_m$ denotes the set of all chord diagrams of degree $m$ 
which have just $n$ circles and no isolated chord.
From (\ref{eqn:ModifiedKontsevichIntegral1}), we have
\begin{eqnarray}\label{eqn:prooftheorem1}
  \bigl[Z_W(\LL)\bigr]^{(4)}
  =\hat{W}_{su(N)}\Biggl(\sum_{m=0}^4\sum_{D\in \bar{\mathfrak D}_m}D
  \hidari \LL,D\migi\Biggr),
\end{eqnarray}
where $\bar{\mathfrak D}_m=\{D\in {\mathfrak D}_m~ | ~\hat{W}_{su(N)}(D)\noeq 0\}$ and we give the table of $\bar{\mathfrak D}_m$ in Appendix A.
   
In (\ref{eqn:prooftheorem1}), we express each chord diagram D in front of $\hidari \LL,D\migi$ as a linear combination of the following CC diagrams
\begin{eqnarray}\label{eqn:CCdiagram}
  &&\Bigl\{\marubc,\maruca,\marudg,\maruec,\maruef,\marufc,\nonumber\\
  &&\hspace{0.3cm}
  \marujl, \marufa, \marufb, \maruki,\maruff,\marufe,\nonumber\\
  &&\hspace{0.3cm}\marufd,\marufi, \marufg, \marufh, \marufk, \marufj \Bigr\}
\end{eqnarray}
regarded as an element in $\hat{\mathcal{A}}$ using Appendix B.

Next, we compute each cofficient of the CC diagram in (\ref{eqn:CCdiagram}). 
For example, the cofficient of \marujl ~is 
\begin{eqnarray*}
  \lefteqn{
  \Bigl(\hbox{the cofficient of} ~\bigl(-\frac{1}{2}\bigr)^2\marujl\Bigr)}
  \hspace{1cm}\\ 
  &=&\sum_{i=1}^n
  \hidari \K_i, \marujk+\marugb+\marugc+2\marugd+2\maruge+3\marugf\migi
  \nonumber\\
  &&+\sum_{i<j}\hidari \{\K_i,\K_j\}, \{\marubb\marubb\}\migi\nonumber\\
  &=&\sum_{i=1}^n\frac{1}{2}\hidari \K_i, \marubb\migi^2
  +\sum_{i<j}\hidari \K_i, \marubb\migi\hidari \K_j, \marubb\migi\nonumber\\
  &=&\frac{1}{2}\biggl\{\sum_{i=1}^n\hidari \K_i, \marubb\migi\biggr\}^2.
\end{eqnarray*}
See Appendix C for the other cofficients of the CC diagrams in (\ref{eqn:CCdiagram}).
Inserting these result into (\ref{eqn:ModifiedKontsevichIntegral2}) and using 
\begin{eqnarray*}
  \lefteqn{\bigl[Z_W(U_0)^{-1}\bigr]^{(4)}}\\
  &&=W_{su(N)}^{(4)}
  \Biggl(\exp\Bigl\{
  \bigl(-\frac{1}{2}\bigr)\marubc\bigl(-\frac{1}{6}\bigr)
  +\bigl(-\frac{1}{2}\bigr)^3\marufa\frac{1}{360}
  +\marufb\bigl(-\frac{1}{360}\bigr)
  \Bigr\}
  \Biggr),
\end{eqnarray*}
we have (\ref{eqn:FormKontsevichIntegral}).

Next we prove the invariace of 
$v_{1},v_{2},v_{3.1},v_{3.2},v_{4.1},v_{4.2},v_{4.3},v_{4.4}$. 
Let $\LL$ be a 1-component link (that is a knot). Then  $v_{1},v_{3.2},v_{4.3},v_{4.4}$ in (\ref{eqn:FormKontsevichIntegral}) vanish since they are defined for more than 2-component link. Since $\hat{Z}_W(L)$ is a link invariant and
\begin{eqnarray}
  &&W_{su(N)}\h\Bigl(-\frac{1}{2}\Bigr)\marubc\m,~~
  W_{su(N)}\h\Bigl(-\frac{1}{2}\Bigr)^2\marudg\m,~~\nonumber\\
  &&W_{su(N)}\h\Bigl(-\frac{1}{2}\Bigr)^3\marufa\m,~~
  W_{su(N)}\h\marufb\m,
\end{eqnarray}
are linearly independent as polynomials of $x,N$, we see that $v_{2},v_{3.1},v_{4.1},v_{4.2}$ are link invariant. We can also prove the invariance of $v_{1},v_{3.2},v_{4.3},v_{4.4}$ in the same way. $\square$

\section{Gauss Diagram Formula} 
In this section, we shall give an explicit Gauss diagram formula for the link invariants $v_{1}$, $v_{2}$, $v_{3.1}$, $v_{3.2}$, $v_{4.1}$, $v_{4.2}$, $v_{4.3}$, $v_{4.4}$ in terms of Gauss diagrams (See Theorem \ref{th:mainresult}). 
Before we state Theorem \ref{th:mainresult}, we shall fix the notation for this purpose.

\subsection{Pairing $\langle \hat{G},\hat{D}\rangle_{\chi}$}
\begin{definition}\hspace{-5pt}{\bf.}\hspace{6pt}(Link Diagram)\label{def:linkdiagram}
Let $X=\cup_{i=1}^nS^1_i$ be $n$-oriented circles and $\vec{y}:X\rightarrow\mathbb{R}^2$ an immersion.
An $n$-component {\it oriented link diagram} $L$ is its image 
$L=\{K_1,\cdots,K_n\}~(K_i=\vec{y}(S^1_i))$ together with the information of overpass or underpass at each crossing. We write the information of each crossing as in Figure \ref{crossing}. We call $\pm1$ assigned to a crossing {\it the signature of the crossing}. We often write $\pm$ instead of $\pm1$ for the signature of the crossing.
\begin{figure}[htbp]
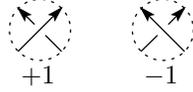

  \begin{center}
    \picff~~~~~~\picfg
  \end{center}
    \caption{the information of the crossing}
	\label{crossing}
\end{figure}
\end{definition}

\begin{definition}\hspace{-5pt}{\bf.}\hspace{6pt}\label{def;ILDiagram}(IL Diagram)
Let $D$ be a chord diagram and $C(D)$ the set of all chords of $D$. By an integer-labelling of $D$, we mean a map $\kappa:C(D)\to \mathbf{Z}$. {\it An Integer-Labeled Chord Diagram} (IL Diagram) is a pair $\{D,\kappa\}$ of a chord diagram $D$ together with an integer-labelling $\kappa$. Two IL diagram $\{D,\kappa\}$, $\{D^\prime,\kappa^\prime\}$ are regarded as equal if $D$, $D^\prime$ are equal as chord diagrams and the homeomorphism $F:D\to D^\prime$ preserves integer-labelling $\kappa^\prime(F(c))=\kappa(c)~(c\in C(D))$. $\square$ 
\end{definition}

We shall define a {\it Gauss Diagram} and {\it ML Diagram} as special cases of IL Diagrams.
\begin{definition}\hspace{-5pt}{\bf.}\hspace{6pt}\label{Gaussdiagram} (Gauss Diagram) 
An IL diagram $\{G,\epsilon\}$ is called a {\it Gauss Diagram} if $\epsilon(c)=\pm 1~(c\in C(G))$. An integer-labelling $\epsilon$ of the Gauss diagram is called a {\it signature-labelling}.

Let $\{L:a_1,\cdots,a_m\}$ be a link diagram $L$ where we select some distinct crossings $a_1,\cdots,a_m$ out of all crossings of $L$. Define a Gauss diagram $P(\{L:a_1,\cdots,a_m\})$ as follows. For each $a_i$, set $\vec{y}^{-1}(a_i)=\{s(a_i),s^\prime(a_i)\}$ as the inverse image of $a_i$. For each crossing $a_i$, we join $s(a_i),s^\prime(a_i)$ by a chord on $X$ and label this chord by the signature of $a_i$ ($i=1,\cdots,m$). We define a Gauss diagram $P(\{L:a_1,\cdots,a_m\})$ to be the result.

Specially, If $\{a_1,\cdots,a_m\}$ are all the crossings of $L$ (this means we select all the crossings of $L$) , we write
$G(L)=P(\{L:a_1,\cdots,a_m\})$ and call it the Gauss diagram of $L$.

For example, see Figure \ref{fig:Gaussdiagram}.
\begin{figure}[h]
\begin{center}
\setlength{\unitlength}{1cm}
\begin{picture}(8,5)(4,1.5)
\put(2.5,2.4){\includegraphics{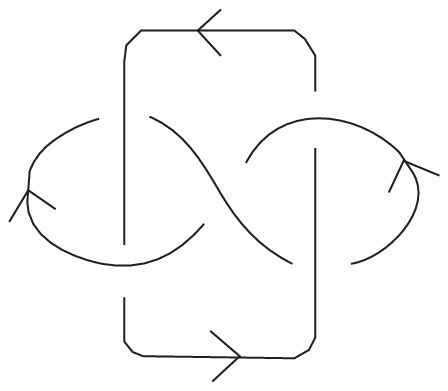} }
\put(3.2,5.3){$a_1$}
\put(3.2,3.3){$a_2$}
\put(8,4){\pichc}
\put(3.3,2){$\{K_1,K_2:a_1,a_2\}$}
\put(6.8,2.3){\picfc} 
\put(9,2.5){$P\bigl(\{K_1,K_2:a_1,a_2\}\bigr)$}
\end{picture} 
\end{center}
\
\begin{center}
\setlength{\unitlength}{1cm}
\begin{picture}(8,5)(4,1.5)
\put(2.5,2.4){\includegraphics{painte.eps} }
\put(8,4){\picdb}
\put(3.8,2){$\{K_1,K_2\}$}
\put(6.8,2.3){\picfc} 
\put(10,2.5){$G\bigl(\{K_1,K_2\}\bigr)$}
\end{picture} 
\end{center}
\caption{}
\label{fig:Gaussdiagram}
\end{figure}
\end{definition}

\begin{definition}\hspace{-5pt}{\bf.}\hspace{6pt}\label{MLdiagram}(ML Diagram)
An IL diagram $\{D,m\}$ is called a {\it Multiplicity-Labeled Diagram} (ML Diagram) if $m(c)=1,2~(c\in C(D))$. In figures, we draw a chord $c$ with $m(c)=1$ by a thin line and a chord $c$ with $m(c)=2$ by a thin line with a letter "2" as follows:
\begin{eqnarray*}
  \picdc~~\hbox{m(c)=1},\hspace{1.5cm}\picdd~~\hbox{m(c)=2}.
\end{eqnarray*}
We give two examples of ML diagrams,
\begin{eqnarray*}
  \marugg,\hspace{1.5cm}\maruhg.
\end{eqnarray*}
\end{definition}

\begin{definition}\hspace{-5pt}{\bf.}\hspace{6pt}\label{def:crossing number}
Let $\hat{G}=\{G,\epsilon\}$ be a Gauss diagram and $\hat{D}=\{D,m\}$ a ML diagram. Let $\psi:D\to G$ be an embedding of $D$ into $G$ which maps the circles of $D$ to those of $G$ preserving the orientations and each chord of D to a chord of $G$. Let $C(G)$ be the set of all chords of $G$. For $\psi$, define a map $\kappa_\psi:C(G)\to\{0,1,2\}$ by
\begin{eqnarray*}
  \kappa_\psi(c)=\left\{
  \begin{array}{cl}
  m(\psi^{-1}(c)) & \mbox{if}~~c\in \psi(D) \\
  0 & \mbox{if}~~c\noin \psi(D) 
  \end{array}\right.
\end{eqnarray*}
Two embedding $\psi,\varphi$ are said to be equal if $\kappa_\psi=\kappa_{\varphi}$. The equivalence class of an embedding $\psi$ is denoted by $[\psi]$.

Let $C(D)$ be the set of all chords of $D$. Define $\mathcal{E}([\psi])$ by
\begin{eqnarray*}
  \mathcal{E}([\psi])=\prod_{c\in C(D)}\bigl\{\epsilon(\psi(c))\bigr\}^{m(c)},
\end{eqnarray*}
where the product is taken over all chords of $D$. Notice this definition is well defined.

Define a pairing of a Gauss diagram and ML diagram $\langle \hat{G},\hat{D}\rangle_{\chi}$ by
\begin{eqnarray*}
  \langle \hat{G},\hat{D}\rangle_{\chi}=\sum_{[\psi]}\mathcal{E}([\psi]),
\end{eqnarray*}
where the sum is taken over all the distinct equivalence classes $[\psi]$.

Let $\hat{G}_i$ be a Gauss diagram and $\hat{D}_i$ a ML diagram.
More generally, for formal linear combinations 
$\displaystyle \sum_i b_i~\hat{G}_i$ and
$\displaystyle \sum_j c_j~\hat{D}_j$
($b_i, c_j\in \mathbf{C}$), set
\begin{eqnarray*}
  \Bigl\langle \sum_i b_i~\hat{G}_i,\sum_j c_j~\hat{D}_j\Bigr\rangle_{\chi}
  =\sum_i\sum_j b_i~c_j~\langle \hat{G}_i,\hat{D}_j\rangle_{\chi},
\end{eqnarray*}

\end{definition}

\begin{example}\hspace{-5pt}{\bf.}\hspace{6pt}
\begin{eqnarray*}
  &&\Bigl\langle \picde,\marubb\Bigr\rangle_{\chi}
  =\epsilon_1\epsilon_2+\epsilon_1\epsilon_3,\nonumber\\
  &&\Bigl\langle \picde,\marudp\Bigr\rangle_{\chi}
  =(\epsilon_1)^2\epsilon_2+\epsilon_1(\epsilon_2)^2
  +(\epsilon_1)^2\epsilon_3+\epsilon_1(\epsilon_3)^2,\nonumber\\
  &&\Bigl\langle \picdp,\marueb\Bigr\rangle_{\chi}
  =(\epsilon_1+\epsilon_2+\epsilon_7)\epsilon_4\epsilon_5,\nonumber\\
  &&\Bigl\langle \picdp,\picdj\Bigr\rangle_{\chi}
  =(\epsilon_1+\epsilon_2+\epsilon_7)
  \bigl((\epsilon_4)^2\epsilon_5+\epsilon_4(\epsilon_5)^2\bigr),
\end{eqnarray*}
where $\epsilon_i=\pm1$ denotes the signature-labelling.
\end{example}

\subsection{Gauss diagram formula}

Next, we shall introduce a concept for splitting the crossings of a link diagram. 
\begin{definition}\hspace{-5pt}{\bf.}\hspace{6pt}(Splitting of the crossings)\label{def:splitting}
Let $\{L:a_1,\cdots,a_m\}$ be a link diagram $L$ where we select some distinct crossings $a_1,\cdots,a_m$ out of all crossings of $L$. 
By a splitting information, we mean a finite sequence $[s_1,\cdots,s_m]$  
$(s_i=\alpha,\beta,\gamma)$, where $\alpha, \beta, \gamma$ are formal letters. 
For example, $[\alpha,\beta,\alpha,\gamma]$ $(m=4)$.

We shall define a link diagram $Q(\{L:a_1,\cdots,a_m\},[s_1,\cdots,s_m])$ as follows.
For $1\le i\le m$, we replace each crossing $a_i$ by
\begin{eqnarray*} 
  \picgc\to\left\{
  \begin{array}{cc}
  \picgd & \mbox{if}~~s_i=\alpha \\
  &\\
  \picge & \mbox{if}~~s_i=\beta \\
  &\\
  \picgf & \mbox{if}~~s_i=\gamma \\
  &\\
  &
  \end{array}
  \right.
\end{eqnarray*}
and give any orientation to the resulting diagram.

Define $Q(\{L:a_1,\cdots,a_m\},[s_1,\cdots,s_m])$ to be the resulting oriented link diagram.
We remark that our caluculation in the sequel does not depend on the paticular choice of the orientation of $Q(\{L:a_1,\cdots,a_m\},[s_1,\cdots,s_m])$.

More generally, for a formal linear combination of splitting information 
$\sum_i b_i~\delta_i$ 
$(b_i\in \mathbf{C},~~\delta_i=[s_1^i,\cdots,s_m^i])$, set
\begin{eqnarray*}
  Q(\{L:a_1,\cdots,a_m\},\sum_i b_i~\delta_i)
  =\sum_i b_i~Q(\{L:a_1,\cdots,a_m\},\delta_i).
\end{eqnarray*}
\end{definition}

\begin{example}\hspace{-5pt}{\bf.}\hspace{6pt}We give a trivial example:
\begin{eqnarray*}
  Q(\{L:a\},[\gamma])=L
\end{eqnarray*}
\end{example}

\begin{example}\hspace{-5pt}{\bf.}\hspace{6pt}
We give two nontrivial examples (see Figure \ref{fig:split1} and Figure \ref{fig:split2}).
As for Figure \ref{fig:split2}, notice the orientation of knots is partly changed.  

\begin{figure}[h]
\setlength{\unitlength}{1cm}
\begin{picture}(8,4.5)(1.7,2)
\put(3,3){\includegraphics{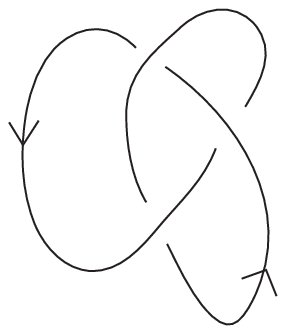} }
\put(5.5,5){$a$}
\put(4.2,2.3){$\{K:a\}$}
\put(7,2){\picfc}
\put(10,3){\includegraphics{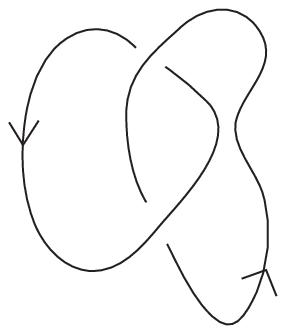}} 
\put(12.5,5){$a$}
\put(10,2.3){$Q(\{K:a\},[\alpha])$}
\end{picture}
\caption{$\{K:a\} \to Q(\{K:a\},[\alpha])$}
\label{fig:split1}
\end{figure}

\begin{figure}[h]
\setlength{\unitlength}{1cm} 
\begin{picture}(8,4.5)(1.7,2)
\put(2.1,2.2){\includegraphics{painte.eps} }
\put(3.0,2){$\{K_1,K_2:a_1,a_2,a_3\}$}
\put(8.5,2.2){\includegraphics{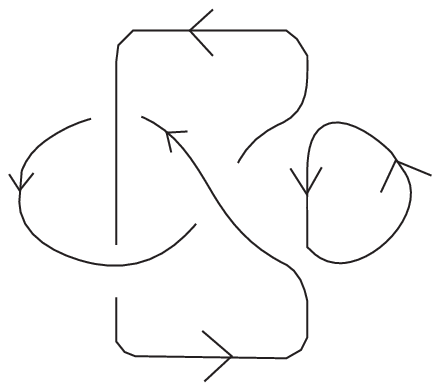}}
\put(9,2){$Q(\{K_1,K_2:a_1,a_2,a_3\},[\beta,\beta,\gamma])$}
\put(7,2){\picfc}
\put(1.7,4.2){$K_1$}
\put(4.2,6.1){$K_2$}
\put(5.4,3){$a_1$} 
\put(5.4,5){$a_2$}
\put(2.8,5.1){$a_3$}
\end{picture}
\caption{$\{K_1,K_2:a_1,a_2,a_3\} 
\to Q(\{K_1,K_2:a_1,a_2,a_3\},[\beta,\beta,\gamma])$}
\label{fig:split2}
\end{figure}

\end{example}

\begin{definition}\hspace{-5pt}{\bf.}\hspace{6pt}Let $L=\{K_1,K_2,\cdots,K_n\}$ be an $n$-component link diagram. Define $S(L)$ to be the formal sum of each component $K_i$ 
\begin{eqnarray*}
   S(L)=\sum_{i=1}^n K_i.
\end{eqnarray*}
For example, see Figure \ref{fig:SofL}.

\begin{figure}[h]
\setlength{\unitlength}{1cm} 
\begin{picture}(8,5)(0,2)
\put(1.7,2.5){\includegraphics{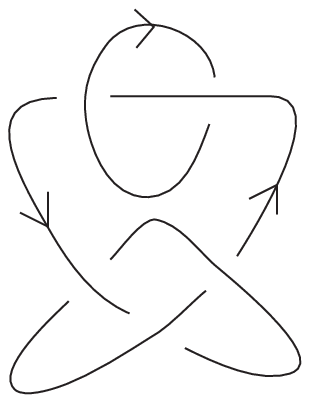} }
\put(1,4){$L=$}
\put(7,2.5){\includegraphics{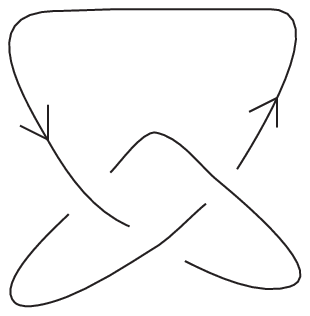}}
\put(10,4){\Large $+$}
\put(11,3){\includegraphics{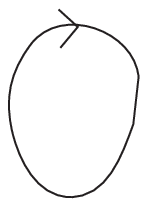}}
\put(4.8,2.3){\picfc}
\put(6,4){$S(L)=$}
\end{picture}
\caption{$L=\{K_1,K_2\}$ and $S(L)=K_1+K_2$}
\label{fig:SofL}
\end{figure}
\end{definition}

\begin{definition}\hspace{-5pt}{\bf.}\hspace{6pt}($L$ and $\alpha(L)$)\label{def:alphaL}
Let $L$ be an $n$-component link diagram. Let $\alpha(L)$ be a trivial link diagram of $n$-separated trivial knots which is obtained by swiching the signature of the crossings of $L$ properly (see Figure \ref{fig:unknot}). There are several ways to obtain $\alpha(L)$ from $L$. So $\alpha(L)$ cannot be uniquely determined from $L$. But the caluculation in the sequel does not depend on the way we choose.
\begin{figure}[h]
\setlength{\unitlength}{1cm}
\begin{picture}(8,5)(1,2)
\put(3,2.3){\includegraphics{paintn.eps} }
\put(9,2.3){\includegraphics{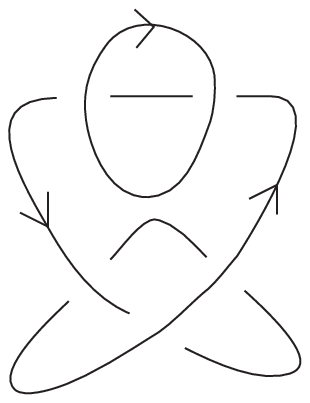}}
\put(4.3,2){$L$} 
\put(7,2){\picfc} 
\put(10.4,2){$\alpha(L)$}
\end{picture}
\caption{$L \to \alpha(L)$}
\label{fig:unknot}
\end{figure}
\end{definition}

\begin{definition}\hspace{-5pt}{\bf.}\hspace{6pt}Let $L_i$ be a link diagram.  
For a formal linear combination $\displaystyle \sum_i b_iL_i$ 
$(b_i\in \mathbf{C})$, we extend the definition of $G,\alpha,S$ by
\begin{eqnarray*}
  &&G(\sum_i b_i~L_i)=\sum_i b_i~G(L_i),~~~~
  \alpha(\sum_i b_i~L_i)=\sum_i b_i~\alpha(L_i),\\
  &&S(\sum_i b_i~L_i)=\sum_i b_i~S(L_i).\\
\end{eqnarray*}
\end{definition}

\begin{theorem}\hspace{-5pt}{\bf.}\hspace{6pt}\label{th:mainresult}(Gauss diagram formula)
Let $K$, $\{K_1,K_2\}$, $\{K_1,K_2,K_3\}$ to be the link diagrams which correspond to links $\K$, $\{\K_1,\K_2\}$, $\{\K_1,\K_2,\K_3\}$ respectively. The link invariants $v_{1},v_{2},v_{3.1},v_{3.2},v_{4.1},v_{4.2},v_{4.3}$ and $v_{4.4}$ have the explicit combinatorial expressions as follows: 
\begin{eqnarray}
  &&\hspace{-0.7cm}\bullet~~v_{1}(\{\K_1,\K_2\})
  =\hidarii G(\{K_1,K_2\}),\marucb\migii_{\chi},\label{eqn:main result12}\\
  &&\hspace{-0.7cm}\bullet~~v_2(\K)
  =-\frac{1}{6}+\hidarii \bar{G}(K),\marubb\migii_{\chi},
  \label{eqn:main result2} 
  \\
  &&\hspace{-0.7cm}\bullet~~v_{3.1}(\K)=
  \hidarii G(K),~~2\marude+\marudd+\frac{1}{2}\marudp\migii_{\chi} 
  -I_{3.1}(K),\label{eqn:main result31} \\
  &&\hspace{-0.7cm}
  \bullet~~v_{3.2}(\{\K_1,\K_2\})=\hidarii G(\{K_1,K_2\}),~~\marueb
  +\marued+\frac{1}{3}\marucb\migii_{\chi}\nonumber\\
  &&\hspace{7cm}-I_{3.2}(\{K_1,K_2\})\label{eqn:main result32} \\
  &&\hspace{-0.7cm}\bullet~~v_{4.1}(\K)= 
  \hidarii \bar{G}(K),~~\maruga+\marugb+2\marugc+4\marugd
  +5\maruge+7\marugf\migii_{\chi}\nonumber\\
  &&\hspace{3cm}+\hidarii \bar{G}(K),
  \frac{1}{6}\marubb+\frac{1}{2}\marugg+2\marugh+2\marugi\migii_{\chi} 
  \nonumber\\
  &&\hspace{3cm}
  -I_{4.1.1}(K)-I_{4.1.2}(K)+\frac{1}{360},\label{eqn:main result41}\\
  &&\hspace{-0.7cm}\bullet~~v_{4.2}(\K)
  =\hidarii \bar{G}(K),~~\marugd+\maruge+\marugf
  +\frac{1}{2}\marugh-\frac{1}{6}\marubb\migii_{\chi}
  \nonumber\\
  &&\hspace{7cm}
  -I_{4.2}(K)-\frac{1}{360},
  \label{eqn:main result42}\\ \nonumber\\ 
  &&\hspace{-0.7cm}\bullet~~v_{4.3}(\{\K_1,\K_2\})\nonumber\\
  &&\hspace{-0.1cm}
  =\hidarii \bar{G}(\{K_1,K_2\}),
  ~~\maruha+\maruhb+2\maruhc+\maruhd\migii_{\chi}
  \nonumber\\
  &&\hspace{-0.1cm}
  +\hidarii \bar{G}(\{K_1,K_2\}),~~
  \maruhe+\maruhf+\frac{1}{2}\maruhg+\frac{1}{2}\maruhh\migii_{\chi} 
  \nonumber\\
  &&\hspace{3cm}-I_{4.3.1}(\{K_1,K_2\})-I_{4.3.2}(\{K_1,K_2\}),
  \label{eqn:main result43}\\
  &&\hspace{-0.7cm}\bullet~~v_{4.4}(\{\K_1,\K_2,\K_3\})
  =\hidarii \bar{G}(\{K_1,K_2,K_3\}),~~\maruhj+\maruhk+\maruhl\migii_{\chi} 
  \nonumber\\
  &&\hspace{3cm}
  -I_{4.4}(\{K_1,K_2,K_3\}),
  \label{eqn:main result44}
\end{eqnarray}
where $\bar{G}(L)=G(L)-G(\alpha(L))$. 
Set $R=G\circ \alpha \circ S \circ Q$ and 
$\bar{P}(\{L:a,b\})=P(\{L:a,b\})-P(\{\alpha(L):a,b\})$.
Here $I_{3.1},I_{3.2},I_{4.1.1},I_{4.1.2}
,I_{4.2},I_{4.3.1},I_{4.3.2},I_{4.4}$ are given as follows:
\begin{eqnarray}
  &&\hspace{-0.7cm}\bullet~~
  I_{3.1}(K)=\sum_{a}\hidarii P(\{K:a\}),\marua\migii_{\chi}
  \hidarii R(\{K:a\},[\gamma]-[\alpha]),
  \marubb\migii_{\chi},\label{eqn:I31}\\
  &&\nonumber\\
  &&\hspace{-0.7cm}\bullet~~I_{3.2}(\{K_1,K_2\})
  =\sum_{a}\hidarii (P(\{K_1,K_2:a\}),\marucb\migii_{\chi}
  \nonumber\\
  &&\hspace{4cm}\times
  \hidarii R(\{K_1,K_2:a\},[\alpha]-[\gamma]),\marubb\migii_{\chi}, 
  \nonumber\\ 
  &&\hspace{-0.7cm}\bullet~~I_{4.1.1}(K)
  =\sum_{(a_1,a_2)}\hidarii \bar{P}(\{K:a_1,a_2\}),\marubb\migii_{\chi}
  \nonumber\\
  &&\hspace{2cm}\times
  \hidarii R(\{K:a_1,a_2\},
  3[\gamma,\gamma]-2[\alpha,\gamma]-2[\gamma,\alpha]+[\beta,\beta])
  ,\marubb\migii_{\chi},
  \nonumber\\
  &&\hspace{-0.7cm}\bullet~~I_{4.1.2}(K)
  =\sum_{(a_1,a_2)}\hidarii \bar{P}(\{K:a_1,a_2\}),\maruba\migii_{\chi}
  \nonumber\\
  &&\hspace{2cm}\times
  \hidarii R(\{K:a_1,a_2\},[\gamma,\gamma]-[\alpha,\gamma]
  -[\gamma,\alpha]+[\alpha,\alpha]),
  \marubb\migii_{\chi}, 
  \nonumber\\
  &&\hspace{-0.7cm}\bullet~~I_{4.2}(K)
  =\sum_{(a_1,a_2)}\hidarii \bar{P}(\{K:a_1,a_2\}),\marubb\migii_{\chi}
  \nonumber\\
  &&\hspace{2cm}\times
  \hidarii R(\{K:a_1,a_2\},[\gamma,\gamma]-[\alpha,\gamma]
  -[\gamma,\alpha]+[\beta,\beta])
  ,\marubb\migii_{\chi}, 
  \nonumber\\
  &&\hspace{-0.7cm}\bullet~~I_{4.3.1}(\{K_1,K_2\})=
  \sum_{(a_1,a_2)}\hidarii \bar{P}(\{K_1,K_2:a_1,a_2\}),
  \maruca\migii_{\chi}
  \nonumber\\
  &&\hspace{2cm}\times
  \hidarii R(\{K_1,K_2:a_1,a_2\},[\gamma,\gamma]-[\beta,\beta])
  ,\marubb\migii_{\chi}, 
  \nonumber\\ 
  &&\nonumber\\ 
  &&\hspace{-0.7cm}\bullet~~I_{4.3.2}(\{K_1,K_2\})
  =\sum_{(a_1,a_2)}\hidarii \bar{P}(\{K_1,K_2:a_1,a_2\}),
  \maruhi\migii_{\chi}
  \nonumber\\ 
  &&\hspace{2cm}\times
  \hidarii R(\{K_1,K_2:a_1,a_2\},[\alpha,\gamma]+[\gamma,\alpha]
  -[\gamma,\gamma]-[\alpha,\alpha])
  ,\marubb\migii_{\chi}, 
  \nonumber\\[1cm]
  &&\hspace{-0.7cm}\bullet~~I_{4.4}(\{K_1,K_2,K_3\})
  =\sum_{(a_1,a_2)}\hidarii \bar{P}(\{K_1,K_2,K_3:a_1,a_2\}),
  \maruhm\migii_{\chi}
  \nonumber\\
  &&\hspace{1.5cm}\times
  \hidarii R(\{K_1,K_2,K_3:a_1,a_2\},[\gamma,\gamma]+[\alpha,\alpha]
  -[\alpha,\gamma]-[\gamma,\alpha])
  ,\marubb\migii_{\chi},
  \nonumber\\\nonumber
\end{eqnarray}
where the sum $\displaystyle \sum_{a}$~
$\Bigl($resp.$\displaystyle \sum_{(a_1,a_2)}\Bigr)$ is taken over all the crossings (resp. all the unordered pairs of the crossings). 
$\square$
\end{theorem}

\noindent{\it Remark.}~~
The Gauss diagram formulas in Theorem \ref{th:mainresult} is expressed
by the pairing $\langle \hat{G},\hat{D}\rangle_{\chi}$ in Definition \ref{def:crossing number}. So it is easy to compute the link invariants $v_{1},v_{2},v_{3.1},v_{3.2},$ $v_{4.1},v_{4.2},v_{4.3},v_{4.4}$ for any link. $\square$ 
\newline

See section \ref{sec:ProofTheorem2} for the proof of Theorem \ref{th:mainresult}. 

\begin{conjecture}\hspace{-5pt}{\bf.}\hspace{6pt}There exist Gauss diagram formulas for any Vassiliev invariants of any degree.
$\square$
\end{conjecture}
 
\section{Homfly Polynomial and Some Caluculations}
\subsection{Relation to Homfly Polynomial}
In this section, we shall discuss the relation between Theorem \ref{th:mainresult} and Homfly polynomial.
\begin{definition}\hspace{-5pt}{\bf.}\hspace{6pt}(Homfly polynomial)
For a link diagram $L$, the Homfly polynomial $P_L(t,z)$ is characterized by the skein relation:
\begin{eqnarray*}
 tP_{L_+}(t,z)-t^{-1}P_{L_-}(t,z)=zP_{L_0}(t,z)
\end{eqnarray*}
\begin{eqnarray}\label{eqn:skein relation}
 P_U&=&1,
\end{eqnarray}
where $U$ denotes a trivial knot. The links $L_+,L_-,L_0$ are given in Figure 
\ref{fig:skeinrelation}. $\square$
\begin{figure}[htbp]
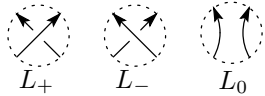

\begin{eqnarray*}
  \picfn~~~~\picfo~~~~\picfm
\end{eqnarray*}
    \caption{skein relation}
	\label{fig:skeinrelation}
\end{figure}
\end{definition}

It is known that the Kontsevich integral and the weight system of $su(N)$ gives the 
Homfly polynomial. More precisely, the following fact holds.
\begin{fact}\hspace{-5pt}{\bf.}\hspace{6pt}\label{fact:Homfly}
Let $\LL=\{\K_1,\cdots,\K_n\}$ be a link.
Define $\hat{P}_{\LL}(x,N)$ by
  \begin{eqnarray}\label{eq:Homfly}
    \hat{P}_{\LL}(x,N)=N^{n-1}\exp\Bigl(-x\frac{N^2-1}{2N}
    w(\LL)\Bigr)~
    \frac{\hat{Z}_W(\LL)}{\hat{Z}_W(U)},
  \end{eqnarray}
where $w(\LL)$ is given by
\begin{eqnarray*}
  w(\LL)=\sum_{1\le i< j\le n}v_{1}(\{\K_i,\K_j\}).
\end{eqnarray*}
Let $L=\{K_1,\cdots,K_n\}$ be the link diagram of $\LL$. Then,
\begin{eqnarray*}
  P_L(e^{\frac{Nx}{2}},e^{\frac{x}{2}}-e^{-\frac{x}{2}})=\hat{P}_{\LL}(x,N)
\end{eqnarray*}
holds. $\square$
\end{fact}

Since $v_{1},v_{2},v_{3.1},v_{3.2},v_{4.1},v_{4.2},v_{4.3},v_{4.4}$ are link invariants and depend only on its link diagrams, we write $v_{1}(\{K_i,K_j\})$ instead of $v_{1}(\{\K_i,\K_j\})$, etc. From Theorem \ref{th:mainresult} and Fact \ref{fact:Homfly}, we 
immediately obtain the following corollary.
\begin{corollary}\hspace{-5pt}{\bf.}\hspace{6pt}\label{co:Homfly}
Up to degree four, the power series expansion of Homfly polynomial with respect to $x$ has the explicit Gauss diagram formula as follows:  
\begin{eqnarray}\label{eqn:Homflyformula}
  &&\hspace{-1cm}
  \Bigl[P_L(e^{\frac{Nx}{2}},e^{\frac{x}{2}}-e^{-\frac{x}{2}})\Bigr]^{(4)}
  \nonumber\\
  &&\hspace{-0.5cm}
  =W_{su(N)}^{(4)}\Biggl(N^{n-1}\biggl\{\exp\Bigl(
  \sum_{D\in\bar{\mathfrak{D}}_K}D~u(D:L)
  \Bigr)\biggr\}\biggl\{\sum_{D\in\mathfrak{D}_L}D~w(D:L)
  \biggr\}\Biggr),
\end{eqnarray}
where $W_{su(N)}^{(4)}(D)=\bigl[W_{su(N)}(D)\bigr]^{(4)}$. The first sum 
$\displaystyle \sum_{D\in\bar{\mathfrak{D}}_K}$
is taken over the following CC diagrams:
\begin{eqnarray*}
  &&\bar{\mathfrak{D}}_K=\Bigl\{\marua,\marubc,\marudg,\marufa,\marufb\Bigr\}.
\end{eqnarray*}
The second sum $\displaystyle \sum_{D\in\mathfrak{D}_L}$ is taken over the same CC diagrams as (\ref{eq:CC diagrams}). Here $w(D:L)$ is the same as (\ref{eq:w(D,L)}) and $u(D:L)$ is given as follows:
\begin{eqnarray*}
  &&\bullet~~u\Bigl(\marua:L\Bigr)=
  -\sum_{1\le i< j\le n}v_1(\{K_i,K_j\}),\nonumber\\
  &&\bullet~~u\Bigl(\marubc:L\Bigr)=(-\frac{1}{2})\bigl\{\frac{1}{6}+\sum_{i=1}^nv_2  (K_i)\bigr\},\nonumber\\
  &&\bullet~~u\Bigl(\marudg:L\Bigr)=(-\frac{1}{2})^2\sum_{i=1}^n
  v_{3.1}(K_i),\nonumber\\
  &&\bullet~~u\Bigl(\marufa:L\Bigr)=(-\frac{1}{2})^3
  \bigl\{-\frac{1}{360}+\sum_{i=1}^nv_{4.1}(K_i)\bigr\},\nonumber\\
  &&\bullet~~u\Bigl(\marufb:L\Bigr)=\sum_{i=1}^n
  \bigl\{\frac{1}{360}+v_{4.2}(K_i)\bigl\}.\nonumber\\
\end{eqnarray*}
\end{corollary}
\rightline{$\square$}

\subsection{Some caluculations}
\begin{figure}[htbp]
	\setlength{\unitlength}{1cm}
\begin{picture}(8,8)(1.5,-1)
\put(3,2.8){\includegraphics{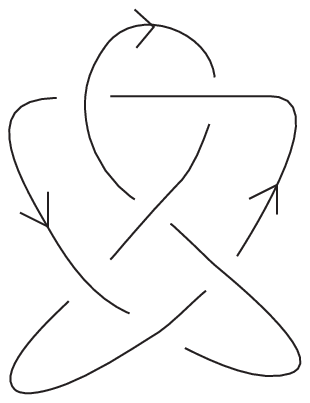} }
\put(5.3,6){\footnotesize a1}
\put(3.3,6){\footnotesize a2}
\put(3.1,3.8){\footnotesize a3}
\put(4.4,3){\footnotesize a4} 
\put(5.5,4){\footnotesize a5}
\put(4.3,5.0){\footnotesize a6}

\put(9,2.8){\includegraphics{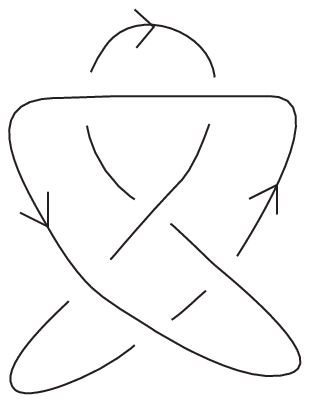} }
\put(11.3,6){\footnotesize a1}
\put(9.3,6){\footnotesize a2}
\put(9.1,3.8){\footnotesize a3}
\put(10.4,3){\footnotesize a4} 
\put(11.5,4){\footnotesize a5}
\put(10.3,5.0){\footnotesize a6}
 
\put(4.2,2.2){$K$}
\put(10.2,2.2){$\alpha(K)$}
\put(3,-0.5){\picfd}
\put(9.3,-0.5){\picfe}
\end{picture}
	\caption{the knot diagram $K$, $\alpha(K)$ 
	and their Gauss diagram $G(K)$, $G(\alpha(K))$}
	\label{fig:62}
\end{figure}

We give an example of Theorem \ref{th:mainresult} 
(Gauss diagram formula). As an example, we compute $v_{2},v_{3.1},v_{4.1},v_{4.2}$ for a knot diagram $K$ given in Figure \ref{fig:62}.
\subsubsection{$v_2(K)$}
Using the Gauss diagram $G(K)$, $G\bigl(\alpha(K)\bigr)$ in Figure \ref{fig:62}, we get
\begin{eqnarray}\label{eqn:example1}
  \hidarii G(K),\marubb\migii_{\chi}=-5,\hspace{1cm}
  \hidarii G(\alpha(K)),\marubb\migii_{\chi}=-1. 
\end{eqnarray}
Inserting these into Gauss diagram formula (\ref{eqn:main result2}), we obtain
\begin{eqnarray}\label{eqn:v2-6}
  v_2(K)&=&-\frac{1}{6}+\hidarii \bar{G}(K),\marubb\migii_{\chi} \nonumber\\
  &=&-\frac{1}{6}+\hidarii G(K),\marubb\migii_{\chi}
  -\hidarii G(\alpha(K)),\marubb\migii_{\chi}\nonumber\\
  &=&-\frac{1}{6}-4.
\end{eqnarray}

\noindent{\it Remark.}
Notice that we cannot replace $\alpha(K)$ by a trivial knot $U$ in the second equation of (\ref{eqn:example1}) since $\hidarii G(\alpha(K)),\marubb\migii_{\chi}$ is not a knot invariant. 
  
\subsubsection{$v_{3.1}(K)$} 
Using the Gauss diagram $G(K)$ in Figure \ref{fig:62}, we get
\begin{eqnarray}\label{eqn: 6-v31}
\begin{array}{cc}
  \hidarii G(K),\marude\migii_{\chi}=5, &
  \hidarii G(K),\marudd\migii_{\chi}=2, \\
  \hidarii G(K),\marudp\migii_{\chi}=-6. &
\end{array}
\end{eqnarray}
\begin{figure}[tbp]
	\setlength{\unitlength}{1cm}
\begin{picture}(8,4)(1.5,2.5)
\put(6,2.9){\includegraphics{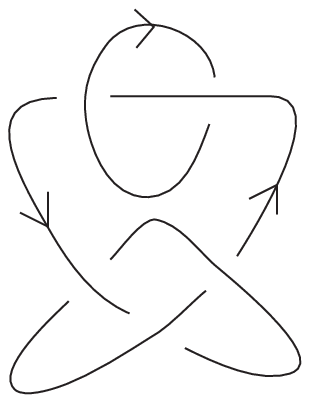} }
\put(8.3,6.1){\footnotesize a1}
\put(6.3,6.1){\footnotesize a2}
\put(6.2,3.9){\footnotesize a3} 
\put(7.4,3.1){\footnotesize a4} 
\put(8.4,4.1){\footnotesize a5}
\put(7.3,5.1){\footnotesize a6}

\end{picture}
	\caption{link diagram $Q(\{K:\mbox{a6}\},[\alpha])$}
	\label{fig:6split}
\end{figure}
Next we shall caluculate $I_{3.1}(K)$. 
Considering Figure \ref{fig:6split}, we have
\begin{eqnarray*}
  &&\hspace{-0.5cm}R(\{K:\mbox{a}_6\},[\gamma]-[\alpha])\\
  &&
  =~~\Biggl\{\piche\Biggr\}~~-~~\Biggl\{\pichd\Biggr\}~~
  -~~\Biggl\{\pichf\Biggr\}.
  \\
\end{eqnarray*}
We can caluculate the other $R(\{K:\mbox{a}_i\},[\gamma]-[\alpha])$ ($i=1,\cdots,5$) in the same way.
Then we have
\begin{eqnarray*}
  &&\hidarii R(\{K:\mbox{a}_i\},[\gamma]-[\alpha]),\marubb\migii_{\chi}
  =-1~~(i=1,\cdots,5),\\
  &&\hidarii R(\{K:\mbox{a}_6\},[\gamma]-[\alpha]),\marubb\migii_{\chi}
  =0.
\end{eqnarray*}
Inserting these into (\ref{eqn:I31}), we obtain
\begin{eqnarray}\label{eqn:exampleI31}
  I_{3.1}(K)=1.
\end{eqnarray}
Inserting (\ref{eqn: 6-v31}) (\ref{eqn:exampleI31}) into Gauss diagram formula (\ref{eqn:main result31}) yields  
\begin{eqnarray}\label{eqn:v31-6}
  v_{3.1}(K)=8.
\end{eqnarray}

\subsubsection{$v_{4.1}(K)$ and $v_{4.2}(K)$}
Using the Gauss diagram $G(K)$, 
$G\bigl(\alpha(K)\bigr)$ in Figure \ref{fig:62}, we get
\begin{eqnarray*}
\begin{array}{lll}
  \hidarii \bar{G}(K),\maruga\migii_{\chi}=0, &
  \hidarii \bar{G}(K),\marugb\migii_{\chi}=0, & 
  \hidarii \bar{G}(K),\marugc\migii_{\chi}=-6, \\
  \hidarii \bar{G}(K),\marugd\migii_{\chi}=4, &
  \hidarii \bar{G}(K),\maruge\migii_{\chi}=2, & 
  \hidarii \bar{G}(K),\marugf\migii_{\chi}=-2, \\
  \hidarii \bar{G}(K),\marubb\migii_{\chi}=-4, &
  \hidarii \bar{G}(K),\marugg\migii_{\chi}=-20, & 
  \hidarii \bar{G}(K),\marugh\migii_{\chi}=12, \\
  \hidarii \bar{G}(K),\marugi\migii_{\chi}=0, & & \\
  & & \\
  I_{4.1.1}(K)=4, & I_{4.1.2}(K)=-2, & I_{4.2}(K)=-2.
\end{array}
\end{eqnarray*}
Inserting these equation into Gauss diagram formula (\ref{eqn:main result41}) and (\ref{eqn:main result42}) yields  
\begin{eqnarray}\label{eqn:v41-6}
  v_{4.1}(K)=\frac{1}{360}+\frac{34}{3}, ~~~~~~
  v_{4.2}(K)=-\frac{1}{360}+\frac{38}{3}.
\end{eqnarray}
\subsubsection{}
Inserting (\ref{eqn:v2-6}),(\ref{eqn:v31-6}),(\ref{eqn:v41-6}) into the right side of (\ref{eqn:Homflyformula}), we get
\begin{eqnarray}\label{eq:example}
  &&\hspace{-1cm}
  \Bigl\{\mbox{the right side of (\ref{eqn:Homflyformula})}\Bigr\}
  \nonumber\\
  &&
  =1+(N^2-1)x^2+N(N^2-1)x^3
  +\frac{-13+6N^2+7N^4}{12}x^4.
\end{eqnarray}
The Homfly polynomial of $K$ is caluculated by the skein relation
(\ref{eqn:skein relation}):  
\begin{eqnarray*}
  P_{K}(t,z)=t^4z^2+t^4-t^2z^4-3t^2z^2-2t^2+z^2+2.
\end{eqnarray*}
We can easiliy check that $\bigl[P_{K}(e^{\frac{Nx}{2}},e^{\frac{x}{2}}-e^{-\frac{x}{2}})]^{(4)}$ coincides with (\ref{eq:example}).

\section{Proof of Theorem \ref{th:mainresult}}\label{sec:ProofTheorem2}
In this section we shall derive the Gauss diagram formula from Kontsevich integral (Proof of Theorem \ref{th:mainresult}).

\subsection{Sketch of the Proof of Theorem \ref{th:mainresult}}
We begin by briefly sketching the proof of Theorem \ref{th:mainresult}.

The integrand of $\langle\!\langle \LL,D \rangle\!\rangle$ in $v_{1},v_{2},v_{3.1},v_{3.2},v_{4.1},v_{4.2},v_{4.3},v_{4.4}$ 
(see Definition \ref{df:Kontsevich} and (\ref{eq:v12}) $\sim$ (\ref{eq:v44}))
has the following form:
\begin{eqnarray}\label{eqn:logtheta}
  d\log(z_{i_kj_k}(t_k))=di\theta_{i_kj_k}(t_k)+d\log r_{i_kj_k}(t_k),
\end{eqnarray}
where $\theta_{i_kj_k}(t_k)$ and $r_{i_kj_k}(t_k)$ are defined 
by the polar form 
\begin{eqnarray*}
z_{i_kj_k}(t_k)=r_{i_kj_k}(t_k)\exp(i\theta_{i_kj_k}(t_k)).
\end{eqnarray*}
We expand the integrand of $\langle\!\langle \LL,D \rangle\!\rangle$ according to (\ref{eqn:logtheta}). For example, if the degree of $D$ is two, the integrand of $\langle\!\langle \LL,D \rangle\!\rangle$ is expanded as follows:
\begin{eqnarray}\label{eqn:expand}
  \prod_{k=1}^2\Bigl\{\epsilon~d\log(z_{i_kj_k}(t_k))\Bigr\}
  &=&\Bigl\{\epsilon~di\theta_{i_1j_1}(t_1)\Bigr\}
  \Bigl\{\epsilon~di\theta_{i_2j_2}(t_2)\Bigr\}\nonumber\\
  &&+\Bigl\{\epsilon~di\theta_{i_1j_1}(t_1)\Bigr\}
  \Bigl\{\epsilon~d\log r_{i_2j_2}(t_2)\Bigr\}\nonumber\\
  &&+\Bigl\{\epsilon~d\log r_{i_1j_1}(t_1)\Bigr\}
  \Bigl\{\epsilon~di\theta_{i_2j_2}(t_2)\Bigr\}\nonumber\\
  &&+\Bigl\{\epsilon~d\log r_{i_1j_1}(t_1)\Bigr\}
  \Bigl\{\epsilon~d\log r_{i_2j_2}(t_2)\Bigr\}.\nonumber\\
\end{eqnarray}

Without loss of generality, we may replace the link $\LL$ with the link $A^b(L)$ in a nice position to caluculate (see Definition \ref{def:flatlimitlink}). 
Then key observation is as follows.
\begin{itemize}
  \item  The integrals which have odd number of $d\log r_{i_kj_k}(t_k)$'s are pure imaginary and do not contribute to the caluculation, since the Kontsevich integral is real valued
  (See Lemma \ref{lem:decomposition}).
  For example, the second and third terms in the right side of 
  (\ref{eqn:expand}) do not contribute to the caluculation. 
  \item The part of $di\theta_{i_kj_k}(t_k)$ integral is localized around the 
  cylinders (crossings) of the link $A^b(L)$, since $\theta_{i_kj_k}(t_k)$ does not vary on the plane $\mathbb{R}\times\{0\}\times\mathbb{R}=\{(x_1,0,x_3)\in\mathbb{R}^3\}$. Thus it is easy to evaluate
  (See Lemma \ref{lem:crossingnumber} and Lemma \ref{lem:split}).
  For example, the first term in the right hand side of (\ref{eqn:expand}) 
  is easily caluculated.
  \item The part of $d\log r_{i_kj_k}(t_k)$ integral is difficult to evaluate. But we can avoid this $d\log r_{i_kj_k}(t_k)$ integral as follows. First, since $r_{i_kj_k}(t_k)$ takes the same value for both signature $\pm$ of the cylinder,  the part of $d\log r_{i_kj_k}(t_k)$ integral does not depend on the 
  signatures of the link $A^b(L)$. In other words, it essentially depends only on  its projection to 
  $\mathbb{R}\times\{0\}\times\mathbb{R}=\{(x_1,0,x_3)\in\mathbb{R}^3\}$ 
  and takes the same value for $L$ and $\alpha(L)$.
  For example, the fourth term in the right hand side of (\ref{eqn:expand}) 
  does not depend on the signatures of the link $A^b(L)$. 
  Second, the modified Kontsevich integral is a link invariant. These two point leads to the final Gauss diagram formula.
\end{itemize} 
\subsection{Preparation for the proof of Theorem \ref{th:mainresult}}
In this section, we shall fix notations for the proof of 
Theorem \ref{th:mainresult}. 
\subsubsection{}
\begin{definition}\hspace{-5pt}{\bf.}\hspace{6pt}(Dotted Diagram)
An IL diagram $\{D,\kappa\}$ is called a {\it Dotted diagram} if $\kappa(c)=0,1~(c\in C(D))$. A chord $c$ is called a {\it normal chord} if $\kappa(c)=1$ and a {\it dotted chord} if $\kappa(c)=0$. In figures, we draw a normal chord $(\kappa(c)=1)$ by a thin line and a dotted chord $(\kappa(c)=0)$ by a dotted line as follows: 
\begin{eqnarray*}
  \picdc~~\hbox{normal chord}~(\kappa(c)=1),
  \hspace{1.5cm}\picee~~\hbox{dotted chord}~(\kappa(c)=0).
\end{eqnarray*}
We give two exmples of dotted diagrams,
\begin{eqnarray*}
  \marueh,\hspace{1.5cm}\marudj.
\end{eqnarray*}
\end{definition}

\begin{definition}\hspace{-5pt}{\bf.}\hspace{6pt}For a complex number $z$, define $\theta(z)$ and $r(z)$ by the polar form $z=r\exp(i\theta)$.
Let $\LL$ be a link and $\hat{D}$ a dotted diagram of degree $m$ which has $l$-normal chords and $(m-l)$-dotted chords, where $l$ is fixed.
We consider $m$-planes $t=t_k, (k=1,\cdots,l,\cdots,m)$ where
$t_{\mathrm{min}}<t_1<\cdots<t_l<t_{\mathrm{max}}$ and $t_{\mathrm{min}}<t_{l+1}<\cdots<t_m<t_{\mathrm{max}}$. For $1\le k\le m$, set $(\pi^{-1})(t_k)=\{s_k^1,\cdots,s_k^{n(t_k)}\}$, where $n(t_k)$ denotes the number of points on the section $t=t_k$ of the link $\LL$. For $1\le k\le l$, set $\theta_{ij}(t_k)=\theta\circ z\bigl\{\vec{x}(s^i_k)-\vec{x}(s^j_k)\bigr\}$. For $(l+1)\le k\le m$, set $r_{ij}(t_k)=r\circ z\bigl\{\vec{x}(s^i_k)-\vec{x}(s^j_k)\bigr\}$.
Define the collection of all pairings by $P=\{(i_1,j_1),\cdots,(i_m,j_m):1\le i_k\le j_k\le n(t_k)~~(k=1,\cdots,m)\}$. For $p\in P$, we shall define a dotted diagram $D_p$ of degree $m$. For all $k$, we join $s_k^{i_k}$ and $s_k^{j_k}$ by normal chords if $1\le k\le l$, and join $s_k^{i_k}$ and $s_k^{j_k}$ by dotted chords if $(l+1)\le k\le m$ on $X$. Define $D_p$ to be the resulting dotted diagram of degree $m$ which has $l$-normal chords and $(m-l)$-dotted chords.
Define $\hidarii \LL,\hat{D}\migii$ by
\begin{eqnarray}\label{eqn:dotted integral}
  &&\hspace{-1.2cm}
  \hidarii \LL,\hat{D}\migii=\frac{1}{(i\pi)^m}
  \int_{{\scriptstyle t_{\mathrm{max}}>t_1>\cdots>t_l>t_{\mathrm{min}}}
  \atop{\scriptstyle t_{\mathrm{max}}>t_{l+1}>\cdots>t_{m}>t_{\mathrm{min}}}}
  \sum_{p\in P}
  \prod_{k=1}^{l}\{\epsilon~di\theta_{i_kj_k}(t_k)\}
  \nonumber\\
  &&\hspace{5cm}\times
  \prod_{k=l+1}^{m}\{\epsilon~d\log r_{i_kj_k}(t_k)\}
  \Theta(D_P,\hat{D}),
\end{eqnarray}
where the sum is taken over all the pairings  $p\in P$.
$\Theta(D_p,\hat{D})$ is defined by
\begin{eqnarray*}
  \Theta(D_p,\hat{D})=
  \left\{\begin{array}{cl}
  1 & \mbox{if}~~D_p=\hat{D} \\
  0 & \mbox{if}~~D_p\ne \hat{D} 
  \end{array}\right..
\end{eqnarray*}

Let $\hat{D}_i$ be a dotted chord diagram.
More generally, for a formal linear combination of dotted diagrams 
$\displaystyle \sum_i b_i~\hat{D}_i$ 
($b_i\in \mathbf{C}$), set
\begin{eqnarray*}
  \hidarii \LL, \sum_i b_i~\hat{D}_i\migii
  =\sum_i b_i~\hidarii \LL, \hat{D}_i\migii. 
\end{eqnarray*}

\noindent{\it Remark}. Roughly speaking, a normal chord represents "$di\theta$" integral, and a dotted chord represents "$d\log r$" integral.
\end{definition}

\begin{lemma}\hspace{-5pt}{\bf.}\hspace{6pt}\label{lem:decomposition}
\begin{eqnarray}
  &&\bullet~~\mbox{Re}~\hidari \{\K_1,\K_2\},\marucb\migi
  =\hidarii \{\K_1,\K_2\},\marucb\migii,
  \label{eqn:decomposition0}\\
  &&\bullet~~
  \mbox{Re}~\hidari \K,\marubb\migi=\hidarii \K,\marubb+\marubd\migii,
  \label{eqn:decomposition1}\\
  &&\bullet~~
  \mbox{Re}~\hidari \K,\marudd\migi
  =\hidarii \K,\marudd+\marudi+\marudj\migii,
  \label{eqn:decomposition2}\\
  &&\bullet~~
  \mbox{Re}~\hidari \K,\marude\migi=\hidarii \K,\marude+\marudk\migii,
  \label{eqn:decomposition3}\\
  &&\bullet~~
  \mbox{Re}~\hidari \{\K_1,\K_2\},\marueb\migi 
  \nonumber\\
  &&\hspace{2cm}
  =\hidarii \{\K_1,\K_2\},\marueb+\marueg+\marueh\migii,
  \label{eqn:decomposition4}\\
  &&\bullet~~
  \mbox{Re}~\hidari \{\K_1,\K_2\},\marued\migi
  \nonumber\\
  &&\hspace{2cm}
  =\hidarii \{\K_1,\K_2\},\marued+\maruei\migii,
  \label{eqn:decomposition5}
\end{eqnarray} 
where Re $\langle\!\langle K,D \rangle\!\rangle$ denotes the real part of complex number $\langle\!\langle K,D \rangle\!\rangle$.
\end{lemma} 
{\bf Proof.}~~
We expand the integrand of $\langle\!\langle K,D \rangle\!\rangle$ according to (\ref{eqn:logtheta}). Only the integrals which have even number of $d\log r_{i_kj_k}(t_k)$'s contribute, since the integrals which have odd number of $d\log r_{i_kj_k}(t_k)$'s are pure imaginary. This proves the above lemma. $\square$

\subsubsection{}
\begin{definition}\hspace{-5pt}{\bf.}\hspace{6pt}(AF Link)\label{def:flatlimitlink}
Let $L$ be a link diagram in $\mathbb{R}\times\{0\}\times\mathbb{R}=\{(x_1,0,x_3)\in\mathbb{R}^3\}$. Without loss of generality, we may assume that the two curves 
around each crossing are given by
\begin{eqnarray}\label{eqn:rectangle}
  \left\{
  \begin{array}{c}
  x_3=\frac{\pi}{2}(-x_1+b) \\
  x_2=0 \\
  \end{array}
  \right.~(-b \le x_1 \le b),~~
  \left\{
  \begin{array}{c}
  x_3=\frac{\pi}{2}(x_1+b) \\
  x_2=0 \\
  \end{array}
  \right.~(-b \le x_1 \le b)
\end{eqnarray} 
with some parallel transformation (see the left side of Figure \ref{fig:cylinder}). 
Here $b$ is sufficiently small. In other words, this assumption is that two curves around the crossing are on the diagonal lines of some sufficiantly small rectangle parallel to $t$-axis.

 For this link diagram $L$, we shall define a link $A^b(L)$ called the {\it Almost-Flat Link} (AF Link) of $L$ as follows. For each crossing of $L$, we replace the two curves (\ref{eqn:rectangle}) by
\begin{eqnarray}\label{eqn:cylinder1}
  \left\{
  \begin{array}{c}
  x_1=b\cos (x_3/b) \\
  x_2=b\sin (x_3/b) \\
  \end{array}
  \right.(0\le x_3\le b\pi),~~
  \left\{
  \begin{array}{c}
  x_1=-b\cos (x_3/b) \\
  x_2=-b\sin (x_3/b) \\
  \end{array}
  \right.(0\le x_3\le b\pi),
\end{eqnarray} 
or
\begin{eqnarray}\label{eqn:cylinder2}
  \left\{
  \begin{array}{c}
  x_1=b\cos (x_3/b) \\
  x_2=-b\sin (x_3/b) \\
  \end{array}
  \right.(0\le x_3\le b\pi),~
  \left\{
  \begin{array}{c}
  x_1=-b\cos (x_3/b) \\
  x_2=b\sin (x_3/b) \\
  \end{array} 
  \right.(0\le x_3\le b\pi)
\end{eqnarray} 
according to the signature of the crossing (see Figure \ref{fig:cylinder}). In other words, we replace the two curves on the rectangle by two curves winding around the cylinder so that projecting two curves winding around the cylinder to $\mathbb{R}\times\{0\}\times\mathbb{R}$ yields the signature of the crossing . For sufficiently small $b$, we define the AF link $A^b(L)$ to be the resulting link.

Since the cylinders of the AF link $A^b(L)$ is one-to-one correspondent to the crossings of the link diagram $L$, we define the signature of each cylinder to be the signature of the corresponding crossing.

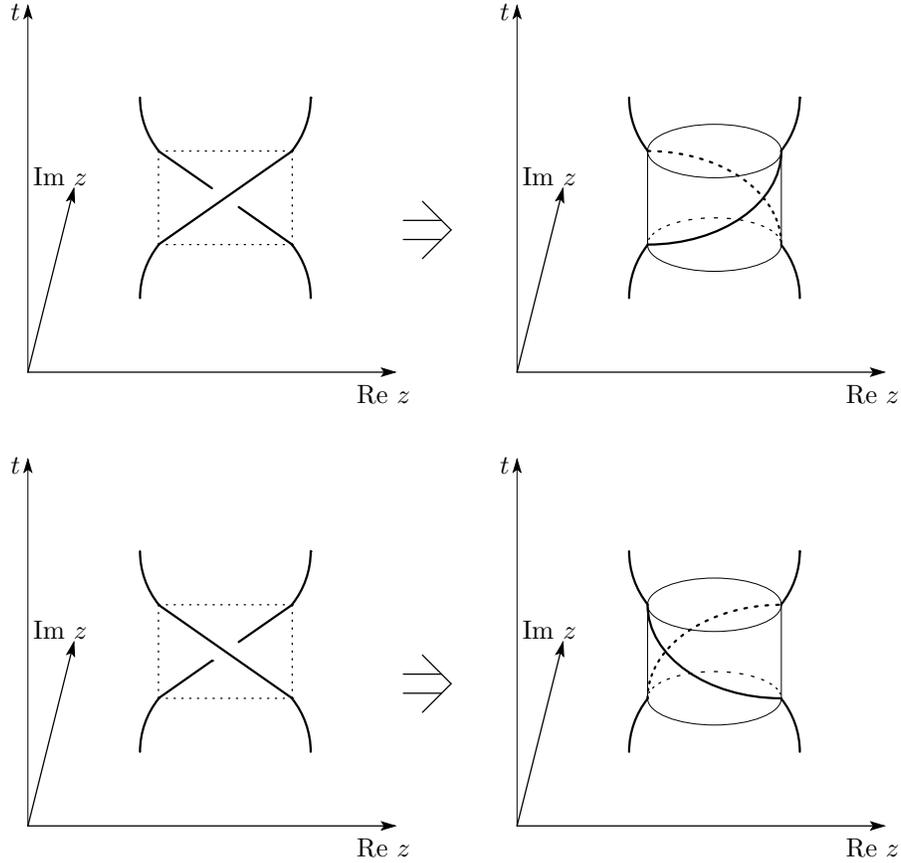
\begin{figure}[htbp]
\setlength{\unitlength}{1cm}
\begin{picture}(15,6)(1.5,0)
\put(2,0){\picfj}
\put(1.9,4.7){$t$}\put(6.5,-0.4){Re $z$}\put(2.2,2.5){Im $z$}
\put(3.5,1){\picfi}
\put(7,0){\picfc}
\put(8.5,0){\picfj}
\put(8.4,4.7){$t$}\put(13,-0.4){Re $z$}\put(8.7,2.5){Im $z$}
\put(10,1){\picfh}
\end{picture}
\begin{picture}(15,6)(1.5,0)
\put(2,0){\picfj}
\put(1.9,4.7){$t$}\put(6.5,-0.4){Re $z$}\put(2.2,2.5){Im $z$}
\put(3.5,1){\picfl}
\put(7,0){\picfc}
\put(8.5,0){\picfj}
\put(8.4,4.7){$t$}\put(13,-0.4){Re $z$}\put(8.7,2.5){Im $z$}
\put(10,1){\picfk}
\end{picture}
\newline
	\caption{$L$ and $A^b(L)$}
	\label{fig:cylinder}
\end{figure} 
\end{definition}

\begin{definition}\hspace{-5pt}{\bf.}\hspace{6pt}(Direction of the Crossing)
Let $L$ be a link diagram in $\mathbb{R}\times\{0\}\times\mathbb{R}=\{(x_1,0,x_3)\in\mathbb{R}^3\}$ as explained in Definition \ref{def:flatlimitlink}. 
We assign a signature $\pm$ to each crossing of $L$ as usual (Definiton \ref{def:linkdiagram}). Moreover, we assign "s" ("s" denotes "same") to each crossing of $L$ if the directions of the two arrows (orientations) are the same with respect to $t$-axis, and assign "d" if they are different ("d" denotes "different"):
\begin{eqnarray*}
  &&\picen~~~~\picef~~\piceg~~~~\piceh~~\picei~~
  ~~~~\picej~~\picek~~~~\picel~~\picem,
  \\ 
\end{eqnarray*}
where "s", "d" are formal letters.
We call "s","d" the {\it Direction of the Crossing}.
\end{definition}
{\it Remark.} Of courese, the concept of direction of the crossing depend on how to choose $t$-axis. So, it is not the proper quantity of link diagrams. 
Although the concept of directions appears in the computaion, it disappears in the final result (see Theorem \ref{th:mainresult}). $\square$
\newline

We shall extend the definition of IL diagram, Gauss diagram, ML diagram and the pairing $\langle \hat{G},\hat{D}\rangle_{\chi}$ to include the concept of direciton (see Definition \ref{def;ILDiagram}, Definition \ref{Gaussdiagram}, Definition \ref{MLdiagram}, Definition \ref{def:crossing number}).
\begin{definition}\hspace{-5pt}{\bf.}\hspace{6pt}({\it IDL Diagram})
Let $D$ be a chord diagram, and let $C(D)$ be the set of all chords of $D$. By a direction-labelling of $D$, we mean a map $f:C(D)\to \{s,d,n\}$, where $s,d,n$ are formal letters. {\it An Integer-Direction-Labeled Chord Diagram} (IDL Diagram) is a triple $\{D,\kappa,f\}$ of a chord diagram $D$ together with an integer-labelling $\kappa$ and a direction-labelling $f$. Two IDL 
diagrams $\{D,\kappa,f\}$, $\{D^\prime,\kappa^\prime,f^\prime\}$ are regarded as equal if $D$, $D^\prime$ are equal as chord diagrams and the homeomorphism $F:D\to D^\prime$ preserves integer-labelling $\kappa^\prime(F(c))=\kappa(c)~(c\in C(D))$ and direction-labelling $f^\prime(F(c))=f(c)~(c\in C(D))$.
\end{definition}

\begin{definition}\hspace{-5pt}{\bf.}\hspace{6pt}(Extended Gauss Diagram)
An IDL diagram $\{G,\epsilon,f\}$ is called a {\it Extended Gauss Diagram} if $\epsilon(c)=\pm 1$ and $f(c)=s,d~(c\in C(D))$.

Let $L$ be a link diagram as explained in Definition \ref{def:flatlimitlink}.
and let $\{L:a_1,\cdots,a_m\}$ be a link diagram $L$ where we select some distinct crossings $a_1,\cdots,a_m$ out of all crossings of $L$. Define a extended Gauss diagram $P^e(\{L:a_1,\cdots,a_m\})$ as follows. For each $a_i$, set $\vec{y}^{-1}(a_i)=\{s(a_i),s^\prime(a_i)\}$ as the inverse image of $a_i$. For each crossing $a_i$, we join $s(a_i),s^\prime(a_i)$ by a chord on $X$ and label this chord by the signature and the direction of $a_i$ ($i=1,\cdots,m$). We define an extended Gauss diagram $P^e(\{L:a_1,\cdots,a_m\})$ to be the result.

Specially, If $\{a_1,\cdots,a_m\}$ are all the crossings of $L$ (this means we select all the crossings of $L$) , we write
$G^e(L)=P^e(\{L:a_1,\cdots,a_m\})$ and call it the extended Gauss diagram of $L$. 
For example, see Fig \ref{fig:extendedGaussdiagram}.
\begin{figure}[htbp]
\setlength{\unitlength}{1cm}
  \begin{picture}(15,4)(0,0.5)
  \put(2.85,0.3){\includegraphics{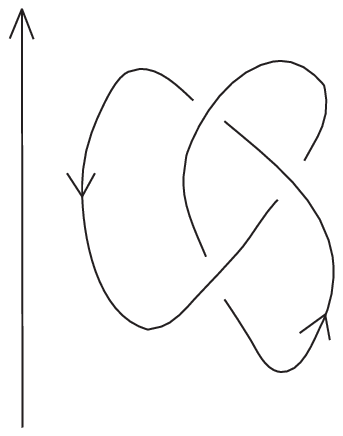}}
  \put(5.9,2.75){a} \put(4.75,3.7){b} \put(4.8,1.4){c}
  \put(4.7,0.5){$K$} \put(2.7,4.3){$t$} 
  \put(8,0.5){\piceo}
  \end{picture}
    \caption{}
	\label{fig:extendedGaussdiagram}
\end{figure}
\end{definition}
\newpage

\begin{definition}\hspace{-5pt}{\bf.}\hspace{5pt}(Extended ML Diagram)
An IDL diagram $\{D,m,\Delta\}$ is called a 
{\it Extended Multiplicity-Labeled Diagram} (Extended ML Diagram) if $m(c)=1,2~$ and $\Delta(c)=s,d,n~(c\in C(D))$. In figures, we draw a chord as follows.
\begin{eqnarray*}
  &&\picdc~~m(c)=1,\Delta(c)=n,
  \hspace{1.5cm}\picdd~~m(c)=2,\Delta(c)=n\nonumber\\
  &&\picea~~m(c)=1,\Delta(c)=s,
  \hspace{1.5cm}\picec~~m(c)=2,\Delta(c)=s
  \nonumber\\
  &&\piceb~~m(c)=1,\Delta(c)=d,
  \hspace{1.5cm}\piced~~m(c)=2,\Delta(c)=d
\end{eqnarray*}
We give two examples of ML diagrams.
\begin{eqnarray*}
  \maruej,\hspace{1.5cm}\marudl
\end{eqnarray*}
\end{definition}

\begin{definition}\hspace{-5pt}{\bf.}\hspace{5pt}Let $\hat{G}^e=\{G,\epsilon,f\}$ be an extended Gauss diagram and $\hat{D}^e=\{D,m,\Delta\}$ be a extended ML diagram. Let $\psi:D\to G$ be an embedding of $D$ into $G$ which maps the circles of $D$ to those of $G$ preserving the orientations and maps each chord of $D$ to a chord of $G$. Let $C(G)$ be the set of all the chords of $G$. For $\psi$, define a map $\kappa^e_\psi:C(G)\to\{0,1,2\}\times\{s,d,n\}$ by
\begin{eqnarray*}
  \kappa^e_\psi(c)=\left\{
  \begin{array}{cl}
  (m(\psi^{-1}(c)),\Delta(\psi^{-1}(c))) & \mbox{if}~~c\in \psi(D) \\
  (~~0~~,~~n~~) & \mbox{if}~~c\noin \psi(D) 
  \end{array}\right.
\end{eqnarray*}
Two embedding $\psi,\varphi$ are said to be equal if $\kappa^e_\psi=\kappa^e_{\varphi}$. The equivalence class of an embedding $\psi$ is denoted by $[\psi]$.
Define $\delta:\{s,d,n\}\times\{s,d\}\to \pm1$ by
\begin{eqnarray*}
  \delta(s:f)=
  \left\{
  \begin{array}{cc}
  1 & (f=s) \\
  0 & (f=d) 
  \end{array}\right.\hspace{-0.2cm},~~
  \delta(d:f)=
  \left\{
  \begin{array}{cc}
  0 & (f=s) \\
  1 & (f=d) 
  \end{array}\right.\hspace{-0.2cm},~~
  \delta(n:f)=
  \left\{
  \begin{array}{cc}
  1 & (f=s) \\
  1 & (f=d) 
  \end{array}\right..
\end{eqnarray*}
Let $C(D)$ be the set of all chords of $D$. Define $\mathcal{E}^e([\psi])$ by
\begin{eqnarray*}
  \mathcal{E}^e([\psi])=\prod_{c\in C(D)}\bigl\{\epsilon(\psi(c))\bigr\}^{m(c)}
  \delta(\Delta(c),f(\psi(c))),
\end{eqnarray*}
where the product is taken over all chords of $D$. Notice this definition is well defined.

Define $\langle \hat{G}^e,\hat{D}^e\rangle_{\chi^e}$ by
\begin{eqnarray*}
  \langle \hat{G}^e,\hat{D}^e\rangle_{\chi^e}=
  \sum_{[\psi]}\mathcal{E}^e([\psi]),
\end{eqnarray*}
where the sum is taken over all the distinct equivalence classes $[\psi]$.

Let $\hat{G}_i^e$ be a extended Gauss diagram and $\hat{D}_i^e$ a extended ML diagram.
More generally, for formal linear combinations $\sum_i b_i~\hat{G}_i^e$ and
$\sum_j c_j~\hat{D}_j^e$
($b_i, c_j\in \mathbf{C}$), set
\begin{eqnarray}
  \Bigl\langle \sum_i b_i~\hat{G}_i^e,\sum_j c_j~\hat{D}_j^e
  \Bigr\rangle_{\chi^e}
  =\sum_i\sum_j b_i~c_j~\langle \hat{G}_i^e,\hat{D}_j^e\rangle_{\chi^e}.
\end{eqnarray}
\newline
{\it Remark}. This definition $\langle \hat{G}^e,\hat{D}^e\rangle_{\chi^e}$ is a natural extension of $\langle \hat{G},\hat{D}\rangle_{\chi}$ in Definition \ref{def:crossing number}. More precisely, if $\Delta(c)=n$ for all $c\in C(D)$, then
\begin{eqnarray*}
   \langle \hat{G}^e,\hat{D}^e\rangle_{\chi^e}=
   \langle \hat{G},\hat{D}\rangle_{\chi},
\end{eqnarray*}
where $\hat{G}=\{G,\epsilon\}$, $\hat{D}=\{D,m\}$.
\end{definition}
\begin{lemma}\hspace{-5pt}{\bf.}\hspace{5pt}\label{lem:crossingnumber}
Let $A^b(K)$ and $A^b(\{K_1,K_2\})$ be the AF links which correspond to link diagrams $K$ and $\{K_1,K_2\}$ respectively. Then, for sufficiently small $b$,
\begin{eqnarray}
  &&\hspace{-0.5cm}\bullet~~ 
  \hidarii A^b(\{K_1,K_2\}),\marucb\migii
  =\hidarii G(\{K_1,K_2\}),\marucb\migii_{\chi},
  \label{eqn:crossingnumber0}\\
  &&\hspace{-0.5cm}\bullet~~ 
  \hidarii A^b(K),\marubb\migii=\hidarii G^e(K),\marubb
  +\frac{1}{2}\marube\migii_{\chi^e}+O(b),\label{eqn:crossingnumber1}\\
  &&\hspace{-0.5cm}\bullet~~ 
  \hidarii A^b(K),\marudd\migii=\hidarii G^e(K),\marudd
  +\frac{1}{2}\marudl~\migii_{\chi^e}+O(b),\label{eqn:crossingnumber311}\\
  &&\hspace{-0.5cm}\bullet~~ 
  \hidarii A^b(K),\marude\migii
  =\hidarii G^e(K),\marude+\frac{1}{2}\marudm~
  +\frac{1}{3!}\marudn\migii_{\chi^e}+O(b),\label{eqn:crossingnumber312}\\
  &&\hspace{-0.5cm}\bullet~~ 
  \hidarii A^b(\{K_1,K_2\}),\marueb\migii
  =\hidarii G(\{K_1,K_2\}),\marueb \migii_{\chi}+O(b),
  \label{eqn:crossingnumber321}\\
  &&\hspace{-0.5cm}\bullet~~ 
  \hidarii A^b(\{K_1,K_2\}),\marued\migii\nonumber\\
  &&\hspace{0cm}=\hidarii G^e(\{K_1,K_2\}),\marued
  +\frac{1}{2}\maruej+\frac{1}{3!}\maruek\migii_{\chi^e}+O(b).
  \nonumber\\
  \label{eqn:crossingnumber322}
\end{eqnarray}
{\bf Proof.} We shall prove (\ref{eqn:crossingnumber1}).  We can prove all the other in the same way. In (\ref{eqn:dotted integral}), $di\theta_{i_kj_k}(t_k)$ integral is localized around the cylinders of the AF knot $A^b(K)$, since $\theta_{i_kj_k}(t_k)$ does not vary on the plane $\mathbb{R}\times\{0\}\times\mathbb{R}=\{(x_1,0,x_3)\in\mathbb{R}^3\}$. Let $a$ be a cylinder of $A^b(K)$ and $I_a$ the small interval on $t$-axis which contains the cylinder $a$. We assume $I_a$ contains only one cylinder $a$ and the other curves in $t\in I_a$ is straight lines parallel to $t$-axis as follows: 

\begin{eqnarray}
  \picfq. \nonumber
\end{eqnarray}
In the above figure, we draw the curves of $A^b(K)$ by thick lines.  Then, we have
\begin{eqnarray}\label{eqn:proofcrossingnumber}
  &&\hspace{-1cm}\hidarii A^b(K),\marubb\migii\nonumber\\
  &&=\sum_{(a,b)}\frac{1}{(i\pi)^2}
  \int_{t_1\in I_a, t_2\in I_b}
  \sum_{p\in P}
  \prod_{k=1}^{2}\{\epsilon~di\theta_{i_kj_k}(t_k)\}
  \Theta(D_P,\marubb)\nonumber\\
  &&~~+\sum_{a}\frac{1}{(i\pi)^2}
  \int_{t_1>t_2\in I_a}
  \sum_{p\in P}
  \prod_{k=1}^{2}\{\epsilon~di\theta_{i_kj_k}(t_k)\}
  \Theta(D_P,\marubb),
\end{eqnarray}
where the first sum is taken over all the unordered pair of cylinders $(a,b)$ and the second sum is taken over all cylinders $a$.
By {\it a pairing on cylinder}, we mean a pairing $p\in P$ for which both $\vec{x}(s^{i_k}_k)$, $\vec{x}(s^{j_k}_k)$ are on the cylinder $(k=1,2)$.
Only the pairings on cylinder contribute to the caluculatin for the following reason. For example, we consider the second term in the right side of (\ref{eqn:proofcrossingnumber}) and the pairing 
$\{(i_1,j_1)(i_2,j_2)\}=\{(1,2)(1,3)\}$. 
Assume $\vec{x}(s^{1}_1)$, $\vec{x}(s^{2}_1)$, $\vec{x}(s^{1}_2)$  are on the cylinder and $\vec{x}(s^{3}_2)$ is on sraight line :
\begin{center}
\setlength{\unitlength}{1cm}
\begin{picture}(15,3)(-2,-0.3)
\put(1.9,2.8){$t$}
\put(2,0){\picgg}
\put(4.05,0){$1$}
\put(5.05,0){$2$}
\put(7.1,0){$3$}
\put(4.55,2){$a$}
\end{picture}
\end{center}
where we draw the curves of $A^b(K)$ by thick lines.
Then 
\begin{eqnarray*}
  &&\biggl\vert 
  \int_{t_1>t_2\in I_a}
  \{\epsilon~di\theta_{12}(t_1)\}\{\epsilon~di\theta_{13}(t_2)\}
  \Theta(D_P,\marubb)\biggr\vert\\
  &<&
  \biggl\{\int_{t_1\in I_a}\bigl\vert
  d\theta_{12}(t_1)\bigr\vert\biggr\}
  \biggl\{\int_{t_2\in I_a}\bigl\vert
  d\theta_{13}(t_2)\bigr\vert\biggr\}\\
  &<& \mbox{(constant)}\times b.
\end{eqnarray*}
So the integral correspond to the pairing which is {\it not} on cylinder is bounded by $b$.

Anyway since only the pairings on cylinder contribute to the caluculation, 
\begin{eqnarray*}
  \frac{1}{(i\pi)}\int_{t_k\in I_a}\{\epsilon~di\theta_{i_kj_k}(t_k)\}
\end{eqnarray*}
gives the signature of the cylinder $a$. Therefore we see that the first term in the right side of (\ref{eqn:proofcrossingnumber}) gives $\hidarii G(K),\marubb\migii_{\chi}$ and the second term gives $\frac{1}{2}\hidarii G^e(K),\marube\migii_{\chi^e}$, considering the restriction of $\Theta(D_P,\marubb)$. $\square$
\end{lemma}

\subsubsection{}
\begin{definition}\hspace{-5pt}{\bf.}\hspace{5pt}Let $\{A^b(L):a_1,\cdots,a_l\}$ be a AF link $A^b(L)$ where we select $l$-cylinders $a_1,\cdots,a_l$ out of all cylinders of $A^b(L)$. Let $\hat{D}$ be a dotted chord diagram of degree $(m+l)$ which has $l$-normal chords and $m$-dotted chords.
We consider $m$-planes $t=t_k,(k=1,\cdots,m)$ where
$t_{\mathrm{max}}>t_1>\cdots>t_m>t_{\mathrm{min}}$. For $1\le k\le m$, set $(\pi^{-1})(t_k)=\{s_k^1,\cdots,s_k^{n(t_k)}\}$, where $n(t_k)$ denotes the number of points on the section $t=t_k$ of $A^b(L)$. For $1\le k\le m$, set $r_{ij}(t_k)=r\circ z\bigl\{\vec{x}(s^i_k)-\vec{x}(s^j_k)\bigr\}$.
Define the collection of all pairings by $P=\{(i_1,j_1),\cdots,(i_m,j_m):1\le i_k\le j_k\le n(t_k)~~(k=1,\cdots,m)\}$.

For each pairing $p\in P$ and the specific cylinders $a_1,\cdots,a_l$, we shall define a dotted diagram $D_{p,a_1\cdots a_{l}}$ of degree $(m+l)$.  For each cylinder $a_k$ $ (1\le k\le l)$, we mark two distinct points $d_k,d_k^\prime$ with the same heights on each curves winding around $a_k$ (see Figure \ref{fig:cylinder2}). We set this height to be $x_3=\frac{b\pi}{2}$ in (\ref{eqn:cylinder1}) (\ref{eqn:cylinder2}). Next set $\hat{s}_k=\vec{x}^{-1}(d_k),\hat{s}_k^\prime=\vec{x}^{-1}(d_k^\prime)$ as the inverse images of $d_k,d_k^\prime$.
For $p\in P$ and $a_1,\cdots,a_l$, join $s_k^{i_k}$ and $s_k^{j_k}$ by dotted chords $(k=1,\cdots,m)$ and join $\hat{s}_k$ and $\hat{s}_k^\prime$ by a normal chord $(k=1,\cdots,l)$ on $X$. Let $D_{p,a_1\cdots a_{l}}$ be the resulting dotted diagram of degree $(m+l)$ which has $l$-normal chords and $m$-dotted chords.

For sufficiently small $b$, define $\hh \{A^b(L):a_1,\cdots,a_{m-l}\},\hat{D}\mm$ by
\begin{eqnarray*}
  &&\hspace{-1cm}\hh \{A^b(L):a_1,\cdots,a_{l}\},\hat{D}\mm
  \\
  &&=\frac{1}{(i\pi)^m}
  \int_{t_1>\cdot\cdot\cdot>t_{m}}\sum_{p\in P}
  \prod_{k=1}^{m}\{\epsilon~d\log r_{i_k j_k}(t_k)\}~
  \Theta(D_{p,a_1\cdots a_{l}},\hat{D}),
\end{eqnarray*}
where the sum is taken over all pairings $p\in P$.
\newline
\newline
{\it Remark}. Roughly speaking, a normal chord represents the cylinder of the AF link, and a dotted chord represents the "$d\log r$" integral.
\begin{figure}[htbp]
\setlength{\unitlength}{1cm}
\begin{picture}(15,5)(-1.8,0)
\put(2,0){\picfj}
\put(1.9,4.7){$t$}\put(6.5,-0.4){Re $z$}\put(2.2,2.5){Im $z$}
\put(3.5,1){\picgb}
\end{picture}	\caption{}
	\label{fig:cylinder2}
\end{figure} 
\end{definition}
\newpage

\begin{lemma}\hspace{-5pt}{\bf.}\hspace{5pt}\label{lem:split}
Let $A^b(K)$ and $A^b(\{K_1,K_2\})$ be the AF links which correspond to link diagrams $K$ and $\{K_1,K_2\}$ respectively. Then, for sufficiently small $b$,
\begin{eqnarray}
  &&\hspace{-0.5cm}\bullet~~ 
  \hidarii A^b(K),\marudi\migii=\sum_{a}
  \hidarii P\bigl(\{K:a\}\bigr),\marua\migii_{\chi} \hh\{A^b(K):a\},\marudi\mm
  \nonumber\\
  &&\hspace{8cm}
  +O(b),\label{eqn:spilit1}\\
  &&\hspace{-0.5cm}\bullet~~ 
  \hidarii A^b(K),\marudj\migii=\sum_{a}
  \hidarii P\bigl(\{K:a\}\bigr),\marua\migii_{\chi} \hh\{A^b(K):a\},\maruia\mm
  \nonumber\\
  &&\hspace{8cm}
  +O(b),\label{eqn:spilit2}\\
  &&\hspace{-0.5cm}\bullet~~ 
  \hidarii A^b(K),\marudk\migii=\sum_{a}
  \hidarii P\bigl(\{K:a\}\bigr),\marua\migii_{\chi} \hh\{A^b(K):a\},\marudk\mm
  \nonumber\\
  &&\hspace{8cm}
  +O(b),\label{eqn:spilit3}\\
  &&\hspace{-0.5cm}\bullet~~ 
  \hidarii A^b(\{K_1,K_2\}),\marueg\migii\nonumber\\
  &&\hspace{0cm}
  =\sum_{a}\hidarii P\bigl(\{K_1,K_2:a\}\bigr),\marua\maru\migii_{\chi} 
  \hh\bigl\{A^b(\{K_1,K_2\}):a\bigr\},\marueg\mm
  \nonumber\\
  &&\hspace{8cm}
  +O(b),\label{eqn:spilit4}\\
  &&\hspace{-0.5cm}\bullet~~ 
  \hidarii A^b(\{K_1,K_2\}),\marueh\migii\nonumber\\
  &&\hspace{0cm}
  =\sum_{a}\hidarii P\bigl(\{K_1,K_2:a\}\bigr),\marucb\migii_{\chi}
  \hh\bigl\{A^b(\{K_1,K_2\}):a\bigr\},\maruid\mm
  \nonumber\\
  &&\hspace{8cm}
  +O(b),\label{eqn:spilit5}\\[1cm]
  &&\hspace{-0.5cm}\bullet~~ 
  \hidarii A^b(\{K_1,K_2\}),\maruei\migii\nonumber\\
  &&\hspace{0cm}
  =\sum_{a}\hidarii P\bigl(\{K_1,K_2:a\}\bigr),\marucb\migii_{\chi}
  \hh\bigl\{A^b(\{K_1,K_2\}):a\bigr\},\maruie\mm
  \nonumber\\
  &&\hspace{8cm}
  +O(b),
  \label{eqn:spilit6}
\end{eqnarray} 
where the sum is taken over all the cylinders $a$ of the AF link. Notice the terms $P\bigl(\{K:a\}\bigr)$, etc make sence, since a cylinder of a AF link is identified with the crossing of the corresponding link diagram.  
\newline
\newline
{\bf Proof.} We shall prove (\ref{eqn:spilit1}). We can prove the other in the same way. In (\ref{eqn:dotted integral}), $di\theta_{i_kj_k}(t_k)$ integral is localized around the cylinders of the AF link $A^b(K)$. So we make the same assumption for the small interval $I_a$ as in the proof of Lemma \ref{lem:crossingnumber}. Considering this position, (\ref{eqn:spilit1}) becomes:
\begin{eqnarray*}
  &&\hspace{-1cm}
  \hidarii A^b(K),\marudi\migii=\sum_{a}\frac{1}{(i\pi)^3}\int_{t_1\in I_a}
  \int_{t_{\mathrm{max}}>t_2>t_3>t_{\mathrm{min}}}
  \sum_{p\in P}
  \{\epsilon~di\theta_{i_1j_1}(t_1)\}\\
  &&\hspace{5cm}\times
  \prod_{k=2}^{3}\{\epsilon~d\log r_{i_kj_k}(t_k)\}
  \Theta(D_P,\marudi).
\end{eqnarray*}
The first sum is taken over all cylinder $a$. Let $P_c$ be a set of all the pairings $p\in P$ where both $\vec{x}(s^{i_1}_1)$ and $\vec{x}(s^{j_1}_1)$ are on the cylinder.
Only the pairings $p\in P_c$ contribute to the caluculation for the same reason as in the proof of Lemma \ref{lem:crossingnumber}. So
\begin{eqnarray*}
  \frac{1}{(i\pi)}\int_{t_1\in I_a}\{\epsilon~di\theta_{i_1j_1}(t_1)\}
\end{eqnarray*}
gives the signature of the cylinder $a$, which is equal to 
$\hidarii P\bigl(\{K:a\}\bigr),\marua\migii_{\chi}$.
The remaining part gives $\hh\{A^b(K):a\},\marudi\mm$. $\square$
\end{lemma}

\subsubsection{}
\begin{definition}\hspace{-5pt}{\bf.}\hspace{5pt}Let $\{L:a_1,\cdots,a_m\}$ be a link diagram $L$ where we select some distinct crossings $a_1,\cdots,a_m$ out of all crossings of $L$. Define a chord diagram $P_0(\{L:a_1,\cdots,a_m\})$ as follows. For each $a_i$, set $\vec{y}^{-1}(a_i)=\{s(a_i),s^\prime(a_i)\}$ as the inverse image of $a_i$. For each crossing $a_i$, we join $s(a_i),s^\prime(a_i)$ by a chord on $X$. We define a chord diagram $P_0(\{L:a_1,\cdots,a_m\})$ to be the result. $\square$
\end{definition} 
  
\noindent{\it Remark.}  $P_0(\{L:a_1,\cdots,a_m\})$ is obtained from $P(\{L:a_1,\cdots,a_m\})$ by dropping the signature labelling.

\begin{definition}\hspace{-5pt}{\bf.}\hspace{5pt}\label{def:NewKnotDiagram}
We shall define knot diagram $K^{[1]}_\pm, K^{[2]}, K^{[3]}_\pm,\cdots,$ etc 
as follows (see also Figure \ref{fig:split}). 
\begin{eqnarray*}
  &&\bullet~~ \mbox{Set}~~ 
  Q(\{K:a\},[\alpha])=\{K^{[1]}_+,K^{[1]}_-\}.\\
  &&\bullet~~ \mbox{If } P_0(\{K_1,K_2:a\})=\marucb,~~~~\mbox{set}~~       
  Q(\{K_1,K_2:a\},[\alpha])=\{K^{[2]}\}.\\
  &&\bullet~~ \mbox{If } P_0(\{K:a_1,a_2\})=\marubb,~~~~\mbox{set}~~       
  Q(\{K:a_1,a_2\},[\alpha,\gamma])=\{K^{[3]}_+,K^{[3]}_-\},\\
  &&\hspace{5.85cm}       
  Q(\{K:a_1,a_2\},[\gamma,\alpha])=\{K^{[4]}_+,K^{[4]}_-\},\\
  &&\hspace{5.05cm}\mbox{and}~~       
  Q(\{K:a_1,a_2\},[\beta,\beta])=\{K^{[5]}_+,K^{[5]}_-\}.\\
  &&\bullet~~ \mbox{If } P_0(\{K:a_1,a_2\})=\maruba,\\[0cm]
  &&\hspace{3.7cm}~~~~\mbox{set}~~ 
  Q(\{K:a_1,a_2\},[\alpha,\gamma])=\{K^{[6]}_+,K^{[6]}_-\},\\
  &&\hspace{4.85cm}       
  Q(\{K:a_1,a_2\},[\gamma,\alpha])=\{K^{[7]}_+,K^{[7]}_-\},\\
  &&\hspace{4.05cm}\mbox{and}~~       
  Q(\{K:a_1,a_2\},[\alpha,\alpha])=\{K^{[8]}_+,K^{[8]}_0,K^{[8]}_-\}.\\
  &&\bullet~~ \mbox{If } P_0(\{K_1,K_2:a_1,a_2\})=\maruca,\\[0.1cm]
  &&\hspace{3.7cm}~~~~\mbox{set}~~     
  Q(\{K_1,K_2:a_1,a_2\},[\beta,\beta])=\{K^{[9]}_+,K^{[9]}_-\},\\
  &&\bullet~~ \mbox{If } P_0(\{K_1,K_2:a_1,a_2\})=\marucc,\\[0.1cm]
  &&\hspace{2.35cm}~~~~\mbox{set}~~        
  Q(\{K_1,K_2:a_1,a_2\},[\gamma,\alpha])=\{K^{[10]}\},\\
  &&\hspace{3.5cm}       
  Q(\{K_1,K_2:a_1,a_2\},[\alpha,\alpha])=\{K^{[11]}_+,K^{[11]}_-\},\\
  &&\hspace{3.5cm}       
  Q(\{K_1,K_2:a_1,a_2\},[\gamma,\gamma])=\{K^{[12]}_+,K^{[12]}_-\},\\
   &&\hspace{2.7cm}\mbox{and}~~       
  Q(\{K_1,K_2:a_1,a_2\},[\alpha,\gamma])
  =\{K^{[13]}_+,K^{[13]}_0,K^{[13]}_-\}.\\
  &&\hspace{0cm} \mbox{(For convenience, assume $Q(\{K_1,K_2:a_1,a_2\},[\gamma,\alpha]))$ is   1-component link.)}\\
  &&\bullet~~ \mbox{If } P_0(\{K_1,K_2,K_3:a_1,a_2\})=\maruhm,\\[0.1cm]
  &&\hspace{2.5cm}
  ~~~~\mbox{set}~~        
  Q(\{K_1,K_2,K_3:a_1,a_2\},[\alpha,\alpha])=\{K^{[14]}\},\\
  &&\hspace{3.65cm}    
  Q(\{K_1,K_2,K_3:a_1,a_2\},[\gamma,\alpha])=\{K^{[15]}_+,K^{[15]}_-\},\\
  &&\hspace{2.85cm}\mbox{and}~~       
  Q(\{K_1,K_2,K_3:a_1,a_2\},[\alpha,\gamma])=\{K^{[16]}_+,K^{[16]}_-\}.\\
  &&\bullet~~ \mbox{If } P_0(\{K_1,K_2:a\})=\{\marua\maru\},\\[0.1cm]
  &&\hspace{2.5cm}
  ~~~~\mbox{set}~~        
  Q(\{K_1,K_2:a\},[\alpha])=\{K^{[17]}_+,K^{[17]}_-,K^{[17]}_0\}.
\end{eqnarray*}
\begin{figure}[htbp]
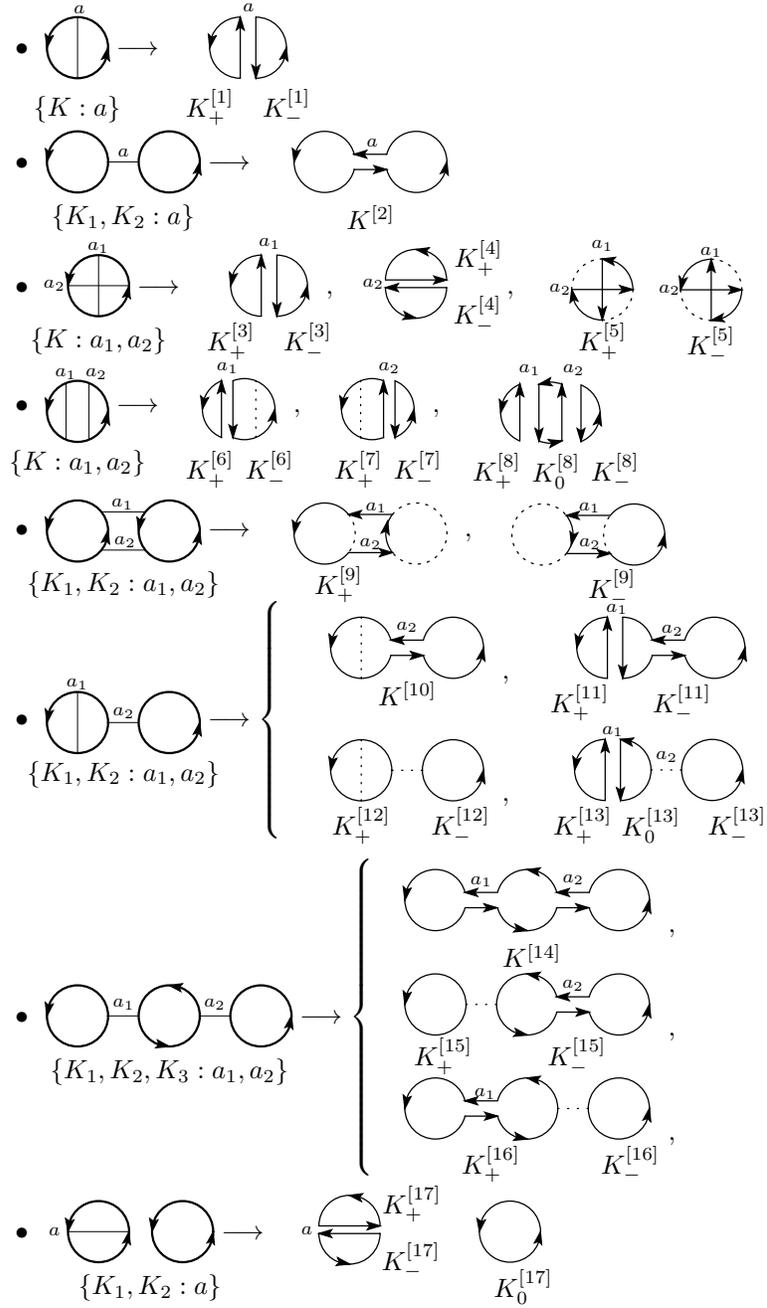
 
\begin{eqnarray}
  &&\bullet~~\picca\longrightarrow~~~~\piccb\nonumber\\[0.5cm]
  &&\bullet~~\piccc\longrightarrow~~~~\piccd\nonumber\\[0.3cm]
  &&\bullet~~\picce\longrightarrow~~~~\piccf~~,~~~~\piccg~~,
  ~~~~\picch~~~~\picci
  \nonumber\\[0.3cm]
  &&\bullet~~\piccj\longrightarrow~~~~\picck~~,~~~~\piccl~~,
  ~~~~~~\piccm
  \nonumber\\[0.5cm]
  &&\bullet~~\piccn\longrightarrow~~~~\picco~~,~~~~\piccp
  \nonumber\\[0.6cm]
  &&\bullet~~\piccq\longrightarrow\left\{
  \begin{array}{cc}
  ~~~~\piccr~~, & ~~~~\piccs \\
  & \\
  ~~~~\picct~~, & ~~~~\piccu \\
  & \\
  \end{array}
  \right.
  \nonumber\\[0cm]
  &&\bullet~~\piccw\longrightarrow\left\{
  \begin{array}{cc}
  ~~\piccx~~, &  \\
  & \\
  ~~\piccy~~, &  \\
  & \\
  ~~\piccz~~, &  \\
  & \\
  \end{array}
  \right.
  \nonumber\\[-0.5cm]
  &&\bullet~~\picgh\picgj\longrightarrow~~~~\picgk~~~~\picgi~~
  \nonumber\\ \nonumber
\end{eqnarray}
    \caption{splitting the crossing to obtain new knot diagrams}
	\label{fig:split}
\end{figure}
\end{definition}

\newpage
\subsection{The main part of the proof of Theorem \ref{th:mainresult}}
This section is the main part of the proof of Theorem \ref{th:mainresult}.
The argument is inductive, that is, we shall use the result of $v_{1}$ and $v_{2}$ for the proof of higher degrees $v_{3.1},v_{3.2},v_{4.1},v_{4.2},v_{4.3},v_{4.4}$.
\subsubsection{$v_{1}$}
\noindent{\bf Proof of} (\ref{eqn:main result12}) {\bf in Theorem \ref{th:mainresult}.} 
We expand (\ref{eq:v12}) according to (\ref{eqn:logtheta}) and obtain:
\begin{eqnarray*}
  v_{1}(A^b(\{K_1,K_2\}))=\hidarii A^b(\{K_1,K_2\}),\marucb\migii
  +\hidarii A^b(\{K_1,K_2\}),\marudo\migii.
\end{eqnarray*}
The second term in the right side of this equation identically vanishis, since 
\linebreak[4]
$v_{1}(A^b(\{K_1,K_2\}))$ and the first term in the right side of this equation are real valued and the second term is pure imaginary. Inserting (\ref{eqn:crossingnumber0}) into the first term yields (\ref{eqn:main result12}). 
$\square$
\newline

From this proof, we obtain the following lemma.
\begin{lemma}\hspace{-5pt}{\bf.}\hspace{5pt}~~\label{lemma:v12}
$\displaystyle
  \hidarii A^b(\{K_1,K_2\}),\marudo\migii=0.~~\square
$
\end{lemma}
Later we shall use this lemma in the proof of higher degrees. 

\subsubsection{$v_2$}
\noindent{\bf Proof of} (\ref{eqn:main result2}) {\bf in Theorem \ref{th:mainresult}.}
For sufficiently small $b$, inserting $(\ref{eqn:decomposition1})$ into
(\ref{eq:v2}) and using $(\ref{eqn:crossingnumber1})$, we have 
\begin{eqnarray}\label{eqn:v2-1}
  &&\hspace{-1.5cm}
  v_2(A^b(K))=\hidarii G(K),\marubb\migii_\chi\nonumber\\
  &&\hspace{1cm}
  +\frac{1}{2}\hidarii G^e(K),\marube\migii_{\chi^e}
  +\hidarii A^b(K),\marubd\migii-\frac{1}{6}m(K)+O(b),
\end{eqnarray}
where we have dropped the imaginary part, since $v_2(A^b(K))$ is real valued. 
We replace $K$ by $\alpha(K)$ in (\ref{eqn:v2-1}) and obtain
\begin{eqnarray}\label{eqn:v2-2}
  &&\hspace{-1.5cm}
  v_2(A^b(\alpha(K)))=\hidarii G(\alpha(K)),\marubb\migii_\chi\nonumber\\
  &&\hspace{1cm}
  +\frac{1}{2}\hidarii G^e(K),\marube\migii_{\chi^e}
  +\hidarii A^b(K),\marubd\migii-\frac{1}{6}m(K)+O(b).
\end{eqnarray} 
The second, third and fourth terms on the right hand side of (\ref{eqn:v2-1}) and (\ref{eqn:v2-2}) take the same value for $K$ and $\alpha(K)$, since they are independent on the signatures.
Since $v_2$ is a knot invariant, $v_2(A^b(\alpha(K)))$ is equal to 
$v_2(U)$, where $U$ denotes a trivial knot. Subtracting (\ref{eqn:v2-2}) from 
(\ref{eqn:v2-1}) yields the Gauss diagram formula (\ref{eqn:main result2}) in Theorem \ref{th:mainresult}:
\begin{eqnarray}\label{eqn:v2-3}
  v_2(A^b(K))=v_2(U)+\hidarii \bar{G}(K),\marubb\migii_{\chi}, 
\end{eqnarray}
where $\bar{G}(K)=G(K)-G(\alpha(K))$, $v_2(U)=-\frac{1}{6}$ and we have dropped $b$-dependent term $O(b)$ since $v_2(A^b(K))$ does not depend on $b$.
$\square$
\newline

From (\ref{eqn:v2-1}) and (\ref{eqn:v2-3}), we obtain the following lemma.
\begin{lemma}\hspace{-5pt}{\bf.}\hspace{5pt}~~\label{lemma:v2-4}
\begin{eqnarray*}
   &&\hspace{-1cm}
   \hidarii A^b(K),\marubd\migii=-\hidarii G(\alpha(K)),\marubb\migii_\chi
  -\frac{1}{2}\hidarii G^e(K),\marube\migii_{\chi^e}\\
  &&\hspace{6cm}
  -\frac{1}{6}(1-m(K))+O(b).
  ~~\square
\end{eqnarray*}
\end{lemma}
Later we shall use this lemma for the proof of higher degrees. Notice that we have obtained the Gauss diagram formula for the difficult integral $\hidarii A^b(K),\marubd\migii$.

\subsubsection{$v_{3.1}$}
To prepare for the proof of (\ref{eqn:main result31}) in Theorem \ref{th:mainresult},
we shall first prove Lemma \ref{lemma:v31-1} and Lemma \ref{lemma:v31-2}.
\begin{lemma}\hspace{-5pt}{\bf.}\hspace{5pt}\label{lemma:v31-1}~~
\begin{eqnarray*}
  \hh\{A^b(K):a\},\maruia+\marudk\mm=\hidarii A^b(K),\marubd\migii-\sum_{s=\pm}
  \hidarii A^b(K^{[1]}_s),\marubd\migii+O(b).
\end{eqnarray*}
\end{lemma}
{\bf Proof.}~~Use the identities: 
\begin{eqnarray*}
  &&\hidarii A^b(K),\marubd\migii
  =\hh\{A^b(K):a\},\maruib+\maruia+\marudk\mm,\nonumber\\
  &&\sum_{s=\pm}\hidarii A^b(K^{[1]}_s),\marubd\migii
  =\hh\{A^b(K):a\},\maruib\mm+O(b),
\end{eqnarray*}
where $K^{[1]}_s$ are given in Definition \ref{def:NewKnotDiagram}. $\square$

\begin{lemma}\hspace{-5pt}{\bf.}\hspace{5pt}\label{lemma:v31-2}~~
$\displaystyle
  \hh\{A^b(K):a\},\marudi+\marudk\mm=O(b).
$  
\end{lemma}
{\bf Proof.}
\begin{eqnarray*}
  \hh\{A^b(K):a\},\marudi+\marudk\mm&
  =&\frac{1}{2}\hh\{A^b(K):a\},\maruic\mm^2\nonumber\\
  &=&\frac{1}{2}\hidarii A^b\bigl(\{K^{[1]}_+,K^{[1]}_-\}\bigr),\marudo\migii^2
  +O(b)\nonumber\\[0.1cm]
  &=&O(b)
\end{eqnarray*}
The last step follows from Lemma \ref{lemma:v12}. $\square$
\newline

\noindent{\bf Proof of} (\ref{eqn:main result31}) {\bf in Theorem \ref{th:mainresult}.}~~
Inserting (\ref{eqn:decomposition2}),(\ref{eqn:decomposition3}) into 
(\ref{eq:v31}) and using (\ref{eqn:spilit1}), (\ref{eqn:spilit2}),(\ref{eqn:spilit3}) yields:
\begin{eqnarray*}
  \lefteqn{v_{3.1}(A^b(K))}\\
  &&=\hidarii A^b(K),\marudd+2\marude\migii\nonumber\\
  &&+\sum_{a}
  \hidarii P(\{K:a\}),\marua\migii_\chi \hh\{A^b(K):a\},
  \marudi+\maruia+2\marudk\mm+O(b)
  \nonumber\\
  &&=\hidarii A^b(K),\marudd+2\marude\migii\nonumber\\
  &&\hspace{-0.5cm}
  +\sum_{a}
  \hidarii P(\{K:a\}),\marua\migii_\chi \biggl\{\hidarii A^b(K),\marubd\migii-
  \sum_{s=\pm}\hidarii A^b(K^{[1]}_s),\marubd\migii\biggr\}+O(b).\\
\end{eqnarray*}
The last step follows from Lemma \ref{lemma:v31-1} and Lemma \ref{lemma:v31-2}.
We insert Lemma \ref{lemma:v2-4} and (\ref{eqn:crossingnumber311}) (\ref{eqn:crossingnumber312}) into this and use
\begin{eqnarray*}
  &&\bullet~~
  \sum_{a}\hidarii P(\{K:a\}),\marua\migii_\chi \hidarii G^e(K)-
  \sum_{s=\pm}G(K^{[1]}_s),\marube\migii_{\chi^e}\nonumber\\
  &&\hspace{6cm}
  =\hidarii G^e(K),\marudm~+\marudn\migii_{\chi^e},
  \nonumber\\
  &&\bullet~~
  \sum_{a}\hidarii P(\{K:a\}),\marua\migii_\chi
  \Bigl\{ \bigl(1-m(K)\bigr)
  -\sum_{s=\pm}\bigl(1-m(K^{[1]}_s)\bigr)\Bigr\}\nonumber\\
  &&\hspace{6cm}
  =-\hidarii G^e(K),\marudn\migii_{\chi^e}.
\end{eqnarray*}
Then we have
\begin{eqnarray*}
  &&\hspace{-0.8cm}
  v_{3.1}(A^b(K))=
  \hidarii G(K),~~2\marude+\marudd+\frac{1}{2}\marudp\migii_{\chi}
  \nonumber\\ 
  &&\hspace{1.2cm}
  -\sum_{a}\hidarii P(\{K:a\}),\marua\migii_\chi
  \hidarii G(\alpha(K))-\sum_{s=\pm}G(\alpha(K_s^{[1]})),
  \marubb\migii_{\chi},
\end{eqnarray*}
where we have dropped $b$-dependent term $O(b)$ since $v_{3.1}(A^b(K))$ does not depend on $b$. We can easily see that this equation is the same as the Gauss diagram formula (\ref{eqn:main result31}) in Theorem \ref{th:mainresult}.
$\square$

\subsubsection{$v_{3.2}$}
To prepare for the proof of (\ref{eqn:main result32}) in Theorem \ref{th:mainresult},
we shall prove Lemma \ref{lemma:v32-1} and Lemma \ref{lemma:v32-2}.
\newpage

\begin{lemma}\hspace{-5pt}{\bf.}\hspace{5pt}\label{lemma:v32-1}~~
If $P_0(\{K_1,K_2:a\})=\marucb$, then
\begin{eqnarray*}
  &&\hspace{-0.8cm}
  \hh\bigl\{A^b(\{K_1,K_2\}):a\bigr\},\maruid+\maruie\mm
   \nonumber\\ 
  &&\hspace{3cm}
  =\hidarii A^b(K^{[2]})
  ,\marubd\migii
  -\sum_{i=1}^2\hidarii A^b(K_i),\marubd\migii+O(b).
\end{eqnarray*}
holds.
\end{lemma}
{\bf Proof.}~~The proof is similar to Lemma \ref{lemma:v31-1}~~.
\begin{eqnarray*}
  &&
  \hh\bigl\{A^b(\{K_1,K_2\}):a\bigr\},\maruid+\maruie\mm\\
  &&\hspace{3cm}
  =\hh\{A^b(K^{[2]}):a\},\maruia+\marudk\mm+O(b)\nonumber\\
  &&\hspace{3cm}
  =\hidarii A^b(K^{[2]}),\marubd\migii
  -\sum_{i=1}^2\hidarii A^b(K_i),\marubd\migii+O(b).~~\square\nonumber
\end{eqnarray*}

\begin{lemma}\hspace{-5pt}{\bf.}\hspace{5pt}\label{lemma:v32-2}~~
$\displaystyle
  \hh\{A^b(\{K_1,K_2\}):a\},\marueg\mm=O(b).
$  
\end{lemma}
{\bf Proof.}
\begin{eqnarray*}
  \lefteqn{\hh\{A^b(\{K_1,K_2\}):a\},\marueg\mm} \nonumber\\
  &&=\hidarii A^b(\{K^{[17]}_+,K^{[17]}_0\}),\marudo\migii
  \hidarii A^b(\{K^{[17]}_-,K^{[17]}_0\}),\marudo\migii+O(b)
  \nonumber\\
  &&=O(b).
\end{eqnarray*}
The last step follows from Lemma \ref{lemma:v12}. $\square$
\newline

\noindent{\bf Proof of} (\ref{eqn:main result32}) {\bf in Theorem \ref{th:mainresult}.}~~
After inserting (\ref{eqn:decomposition4})
,(\ref{eqn:decomposition5}) into (\ref{eq:v32}), we use (\ref{eqn:spilit4}),(\ref{eqn:spilit5}),(\ref{eqn:spilit6}). Then we obtain
\begin{eqnarray*}
  \lefteqn{v_{3.2}(A^b(\{K_1,K_2\}))}\nonumber\\
  &&=\hidarii A^b(\{K_1,K_2\}),\marueb+2\marued\migii\nonumber\\
  &&~~+\sum_{a}
  \hidarii P({K_1,K_2:a}),\marucb\migii_{\chi}\\
  &&\hspace{3cm}\times
  \hh \{A^b(\{K_1,K_2\}):a\},\maruid+\maruie\mm\nonumber\\
  &&~~+\sum_{a}
  \hidarii P({K_1,K_2:a}),\marua\maru\migii_{\chi}\\
  &&\hspace{3cm}\times
  \hh \{A^b(\{K_1,K_2\}):a\},\marueg\mm+O(b)
  \nonumber\\
  &&=\hidarii A^b(\{K_1,K_2\}),\marueb+2\marued\migii\nonumber\\
  &&~~+\sum_{a}
  \hidarii P({K_1,K_2:a}),\marucb\migii_{\chi}\\
  &&\hspace{3cm}\times
  \biggl\{\hidarii A^b(K^{[2]}),\marubd\migii
  -\sum_{i=1}^2\hidarii A^b(K_i),\marubd\migii\biggr\}+O(b).
\end{eqnarray*}
The last step follows from Lemma \ref{lemma:v32-1} and \ref{lemma:v32-2}.
We insert Lemma \ref{lemma:v2-4} and (\ref{eqn:crossingnumber321}) (\ref{eqn:crossingnumber322}) into this and use
\begin{eqnarray*}
  &&\bullet~~
  \sum_{a}\hidarii P({K_1,K_2:a}),\marucb\migii_{\chi}
  \hidarii G^e(K^{[2]})-\sum_{i=1,2}G^e(K_i),\marube\migii_{\chi^e}\nonumber\\
  &&\hspace{3cm}=\hidarii G^e(K_1,K_2),\maruej\migii_{\chi^e},
  \nonumber\\
  &&\bullet~~
  \sum_{a}\hidarii P({K_1,K_2:a}),\marucb\migii_{\chi}
  \Bigl\{ \bigl(1-m(K^{[2]})\bigr)
  -\sum_{i=1,2}\bigl(1-m(K_i)\bigr)\Bigr\}\nonumber\\
  &&\hspace{3cm}=-\hidarii G^e(\{K_1,K_2\}),\marucb+\maruel\migii_{\chi^e}.
\end{eqnarray*}
Then we have 
\begin{eqnarray*}
  \lefteqn{v_{3.2}(A^b(\{K_1,K_2\}))}\\
  &&=\hidarii G(K_1,K_2),~~\marueb
  +\marued+\frac{1}{3}\marucb\migii_{\chi}\nonumber\\
  &&-\sum_{a}\hidarii P({K_1,K_2:a}),\marucb\migii_{\chi}
  \hidarii G(\alpha(K^{[2]}))
  -\sum_{i=1,2}G(\alpha(K_i)),\marubb\migii_{\chi}.\nonumber
\end{eqnarray*}
This is the same as the Gauss diagram formula (\ref{eqn:main result32}) in Theorem \ref{th:mainresult}.
$\square$

\subsubsection{$v_{4.1}, v_{4.2}, v_{4.3}, v_{4.4}$}
The computaion of degree four $v_{4.1},v_{4.2},v_{4.3},v_{4.4}$ are long but straightforward. We use the same argument as degree two and three. 
\newline 

\noindent{\bf Proof of} (\ref{eqn:main result41}) {\bf in Theorem \ref{th:mainresult}.} 
We calculate $v_{4.1}$ in the same way as the lower degree, and obtain:
\begin{eqnarray*}
  \lefteqn{v_{4.1}(A^b(K))}\\
  &&=\hidarii A^b(K),~~\maruga+\marugb+2\marugc+4\marugd
  +5\maruge+7\marugf\migii \nonumber\\
  &&+\sum_{(a_1,a_2)}\hidarii P({K:a_1,a_2}),\marubb\migii_{\chi}
  \nonumber\\
  &&\hspace{-0.5cm}\times
  \biggl\{3\hidarii A^b(K),\marubd\migii
  -2\sum_{n=3,4}\sum_{s=\pm}\hidarii A^b(K^{[n]}_s),\marubd\migii
  +\sum_{s=\pm}\hidarii A^b(K^{[5]}_s),\marubd\migii\biggr\}
  \nonumber\\
  &&+\sum_{(a_1,a_2)}\hidarii P({K:a_1,a_2}),\maruba\migii_{\chi}
  \nonumber\\
  &&\hspace{-0.3cm}\times
  \biggl\{\hidarii A^b(K),\marubd\migii
  -\sum_{n=6,7}\sum_{s=\pm}\hidarii A^b(K^{[n]}_s),\marubd\migii
  +\sum_{s=\pm,0}\hidarii A^b(K^{[8]}_s),\marubd\migii
  \biggr\}\nonumber\\
  &&+(\mbox{signature independent terms})+O(b),
\end{eqnarray*}
where we have used Lemma \ref{lemma:v12}. In the above equation, 
"(signature independent terms)" means the terms which take the same value for $K$ and $\alpha(K)$.  
After inserting Lemma \ref{lemma:v2-4} into this, we use 
\begin{eqnarray*}
  &&\bullet~~
  \sum_{(a_1,a_2)}\hidarii P({K:a_1,a_2}),\marubb\migii_{\chi}\nonumber\\
  &&\hspace{2cm}\times
  \hidarii 3~G(K)-2\sum_{n=3,4}\sum_{s=\pm}G(K^{[n]}_s)
  +\sum_{s=\pm}G(K^{[5]}_s),\marube\migii_{\chi} \nonumber\\
  &&\hspace{4cm}=\hidarii G(K),3\marujd+\maruja~+2\marujc+\marugi\migii_{\chi},
  \nonumber\\
  &&\bullet~~
  \sum_{(a_1,a_2)}\hidarii P({K:a_1,a_2}),\marubb\migii_{\chi}\nonumber\\
  &&\hspace{1.5cm}\times
  \Bigl\{ 3\bigl(1-m(K)\bigr)
  -2\sum_{n=3,4}\sum_{s=\pm}\bigl(1-m(K^{[n]}_s)\bigr)
  +\sum_{s=\pm}\bigl(1-m(K^{[5]}_s)\bigr)\Bigr\} \nonumber\\
  &&\hspace{4cm}=\hidarii G(K),\marubb-3\marujd\migii_{\chi},
\end{eqnarray*}
and
\begin{eqnarray*}
  &&\bullet~~
  \sum_{(a_1,a_2)}\hidarii P({K:a_1,a_2}),\maruba\migii_{\chi}\nonumber\\
  &&\hspace{2cm}\times
  \hidarii G(K)-\sum_{n=6,7}\sum_{s=\pm}G(K^{[n]}_s)
  +\sum_{s=\pm,0}G(K^{[8]}_s),\marube\migii_{\chi}
  \nonumber\\
  &&\hspace{6cm}
  =\hidarii G(K),\marujb\migii_{\chi},\nonumber\\
  &&\bullet~~
  \Bigl\{ \bigl(1-m(K)\bigr)
  -\sum_{n=6,7}\sum_{s=\pm}\bigl(1-m(K^{[n]}_s)\bigr)
  +\sum_{s=\pm,0}\bigl(1-m(K^{[8]}_s)\bigr)\Bigr\}=0.
\end{eqnarray*}
Lastly we cancel the signature independent terms by using $\alpha(K)$ in the same argument as in the proof of $v_2$. Then we have
\begin{eqnarray*}
  \lefteqn{v_{4.1}(A^b(K))}\\
  &&=\frac{1}{360}+\hidarii \bar{G}(K),~\maruga+\marugb+2\marugc+4\marugd
  +5\maruge+7\marugf\migii_{\chi} \nonumber\\
  &&\hspace{3cm}+\hidarii \bar{G}(K),~
  \frac{1}{6}\marubb+\frac{1}{2}\marugg+2\marugh+2\marugi\migii_{\chi} 
  \nonumber\\
  &&-\sum_{(a_1,a_2)}\hidarii \bar{P}(\{K:a_1,a_2\}),\marubb\migii_{\chi}
  \nonumber\\
  &&\hspace{1cm}\times
  \hidarii 3~G(\alpha(K))-2\sum_{n=3,4}\sum_{s=\pm}G(\alpha(K^{[n]}_s))
  +\sum_{s=\pm}G(\alpha(K^{[5]}_s)),\marubb\migii_{\chi} 
  \nonumber\\
  &&-\sum_{(a_1,a_2)}\hidarii \bar{P}(\{K:a_1,a_2\}),\maruba\migii_{\chi}
  \nonumber\\
  &&\hspace{1cm}\times
  \hidarii G(\alpha(K))-\sum_{n=6,7}\sum_{s=\pm}G(\alpha(K^{[n]}_s))
  +\sum_{s=\pm,0}G(\alpha(K^{[8]}_s)),\marubb\migii_{\chi}. 
  \nonumber\\
\end{eqnarray*}
This is the same as the Gauss diagram formula (\ref{eqn:main result41}) in Theorem \ref{th:mainresult}.  $\square$
\newline 

\noindent{\bf Proof of} (\ref{eqn:main result42}) {\bf in Theorem \ref{th:mainresult}.} We calculate $v_{4.2}$ in the same way as the lower degree and obtain:
\begin{eqnarray*}
  \lefteqn{v_{4.2}(A^b(K))}\\
  &&=
  \hidarii A^b(K),~~\marugd+\maruge+\marugf\migii \nonumber\\
  &&+\sum_{(a_1,a_2)}\hidarii P({K:a_1,a_2}),\marubb\migii_{\chi}
  \nonumber\\
  &&\hspace{0cm}\times
  \biggl\{\hidarii A^b(K),\marubd\migii
  -\sum_{n=3,4}\sum_{s=\pm}\hidarii A^b(K^{[n]}_s),\marubd\migii
  +\sum_{s=\pm}\hidarii A^b(K^{[5]}_s),\marubd\migii\biggr\}
  \nonumber\\
  &&+(\mbox{signature independent term})+O(b),
\end{eqnarray*}
where we have used Lemma \ref{lemma:v12}. After inserting Lemma \ref{lemma:v2-4} into this, we use 
\begin{eqnarray*}
  &&\bullet~~
  \sum_{(a_1,a_2)}\hidarii P({K:a_1,a_2}),\marubb\migii_{\chi}
  \nonumber\\
  &&\hspace{2cm}\times
  \hidarii G(K)-\sum_{n=3,4}\sum_{s=\pm}G(K^{[n]}_s)
  +\sum_{s=\pm}G(K^{[5]}_s),\marube\migii_{\chi} \nonumber\\
  &&\hspace{4cm}=\hidarii G(K),\marugi+\marujd \migii_{\chi},\nonumber\\
  &&\bullet~~
  \sum_{(a_1,a_2)}\hidarii P({K:a_1,a_2}),\marubb\migii_{\chi}
  \nonumber\\
  &&\hspace{1.5cm}\times
  \Bigl\{ \bigl(1-m(K)\bigr)
  -\sum_{n=3,4}\sum_{s=\pm}\bigl(1-m(K^{[n]}_s)\bigr)
  +\sum_{s=\pm}\bigl(1-m(K^{[5]}_s)\bigr)\Bigr\} \nonumber\\
  &&\hspace{4cm}=\hidarii G(K),-\marubb+\maruje-\marujd\migii_{\chi}.
\end{eqnarray*}
Lastly we cancel the signature independent terms by using $\alpha(K)$ in the same argument as in the proof of $v_2$. Then we have
\begin{eqnarray*}
  \lefteqn{v_{4.2}(A^b(K))}\\
  &&=-\frac{1}{360}+\hidarii \bar{G}(K),~~\marugd+\maruge+\marugf
  +\frac{1}{2}\marugh-\frac{1}{6}\marubb\migii_{\chi}
  \nonumber\\
  &&~~-\sum_{(a_1,a_2)}\hidarii \bar{P}(\{K:a_1,a_2\}),\marubb\migii_{\chi}
  \nonumber\\
  &&\hspace{1.5cm}\times
  \hidarii G(\alpha(K))-\sum_{n=3,4}\sum_{s=\pm}G(\alpha(K^{[n]}_s))
  +\sum_{s=\pm}G(\alpha(K^{[5]}_s)),\marubb\migii_{\chi}. 
  \nonumber\\
\end{eqnarray*}
This is the same as the Gauss diagram formula (\ref{eqn:main result42}) in Theorem \ref{th:mainresult}.  $\square$
\newline

\noindent
{\bf Proof of} (\ref{eqn:main result43}) {\bf in Theorem \ref{th:mainresult}.}
We calculate $v_{4.3}$ in the same way as the lower degree, and we obtain 
\begin{eqnarray*}
  \lefteqn{v_{4.3}(A^b(\{K_1,K_2\}))}\\
  &&=\hidarii A^b(\{K_1,K_2\}),\maruha+\maruhb
  +2\maruhc \nonumber\\
  &&\hspace{4cm}+\maruhd+\maruhe+\maruhf\migii\nonumber\\
  &&+\sum_{(a_1,a_2)}\hidarii P({K_1,K_2:a_1,a_2}),\maruca\migii_{\chi}
   \nonumber\\
  &&\hspace{1.5cm}\times
  \biggl\{\sum_{i=1,2}\hidarii A^b(K_i),\marubd\migii
  -\sum_{s=\pm}\hidarii A^b(K^{[9]}_s),\marubd\migii\biggr\}\nonumber\\
  &&+\sum_{(a_1,a_2)}\hidarii P({K_1,K_2:a_1,a_2}),\marucc\migii_{\chi}
  \nonumber\\
  &&
  \times\biggl\{\hidarii A^b(K^{[10]}),\marubd\migii
  -\sum_{n=11,12}\sum_{s=\pm}\hidarii A^b(K^{[n]}_s),\marubd\migii
  \nonumber\\
  &&\hspace{6cm}
  +\sum_{s=\pm,0}\hidarii A^b(K^{[13]}_s),\marubd\migii\biggr\},\nonumber\\
\end{eqnarray*}
where we have used Lemma \ref{lemma:v12}. After inserting Lemma \ref{lemma:v2-4} into this, we use 
\begin{eqnarray*}
  &&\hspace{-1cm}\bullet~
  \sum_{(a_1,a_2)}\hidarii P({K_1,K_2:a_1,a_2}),\maruca\migii_{\chi}
  \hidarii \sum_{i=1,2}G(K_i)
  -\sum_{s=\pm}G(K^{[9]}_s),\marube\migii_{\chi} \nonumber\\
  &&\hspace{3cm}=\hidarii G(\{K_1,K_2\}),\marujf-\marujg\migii_{\chi},
  \nonumber\\
  &&\hspace{-1cm}\bullet~
  \sum_{(a_1,a_2)}\hidarii P({K_1,K_2:a_1,a_2}),\maruca\migii_{\chi}
  \nonumber\\
  &&\hspace{2cm}\times
  \Bigl\{\sum_{i=1,2}\bigl(1-m(K_i)\bigr)
  -\sum_{s=\pm}\bigl(1-m(K^{[9]}_s)\bigr)\Bigr\} \nonumber\\
  &&\hspace{4cm}=\hidarii G(\{K_1,K_2\}),\marujh\migii_{\chi},
\end{eqnarray*}
and
\begin{eqnarray*}
  &&\bullet~
  \sum_{(a_1,a_2)}\hidarii P({K_1,K_2:a_1,a_2}),\marucc\migii_{\chi}
  \nonumber\\
  &&\hspace{1.5cm}\times
  \hidarii G(K^{[10]})-\sum_{n=11,12}\sum_{s=\pm}G(K^{[n]}_s)
  +\sum_{s=\pm,0}G(K^{[13]}_s),\marube\migii_{\chi}\nonumber\\
  &&\hspace{5cm}=\hidarii G(\{K_1,K_2\}),\maruji\migii_{\chi},\nonumber\\
  &&\bullet~
  \bigl(1-m(K^{[10]})\bigr)
  -\sum_{n=11,12}\sum_{s=\pm}\bigl(1-m(K^{[n]}_s)\bigr)
  +\sum_{s=\pm,0}\bigl(1-m(K^{[13]}_s)\bigr)=0.
\end{eqnarray*}
Lastly we cancel the signature independent terms by using $\alpha(K)$ in the same argument as in the proof of $v_2$. Then we have
\begin{eqnarray*}
  \lefteqn{v_{4.3}(A^b(\{K_1,K_2\}))}\\
  &&=\hidarii \bar{G}(\{K_1,K_2\}),
  ~\maruha+\maruhb+2\maruhc+\maruhd\migii_{\chi}
  \nonumber\\
  &&\hspace{-0.1cm}
  +\hidarii \bar{G}(\{K_1,K_2\}),
  \maruhe+\maruhf+\frac{1}{2}\maruhg+\frac{1}{2}\maruhh\migii_{\chi} 
  \nonumber\\
  &&-\sum_{(a_1,a_2)}\hidarii \bar{P}(\{K_1,K_2:a_1,a_2\}),\maruca\migii_{\chi}
  \nonumber\\
  &&\hspace{4cm}\times
  \hidarii \sum_{i=1,2}G(\alpha(K_i))
  -\sum_{s=\pm}G(\alpha(K^{[9]}_s)),\marubb\migii_{\chi} 
  \nonumber\\
  &&-\sum_{(a_1,a_2)}\hidarii \bar{P}(\{K_1,K_2:a_1,a_2\}),\maruhi\migii_{\chi}
  \nonumber\\
  &&\hspace{0cm}\times
  \hidarii G(\alpha(K^{[10]}))-\sum_{n=11,12}\sum_{s=\pm}G(\alpha(K^{[n]}_s))
  +\sum_{s=\pm,0}G(\alpha(K^{[13]}_s)),\marubb\migii_{\chi} 
  \nonumber\\
\end{eqnarray*}
This is the same as the Gauss diagram formula (\ref{eqn:main result43}) in Theorem \ref{th:mainresult}.  $\square$
\newline

\noindent
{\bf Proof of} (\ref{eqn:main result44}) {\bf in Theorem \ref{th:mainresult}.}
We calculate $v_{4.4}$ in the same way as the lower degree, and we obtain
\begin{eqnarray*}
  \lefteqn{v_{4.4}(A^b(\{K_1,K_2,K_3\}))}\\
  &&=\hidarii A^b(\{K_1,K_2,K_3\}),
  ~~\maruhj+\maruhk+\maruhl\migii\nonumber\\[0.2cm]
  &&\hspace{0cm}
  +\sum_{(a_1,a_2)}\hidarii P({K_1,K_2,K_3:a_1,a_2}),\maruhm\migii_{\chi}
  \nonumber\\
  &&\hspace{1cm}\times\biggl\{\hidarii A^b(K^{[14]}),\marubd\migii
  -\sum_{n=15,16}\sum_{s=\pm}\hidarii A^b(K^{[n]}_s),\marubd\migii
  \nonumber\\
  &&\hspace{6cm}
  +\sum_{i=1,2,3}\hidarii A^b(K_i),\marubd\migii\biggr\},\nonumber\\
\end{eqnarray*}
where we have used Lemma \ref{lemma:v12}. After inserting Lemma \ref{lemma:v2-4} into this, we use 
\begin{eqnarray*}
  &&\bullet~~
  \sum_{(a_1,a_2)}\hidarii P({K_1,K_2,K_3:a_1,a_2}),\maruhm\migii_{\chi}
  \nonumber\\
  &&\hspace{1cm}\times
  \hidarii G(K^{[14]})-\sum_{n=15,16}\sum_{s=\pm}G(K^{[n]}_s)
  +\sum_{i=1,2,3}G(K_i),\marube\migii_{\chi}\nonumber\\
  &&\hspace{3cm}=\hidarii G(\{K_1,K_2,K_3\}),\marujj\migii_{\chi},
  \nonumber\\[0.5cm]
  &&\bullet~~\bigl(1-m(K^{[14]})\bigr)
  -\sum_{n=15,16}\sum_{s=\pm}\bigl(1-m(K^{[n]}_s)\bigr)
  +\sum_{i=1,2,3}\bigl(1-m(K_i)\bigr)=0.
\end{eqnarray*}
Lastly we cancel the signature independent terms by using $\alpha(K)$ in the same argument as in the proof of $v_2$. Then we have 
\begin{eqnarray*}
  \lefteqn{v_{4.4}(A^b(\{K_1,K_2,K_3\}))}\\
  &&=\hidarii \bar{G}(\{K_1,K_2,K_3\}),~~\maruhj+\maruhk+\maruhl\migii_{\chi} 
  \nonumber\\[0.2cm]
  &&-\sum_{(a_1,a_2)}\hidarii \bar{P}(\{K_1,K_2,K_3:a_1,a_2\})
  ,\maruhm\migii_{\chi}
  \nonumber\\
  &&\hspace{1cm}\times
  \hidarii G(\alpha(K^{[14]})) 
  -\sum_{n=15,16}\sum_{s=\pm}G(\alpha(K^{[n]}_s))
  +\sum_{i=1}^3G(\alpha(K_i)),\marubb\migii_{\chi}. 
  \nonumber\\
\end{eqnarray*}
This is the same as the Gauss diagram formula (\ref{eqn:main result44}) in Theorem \ref{th:mainresult}.  $\square$
\newline

\noindent{\it Remark.}
Notice the concept of direction "s", "d" disappear in the final Gauss diagram formula, as we have expected.

\section{Consistency Check}
We write $\hat{P}_L^{(4)}$ for the right hand side of (\ref{eqn:Homflyformula}):\begin{eqnarray} \label{eqn:Consistencycheck1}
  \hat{P}_{L}^{(4)}
  =W_{su(N)}^{(4)}\Biggl(N^{n-1}\biggl\{\exp\Bigl(
  \sum_{D\in\mathfrak{D}_K}D~u(D:L)
  \Bigr)\biggr\}\biggl\{\sum_{D\in\mathfrak{D}_L}D~w(D:L)
  \biggr\}\Biggr).
\end{eqnarray}
From Corollary \ref{co:Homfly}, it is trivial that $\hat{P}_L^{(4)}$ satisfies the Homfly skein relation (\ref{eqn:skein relation}) up to degree four. 
But as a consitency check of the Gauss diagram formula, we will prove directly that $\hat{P}_L^{(4)}$ really satisfies the HOMFLY skein relation up to degree four {\it by using the Gauss diagram formula in Theorem} \ref{th:mainresult}:  
\begin{eqnarray}\label{eqn:Consistencycheck2}
  \Bigl[\exp(\frac{Nx}{2})\hat{P}_{L_+}^{(4)}
  -\exp(-\frac{Nx}{2})\hat{P}_{L_-}^{(4)}
  -(e^{\frac{x}{2}}-e^{-\frac{x}{2}})\hat{P}_{L_0}^{(4)}~\Bigr]^{(4)}=0.
\end{eqnarray}
 
\noindent{\bf Proof.}
There are two cases:
\begin{itemize}
  \item[(1)] 
  $L_+$ and $L_-$ have $n$-component, while $L_0$ has $(n+1)$-components.
  \item[(2)] 
  $L_+$ and $L_-$ have $(n+1)$-component, while $L_0$ has $n$-components.
\end{itemize}
First, we consider the case (1). Set 
\begin{eqnarray*}
   L_+&=&\{K_1,\cdot\cdot\cdot,K_{n-1},K_n^+\}\nonumber\\
   L_-&=&\{K_1,\cdot\cdot\cdot,K_{n-1},K_n^-\}\nonumber\\
   L_0&=&\{K_1,\cdot\cdot\cdot,K_{n-1},K_n^0,K_{n+1}^0\}.\nonumber\\
\end{eqnarray*}
The $(n-1)$-components $K_1,\cdot\cdot\cdot,K_{n-1}$ are common in $L_+,L_-,L_0$, while $K_n^+$, $K_n^-$, $K_n^0$, $K_{n+1}^0$ are the same except inside the dashed circle: 
\begin{center}
\setlength{\unitlength}{1cm}
\begin{picture}(10,2)(1.5,0.5)
\put(2,1){\picfr}
\put(5,1){\picfs}
\put(8,1.1){\picft}
\end{picture}
\end{center}
Inserting (\ref{eqn:Consistencycheck1}) into the left side of (\ref{eqn:Consistencycheck2}) and using Appendix D, we have
\begin{eqnarray*}
  &&\hspace{-1cm}
  \Bigl[\exp(\frac{Nx}{2})\hat{P}_{L_+}^{(4)}
  -\exp(-\frac{Nx}{2})\hat{P}_{L_-}^{(4)}
  -(e^{\frac{x}{2}}-e^{-\frac{x}{2}})\hat{P}_{L_0}^{(4)}~\Bigr]^{(4)}\\
  &&\hspace{-1cm}
  =-\frac{1}{4}(N^2-1)x^2~V_1+\frac{1}{8}N(N^2-1)x^3~V_2
  -\frac{(N^2-1)}{8N}x^3~V_3-\frac{N^2(N^2-1)}{16}x^4~V_4\\
  &&+\frac{(N^2-1)(N^2+2)}{16}x^4~V_5
  +\frac{(N^2-1)}{16}x^4~V_6-\frac{(N^2-1)}{16N^2}x^4~V_7,
\end{eqnarray*}
where,
\begin{eqnarray*}
  &&\bullet~~V_1=\Bigl\{\sum_{s=\pm}s~v_2(K_n^s)\Bigr\}
  -2v_{1}\bigl(\{K_n^0,K_{n+1}^0\}\bigr),\nonumber\\
  &&\bullet~~V_2=\Bigl\{\sum_{s=\pm}s~v_{3.1}(K_n^s)\Bigr\}
  -\Bigl\{\sum_{s=\pm}v_2(K_n^s)\Bigr\}
  +2\Bigl\{\sum_{i=n}^{n+1}v_2(K_i^0)\Bigr\}
  \nonumber\\
  &&\hspace{6cm}
  -\Bigl\{v_{1}(K_n^0,K_{n+1}^0)\Bigr\}^2
  +\frac{1}{3},\nonumber\\
  &&\bullet~~V_3=\sum_{i=1}^{n-1}\biggl[\Bigl\{\sum_{s=\pm}
  s~v_{3.2}\bigl(\{K_i,K_n^s\}\bigr)\Bigr\}
  -2v_{1}(K_i,K_n^0)v_{1}\bigl(\{K_i,K_{n+1}^0\}\bigr)\biggr],\nonumber\\
  &&\bullet~~V_4=\Bigl\{\sum_{s=\pm}s~v_{4.1}(K_n^s)\Bigr\}
  -\Bigl\{\sum_{s=\pm}v_{3.1}(K_n^s)\Bigr\}
  +2\Bigl\{\sum_{i=n}^{n+1}v_{3.1}(K_i^0)\Bigr\}\nonumber\\
  &&\hspace{1cm}
  +v_{3.2}\bigl(\{K_n^0,K_{n+1}^0\}\bigr)
  -\frac{3}{2}v_{1}\bigl(\{K_n^0,K_{n+1}^0\}\bigr)
  \Bigl\{\sum_{s=\pm}v_2(K_n^s)-2\sum_{i=n}^{n+1}v_2(K_i^0)\Bigr\}\nonumber\\
  &&\hspace{5cm}-\frac{1}{3}\Bigl\{v_{1}(K_n^0,K_{n+1}^0)\Bigr\}^3
  +\frac{7}{6}v_{1}\bigl(\{K_n^0,K_{n+1}^0\}\bigr),\nonumber\\
  &&\bullet~~V_5=\Bigl\{\sum_{s=\pm}s~v_{4.2}(K_n^s)\Bigr\}
  +v_{3.2}\bigl(\{K_n^0,K_{n+1}^0\}\bigr)
  +\frac{1}{6}v_{1}\bigl(\{K_n^0,K_{n+1}^0\}\bigr)\nonumber\\
  &&\hspace{4cm}-\frac{1}{2}v_{1}\bigl(\{K_n^0,K_{n+1}^0\}\bigr)
  \Bigl\{\sum_{s=\pm}v_2(K_n^s)-2\sum_{i=n}^{n+1}v_2(K_i^0)\Bigr\},
  \nonumber\\
  &&\bullet~~V_6=\sum_{i=1}^{n-1}\biggl[\Bigl\{\sum_{s=\pm}
  s~v_{4.3}\bigl(\{K_i,K_n^s\}\bigr)\Bigr\}
  -\Bigl\{\sum_{s=\pm}v_{3.2}\bigl(\{K_i,K_n^s\}\bigr)\Bigr\}\nonumber\\
  &&\hspace{0.5cm}
  +2\Bigl\{\sum_{j=n}^{n+1}v_{3.2}\bigl(\{K_i,K_j^0\}\bigr)\Bigr\}
  -2v_{1}(K_n^0,K_{n+1}^0)v_{1}\bigl(\{K_i,K_{n}^0\}\bigr)
  v_{1}\bigl(\{K_i,K_{n+1}^0\}\bigr)\biggr],\nonumber\\
  &&\bullet~~V_7=\sum_{1\le i<j\le n-1}\biggl[\Bigl\{\sum_{s=\pm}
  s~v_{4.4}\bigl(\{K_i,K_j,K_n^s\}\bigr)\Bigr\}
  \nonumber\\
  &&\hspace{4cm}
  -2v_{1}\bigl(\{K_i,K_j\}\bigr)
  \Bigl\{v_{1}\bigl(\{K_i,K_{n}^0\}\bigr)v_{1}\bigl(\{K_j,K_{n+1}^0\}\bigr)
  \nonumber\\
  &&\hspace{6cm}
  +v_{1}\bigl(\{K_i,K_{n+1}^0\}\bigr)v_{1}\bigl(\{K_j,K_{n}^0\}\bigr)
  \Bigr\}\biggr].\nonumber\\
  &&
\end{eqnarray*}
Inserting the Guass diagram formula (Theorem \ref{th:mainresult}) into each $V_i$, we find out all of these equations vanishes identically $V_i=0$ ($i=1,\cdots,7$). This shows the skein relation (\ref{eqn:Consistencycheck2}) holds in case of (1).

Next we consider the case (2). Set 
\begin{eqnarray*}
   L_+&=&\{K_1,\cdot\cdot\cdot,K_{n-1},K_n^+,K_{n+1}^+\}\nonumber\\
   L_-&=&\{K_1,\cdot\cdot\cdot,K_{n-1},K_n^-,K_{n+1}^-\}\nonumber\\
   L_0&=&\{K_1,\cdot\cdot\cdot,K_{n-1},K_n^0\}.\nonumber\\
\end{eqnarray*}
The $(n-1)$-components $K_1,\cdot\cdot\cdot,K_{n-1}$ are common in $L_+,L_-,L_0$, while $K_n^+$, $K_{n+1}^+$, $K_n^-$, $K_{n+1}^-$, $K_n^0$ are the same except inside the dashed circle: 
\begin{center}
    \setlength{\unitlength}{1cm}
    \begin{picture}(10,2.5)(1.5,0.5)
    \put(2,1){\picfu}
    \put(5,1){\picfw}
    \put(8,1.1){\picfx}
    \end{picture}
\end{center}
Inserting (\ref{eqn:Consistencycheck1}) into the left side of (\ref{eqn:Consistencycheck2}) and using Appendix D, we have
\begin{eqnarray*}
  &&\hspace{-2cm}
  \Bigl[\exp(\frac{Nx}{2})\hat{P}_{L_+}^{(4)}
  -\exp(-\frac{Nx}{2})\hat{P}_{L_-}^{(4)}
  -(e^{\frac{x}{2}}-e^{-\frac{x}{2}})\hat{P}_{L_0}^{(4)}~\Bigr]^{(4)}\\
  &=&-\frac{(N^2-1)}{8}x^3~V_8+\frac{N(N^2-1)}{16}x^4~V_9
  -\frac{(N^2-1)}{16N}x^4~V_{10},
\end{eqnarray*}
where,

\begin{eqnarray*}
  &&\bullet~~V_{8}=
  \Bigl\{\sum_{s=\pm}s~v_{3.2}\bigl(\{K_n^s,K_{n+1}^s\}\bigr)\Bigr\}
  +2\Bigl\{\sum_{i=n}^{n+1}v_2(K_i^+)\Bigr\}
  -2v_2(K_n^0)-\frac{1}{3},\nonumber\\
  &&\bullet~~V_{9}=
  \Bigl\{\sum_{s=\pm}s~v_{4.3}\bigl(\{K_n^s,K_{n+1}^s\}\bigr)\Bigr\}
  +2\Bigl\{\sum_{i=n}^{n+1}v_{3.1}(K_i^+)\Bigr\}-2v_{3.1}(K_n^0)\nonumber\\
  &&\hspace{2cm}+\frac{p}{2}
  \Bigl\{\sum_{s=\pm}s~v_{3.2}\bigl(\{K_n^s,K_{n+1}^s\}\bigr)\Bigr\}
  -\frac{1}{2}
  \Bigl\{\sum_{s=\pm}v_{3.2}\bigl(\{K_n^s,K_{n+1}^s\}\bigr)\Bigr\},
  \nonumber\\
  &&\bullet~~V_{10}=\sum_{i=1}^{n-1}\biggl[
  \Bigl\{\sum_{s=\pm}s~v_{4.4}\bigl(\{K_i,K_n^s,K_{n+1}^s\}\bigr)\Bigr\}
  +2\Bigl\{\sum_{j=n}^{n+1}v_{3.2}\bigl(\{K_i,K_j^+\}\bigr)\Bigr\}
  \nonumber\\
  &&\hspace{8cm}
  -2v_{3.2}\bigl(\{K_i,K_{n}^0\}\bigr)\biggr],\nonumber\\
  &&
\end{eqnarray*}
where $p=v_{1}\bigl(\{K_n^+,K_{n+1}^+\}\bigr)-1$. 
Inserting the Guass diagram formula (Theorem \ref{th:mainresult}) into each $V_i$, we find out all of these equations vanishes identically $V_i=0$ ($i=8,9,10$). This shows the skein relation (\ref{eqn:Consistencycheck2}) holds in case of (2). $\square$

\vspace{1cm}
\noindent{\bf Acknowledgements}

The author would like to express his hearty thanks to his thesis advisor Professor T.Eguchi for his advise and constant encouragement, and also to Professor T.Kohno for his kind suggestion for further discussion.

\vspace{0.5cm}
\noindent{\bf Appendix A}
\vspace{0.3cm}

We list the chord diagrams of $n$-circles without any isolated chord up to degree four (see the proof of Theorem \ref{th:Kontsevich}).

For convenience, we omit the circles which have no chord. For example,
we write $\bigl(\marubb\bigr)$ instead of $\bigl(\marubb\maru\cdots\maru\bigr)$
, etc.
\begin{eqnarray*} 
  &&\bullet~~\bar{\mathfrak D}_2=\Bigl\{\marubb,\maruca\Bigr\}
  \\
  &&\bullet~~\bar{\mathfrak D}_3=\Bigl\{\marudd,\marude,
  \marueb,\maruec,\marued,\marufc\Bigr\}
  \\
  &&\bullet~~\bar{\mathfrak D}_4=\Bigl\{
  \marujk,\marugb,\maruga,\marugc,\marugd,\maruge,\marugf,
  \\
  &&\hspace{2cm}\marukg,\maruha,\maruhb,\maruhc,\\
  &&\hspace{2cm}\marukn,\maruhd,\marufd,\marukq,\\
  &&\hspace{2cm}\maruhe,\maruhf,(\marubb\marubb),(\marubb\maruca),\\
  &&\hspace{2cm}\maruhk,\marufg,\maruhj,
  \marufh,\maruhl,\marufk,\marufj \Bigr\}
\end{eqnarray*}

\vspace{0.5cm}
\noindent{\bf Appendix B}
\vspace{0.3cm}

We expand chord diagrams into CC diagrams, using AS, IHX, STU relation as follows (see the proof of Theorem \ref{th:Kontsevich}). 
\begin{eqnarray*}
  &&\bullet~
  \marubb=\maruba+(-\frac{1}{2})\marubc\nonumber\\
  &&\bullet~
  \marudd=\marudh+2(-\frac{1}{2})\marudf+(-\frac{1}{2})^2\marudg\nonumber\\
  &&\bullet~
  \marude=\marudh+3(-\frac{1}{2})\marudf+2(-\frac{1}{2})^2\marudg\nonumber\\
  &&\bullet~
  \marueb=\maruea+(-\frac{1}{2})\maruef\nonumber\\
  &&\bullet~
  \marued=\maruec+(-\frac{1}{2})\maruef\nonumber\\
\end{eqnarray*}

\begin{eqnarray*}
  &&\bullet~
  \marujk=\marukf+2(-\frac{1}{2})\maruke+(-\frac{1}{2})^2\marujl\nonumber\\
  &&\bullet~
  \marugb=\marukf+3(-\frac{1}{2})\maruke+2(-\frac{1}{2})^2\marukd
  +(-\frac{1}{2})^2\marujl+(-\frac{1}{2})^3\marufa\nonumber\\
  &&\bullet~
  \maruga=\marukf+3(-\frac{1}{2})\maruke+3(-\frac{1}{2})^2\marukd
  +(-\frac{1}{2})^3\marufa\nonumber\\
  &&\bullet~
  \marugc=\marukf+4(-\frac{1}{2})\maruke+4(-\frac{1}{2})^2\marukd
  +(-\frac{1}{2})^2\marujl+2(-\frac{1}{2})^3\marufa\nonumber\\
  &&\bullet~
  \marugd=\marukf+4(-\frac{1}{2})\maruke+4(-\frac{1}{2})^2\marukd
  +2(-\frac{1}{2})^2\marujl
  \\&&\hspace{8cm}
  +4(-\frac{1}{2})^3\marufa+\marufb\nonumber\\
  &&\bullet~
  \maruge=\marukf+5(-\frac{1}{2})\maruke+6(-\frac{1}{2})^2\marukd
  +2(-\frac{1}{2})^2\marujl
  \\&&\hspace{8cm}
  +5(-\frac{1}{2})^3\marufa+\marufb\nonumber\\
  &&\bullet~
  \marugf=\marukf+6(-\frac{1}{2})\maruke+8(-\frac{1}{2})^2\marukd
  +3(-\frac{1}{2})^2\marujl
  \\&&\hspace{8cm}
  +7(-\frac{1}{2})^3\marufa+\marufb
\end{eqnarray*}

\begin{eqnarray*}
  &&\bullet~
  \marukg=\marukh+(-\frac{1}{2})\maruki\nonumber\\
  &&\bullet~
  \maruha=\marukh+(-\frac{1}{2})\marukk+(-\frac{1}{2})\maruki
  \\&&\hspace{8cm}
  +(-\frac{1}{2})^2\maruff\nonumber\\
  &&\bullet~
  \maruhb=\marukh+2(-\frac{1}{2})\marukk+(-\frac{1}{2})^2\maruff\nonumber\\
  &&\bullet~
  \maruhc=\marukh+2(-\frac{1}{2})\marukk+(-\frac{1}{2})\maruki
  \\&&\hspace{8cm}
  +2(-\frac{1}{2})^2\maruff\nonumber\\
  &&\bullet~
  \marukn=\marukm+(-\frac{1}{2})\marufe\nonumber\\
  &&\bullet~
  \maruhd=\marukm+(-\frac{1}{2})\marufe+(-\frac{1}{2})\marukk
  \\&&\hspace{8cm}
  +(-\frac{1}{2})^2\maruff\nonumber\\
  &&\bullet~
  \marukq=\marufd+(-\frac{1}{2})\marufe\nonumber\\
  &&\bullet~
  \maruhe=\marufd+2(-\frac{1}{2})\marufe+(-\frac{1}{2})^2\maruff\nonumber\\
  &&\bullet~
  \maruhf=\marukp+2(-\frac{1}{2})\marukk+(-\frac{1}{2})^2\maruff\nonumber\\
  &&\bullet~
  \maruhk=\marukr+(-\frac{1}{2})\marufi\nonumber\\
  &&\bullet~
  \maruhj=\marufg+(-\frac{1}{2})\marufi\nonumber\\
  &&\bullet~
  \maruhl=\marufh+(-\frac{1}{2})\marufi\nonumber\\
\end{eqnarray*}

Notice we set the diagrams 
\begin{eqnarray*}
  &&\marua, \maruba, \marudh, \marudf, \maruea, \\
  &&\marukf, \maruke, \marukd, \marukh, \marukk, \marukm, \marukr, 
  \mbox{etc}, 
\end{eqnarray*}
 to be 0 by framing independence.

\vspace{0.5cm}
\noindent{\bf Appendix C}
\vspace{0.3cm}

We compute each cofficient of the CC diagram in (\ref{eqn:CCdiagram}) 
(see the proof of Theorem \ref{th:Kontsevich}).
\begin{eqnarray*}
  &&\hspace{-0.5cm}\bullet 
  \Bigl(\hbox{the cofficient of} ~\bigl(-\frac{1}{2}\bigr)\marubc\Bigr) 
  =\sum_{i=1}^n
  \hidari \K_i, \marubb\migi,
  \nonumber\\
  &&\hspace{-0.5cm}\bullet\Bigl(\hbox{the cofficient of} ~\maruca\Bigr) 
  =\sum_{i<j}
  \hidari \{\K_i,\K_j\}, \maruca\migi
  \nonumber\\
  &&\hspace{4.4cm}
  =\sum_{i<j}\frac{1}{2}\hidari \{\K_i,\K_j\}, \marucb\migi^2,\\
  &&\hspace{-0.5cm}\bullet 
  \Bigl(\hbox{the cofficient of} ~\bigl(-\frac{1}{2}\bigr)^2\marudg\Bigr) 
  =\sum_{i=1}^n
  \hidari \K_i, \marudd+2\marude\migi,
  \nonumber\\
  &&\hspace{-0.5cm}\bullet\Bigl(\hbox{the cofficient of} ~\maruec\Bigr) 
  =\sum_{i<j}
  \hidari \{\K_i,\K_j\}, \maruec+\marued\migi
  \nonumber\\
  &&\hspace{4.4cm}
  =\sum_{i<j}\frac{1}{3!}\hidari \{\K_i,\K_j\}, \marucb\migi^3,
  \nonumber\\
  &&\hspace{-0.5cm}\bullet 
  \Bigl(\hbox{the cofficient of} ~\bigl(-\frac{1}{2}\bigr)\maruef\Bigr) 
  \nonumber\\
  &&\hspace{3cm}
  =\sum_{i<j}\hidari \{\K_i,\K_j\}, \marueb+\marued\migi,\\
  &&\hspace{-0.5cm}\bullet
  \Bigl(\hbox{the cofficient of} ~\marufc\Bigr) 
  =\sum_{i<j<k}
  \hidari \{\K_i,\K_j,\K_k\}, \marufc\migi
  \nonumber\\
  &&\hspace{0cm}=
  \sum_{
  1\le i< j<k\le n}
  \hidari \{\K_i,\K_j\}, \marucb\migi \hidari \{\K_j,\K_k\}, \marucb\migi
  \nonumber\\
  &&\hspace{5cm}\times
  \hidari \{\K_k,\K_i\}, \marucb\migi\\
  &&\hspace{-0.5cm}\bullet 
  \Bigl(\hbox{the cofficient of} ~\bigl(-\frac{1}{2}\bigr)^2\marujl\Bigr) 
  \nonumber\\
  &&\hspace{1.5cm}
  =\sum_{i=1}^n
  \hidari \K_i, \marujk+\marugb+\marugc+2\marugd+2\maruge+3\marugf\migi
  \nonumber\\
  &&\hspace{5cm}
  +\sum_{i<j}\hidari \{\K_i,\K_j\}, \{\marubb\marubb\}\migi\nonumber\\
  &&\hspace{1.5cm}
  =\frac{1}{2}\biggl\{\sum_{i=1}^n\hidari \K_i, \marubb\migi\biggr\}^2,
  \nonumber\\
  &&\hspace{-0.5cm}\bullet
  \Bigl(\hbox{the cofficient of} ~\bigl(-\frac{1}{2}\bigr)^3\marufa\Bigr) 
  \nonumber\\
  &&\hspace{1.5cm}
  =\sum_{i=1}^n\hidari \K_i,\maruga+\marugb
  +2\marugc+4\marugd+5\maruge+7\marugf\migi, \nonumber\\
  &&\hspace{-0.5cm}\bullet
  \Bigl(\hbox{the cofficient of} ~\marufb\Bigr) 
  =\sum_{i=1}^n\hidari \K_i,\marugd+\maruge+\marugf\migi, \nonumber\\
  &&\hspace{-0.5cm}\bullet
  \Bigl(\hbox{the cofficient of} ~\bigl(-\frac{1}{2}\bigr)\maruki\Bigr) 
  \nonumber\\
  &&\hspace{1.5cm}
  =\sum_{i<j} 
  \hidari \{\K_i,\K_j\}, \marukg+\maruha+\maruhc\migi
  \nonumber\\
  &&\hspace{2cm}
  +\sum_{i<j<k}\hidari \{\K_i,\K_j,\K_k\}, \{\marubb\maruca\}\migi\nonumber\\
  &&\hspace{1.5cm}=\biggl\{\sum_{i=1}^n\hidari \K_i, \marubb\migi\biggr\}
  \biggl\{\sum_{i<j}\frac{1}{2}\hidari \{\K_i,\K_j\}, \marucb\migi\biggr\}^2,\\
  &&\hspace{-0.5cm}\bullet
  \Bigl(\hbox{the cofficient of} ~\bigl(-\frac{1}{2}\bigr)^2\maruff\Bigr) 
  \nonumber\\
  &&\hspace{1.5cm}
  =\sum_{i<j} \hidari \{\K_i,\K_j\},\maruha+\maruhb+2\maruhc\nonumber\\
  &&\hspace{4.5cm}+\maruhd+\maruhe+\maruhf\migi,
  \nonumber\\
  &&\hspace{-0.5cm}\bullet
  \Bigl(\hbox{the cofficient of} ~\bigl(-\frac{1}{2}\bigr)\marufe\Bigr) 
  \nonumber\\
  &&\hspace{0cm}
  =\sum_{i<j}
  \hidari \{\K_i,\K_j\}, \marukn+\maruhd+\marukq+2\maruhe\migi
  \nonumber\\
  &&\hspace{0cm}=\sum_{i<j}\frac{1}{2}\hidari \{\K_i,\K_j\}, \marucb\migi
  \hidari \{\K_i,\K_j\}, \marueb+\marued\migi,\\
  &&\hspace{-0.5cm}\bullet
  \Bigl(\hbox{the cofficient of} ~\marufd\Bigr) 
  \nonumber\\
  &&\hspace{1cm}
  =\sum_{i<j}
  \hidari \{\K_i,\K_j\}, \marufd+\marukq+\maruhe\migi
  \nonumber\\
  &&\hspace{1cm}=\sum_{i<j}\frac{1}{4!}\hidari \{\K_i,\K_j\}, \marucb\migi^4,\\
  &&\hspace{-0.5cm}\bullet
  \Bigl(\hbox{the cofficient of} ~(-\frac{1}{2})\marufi\Bigr) 
  \nonumber\\
  &&\hspace{2cm}
  =\sum_{i<j<k}
  \hidari \{\K_i,\K_j.\K_k\}, \maruhj+\maruhk+\maruhl\migi
  \nonumber\\
  &&\hspace{-0.5cm}\bullet
  \Bigl(\hbox{the cofficient of} ~\marufg\Bigr) 
  =\sum_{i<j<k}
  \hidari \{\K_i,\K_j.\K_k\}, \marufg+\maruhj\migi
  \nonumber\\
  &&\hspace{0.5cm}=\sum_{
  {{\scriptstyle 1\le i< j<k\le n}\atop
  {\scriptstyle 1\le j<i<k\le n}}\atop
  {\scriptstyle 1\le j<k<i\le n}
  }
  \frac{1}{2}\hidari \{\K_i,\K_j\}, \marucb\migi^2
  \frac{1}{2}\hidari \{\K_i,\K_k\}, \marucb\migi^2\\
  &&\hspace{-0.5cm}\bullet
  \Bigl(\hbox{the cofficient of} ~\marufh\Bigr) 
  =\sum_{i<j<k}
  \hidari \{\K_i,\K_j.\K_k\}, \marufh+\maruhl\migi
  \nonumber\\
  &&\hspace{0.5cm}=
  \sum_{
  {{\scriptstyle 1\le i< j<k\le n}\atop
  {\scriptstyle 1\le j<i<k\le n}}\atop
  {\scriptstyle 1\le j<k<i\le n}
  }\hidari \{\K_i,\K_j\}, \marucb\migi \hidari \{\K_i,\K_k\}, \marucb\migi
  \nonumber\\[-0.8cm]
  &&\hspace{7cm}\times
  \frac{1}{2}\hidari \{\K_j,\K_k\}, \marucb\migi^2\\
  &&\hspace{-0.5cm}\bullet
  \Bigl(\hbox{the cofficient of} ~\marufk\Bigr) 
  =\sum_{i<j<k<l}
  \hidari \{\K_i,\K_j.\K_k,\K_l\}, \marufk\migi
  \nonumber\\
  &&\hspace{0cm}=
  \sum_{
  {{\scriptstyle 1\le i< j<k<l\le n}\atop
  {\scriptstyle 1\le i<k<j<l\le n}}\atop
  {\scriptstyle 1\le i<k<l<j\le n}}
  \frac{1}{2}\hidari \{\K_i,\K_j\}, \marucb\migi^2
  \frac{1}{2}\hidari \{\K_k,\K_l\}, \marucb\migi^2\\
  &&\hspace{-0.5cm}\bullet
  \Bigl(\hbox{the cofficient of} ~\marufj\Bigr) 
  =\sum_{i<j<k<l}
  \hidari \{\K_i,\K_j.\K_k,\K_l\}, \marufj\migi
  \nonumber\\
  &&\hspace{1cm}=
  \sum_{
  {{\scriptstyle 1\le i< k<j<l\le n}\atop
  {\scriptstyle 1\le i<j<k<l\le n}}\atop
  {\scriptstyle 1\le i<j<l<k\le n}}
  \hidari \{\K_i,\K_j\}, \marucb\migi
  \hidari \{\K_j,\K_k\}, \marucb\migi\\
  &&\hspace{1.5cm}\times
  \hidari \{\K_k,\K_l\}, \marucb\migi
  \hidari \{\K_l,\K_i\}, \marucb\migi.
\end{eqnarray*}

\\
\noindent{\bf Appendix D}
\vspace{0.3cm}

We list the acutual table of the weight system. The computaion is straightforward, using Definition \ref{Def:weightsystem}. 
\begin{eqnarray*}   
  &&\bullet~~
  W_{su(N)}\h\Bigl(-\frac{1}{2}\Bigr)\marubc\m=-x^2\frac{N^2-1}{4},\\ 
  &&\bullet~~
  W_{su(N)}\h\Bigl(-\frac{1}{2}\Bigr)^2\marudg\m
  =x^3\frac{N(N^2-1)}{8},\\ 
  &&\bullet~~
  W_{su(N)}\h\Bigl(-\frac{1}{2}\Bigr)^3\marufa\m
  =-x^4\frac{N^2(N^2-1)}{16}, \\
  &&\bullet~~W_{su(N)}\h\marufb\m
  =x^4\frac{(N^2-1)(N^2+2)}{16}, \\
  &&\bullet~~W_{su(N)}\h\maruca\m 
  =x^2\frac{(N^2-1)}{4N^2}, \\
  &&\bullet~~W_{su(N)}\h\maruec\m
  =x^3\frac{(N^2-1)(N^2-2)}{8N^3}, \\
  &&\bullet~~W_{su(N)}\h\Bigl(-\frac{1}{2}\Bigr)\maruef\m
  =-x^3\frac{N^2-1}{8N}, \\
  &&\bullet~~W_{su(N)}\h\marufc\m
  =x^3\frac{N^2-1}{8N^3}, \\
  &&\bullet~~W_{su(N)}\h\marufd\m
  =x^4\frac{(N^2-1)(N^4-3N^2+3)}{16N^4},\nonumber\\
  &&\bullet~~W_{su(N)}\h\Bigl(-\frac{1}{2}\Bigr)\marufe\m
  =-x^4\frac{(N^2-1)(N^2-2)}{16N^2},\nonumber\\
  &&\bullet~~W_{su(N)}\h\Bigl(-\frac{1}{2}\Bigr)^2\maruff\m
  =x^4\frac{(N^2-1)}{16}, \\
  &&\bullet~~W_{su(N)}\h\marufg\m
  =x^4\frac{(N^2-1)^2}{16N^4}, \\
  &&\bullet~~W_{su(N)}\h\marufh\m
  =x^4\frac{(N^2-1)(N^2-2)}{16N^4}, \\
  &&\bullet~~W_{su(N)}\h\Bigl(-\frac{1}{2}\Bigr)\marufi\m
  =-x^4\frac{(N^2-1)}{16N^2}, \\
  &&\bullet~~W_{su(N)}\h\marufk\m
  =x^4\frac{(N^2-1)^2}{16N^4}, \\
  &&\bullet~~W_{su(N)}\h\marufj\m
  =x^4\frac{(N^2-1)}{16N^4}.
\end{eqnarray*}  

\newpage


\begin{thebibliography}{99}
  \bibitem{BN}D.Bar-Natan: On the Vassiliev knot invariants, Topology {\bf 34} (1995), 423-472.
  \bibitem{HS}A.C Hirshfeld and U.Sassenberg: Derivation of the total twist from Chern-Simons theory, Jounal of Knot Theory and Its Ramifications, Vol.5, No.4 (1996), 489-515.
  \bibitem{HS2}A.C Hirshfeld and U.Sassenberg: Explicit formulation of a third order finite knot invariant derived from Chern-Simons theory, Jounal of Knot Theory and Its Ramifications, Vol.5, No.6 (1996), 805-847.
  \bibitem{HS3}A.C Hirshfeld, U.Sassenberg and T.Kloker: A combinatorial link invariant of finite type derived from Chern-Simons theory, Jounal of Knot Theory and Its Ramifications, Vol.6, No.2 (1997), 243-280.
  \bibitem{Ko}M.Kontsevich: Vassiliev's knot invariants, Adv. Sov. Math {\bf 16} (1993), 137-150.
  \bibitem{La}J.M.F. Labastida and E.Perez: Combinatorial formula for Vassiliev Invariants from Chern-Simons Gauge theory, preprint hep-th/9807155.
  \bibitem{La2}J.M.F. Labastida: Chern-Simons Gauge theory: Ten Years After, hep-th/9905057.
  \bibitem{LMO}T.T.Q.Le, J.Murakami and T.Ohtsuki: On a universal perturbative invariant of 3-manifolds, Topology {\bf 37}, No.3 (1998), 539-574.
  \bibitem{Oh}T.Ohtsuki: Combinatorial quantum method in 3-dimensional topology, MSJ Memoirs {\bf 3}, Math. Soc. Japan, 1999.  
  \bibitem{Po}M.Polyak and O.Viro: Gauss diagram formulas for Vassiliev invariant, International Mathematics Research Notices {\bf 11} (1994), 445-452.
  \bibitem{Wi}E.Witten: Quantum field theory and the Jones polynomial, Comm. Math. Phys. {\bf 121} (1989), 351-399.
\end{thebibliography}
\end{document}